\documentclass[a5paper,11pt]{book}
\usepackage{amssymb}
\usepackage{amsmath}
\usepackage{theorem}
\usepackage{xspace}
\usepackage{cite}
\usepackage{endnotes}
\usepackage{mathrsfs}
\usepackage[cp1251]{inputenc}
\usepackage[russian,english]{babel}

\renewcommand{\notesname}{Примечания}   
\makeatletter
\renewcommand{\enoteheading}{\chapter*{\notesname
  \@mkboth{\notesname}{\notesname}
  \addcontentsline{toc}{chapter}{Примечания} }
  \leavevmode\par\vskip-\baselineskip}
\renewcommand{\@makeenmark}{\hbox{$^{\@theenmark )}$}}
\renewcommand{\enoteformat}{\rightskip\z@ \leftskip\z@ \parindent=1.8em
     \leavevmode\llap{\hbox{$^{\@theenmark )}$}}}
\makeatother

\renewcommand{\thechapter}{\Roman{chapter}}

\renewcommand{\chaptermark}[1]{\markboth{\chaptername\ \thechapter. \ #1}{}}

\makeatletter
\renewcommand{\@makechapterhead}[1]{  \vspace*{20\p@}  {\parindent \z@ \raggedright \normalfont
    \ifnum \c@secnumdepth >\m@ne
      \if@mainmatter
        \large \bfseries \@chapapp\space \thechapter        \par\nobreak
        \vskip 10\p@      \fi
    \fi
    \interlinepenalty\@M
    \LARGE \bfseries #1\par\nobreak    \vskip 40\p@
  }}
\renewcommand{\@makeschapterhead}[1]{  \vspace*{20\p@}  {\parindent \z@ \raggedright
    \normalfont
    \interlinepenalty\@M
    \LARGE \bfseries  #1\par\nobreak
    \vskip 40\p@
  }}
\makeatother
\makeatletter
\def\@chapter[#1]#2{\ifnum \c@secnumdepth >\m@ne
                       \if@mainmatter
                         \refstepcounter{chapter}                         \typeout{\@chapapp\space\thechapter.}                         \addcontentsline{toc}{chapter}                                   {\protect \chaptername\ \thechapter.\ \ #1}                       \else
                         \addcontentsline{toc}{chapter}{#1}                       \fi
                    \else
                      \addcontentsline{toc}{chapter}{#1}                    \fi
                    \chaptermark{#1}                    \addtocontents{lof}{\protect\addvspace{10\p@}}                    \addtocontents{lot}{\protect\addvspace{10\p@}}                    \if@twocolumn
                      \@topnewpage[\@makechapterhead{#2}]                    \else
                      \@makechapterhead{#2}                      \@afterheading
                    \fi}
\makeatother
\makeatletter
\def\@sect#1#2#3#4#5#6[#7]#8{  \ifnum #2>\c@secnumdepth
    \let\@svsec\@empty
  \else
    \refstepcounter{#1}    \protected@edef\@svsec{\@seccntformat{#1}\relax}  \fi
  \@tempskipa #5\relax
  \ifdim \@tempskipa>\z@
    \begingroup
      #6{        \@hangfrom{\hskip #3\relax\@svsec}          \interlinepenalty \@M #8\@@par}    \endgroup
    \csname #1mark\endcsname{#7}    \addcontentsline{toc}{#1}{      \ifnum #2>\c@secnumdepth \else
        \protect\numberline{\S \csname the#1\endcsname}      \fi
      #7}  \else
    \def\@svsechd{      #6{\hskip #3\relax
      \@svsec #8}      \csname #1mark\endcsname{#7}      \addcontentsline{toc}{#1}{        \ifnum #2>\c@secnumdepth \else
          \protect\numberline{\csname the#1\endcsname}        \fi
        #7}}  \fi
  \@xsect{#5}}
\makeatother
\newtheorem{Pa}{Paper}[section]

\newtheorem{corollary}[Pa]{Следствие}  

\newtheorem{definition}[Pa]{Определение}  

\newtheorem{lemma}[Pa]{Лемма}  

\begin{document}

\selectlanguage{russian}

\setcounter{endnote}{8}

\renewcommand{\contentsname}{Содержание}
\renewcommand{\bibname}{Литература}

\thispagestyle{empty}

\newpage

\begin{center}
{\large БЕЛОРУССКИЙ ГОСУДАРСТВЕННЫЙ УНИВЕРСИТЕТ}
\quad \\
\quad \\
\quad \\
\bigskip {\bf \LARGE Александр Киселев }\bigskip\\
\quad \\
\quad \\
\quad \\
{\bf \Huge Недостижимость}
\quad \\
\quad \\
{\bf \Huge и}
\quad \\
\quad \\
{\bf \Huge субнедостижимость}\\
\quad \\
\quad \\
{\Large В двух частях} \\
\smallskip
{\Large Часть II} \\
\quad \\
\quad \\

\quad \\
\quad \\
\quad \\
\quad \\
\quad \\
\quad \\
{Минск} \\
\smallskip
{``Издательский центр БГУ''} \\
\smallskip
{2011}

\end{center}

\newpage

\thispagestyle{empty}

\noindent УДК 510.227
\\

{\small \textbf{Киселев, А. А.} Недостижимость и
субнедостижимость. В 2 ч. Ч. 2 / Александр Киселев. -- Минск: Изд.
центр БГУ, 2011. -- \pageref{end} с. -- ISBN 978-985-476-935-6.}
\\

{\footnotesize Данная работа представляет собой перевод с
английского языка монографии Киселева А. А. под тем же названием,
 содержащей доказательство (в $ZF$)
несуществования больших кардиналов, 1-е издание которой вышло в
свет в 2000 году.  Часть II содержит приложения аппарата
субнедостижимых кардиналов и его основных средств --- теорий
редуцированных формульных спектров и матриц, диссеминаторов и
других, которые используются в этом доказательстве и представлены
здесь в их более прозрачной и детализированной форме. Большое
внимание уделяется более глубокой разработке и культивированию
базовых идей, служащих основаниями для основных конструкций и
рассуждений. Доказательство теоремы о несуществовании больших
кардиналов представлено в его подробном виде. Излагается несколько
простых следствий этой теоремы и других хорошо известных
классических результатов.

Предназначено для специалистов по теории множеств и математической
логике, а также для преподавателей и студентов факультетов
математического профиля.

Библиогр.: 47 назв.
\\
Перевод осуществлён по изданию: \textbf{Kiselev, Alexander.}
Inaccessibility and Subinaccessibility. In 2 pt. Pt 2 / Alexander
Kiselev. -- 2nd ed., enrich. and improv. Minsk : Publ. center of
BSU, 2010. }
\quad \\
\begin{center}
{Р\;е\;ц\;е\;н\;з\;е\;н\;т\;ы}
\\
профессор {\em П. П. Забрейко;}
\\
профессор {\em А. В. Лебедев}
\quad \\
\quad \\
\noindent Математическая классификация тем (2000): \\
03E05, 03E15, 03E35, 03E55, 03E60
\quad \\
\end{center}

\noindent {\scriptsize {\bf ISBN 978-985-476-935-6 (ч. 2)} \hfill
$\copyright$ Киселев А. А., 2011}

\noindent {\scriptsize {\bf ISBN 978-985-476-597-6 (pt. 2)} \hfill
$\copyright$ Kiselev Alexander, 2010}

\newpage

\thispagestyle{empty}

\hfill \textit{Посвящается моей матери Анне}

\newpage

\thispagestyle{empty}

\chapter*{Благодарности}
 \hspace*{1em}  Автор высказывает свои первые слова глубокой благодарности Хан\-не
 Калиендо за понимание значительности темы и за сердечную
воодушевляющую помощь в продвижении работы.

Особая признательность выражается проф. С. Р. Когаловскому,
который научил автора Теории Иерархий, и  проф. Акихиро Канамори
за их бесценную воодушевляющую поддержку, придавшую необходимый
импульс завершению работы.

Автор хотел бы также выразить особую благодарность проф. А. В.
Лебедеву и проф. П. П. Забрейко за многолетнюю поддержку его
работ; большой интеллектуальный и моральный  долг высказывается им
обоим за их практическую и духовную помощь.

 Глубокую благодарность автор высказывает А.А.
 Лапцевичу, начальнику Минского  государственного высшего авиационного колледжа;
 С. В. Сизикову, заместителю начальника; А.И. Рипинскому, декану факультета гражданской авиации;
 А. И. Кириленко, заведающего кафедрой естественнонаучных дисциплин,
 за создание условий для плодотворной работы над этим изданием.

Глубокая благодарность выражается  проф. А. В. Тузикову и Ю.
Прокопчуку, которые оказывали автору интенсивную экспертную помощь
в наборе текстов.

Проф. В. М. Романчак  оказывал автору материальную и моральную
помощь в продолжении самого сложного периода исследований темы и
автор направляет ему много слов глубокой благодарности.

Большая благодарность высказывается также Людмиле Лаптёнок,
которая осуществила много первоначальных трудных наборов
предыдущих текстов автора, подготовивших  эту работу.

Но самая глубокая сердечная благодарность высказывается Надежде
Забродиной за многолетние воодушевление и поддержку, без которых
эта работа была бы значительно затруднена.

Группа других специалистов, которые в течении долгих лет
поддерживали и воодушевляли автора, слишком велика, чтобы их
перечислить, и автор выражает им всем свою великую благодарность.

\newpage
{} \thispagestyle{empty}

\tableofcontents

\newpage {} \thispagestyle{empty}

\newpage

\chapter*{Введение}

\addcontentsline{toc}{chapter}{Введение} %
\markboth{Введение}{Введение}

\setcounter{equation}{0}

\hspace*{1em} Эта работа представляет собой  непосредственное
продолжение предыдущей части~I~\cite{Kiselev11} и составляет с ней
единый текст. Автору представилось естественным организовать эту
работу таким образом, что в ней продолжаются все
нумерации~\cite{Kiselev11}, в том числе нумерации глав,
параграфов, определений, формул, утверждений и даже нумерации
примечаний и библиографических ссылок.
\\
Также формулы, понятия или символы, используемые в этой работе без
пояснений, уже были введены и использовались в работах
автора~\cite{Kiselev9}, ~\cite{Kiselev10}, в
части~I~\cite{Kiselev11}, или считаются общепринятыми, или
используются в замечательной книге Йеха~\cite{Jech} ``Lectures in
Set Theory with Particular Emphasis on the Method of Forcing''
(русский перевод `` Теория множеств и метод форсинга'', М., 1973),
содержащей многие базовые понятия и сведения и гораздо более того;
поэтому они будут считаться известными и часто будут
использоваться без пояснений.
\\
Таким образом, читателю было бы удобнее предварительно
ознакомиться с~\cite{Kiselev11} и с основными понятиями и
обозначениями этой работы хотя бы в общих чертах.
\\
Во всяком случае, было бы полезно ознакомиться предварительно с
планом всей этой работы и с кратким изложением предстоящего
развития основных её идей в том виде, как он изложен
в~\cite{Kiselev11} на стр.~11--20.
\\
Руководствуясь этими соображениями, читателю следует помнить, что
все ссылки на предшествующие параграфы с номерами меньшими 7
относятся к части~I~\cite{Kiselev11}, аналогично для нумерации
утверждений и т.д.

Что же касается содержания этой работы и технической стороны дела
в целом, следует отметить, что она следует её предыдущему
изданию~\cite{Kiselev8}  2000 года, но более систематическим
образом.

Также следует отметить, что в этом издании~\cite{Kiselev8} и
предшествующих работах автор стремился избежать использования
личных авторских понятий и символики (кроме самых необходимых),
так как он был обеспокоен трудностями и неприятием, которые они
могут вызвать у читателя.

Однако восприятие этого предыдущего издания~\cite{Kiselev8}
читателями показало, что подобные опасения бесполезны, поэтому
подобное использование неизбежно в любом случае.

Поэтому в настоящей работе автор предпринял иной подход и посчитал
более естественным вовлечь всю систему его личных авторских
понятий и определений в целом, которые он разрабатывал с 1976
года, так как она обладает технической и концептуальной
выразительностью и ведёт напосредственно к сути дела, и, таким
образом, было бы слишком исскуственным избегать её использования
(см. примечание 3~\cite{Kiselev11} в качестве примера). Как это
обычно бывает, некоторые утверждения получили свои усиления;
некоторые фрагменты, считавшиеся очевидными в предыдущих работах
автора, получили свою детализацию; некоторые аргументы изменили
своё расположение на более удобное; также иногда использованы
некоторые удобные переобозначения.

Но что касается основных конструкций и рассуждений следует
отметить, что данная работа следует изданию~\cite{Kiselev8} 2000
года, но в более  прояснённом виде. Основной результат этой работы
таков: система
\[
    ZF+\exists k \hspace{2mm} (k \; \mbox{\it слабо недосижимый кардинал})
\]
несовместна; все рассуждения проводятся в этой теории. Все слабо
недостижимые кардиналы становятся сильно недостижимыми в
конструктивном классе \ $L$ \ и поэтому рассуждения переносятся в
стандартную счётную основную исходную  модель
\[
    \mathfrak{M}=(L_{\chi ^{0}},\; \in , \; =)
\]
теории
\[
    ZF+V=L+\exists k \hspace{2mm} (k \; \mbox{\it слабо недостижымый кардинал}),
\]
и далее \ $k$ \ это наименьший недостижимый кардинал в \
$\mathfrak{M}$. \ В действительности в этой теории используются
только формулы некоторой ограниченной длины; более того, счётность
этой модели нужна только для некоторого технического удобства и
можно обойтись без неё (см. ``Введение''~\cite{Kiselev11}). В этой
модели \ $\mathfrak{M}$ \ строятся так называемые матричные
функции, обладающие несовместными свойствами монотонности и
немонотонности; этим противоречием устанавливается
\\

\noindent {\bf Основная теорема } ($ZF$)
\\
\hspace*{1cm} {\it Не существует слабо недостижимых кардиналов.}
\\

\noindent Это влечёт несуществование сильно недостижимых
кардиналов и несуществование всех других больших кардиналов. Эти
матричные функции конструируются и рассматриваются посредством
элементарного языка из формульных классов (см. определение
2.1~\cite{Kiselev11}) некоторого фиксированного уровня
 $>3$ над стандартной моделью
\[
    (L_k, \in, =)
\]
и затем все конструкции осуществляются посредством этого языка
(если обратное не устанавливается контекстом).
\\

\noindent В дополнение к этому в \S 12 приводятся некоторые
простые следствия основной теоремы и некоторых хорошо известных
результатов.

\setcounter{chapter}{1}

\chapter{Специальная теория: матричные функции}

\markboth{\chaptername\ \thechapter. \ Специальная
теория}{\chaptername\
\thechapter. \ Специальная теория} %

\setcounter{section}{6}

\section{Матричные \ $\protect\delta$\,-функции}

\setcounter{equation}{0}

\hspace*{1em} Здесь мы собираемся осуществить дальнейшее развитие
идеи доказательства основной теоремы и модифицировать простейшие
матричные функции \ $S^{<\alpha_1}_{\chi f} $ \ (см. определение
5.14~\cite{Kiselev11}) таким образом, чтобы их новый специальные
варианты -- так называемые \ $\alpha $-функции -- доставили
требуемое противоречие: они будут обладать свойством \
$\underline{\lessdot}$-монотонности и в то же время будут лишены
этого свойства.

Напомним, что простейшие матричные функции, которые
рассматривались в \S~5~\cite{Kiselev11}, обладают свойством
монотонности, но оказалось, что непосредственное доказательство
требуемого противоречия -- доказательство их немонотонности --
препятствуется следующим обстоятельством: некоторые существенные
свойства нижних уровней универсума не распротраняются до
кардиналов скачка матриц на их носителях, которые являются
значениями таких матричных функций.
\\
С целью разрушить это препятствие мы снабдим такие матрицы их
соответствующими диссеминаторами, и в результате простейшие
матричные функции будут преобразованы в их более сложные формы, \
$\alpha$-функции.

Однако непосредственное формирование этих функций представляется
значительно усложнённым и некоторые важные их особенности
немотивированными.
\\
Поэтому, чтобы представить их введение более прозрачным образом,
мы предварительно предпримем \emph{второе приближение} к идее
доказательства основной теоремы и обратимся к их более простой
форме, то есть к \ $\delta $-функциям.
\\
С этой цеью мы применим результаты \S 6~\cite{Kiselev11} для \
$m=n+1$ \ и фиксированного уровня \ $n>3$, \ но понятие
диссеминатора должно быть уточнено; все диссеминаторы в дальнйшем
будут уровня \ $n+1$ (см. определение 6.9~\cite{Kiselev11}).

\begin{definition}
\label{7.1.} \hfill {} \newline \hspace*{1em} Пусть
\begin{equation*}
\gamma < \alpha <\alpha _{1}\leq k.
\end{equation*}
\noindent \emph{1)}\quad Мы обозначаем через \
$\mathbf{K}_{n}^{\forall <\alpha _{1}}(\gamma ,\alpha )$ \
формулу:

\vspace{6pt}
\begin{equation*}
SIN_{n-1}^{<\alpha _{1}}(\gamma )\wedge \forall \gamma ^{\prime }\leq \gamma
\ (SIN_{n}^{<\alpha _{1}}(\gamma ^{\prime })\longrightarrow SIN_{n}^{<\alpha
}(\gamma ^{\prime }))~.
\end{equation*}
\vspace{0pt}

\noindent Если эта формула выполняется константами \ $\gamma $, \
$\alpha $, \ $\alpha _{1}$, \ то мы будем говорить, что \ $\alpha
$ \ сохраняет \ $SIN_{n}^{< \alpha _{1}}$-кардиналы \ $\leq \gamma
$ \ ниже \ $\alpha _{1}$.
\\
Если \ $S$ \ -- это матрица на носителе \ $\alpha $ \ и её
кардинал предскачка \ $\alpha _{\chi }^{\Downarrow }$ \ после \
$\chi$ \ сохраняет такие  кардиналы, то мы также будем говорить,
что \ $S$ \ на \ $\alpha $ \ сохраняет эти кардиналы ниже \
$\alpha _{1}$.
\newline

\noindent \emph{2)}\quad Мы обозначаем через \
$\mathbf{K}_{n+1}^{\exists }(\chi ,\delta ,\gamma ,\alpha ,\rho
,S)$ \ следующую \ $\Pi _{n-2}$-формулу:

\begin{equation*}
\sigma (\chi ,\alpha ,S)\wedge Lj^{<\alpha }(\chi )\wedge \chi
<\delta <\gamma <\alpha \wedge S\vartriangleleft \rho \leq \chi
^{+}\wedge \rho = \widehat{\rho }\wedge
\end{equation*}
\begin{equation*}
\wedge SIN_{n}^{<\alpha _{\chi }^{\Downarrow }}(\delta )\wedge
SIN_{n+1}^{<\alpha _{\chi }^{\Downarrow }}\left[ <\rho \right] (\delta ).
\end{equation*}
\vspace{0pt}

\noindent Здесь, напомним, \ $\Pi_{n-2}$-формула \ $\sigma(\chi,
\alpha, S)$ \ означает, что \ $S$ \ это сингулярная матрица на её
носителе \ $\alpha$, \ редуцированная к кардиналу \ $\chi$ \ (см.
определение 5.7~\cite{Kiselev11}); \ $\delta$ \ это диссеминатор
для \ $S$ \ на \ $\alpha$ \ с базой данных \ $\rho$ \ уровня \
$n+1$ \ (определение 6.9~\cite{Kiselev11}); верхние индексы \ $<
\alpha_{\chi}^{\Downarrow}$ \ означают ограничение формульных
кванторов кардиналом предскачка \ $\alpha_{\chi}^{\Downarrow}$ \
(см. также определения 2.3, 5.9~\cite{Kiselev11}); \
$\widehat{\rho}$ \ это замыкание \ $\rho$ \ относительно функции
пары; и \ $Lj^{<\alpha}(\chi)$ \ это \ $\Delta_1$-свойство
насыщенности кардинала \ $\chi$ \ ниже \ $\alpha$ \ (см.
определение~6.9~4) \cite{Kiselev11}):
\[
    \chi < \alpha \wedge SIN_{n-1}^{<\alpha}(\chi) \wedge
    \Sigma rng\big(\widetilde{\mathbf{S}}_n^{sin \vartriangleleft \chi}\big)
    \in B_{\chi} \wedge \sup dom \big(\widetilde{\mathbf{S}}_n^{sin \vartriangleleft
    \chi}\big) = \chi.
\]
Мы обозначаем через \ $\mathbf{K}^{<\alpha _{1}}(\chi ,\delta
,\gamma ,\alpha ,\rho ,S)$ \ формулу:

\vspace{6pt}
\begin{equation*}
\mathbf{K}_{n}^{\forall <\alpha _{1}}(\gamma ,\alpha _{\chi }^{\Downarrow
})\wedge \mathbf{K}_{n+1}^{\exists \vartriangleleft \alpha _{1}}(\chi
,\delta ,\gamma ,\alpha ,\rho ,S)\wedge \alpha <\alpha _{1}~.
\end{equation*}
\vspace{0pt}

\noindent \emph{3)}\quad Если эта формула выполняется константами
\ $\chi$, $\delta $, $\gamma $, $\alpha $, $\rho$, $S$, $\alpha
_{1}$, \ то мы будем говорить, что \ $\chi $, $\delta $, $\alpha
$, $\rho $, $S$ \ сильно допустимы для \ $\gamma $ \ ниже \
$\alpha _{1}$.
\newline Если некоторые из них фиксированы или
подразумеваются контекстом, то мы будем говорить, что остальные
также сильно недостижимы для них (и для \ $\gamma $) \ ниже \
$\alpha _{1}$. \newline

\noindent \emph{4)}\quad Матрица \ $S$ \ называется сильно
диссеминаторной матрицей или, короче, \ $\delta $-матрицей, сильно
допустимой на носителе \ $\alpha $ \ для \ \mbox{$\gamma =\gamma
_{\tau }^{<\alpha _{1}}$} \ ниже \ $\alpha _{1}$, \ если она
обладает некоторым диссеминатором \ $\delta <\gamma $ \ с базой \
$\rho $, \ сильно допустимым для них (также ниже \ $\alpha _{1}$).
\newline В каждом подобном случае \ $\delta $-матрица обозначается общим символом
 \ $\delta S$ \ или \ $S$.
\\
Если \ $\alpha_1=k$, \ или \ $\alpha_1$ \ указывается в контексте,
то верхние индексы \ $<\alpha_1$, $\vartriangleleft \alpha_1$ \
здесь и другие упоминания о \ $\alpha_1$ \ опускаются.
\hspace*{\fill} $\dashv$
\end{definition}

\noindent Далее вплоть до конца $\S$~7 понятия допустимости и \
$\delta$-матрицы будут рассматриваться как \textit{сильные}
понятия, поэтому термин ``сильно'' будет опускаться. Все матрицы
будут считаться \ $\delta$-матрицами; в качестве редуцирующего
кардинала \ $\chi $ \ далее будет использоваться полный кардинал \
$\chi ^{\ast }$ \ (см. определение~5.4~\cite{Kiselev11}) -- если
констекст не указывает на другой случай.
\\
Здесь следует обратить внимание на понятие насыщенности кардинала
\ $\chi$ \ ниже \ $\alpha$, \ то есть на \ $\Delta_1$-свойство \
$Lj^{<\alpha}(\chi)$; \ из леммы 5.5~\cite{Kiselev11} следует, что
\ $\chi^{\ast}$ \ это кардинал, насыщенный ниже всякого \ $\alpha
> \chi^{\ast}$, $\alpha \in SIN_{n-2}$.
\\
Символ \ $\chi ^{\ast }$ \ в обозначениях и написаниях формул
будет часто опускаться для некоторой краткости.
\\
Далее каждый ограничивающий кардинал \ $\alpha _{1}$ \ будет
принадлежать \ $SIN_{n-2}$ \ и выполнять условие

\begin{equation}\label{e7.1}
\chi ^{\ast }<\alpha _{1} \leq k \wedge A_{n}^{\vartriangleleft
\alpha_{1}} ( \chi ^{\ast } ) =\left\| u_{n}^{\vartriangleleft
\alpha_{1}} ( \underline{l} ) \right\|,
\end{equation}

\noindent или \ $\alpha_1 = k$ \ (если иное не предусмотрено
контекстом).

Кардинал \ $\alpha_1 \le k$ \ здесь с этим свойством будет
называться {\em эквиинформативным} (равно информативным) с
кардиналом \ $\chi^{\ast}$.
\\
Этот термин вводится здесь вот почему: в этом случае никакое \
$\Sigma_n$-утверждение \ $\varphi(\underline{l})$ \ не имеет
\textit{ординалов скачка} после \ $\chi^{\ast}$ \ ниже \
$\alpha_1$ \ (см. определение 2.4~\cite{Kiselev11}). \ Нетрудно
видеть, что это равносильно следующему: для всякого генерического
расширения \ $\mathfrak{M}[l]$ \ каждое \ $\Pi_n$-утверждение \
$\varphi(l)$, \ которое выполняется в \ $\mathfrak{M}[l]$ \ ниже \
$\chi^{\ast}$, \ выполняется также в этом расширении и ниже \
$\alpha_1$ \ благодаря (\ref{e7.1}) и \ $\alpha_1 \in \Pi_{n-2}$;
\ поэтому каждое \ $\Pi_n$-утверждение \ $\varphi(l)$ \
выполняется или нет в обоих случаях одновременно для каждого
расширения \ $\mathfrak{M}[l]$ \ (см. также примечание
7~\cite{Kiselev11} для иллюстрации важности этого понятия).

Следует обратить внимание на важный пример такого кардинала:  на
кардинал предскачка \ $\alpha_{\chi^{\ast}}^{\Downarrow}$ \ после
\ $\chi^{\ast}$ \ для всякого носителя матрицы \
$\alpha>\chi^{\ast}$. \ Кроме того, будет всегда полагаться для \
$\chi^{\ast}$ \  и \ $\alpha_1$, \ что
\[
    \forall \gamma < \alpha_1 \exists \gamma^{\prime} \in \left[\gamma,
    \alpha_1\right[ \ SIN_n^{<\alpha_1}(\gamma^{\prime}) \wedge
    cf(\alpha_1) \ge \chi^{\ast +}
\]
для удобства некоторых формульных преобразований.
\\
Символы граничения \ $ <\alpha _{1}$, \ $\vartriangleleft \alpha
_{1}$ \ будут опускаться, как всегда, если \ $\alpha _{1}=k$, \
или \ $\alpha_1$ \ подразумевается контекстом.

\begin{definition}
\label{7.2.} \hfill {} \newline \hspace*{1em} Пусть \ $\chi ^{\ast
}<\alpha _{1}$. \newline \quad \newline \emph{1)}\quad Мы называем
матричной \ $\delta $-функцией уровня \ $n$ \ ниже \ $\alpha _{1}$
\ редуцированной к \ $\chi ^{\ast }$ \ следующую \ функцию
\begin{equation*}
\delta S_{f}^{<\alpha _{1}}=(\delta S_{\tau }^{<\alpha
_{1}})_{\tau }~,
\end{equation*}
\noindent принимающую значения для\ $\tau$:
\begin{equation*}
\delta S_{\tau }^{<\alpha _{1}}=\min_{\underline{\lessdot }} \bigl
\{ S \vartriangleleft \chi^{\ast +} : \exists \delta ,\alpha ,\rho
< \gamma_{\tau+1}^{<\alpha_1} ~ \mathbf{K}^{<\alpha _{1}}(\delta
,\gamma _{\tau }^{<\alpha _{1}},\alpha ,\rho ,S) \bigr \};
\end{equation*}

\noindent \emph{2)}\quad следующие сопровождающие ординальные
функции определяются ниже \ $\alpha _{1}$:
\begin{equation*}
\check{\delta}_{f}^{<\alpha _{1}}=(\check{\delta}_{\tau }^{<\alpha
_{1}})_{\tau };\quad \rho _{f}^{<\alpha _{1}}=(\rho _{\tau
}^{<\alpha _{1}})_{\tau };\quad \alpha _{f}^{<\alpha _{1}}=(\alpha
_{\tau }^{<\alpha _{1}})_{\tau },
\end{equation*}

\noindent принимающие значения:
\begin{equation*}
\check{\delta}_{\tau }^{<\alpha _{1}} = \min \{\delta <
\gamma_{\tau}^{<\alpha_1}: \exists \alpha ,\rho <
\gamma_{\tau+1}^{<\alpha_1} ~\mathbf{K}^{<\alpha _{1}}(\delta
,\gamma _{\tau }^{<\alpha _{1}},\alpha ,\rho ,\delta S_{\tau
}^{<\alpha _{1}})\};\quad
\end{equation*}
\vspace{-6pt}
\begin{equation*}
\rho _{\tau }^{<\alpha _{1}} = \min \{\rho < \chi^{\ast +} :
\exists \alpha < \gamma_{\tau+1}^{<\alpha_1} ~ \mathbf{K}^{<\alpha
_{1}}(\check{\delta}_{\tau }^{<\alpha _{1}},\gamma _{\tau
}^{<\alpha _{1}},\alpha ,\rho ,\delta S_{\tau }^{<\alpha _{1}})\};
\end{equation*}
\vspace{-6pt}
\begin{equation*}
\alpha _{\tau }^{<\alpha _{1}} = \min \{\alpha <
\gamma_{\tau+1}^{<\alpha_1} : \mathbf{K} ^{<\alpha
_{1}}(\check{\delta}_{\tau }^{<\alpha _{1}},\gamma _{\tau
}^{<\alpha _{1}},\alpha ,\rho _{\tau }^{<\alpha _{1}},\delta
S_{\tau }^{<\alpha _{1}})\}. \qquad
\end{equation*}

\noindent Для каждой матрицы \ $\delta S_{\tau }^{<\alpha _{1}}$ \
эти функции определяют её производящий диссеминатор \
$\check{\delta}_{\tau }^{<\alpha _{1}}<\gamma _{\tau }^{<\alpha
_{1}}$ \ вместе с его базой \ $ \rho_{\tau }^{<\alpha _{1}}$ \ и
её носителем \ \mbox{$\alpha _{\tau }^{<\alpha _{1}}$}.
\hspace*{\fill} $\dashv$
\end{definition}

\noindent Используя лемму 6.8~\cite{Kiselev11} легко видеть, что
здесь \ $\check{\delta}_{\tau}^{<\alpha _{1}}$ \ -- это
минимальный диссеминатор с базой
\[
    \rho_{\tau}^{<\alpha _{1}} = \widehat{\rho_1}, \
    \rho_1 = Od(\delta S_{\tau}^{<\alpha_{1}}),
\]
то есть замыканием ординала \ $Od(\delta S_{\tau}^{<\alpha _{1}})$
\ относительно функции пары; поэтому такой диссеминатор будет
называться  {\em производящим собственным диссеминатором матрицы \
$\delta S_{\tau}^{<\alpha _{1}}$} \ на \ $\alpha_{\tau}^{<\alpha
_{1}}$ \ ниже \ $\alpha_1$ \ и обозначаться через \
$\check{\delta}_{\tau}^{S < \alpha_1}$ \ (см. также определение
6.9~2)~\cite{Kiselev11}), а его база \ $\rho_{\tau}^{<\alpha_1}$ \
будет обозначаться через \ $\rho_{\tau}^{S <\alpha_1}$.

Нетрудно получить следующие леммы из этих определений и лемм 5.15,
5.16~\cite{Kiselev11}:

\begin{lemma}
\label{7.3.} \hfill {} \newline \hspace*{1em} Для \ $\alpha
_{1}<k$ \ формулы \ $\mathbf{K}_n^{\forall < \alpha_1}$, \
$\mathbf{K}^{<\alpha _{1}}$ \ принадлежат \ $\Delta _{1}$ \ и
поэтому все функции
\begin{equation*}
\delta S_{f}^{<\alpha _{1}},\quad \check{\delta}_{f}^{<\alpha
_{1}},\quad \rho _{f}^{<\alpha _{1}},\quad \alpha _{f}^{<\alpha
_{1}}
\end{equation*}
 \ $\Delta _{1}$-определимы через \ $\chi ^{\ast },\alpha _{1}$.
\\
Для \ $\alpha _{1}=k$ \ формулы \ $\mathbf{K}_n^{\forall}$, \
$\mathbf{K}$ \ принадлежат \ $\Sigma _{n}$ \ и эти функции \
$\Delta _{n+1}$-определимы. \hspace*{\fill} $\dashv$
\end{lemma}

\begin{lemma}
\label{7.4.} \emph{(О абсолютности  \ $\delta $-функций)} \hfill
{}
\newline
\hspace*{1em} Пусть \ $\chi ^{\ast }<\gamma _{\tau +1}^{<\alpha
_{1}}<\alpha _{2}<\alpha _{1}\leq k$, \quad $\alpha _{2}\in
SIN_{n-2}^{<\alpha _{1}}$ \ и
\begin{equation*}
(\gamma _{\tau }^{<\alpha _{1}}+1)\cap SIN_{n}^{<\alpha _{2}}=(\gamma _{\tau
}^{<\alpha _{1}}+1)\cap SIN_{n}^{<\alpha _{1}},
\end{equation*}
тогда на множестве
\begin{equation*}
\{\tau ^{\prime }:\ \ \ \chi ^{\ast }\leq \gamma _{\tau ^{\prime }}^{<\alpha
_{2}}\leq \gamma _{\tau }^{<\alpha _{1}}\}
\end{equation*}
функции
\begin{equation*}
\delta S_{f}^{<\alpha _{2}},\quad \check{\delta}_{f}^{<\alpha
_{2}},\quad \rho _{f}^{<\alpha _{2}},\quad \alpha _{f}^{<\alpha
_{2}}
\end{equation*}
тождественно совпадают соответственно с функциями
\begin{equation*}
\delta S_{f}^{<\alpha _{1}},\quad \check{\delta}_{f}^{<\alpha
_{1}},\quad \rho _{f}^{<\alpha _{1}},\quad \alpha _{f}^{<\alpha
_{1}}.
\end{equation*}
\hspace*{\fill} $\dashv$
\end{lemma}

Следующая лемма и её доказательство демонстрируют идею, которая
будет применяться далее в различных типичных ситуациях:

\begin{lemma}
\label{7.5.} \emph{(О диссеминаторе)} \newline \emph{1)}\quad
Пусть
\begin{itemize}
\item[(i)]  \ $]\tau _{1},\tau _{2}[ \; \subseteq dom\bigl( \delta
S_{f}^{<\alpha _{1}}\bigr) ~,\quad \gamma _{\tau _{2}}\in
SIN_{n}^{<\alpha _{1}}$; \medskip

\item[(ii)]  \ $\tau _{3}\in dom\bigl( \delta S_{f}^{<\alpha _{1}}\bigr)
~,~~\tau _{2} \leq \tau _{3}$; \medskip

\item[(iii)]  \ $\check{\delta}_{\tau _{3}}^{<\alpha _{1}}<\gamma _{\tau
_{2}}^{<\alpha _{1}}$~. \medskip
\end{itemize}
Тогда
\begin{equation*}
\check{\delta}_{\tau _{3}}^{<\alpha _{1}}\leq \gamma _{\tau
_{1}}^{<\alpha _{1}}~.
\end{equation*}

\noindent \emph{2)}\quad Пусть \ $\delta$-матрица \ $S$ \ на её
носителе \ $\alpha$ \ допустима для \ $\gamma_{\tau}^{<\alpha_1}$
\ вместе со своим диссеминатором \ $\delta$ \ и базой \ $\rho$ \
ниже \ $\alpha_1$, \ тогда:

\begin{itemize}
\item[(i)]  \ \ $\{\tau ^{\prime }: \ \delta
< \gamma _{\tau ^{\prime }}^{<\alpha _{1}}\leq \gamma _{\tau
}^{<\alpha _{1}}\}\subseteq dom\bigl( \delta S_{f}^{<\alpha
_{1}}\bigr);$ \medskip

\item[(ii)] эта матрица \ $S$ \ вместе с теми же \
$\delta$, \ $\rho$ \ обладает минимальным носителем \
$\alpha^{\prime} \in \; ]\gamma _{\tau }^{<\alpha _{1}},\gamma
_{\tau +1}^{<\alpha _{1}}[$~, допустимым ниже \ $\alpha_1$. \
\end{itemize}
\end{lemma}

\noindent \textit{Доказательство.} \ 1) Верхние индексы \ $<
\alpha_{1}$, $\vartriangleleft \alpha_{1}$ \ будут опускаться.
Рассмотрим матрицу \ $S^3 = \delta S_{\tau_3}$ \ и \
$\check{\delta}^3 = \check{\delta}_{\tau_3}$, \ $\rho^3 =
\rho_{\tau_3}$. \ Предположим, что 1) неверно, тогда по $(iii)$
\begin{equation*}
\gamma_{\tau_{1}} < \check{\delta}^3< \gamma_{\tau_{2}} \quad
\mbox{\it и } \quad \check{\delta}^3 =  \gamma_{\tau_{4}}
\end{equation*}
для некоторого \ $\tau_{4} \in \; ] \tau_{1},\tau_{2} [ $. \
Расмотрим ситуацию ниже \ $\alpha^{3} =
\alpha_{\tau_{3}}^{\Downarrow}$,\ стоя на \ $\alpha^{3} =
\alpha_{\tau_{3}}^{\Downarrow}$\;. \ Из $(i)$ и леммы~7.4 следует,
что
\begin{equation*}
\delta S_{f}^{<\alpha^{3}} \equiv \delta S_{f} \mbox{\it \ \ на \
\ } ]\tau_1, \tau_2[
\end{equation*}
и матрица \ $S^4 = \delta S_{\tau_{4}}^{<\alpha^{3}}=\delta
S_{\tau_{4}} $ \ на носителе \
$\alpha_{\tau_{4}}^{<\alpha^{3}}=\alpha_{\tau_{4}} $ \ имеет
диссеминатор
\begin{equation*}
\check{\delta}^{4} = \check{\delta}_{\tau_{4}}^{<\alpha^{3}} =
\check{\delta} _{\tau_{4}} < \gamma_{\tau_{4}} = \check{\delta}^3
\mbox{\it \ \ с базой } \ \rho^4 = \rho_{\tau_4}^{<\alpha^3}.
\end{equation*}

\noindent Теперь аргумент из доказательства леммы
6.6~\cite{Kiselev11} нужно  повторить следующим образом. Из \ $
\check{\delta}^{4} < \check{\delta}^{3} $ \ вытекает, что
\begin{equation*}
\rho^4<\rho^3 \mbox{\it \ и поэтому \ } \check{\delta}^{4} \notin
SIN_{n+1}^{< \alpha^{3}} [<\rho_{\tau_{3}}],
\end{equation*}
тогда по лемме 6.6~\cite{Kiselev11} (для \ $m=n+1$) \ существует
некоторое \ \mbox{$\Sigma_{n}$-предложение} \ $\varphi(\alpha,
\overrightarrow{a}) $ \ с кортежем \ $\overrightarrow{a} $ \
констант \ $< \rho_{\tau_{3}} $ \ и некоторый ординал \
$\alpha_{0} \in [ \check{\delta}^4, \alpha^3 [ $ \ такие, что
\begin{equation*}
\forall \alpha < \alpha_{0} \ \ \varphi^{\triangleleft \alpha^{3}}(\alpha,
\overrightarrow{a})\wedge \neg  \varphi^{\triangleleft
\alpha^{3}}(\alpha_{0}, \overrightarrow{a})~.
\end{equation*}
Диссеминатор \ $\check{\delta}^{3} $ \ ограничивает предложение \
$ \exists \alpha \neg \varphi(\alpha, \overrightarrow{a}) $ \ ниже
\ $\alpha^{3} $, \ поэтому \ $\alpha_{0} \in \; ]
\check{\delta}^{4}, \check{\delta} ^{3} [ $. \ Очевидно, \
$\Pi_{n+1}$-предложение
\begin{multline*}
    \forall \alpha, \gamma \Bigl( \neg \varphi(\alpha,
    \overrightarrow{a}) \longrightarrow \exists \gamma_1 \bigl (
    \gamma < \gamma_1 \wedge SIN_{n-1}(\gamma_1) \wedge         \
\\
    \wedge \exists \delta < \alpha \ \ \exists \alpha^{\prime} \
    \mathbf{K} (\delta, \gamma_{1}, \alpha^{\prime}, \rho^4, S^4) \bigr)
    \Bigr)
\end{multline*}

\noindent выполняется ниже \ $\check{\delta}^{3} $ \ и, значит, \
$\check{ \delta}^{3} $ \ продолжает его до \ $\alpha^{3} $, \
потому что
\[
    S^4 \vartriangleleft \rho_4 < \rho_3.
\]
Следовательно, для каждого \ $\gamma_{\tau}^{< \alpha^3}
> \check{\delta}^3$ \ появляется \ $\delta$-матрица \ $S^4$ \
допустимая на некотором носителе
\[
    \alpha \in [ \gamma_{\tau}^{< \alpha^3}, \alpha^3 [
    \mbox{\it \ \ для  \ \ } \gamma_{\tau}^{< \alpha^3}
\]
вместе со своим диссеминатором \ $\check{\delta}^4 <
\check{\delta} ^3 $ \ и его базой \ $\rho^4$.
\\
Отсюда следует, что ниже \ $\alpha^3$ \ определяется
\textit{минимальный} кардинал \ $\check{\delta}^m$ \ и
\textit{минимальная} база \ $ \rho^m$ с этим свойством, то есть
выполняющие следующее утверждение \textit{ниже} \ $\alpha^3$:
\[
    \exists \gamma^m \forall \gamma > \gamma^m \bigl(
    SIN_{n-1}(\gamma) \rightarrow \exists \alpha^{\prime},
    S ~\mathbf{K}(\check{\delta}^m,
    \gamma, \alpha^{\prime}, \rho^m, S) \bigr),
\]
\begin{sloppypar}
\noindent то есть существует \ $\gamma^m < \alpha^3$ \ такой, что
для каждого \ \mbox{$\gamma_{\tau}^{< \alpha^3} \in \; ] \gamma^m,
\alpha^3 [$} \ существует некоторая \ $\delta$-матрица \ $ S$, \
допустимая на некотором носителе \ \mbox{$\alpha \in [
\gamma_{\tau}^{< \alpha^3}, \alpha^3 [$} \ для \ $\gamma_{\tau}^{<
\alpha^3}$ \ ниже \ $\alpha^3$ \ вместе со своим производящим
диссеминатором \ $\check{\delta}^m < \gamma^m$ \ с базой\
$\rho^m$.
\end{sloppypar}

\noindent Очевидно, \ $\check{\delta}^m < \check{\delta}^3$. \ Так
как минимальное значение \ $ \rho^m$ \ определяется ниже \
$\alpha^3$, \ то по лемме 4.6~\cite{Kiselev11} о спектральном типе
это влечёт
\begin{equation*}
\rho^m < OT(\delta S_{\tau_3}) \leq Od(\delta S_{\tau_3}).
\end{equation*}
Но тогда это влечёт противоречие: существует \ $\delta$-матрица \
$ S^m$ \ на некотором носителе \ $\alpha^m \in \; ]
\gamma_{\tau_3}, \alpha^3 [ $, \ допустимая для \
$\gamma_{\tau_3}$ \ вместе с диссеминатором \ $\check{\delta}^m <
\gamma_{\tau_3}$ \ с базой \ $\rho^m$ \ и по условию \
$\mathbf{K}^{\exists}_{n+1}$

\begin{equation*}
S^m \vartriangleleft \rho^m < OT(\delta S_{\tau_3})  \leq
Od(\delta S_{\tau_3}),
\end{equation*}
\vspace{0pt}

\noindent хотя \ $\delta S_{\tau_3}$ \ является  \
$\underline{\lessdot}$-минимальной по определению~7.2.
\\
\quad \\
Утвеждение 2)~$(i)$ повторяет лемму 5.17~2)~$(i)$~\cite{Kiselev11}
и следует из определения~7.2 непосредственно; а утверждение
2)~$(ii)$ можно легко установить посредством аргумента
доказательства леммы 5.17~2)~$(ii)$~\cite{Kiselev11} \ для матрицы
\ $S$ \ вместо \ $S_{\chi \tau}^{< \alpha_{1}} $ \ и для формулы \
$\mathbf{K}$ \ вместо \ $\sigma$; \ мы вернёмся к этому аргументу
в $\S$~8 в более важном случае. \hspace*{\fill} $\dashv$
\\
\quad \\

\noindent Нерелятивизированная функция \ $\delta S_{f} $ \
действительно определена на финальном подинтервале недостижимого
кардинала \ $k$, \ как это показывает

\begin{lemma}
\label{7.6.} \emph{(Об определённости \ $\delta $-функции)}
\newline \hspace*{1em} Существует ординал \ $\delta <k$ \
такой, что \ $\delta S_{f}$ \ определена на множестве
\begin{equation*}
T = \{\tau :\delta <\gamma _{\tau } < k \}.
\end{equation*}
Минимальный из таких ординалов \ $\delta $ \ обозначается через \
$\delta ^{\ast }$, \ следующий за ним в \ $SIN_{n}$ \ кардинал --
через \ $\delta ^{\ast 1}$; \ также вводятся следующие
соответствующие ординалы:
\[
    \tau_1^{\ast }=\tau (\delta ^{\ast }),\quad \tau ^{\ast 1}=\tau
    (\delta ^{\ast 1}),
\]
\[
    \quad \mbox{\it \ \ так что \ \ }
    \delta^{\ast}=\gamma_{\tau_1^{\ast}}, \quad \delta^{\ast
    1}=\gamma_{\tau^{\ast 1}},
\]
\[
    \quad \mbox{\it \ \ и \ \ }
    \alpha ^{\ast
    1}=\alpha _{\tau ^{\ast 1}}^{\Downarrow },\quad \rho ^{\ast
    1}=\rho _{\tau _{{}}^{\ast 1}}.
\]
\end{lemma}

\noindent \textit{Доказательство} \ состоит в непосредственном
применении леммы~6.14~\cite{Kiselev11} для $\alpha_{1} = k, \
m=n+1, \ \chi = \chi^{\ast} $. \hspace*{\fill} $\dashv$

\begin{lemma}
\label{7.7.}
\begin{equation*}
\delta ^{\ast }\in SIN_{n}\cap SIN_{n+1}^{<\alpha ^{\ast 1}}\left[
<\rho ^{\ast 1}\right].
\end{equation*}
\end{lemma}

\noindent \textit{Доказательство.} \ Рассмотрим диссеминатор \
$\check{\delta}_{\tau^{\ast 1 }} $ \ с базой \ $\rho^{\ast 1} $ \
матрицы \ $\delta S_{\tau^{\ast 1 }} $ \ на носителе \
$\alpha_{\tau^{\ast 1 }} $. \ Так как
\begin{equation*}
\delta^{\ast 1} \in SIN_{n}, \quad \check{\delta}_{\tau^{\ast 1 }}
< \delta^{\ast 1}
\end{equation*}
и
\begin{equation*}
\check{\delta}_{\tau^{\ast 1 }} \in SIN_{n}^{ <
\alpha^{\ast1}}\cap SIN_{n+1}^{ < \alpha^{\ast 1}}[< \rho^{\ast
1}],
\end{equation*}
\vspace{0pt}

\noindent то лемма 3.8~\cite{Kiselev11} влечёт \
$\check{\delta}_{\tau^{\ast 1 }} \in SIN_{n}$ \ и по
леммам~7.5~2), 7.6 \ $\check{\delta}_{\tau^{\ast 1 }} =
\delta^{\ast} $. \hspace*{\fill}$\dashv$
\\

\begin{definition}
\label{7.8.} \hfill {}

1. Функция \ $\delta S_{\tau}^{<\alpha_1}$ \ называется монотонной
на интервале \ $[\tau_1, \tau_2[$, \ или на соответствующем
интервале \ $[\gamma_{\tau_1}^{<\alpha_1},
\gamma_{\tau_2}^{<\alpha_1}[$ \ ниже \ $\alpha_1$, \ если \
$\tau_1+1<\tau_2$, \ $]\tau_1, \tau_2[\; \subseteq dom(\alpha
S_f^{<\alpha_1})$ \ и
\[
    \forall \tau^{\prime}, \tau^{\prime\prime}(\tau_1 < \tau^{\prime}
    < \tau^{\prime\prime} < \tau_2 \longrightarrow \delta
    S_{\tau^{\prime}}^{<\alpha_1} \; \underline{\lessdot } \; \delta
    S_{\tau^{\prime\prime}}^{<\alpha_1}).
\]

2. Соответственно этому функция \ $\delta S_f$ \ называется
(тотально) монотонной, если для \ $\tau_1^{\ast} =
\tau(\delta^{\ast})$:
\[
    \forall \tau^{\prime}, \tau^{\prime\prime}(\tau_1^{\ast} < \tau^{\prime}
    < \tau^{\prime\prime} < k \longrightarrow \delta
    S_{\tau^{\prime}} \; \underline{\lessdot } \; \delta S_{\tau^{\prime\prime}} ).
\]
\hspace*{\fill}$\dashv$
\end{definition}

Некоторые простые фрагменты монотонности матричной функции \
$\delta S_f$ \ следуют из определения~7.2 и леммы 7.5~2)~$(ii)$
немедленно:

\begin{lemma}
\label{7.9.} \emph{(О монотонности \ $\delta $-функции)}
\newline \hspace*{1em} Пусть
\[
    \tau_1 < \tau_2 \mbox{\it \ и \ }
    \check{\delta}_{\tau_2}^{<\alpha_1} <
    \gamma_{\tau_1}^{<\alpha_1}.
\]

Тогда
\[
    \delta S_{\tau_1}^{<\alpha_1} \; \underline{\lessdot } \; \delta
    S_{\tau_2}^{<\alpha_1}.
\]
\hspace*{\fill}$\dashv$
\end{lemma}

Обсудим возникшую ситуацию.
\\
Мы видели выше, что простейшая матричная функция \ $S_{f}$ \
является \ $ \underline{ \lessdot }$-монотонной, но для каждого \
$\tau>\tau ^{\ast }$ \ кардинал предскачка \ $\alpha ^{\Downarrow
}$ \ матрицы \ $S_{\tau} $ \ на её соответствующем носителе \
$\alpha \in \; ]\gamma_{\tau}, \gamma_{\tau+1} [ $ \ не сохраняет
субнедостижимость уровней  \ $\geq n$ \ кардиналов \ $\leq \gamma
_{\tau }$, \ и некоторые другие важные свойства нижних уровней
универсума также нарушаются при релятивизации к \ $ \alpha
^{\Downarrow }$ \ (см. леммы 5.17, 5.18 и их обсуждение в конце
$\S$5~\cite{Kiselev11}).

Чтобы преодолеть это препятствие мы снабдили значения этой
функции, то есть матрицы \ $S_{\tau }$, \ диссеминаторами уровня \
$n+1$ \ и потребовали сохранения субнедостижимости уровня \ $n$ \
для кардиналов \ $\leq \gamma _{\tau }$, \ то есть мы перешли к \
$\delta $-функции \ $\delta S_{f}$.
\\

\noindent Но теперь это вызывает новое осложнение: теперь с
помощью лемм \ref{7.3.}-\ref{7.7.} нетрудно видеть, что после этой
модификации \ $\delta $-функция \ $\delta S_{f}$ лишается своего
свойства
 \ $\underline{\lessdot }$-монотонности на
$[\tau_1^{\ast}, k[\;$, \ и именно благодаря тому, что во многих
случаях кардиналы предскачков\ $\alpha^{\Downarrow } $ \ носителей
\ $\delta$-матриц\ $\alpha$, \ наоборот, \textit{вызывают
возникновение субнедостижимости уровня} \ $n$ \ некоторых
кардиналов \ $\leq \gamma _{\tau }$, \ которые становятся
субнедостижимыми (релятивизированно к \ $\alpha^{\Downarrow} $), \
не являясь таковыми до этого (Киселев~\cite{Kiselev4}).
\newline
Выход из этой новой затруднительной ситуации обнаруживается с
помощью следующего явления, которое в дальнейшем предоставит
решение всей проблемы:
\\
А именно, можно видеть, что \ $\delta$-матрица \ $\delta S_{\tau
_{0}}$, \ нарушающая \ $ \underline{\lessdot}$-монотонность на \
$[\tau_1^{\ast}, k[$ \ впервые, то есть для
\begin{equation*}
\tau _{0}=\sup \{ \tau :\delta S_{f}\quad \mbox{\it является}\quad
\underline{\lessdot} \mbox{\it -монотонной\ на} \quad \left] \tau
_{1}^{\ast},\tau \right[ \; \} ~,
\end{equation*}
помещается на некотором носителе \ $ \alpha_{\tau_0} \in \; ]
\delta ^{\ast },\delta ^{\ast 1} [$ \ и \ $\delta S_{\tau
_{0}}\vartriangleleft \rho ^{\ast 1}$ \ по лемме
3.2~\cite{Kiselev11}.
\\
Поэтому из лемм 7.7,~6.3~\cite{Kiselev11} (для \ $m=n+1$, \
$\alpha _{1}=\alpha ^{\ast 1}$) следует, что диссеминатор \
$\check{\delta}_{\tau ^{\ast 1}}$ \ переносит в точности ту же
самую ситуацию, но ниже \ $\alpha ^{0}=\alpha _{\tau
_{0}}^{\Downarrow }$, \ то есть:
\\
\quad \\
\textit{класс \ $SIN_{n}^{<\alpha ^{0}}$ \ содержит некоторые
кардиналы \ $ \gamma _{\tau _{1}}^{<\alpha ^{0}}<\gamma _{\tau
_{2}}^{<\alpha ^{0}}$ такие, что
\begin{equation*}
\left] \tau_{1}, \tau _{2}\right[ \subseteq  dom ( \delta
S_{f}^{<\alpha ^{0}} ),
\end{equation*}
и снова та же самая матрица
\begin{equation*}
\delta S_{\tau ^{0}}=\delta  S_{\tau _{0}^{\prime }}^{<\alpha ^{0}}
\end{equation*}
нарушает монотонность \ $\delta S_{f}^{<\alpha ^{0}}$ \ на  \
$\left] \tau _{1},\tau _{2}\right[ $ \ впервые для некоторого
ординала \ $\tau_{0}^{\prime} \in ]\tau_{1},\tau_{2} [\;$,} \ {\sl
но ниже \ $\alpha^0$.}
\\

Таким образом, здесь мы приходим к третьему и последнему
приближению к главной идее:
\newline \quad \newline \textit {Следующие требования необходимо
наложить на
 \ $\delta$-матрицы:
\newline
1) они должны обладать свойством ``самоисключения" в подобных
ситуациях нарушения монотонности (будем называть его ``свойством
автоэкзорцизивности''); матрицы с этим свойством ( матрицы {\sl
``единичной характеристики''}) должны иметь приоритет над другими
матрицами ( матрицами {\sl ``нулевой характеристики''}
соответственно) в ходе определения матричной функции;
\newline 2) ещё одно требование должно быть наложено на матрицы {\sl
нулевой} характеристики, препятствуя их образованию: их базы
данных должны существенно возрастать, когда предшествующая часть
матричной функции, которая уже определена, содержит нарушение
монотонности, с целью скорректировать это нарушение --
использование матрицы {\sl нулевой} характеристики;
\\
на этих основаниях матричная \ $\delta$-функция должна получить
несовместные свойства монотонности и немонотонности одновременно.
}
\\

Очевидно, все эти соображения требуют рекурсивного определения
матричной функции, определяя её значения в зависимости от свойств
её предшествующих значений.
\\
Мы начинаем осуществлять эту идею в следующем параграфе.

\newpage

\section{Матричные \ $\protect\alpha$\,-функции}
\setcounter{equation}{0}

\hspace*{1em} Для предстоящего рекурсивного определения необходимо
усложнить предыдущую формулу \ $\mathbf{K}_{n+1}^{\exists }$ \
(определение \ref{7.1.}). Но предварительно следует ввести
некоторые субформулы для большей ясности конструкции этой формулы,
где переменная \ $X_1$ \ играет роль матричной функции \ $\alpha
S_{f}^{< \alpha}$,\ а переменная \ $X_2$ \ играет роль
характеристической функции \ $a_f^{< \alpha}$, \ которые будут обе
определены ниже \ $\alpha$; \ последняя функция будет приписывать
соответствующие характеристики (единичные или нулевые)
редуцированным матрицам, служащим  значениями функции \ $\alpha
S_f^{< \alpha}$; \ эти характеристики матриц на их носителях будут
принимать значения единичное \ $a=1$ \ или нулевое \ $a=0$ \
соответственно принципу, обрисованному выше.
\\
В ходе их введения эти формулы будут сопровождаться пояснениями их
смысла, а после результирующего определения~8.2 мы опишем как оно
действует в целом, без вхождения в детали.
\newline

Все эти формулы были использованы в предыдущих работах автора
\cite{Kiselev1,Kiselev2,Kiselev3,Kiselev4,Kiselev5,Kiselev6,Kiselev7,Kiselev8},
но некоторые из них были рассеяны в тексте в их определённых
формах (некоторые в неформализованном, некоторые другие в
семантическом виде), и здесь они излагаются систематическим
образом; также применяется несколько удобных переобозначений.

В этих формулах используются различные кардиналы из классов \
$SIN_{n}$, $SIN_{n-1}$, $SIN_{n-2}$ \ субнедостижимости.
Необходимо иметь ввиду, что после\ $<$- \ или \
$\vartriangleleft$-ограничения этих формул некоторым кардиналом \
$\alpha$ \ (см. определение~2.3~\cite{Kiselev11}) возникают
субнедостижимые классы того же уровня, но релятивизированные к
этому \ $\alpha$; \ например, \ $SIN_{n}$-субнедостижимость \
преобразуется в \ $SIN_{n}^{<\alpha}$-субнедостижимость, \ но уже
ниже \ $\alpha$; \ поэтому все формулы после этого повествуют о
соответствующей ситуации ниже \ $\alpha$.
\\
Такие трансформации обеспечиваются определениями и леммами
3.3-3.8~\cite{Kiselev11}.
\\

\begin{definition}
\label{8.1.} \hfill {} \newline \hspace*{1em} Вводятся следующие
вспомогательные формулы:
\\ \quad
\\
\emph{I.} Интервалы определённости матричной функции:
\\ \quad
\\
\emph{1.0}\quad $A_0(\chi ,\tau _{1},\tau _{2},X_1)$:
\vspace{-6pt}
\begin{multline*}
    \tau_{1}+1<\tau _{2} ~\wedge~ \big(X_1\mbox{\it \ is a function on
    }\left] \tau _{1},\tau _{2}\right[ \; \big)\wedge
\\
    \wedge \tau _{1} =\min \big\{\tau :\left] \tau ,\tau _{2}\right[
    \subseteq dom(X_1)\big\}\wedge
\\
    \wedge \exists \gamma^{1} \big( \chi \leq \gamma ^{1}=\gamma
    _{\tau _{1}} \wedge SIN_{n}(\gamma ^{1}) \big);
\end{multline*}

\noindent очевидно, эта формула означает, что интервал \ $[\tau
_{1},\tau _{2}[$ \ занимает особое место по отношению к матричной
функции \ $X_1$:\ эта функция определена на \ $]\tau _{1},\tau
_{2}[$ \ и \ $\tau_{1}$ \ есть минимальный ординал с этим
свойством, сверх того кардинал \ $ \gamma _{\tau _{1}}$ \
принадлежит \ $SIN_{n}$;
\\
благодаря этой минимальности \ $X_1$ \ никогда не определена для
этого ординала \ $\tau_1$.
\quad \medskip %

\noindent \emph{1.1}\quad $A_{1}(\chi ,\tau _{1},\tau _{2},X_1)$:
\[
    A_0(\chi ,\tau _{1},\tau _{2},X_1)\wedge \exists \gamma^2
    \big( \gamma^2 = \gamma_{\tau_2} \wedge SIN_n (\gamma^2)\big);
\]

\noindent такой интервал \ $[\tau_1, \tau_2[$, \ а также интервал
\ $[\gamma_{\tau_1}, \gamma_{\tau_2}[$, \ будут называться
интервалами определённости матричной функции \ $X_1$ \
 {\sl максимальными влево} (в \
$dom(X_1)$), \ максимальными в том смысле, что не существует
интервала \ $]\tau^{\prime}, \tau_2[$ \ в \ $dom(X_1)$ \ с меньшим
левым концом \ $\tau^{\prime} < \tau_1$; \ сверх того по-прежнему
требуется \ $\gamma_{\tau_1} \in SIN_n$, \ и ещё $\gamma_{\tau_2}
\in SIN_n$.
\\

\noindent \emph{1.2}\quad $A_{1.1}^M(\chi ,\tau _{1},\tau
_{2},X_1)$: \vspace{-6pt}
\[
    \quad A_{1}(\chi, \tau_{1}, \tau_{2}, X_1)\wedge \tau _{2}=
    \sup \big\{ \tau : A_{1}(\chi, \tau_{1}, \tau_{2}, X_1) \big\};
\]

\begin{sloppypar}
\noindent здесь интервал \ $]\tau _{1},\tau _{2}[$ \ ( включённый
в \ $dom(X_1)$) \ -- максимальный в том смысле, что он не
включается ни в какой другой интервал \ \mbox{$]\tau _{1}^{\prime
},\tau _{2}^{\prime }[\;\subseteq dom(X_1)$} \ такой, что \ $
\gamma _{\tau_{2}^{\prime}}\in SIN_{n}$; \ помимо этого
по-прежнему требуется, чтобы \ $\gamma _{\tau _{1}} \in SIN_{n}$,
\  \ $\gamma _{\tau _{2}} \in SIN_{n}$; \ по этой причине такой
интервал \ $[\tau_1, \tau_2[$ \ и соответствующий интервал
 \ $[\gamma_{\tau_1}, \gamma_{\tau_2}[$
\ будут называться {\sl максимальными} интервалами определённости
матричной  функции \ $X_1$ .\
\\
\quad \medskip %
\end{sloppypar}

\noindent \emph{1.3}\quad $A_{1.2}(\tau _{1},\tau _{2},\eta )$:
\vspace{-6pt}
\begin{multline*}
    \exists \gamma ^{1},\gamma ^{2} \Big(\gamma ^{1}=
    \gamma _{\tau_{1}}\wedge \gamma ^{2}=\gamma _{\tau _{2}}
    \wedge
\\
    \wedge \eta =OT\big(\big\{\gamma :\gamma ^{1}<
    \gamma <\gamma^{2}\wedge SIN_{n}(\gamma )\big\}\big)\Big);
\end{multline*}

\noindent здесь, напомним, \ $OT$ \ обозначает порядковый тип
указанного множества, поэтому мы будем называть такой ординал \
$\eta $ \ типом интервала \ $[\tau _{1},\tau _{2}[$ \ и также
типом соответствующего интервала \ $[\gamma _{\tau _{1}},\gamma
_{\tau _{2}}[$~.
\newline
\quad \medskip %

\noindent \emph{1.4}\quad $A_{2}(\chi ,\tau _{1},\tau _{2},\tau
_{3},X_1)$: \vspace{-6pt}
\begin{multline*}
    \ A_{1}(\chi ,\tau _{1},\tau _{3},X_1)\wedge
    \tau _{1}+1<\tau_{2}<\tau _{3} \wedge \tau_{2} =
\\
\qquad = \sup \Big\{\tau <\tau _{3}:\forall \tau ^{\prime },
    \tau ^{\prime\prime }\big(\tau _{1}<\tau ^{\prime }
    <\tau ^{\prime \prime }<\tau \longrightarrow
    X(\tau ^{\prime })\underline{\lessdot } X
    (\tau^{\prime \prime })\big)\Big\};
\end{multline*}

\noindent здесь \ $\tau _{2}$ \ -- это минимальный индекс, на
котором нарушается \ $ \underline{\lessdot }$-монотонность
матричной функции \ $X_1$ \ на интервале \ $\left] \tau _{1},\tau
_{3}\right[ $ .\
\newline
\quad \medskip %

\noindent \emph{1.5}\quad $A_{3}(\chi ,\tau _{1},\tau _{1}^{\prime
},\tau _{2},\tau _{3},X_1,X_2)$: \newline

\noindent $\ A_{2}(\chi ,\tau _{1},\tau _{2},\tau _{3},X_1 )\wedge
\tau _{1}<\tau _{1}^{\prime }<\tau _{2}\wedge \big(X_2\mbox{\it \
это функция } ]\tau _{1},\tau _{3}[\ \big) \wedge $
\newline

\noindent \hfill $\ \wedge \tau _{1}^{\prime }=\min \big\{\tau \in
\;]\tau _{1},\tau _{2}[\;:X_1(\tau )\gtrdot X_1(\tau _{2})\wedge
X_2(\tau )=1\big\}$;
\newline

\noindent здесь указывается, что \ $\underline{\lessdot }$
-монотонность матричной функции \ $X_1$ \ на \ $\left] \tau
_{1},\tau _{3}\right[ $ \ впервые нарушается на индексе \ $\tau
_{2}$ \ и именно из-за матрицы единичной характеристики \
$X_1(\tau _{1}^{\prime })\gtrdot X_1(\tau _{2})$ \ для \ $\tau
_{1}^{\prime }\in \left] \tau _{1},\tau _{2}\right[ $ .\
\newline
\quad \medskip %

\noindent \emph{1.6.a}\quad $A_{4}^b(\chi ,\tau _{1},\tau
_{1}^{\prime },\tau _{2},\tau _{3},\eta ,X_1,X_2)$:
\newline

\qquad \qquad $A_{3}(\chi ,\tau _{1},\tau _{1}^{\prime },\tau
_{2},\tau _{3},X_1,X_2) \wedge A_{1.2}(\tau _{1},\tau _{3},\eta
);$
\\
\quad \\
\noindent \emph{1.6.a(i)}\quad $A_{4}^b(\chi ,\tau _{1}, \tau
_{2}, \eta, X_1, X_2)$:
\[
    \qquad \exists \tau_1^{\prime}, \tau_2^{\prime} \le \tau_2 \
    A_{4}^b(\chi, \tau_{1}, \tau_{1}^{\prime}, \tau_{2}^{\prime},
    \tau_{2}, \eta, X_1, X_2);
\]

\noindent \emph{1.6.b}\quad $A_{4}^{M b}(\chi ,\tau _{1},\tau
_{1}^{\prime },\tau _{2},\tau _{3},\eta ,X_1,X_2)$:
\[
    A_{4}^b(\chi, \tau_{1},\tau _{1}^{\prime },\tau _{2},\tau
    _{3},\eta,X_1,X_2) \wedge
    A_{1.1}^M(\chi ,\tau _{1},\tau _{3},X_1);
\]

\noindent \emph{1.6.b(i)}\quad $A_{4}^{M b}(\chi ,\tau _{1}, \tau
_{2}, \eta, X_1, X_2)$:
\[
    \qquad \exists \tau_1^{\prime}, \tau_2^{\prime} \le \tau_2 \
    A_{4}^{M b}(\chi, \tau_{1}, \tau_{1}^{\prime}, \tau_{2}^{\prime},
    \tau_{2}, \eta, X_1, X_2);
\]
в дальнейшем всякий интервал \ $[\tau _{1},\tau _{3}[$, \
обладающий свойством \ $A_{4}^b$ \ для некоторых \ $\tau
_{1}^{\prime }$, \ $\tau _{2}$, \ $\eta $ \ и соответствующий
интервал \ $[\gamma _{\tau _{1}},\gamma _{\tau _{3}}[$ \ будут
называться блоками типа \ $\eta$, \ а если они вдобавок выполняют
 \ $A_{1.1}^M(\chi ,\tau _{1},\tau
_{3},X_1)$ \ -- \ то {\sl максимальными блоками} этого типа.
\\
Такие блоки будут рассматриваться далее как неприемлимые по
причине их фатального дефекта: нарушения монотонности матричной
функции. По этой причине мы наложим на такие блоки некоторые
обременительные требования с целью избежать их формирования в ходе
определения матричной функции (см. условие \ $\mathbf{K} ^{0}$ \
ниже).
\end{definition}

Приостановим ненадолго это определение \ref{8.1.} чтобы пояснить
смысл и направление развития его последующей части.

Формулы и понятия, введённые выше и последующие, будут
использоваться  в результирующем определении~\ref{8.2.}  в их
релятивизированных формах, то есть их индивидные переменные и
константы будут \ $<$- \ или \ $\vartriangleleft$-ограничиваться
некоторым соответствующим кардиналом \ $\alpha_1$. \ В подобных
случаях используются их настоящие формулировки, но с добавленной
ремаркой ``ниже \ $\alpha_1$''; \ соответственно их обозначения
снабжаются верхним индексом \ $<\alpha_1$ \ или \
$\vartriangleleft \alpha_1$.
\\
Таким образом,
\[
    A_1^{\vartriangleleft\alpha_1}(\chi, \tau_1, \tau_2, X_1)
\]
это формула:

\vspace{-6pt}
\begin{eqnarray*}
    \tau_{1}+1<\tau _{2} ~\wedge~ \big(X_1\mbox{\it \ это функция на
    }\left] \tau _{1},\tau _{2}\right[ \; \big)\wedge \qquad\qquad
\\
    \wedge \tau _{1} =\min \big\{\tau :\left] \tau ,\tau _{2}\right[
    \subseteq dom(X_1)\big\}\wedge \qquad\qquad\qquad\qquad
\\
    \quad \wedge \exists \gamma^{1}, \gamma^{2}
    \big( \chi \leq \gamma ^{1} \wedge
    \gamma ^{1}=\gamma _{\tau_{1}}^{<\alpha_1} \wedge
    \gamma^{2}=\gamma _{\tau _{2}}^{<\alpha_1} \wedge
\\
    \wedge SIN_{n}^{<\alpha_1}(\gamma ^{1})
    \wedge SIN_{n}^{<\alpha_1}(\gamma ^{2})     \big),
\end{eqnarray*}
\vspace{0pt}

\noindent которая означает, что \ $\left]\tau_1, \tau_2 \right[$ \
есть интервал из области определения функции \ $X_1$ \ с
минимальным левым концом \ $\tau_1$, \ и сверх того
соответствующие кардиналы \ $\gamma_{\tau_1}^{<\alpha_1}$,
$\gamma_{\tau_2}^{<\alpha_1}$ \ являются \
$SIN_n^{<\alpha_1}$-кардиналами -- и всё это ниже \ $\alpha_1$.

\noindent Соответственно этому,
\[
    A_4^{b \vartriangleleft\alpha_1}(\chi, \tau_1, \tau_1^{\prime},
    \tau_2, \tau_3, \eta, X_1, X_2)
\]
это формула:
\[
    A_{3}^{\vartriangleleft\alpha_1}(\chi ,\tau _{1},\tau_{1}^{\prime },
    \tau _{2},\tau _{3},X_1,X_2) \wedge
    A_{1.2}^{\vartriangleleft\alpha_1}(\tau_{1},\tau_{3},\eta )
\]
которая означает, что \ $[ \tau_1, \tau_3 [$ \ и \ $[
\gamma_{\tau_1}^{<\alpha_1}, \gamma_{\tau_3}^{<\alpha_1} [$ \ это
\textit{блоки} ниже \ $\alpha_1$ типа \ $\eta$, \ то есть интервал
 \ $]\tau_1, \tau_3[$ \ является максимальным
влево в \ $dom(X_1)$ \ и кардиналы \
$\gamma_{\tau_1}^{<\alpha_1}$, \ $\gamma_{\tau_3}^{<\alpha_1}$ \
оба содержатся в \ $SIN_n^{<\alpha_1}$, \  и \
$\underline{\lessdot}$-монотонность \ $X_1$ \ на \ $]\tau_1,
\tau_3[$ \ нарушается впервые на индексе \ $\tau_2 \in \; ]\tau_1,
\tau_3[$ \ и именно по причине матрицы \ $X_1(\tau_1^{\prime})
\gtrdot X_1(\tau_2)$ \ \textit{единичной} характеристики для
некоторого \ $\tau_1^{\prime} \in \; ]\tau_1, \tau_2[$ \ -- \ и
всё это ниже \ $\alpha_1$.

\noindent Нетрудно видеть, что все эти и последующие подобные
ограниченные формулы содержаться в классе \ $\Delta_1^1$ \ для
каждого \ $\alpha_1 > \chi$, \ $\alpha_1 < k$, \ $\alpha_1 \in
SIN_{n-2}$.

Но чтобы ввести последующие понятия наиболее прозрачным образом,
удобно предварительно прояснить принцип, регулирующий приписывание
характеристик матрицам на их носителях и взаимодействие этих
характеристик между собой, так как характеристическая функция
играет ведущую роль в рекурсивном определении~\ref{8.2.} матричной
функции ниже.
\\
Итак, матрица \ $S$ \ на её носителе \ $\alpha$
--- и сам этот носитель \ $\alpha$ \  --- будет получать
\textit{нулевую} характеристику \ $a=0$, \ если она участвует в
нарушении монотонности матричной функции в следующем смысле:
\\
существует интервал определённости матричной функции

\[
    [ \gamma_{\tau_1}^{<\alpha_{\chi}^{\Downarrow}},
    \gamma_{\tau_3}^{<\alpha_{\chi}^{\Downarrow}} [
\]
ниже кардинала предскачка \ $\alpha_{\chi}^{\Downarrow}$ \ после
 \ $\chi$ \ этого носителя \ $\alpha$, \ где
встречается \textit{та же самая матрица} \ $S$ \ как значение
матричной функции \ $X_1$, \ но ниже\
$\alpha_{\chi}^{\Downarrow}$:
\[
    X_1(\tau_2) = S
\]
для индекса \ $\tau_2 \in ]\tau_1, \tau_3[$, \ который является
минимальным, нарушающим монотонность функции  \ $X_1$ \ на  \ $\
]\tau_1, \tau_3[$ \ ниже \ $\alpha_{\chi}^{\Downarrow}$, \ то есть
когда выполняется
\[
    A_2^{\vartriangleleft \alpha_{\chi}^{\Downarrow}} ( \chi, \tau_1, \tau_2, \tau_3,
    X_1).
\]
И здесь наступает последнее уточнение этого понятия: вдобавок не
должно существовать допустимых матриц для \
$\gamma_{\tau_1}^{<\alpha_{\chi}^{\Downarrow}}$ \
 и все значения матричной функции \ $X_1$ \ на интервале \ $]\tau_1,
\tau_2]$ \   ниже \ $\alpha_{\chi}^{\Downarrow}$ \ должны быть
  \textit{единичной} характеристики:
\\
\quad \\
\hspace*{8em} $ \forall \tau ( \tau_1 < \tau \le \tau_2
\rightarrow X_2(\tau)=1)$. \label{c9}
\endnote{
\ стр. \pageref{c9}. \ Это последнее уточнение не является
необходимым и доказательство основной теоремы можно провести без
него (ценой некоторых незначительных усложнений), однако мы примем
его, чтобы несколько сократить предстоящие рассуждения.
\\
\quad \\
} %
\\
\quad \\
В противном случае \ $S$ \ на \ $\alpha$ \ и сам  \ $\alpha$ \
будут получать \textit{единичную} характеристику \ $a=1$.
\\
И по ходу того, как матричная функция будет получать своё
рекурсивное определение~\ref{8.2.}, матрицы единичной
характеристики буду систематически получать \emph{приоритет} над
матрицами нулевой характеристики --- чтобы избежать нарушение
монотонности этой функции.
\\
Естественно понимать понятие ``приоритет'' в следующе смысле:
когда определяется некоторое значение \ $X_1(\tau)$ \ матричной
функции \ $X_1$ \ и на такое значение представляются матрицы \
$S^0$, \ $S^1$ \ нулевой и единичной характеристики, тогда именно
матрица \ $S^1$ \ должна быть назначена значением \ $X_1(\tau)$.

Но не исключаются некоторые определённые случаи, когда матрицы
нулевой характеристики будут отвергаться  по некоторым другим
причинам, когда \textit{матрица нулевой характеристики} \ $S$ \ на
её носителе \ $\alpha$ \ будет запрещена к представлению на
значение матричной функции; в каждом таком случае мы будем
говорить, что \ $S$ \ на \ $\alpha$ \ {\em подавляется} .
\\
Мы используем термин `` подавление'', а не ``отсутствие
приоритета'', так как подобное подавление будет применяться только
в специальных особых случаях в зависимости от расположения этого
носителя \ $\alpha$.

Итак, мы переходим к описанию случаев, когда действует подобное
подавление нулевой характеристики.
\\
Для этото следует отметить, что приведённые в этом
определении~\ref{8.1.} формулы должны использоваться следующим
особым образом:
\\
До сих пор в этих формулах 1.0--1.6 b $(i)$ символы \ $X_1$, $X_2$
\ означали функции, определённые на ординалах.
\\
Но для рекурсивного определения~\ref{8.2.} матричной функции
необходимо использовать функции, определённые на  {\em парах}
ординалов. Поэтому введём для такой функции  \ $X$ \ другую
функцию
\[
    X[\alpha] = \{(\tau,\eta) : ((\alpha,\tau),\eta) \in X \},
\]
так что
\[
    X(\alpha,\tau) = X[\alpha](\tau)
\]
для каждой пары \ $(\alpha,\tau) \in dom(X)$.
\\
Соответственно этому указанные формулы в определении~\ref{8.2.} и
последующие формулы будут часто использоваться для \ $X_1$, $X_2$
\ как функций
\[
    X_1[\alpha^0], ~
    X_2[\alpha^0],
\]
где \ $\alpha^0$ \ это некоторый ординал.
\\
Теперь вернёмся к определению \ref{8.1.} чтобы сформулировать так
называемое ``условие подаления'', оно возникает в связи с
покрытиями кардиналов блоками специального вида и для этого
необходима следующая группа условий:
\\
\quad \\

\em

\noindent \emph{II.} Условия подавления
\\
\quad \\
\emph{2.1a.}\quad $A_{5.1}^{sc}(\chi, \gamma^m, \gamma, X_1,
X_2)$:
\[
    \gamma^m < \gamma \wedge \forall \gamma^{\prime} \in
    [\gamma^m, \gamma[ \ \exists \tau_1, \tau_2, \eta \ \big (
    \gamma^m \le \gamma_{\tau_1} \le \gamma^{\prime} <
    \gamma_{\tau_2} \le \gamma \wedge
\]
\[
    \qquad \qquad \qquad \qquad \qquad \qquad \qquad
    \wedge A_4^{Mb}(\chi, \tau_1, \tau_2, \eta, X_1, X_2) \big)
    \wedge
\]
\[
    \wedge \forall \gamma^{m \prime} < \gamma^m \neg \forall \gamma^{\prime\prime} \in
    [\gamma^{m \prime}, \gamma[ \ \exists \tau_1^{\prime}, \tau_2^{\prime},
    \eta^{\prime} \big ( \gamma^{m \prime} \le \gamma_{\tau_1^{\prime}} \le \gamma^{\prime\prime}
    < \gamma_{\tau_2^{\prime}} \le \gamma \wedge
\]
\[
    \qquad \qquad \qquad \qquad \qquad \qquad \qquad
    \wedge A_4^{Mb}(\chi, \tau_1^{\prime}, \tau_2^{\prime}, \eta^{\prime}, X_1, X_2)
    \big);
\]
здесь указывается, что интервал  \ $[\gamma^m, \gamma[$ \ с правым
концом \ $\gamma$ \ является объединением максимальных блоков и
что его левый конец \ $\gamma^m$ \ -- минимальный с этим
свойством; такое семейство блоков будет называться  {\sl покрытием
} кардинала \ $\gamma$; \ легко видеть, что при этом условии \
$\gamma^m$, $\gamma$ \ являются \ $SIN_n$-кардиналами;
\\
если опустить здесь правый конец \ $\gamma$, \ то получается
следующее условие:
\\
\quad \\
\emph{2.1b.}\quad $A_{5.1}^{sc}(\chi, \gamma^m, X_1, X_2)$:
\[
    \forall \gamma^{\prime} \ge \gamma^m \ \exists \tau_1, \tau_2, \eta \ \big (
    \gamma^m \le \gamma_{\tau_1} \le \gamma^{\prime} <
    \gamma_{\tau_2} \wedge \qquad \qquad \qquad
\]
\[
    \qquad \qquad \qquad \qquad \qquad \qquad \qquad
    \wedge A_4^{Mb}(\chi, \tau_1, \tau_2, \eta, X_1, X_2) \big)
    \wedge
\]
\[
    \wedge \forall \gamma^{m \prime} < \gamma^m \neg \forall \gamma^{\prime\prime}
    \ge \gamma^{m \prime} \ \exists \tau_1^{\prime}, \tau_2^{\prime},
    \eta^{\prime} \ \big ( \gamma^{m \prime} \le \gamma_{\tau_1^{\prime}} \le \gamma^{\prime\prime}
    < \gamma_{\tau_2^{\prime}} \wedge
\]
\[
    \qquad \qquad \qquad \qquad \qquad \qquad \qquad
    \wedge A_4^{Mb}(\chi, \tau_1^{\prime}, \tau_2^{\prime}, \eta^{\prime}, X_1, X_2)
    \big).
\]
Чтобы сформировать условия подавление продуктивным образом следует
наложить следующие специальные условия на ординалы \ $\gamma^m <
\gamma^{\ast} < \gamma$, $\eta^{\ast}$:
\\
\quad \\
\emph{2.2.}\quad $A_{5.2}^{sc}(\chi, \gamma^m, \gamma^{\ast},
\eta^{\ast}, X_1, X_2)$:
\[
    A_{5.1}^{sc}(\chi, \gamma^m, \gamma^{\ast}, X_1, X_2) \wedge
    \qquad \qquad \qquad \qquad \qquad \qquad \qquad \qquad \qquad
\]
\[
    \wedge \forall \tau_1, \tau_2, \eta \ \big ( \gamma_{\tau_1} <
    \gamma_{\tau_2} \le \gamma^{\ast} \wedge
    A_4^{Mb}(\chi, \tau_1, \tau_2, \eta, X_1, X_2) \rightarrow
    \eta < \eta^{\ast} \big) \wedge
\]
\[
    \wedge \forall \eta < \eta^{\ast} \ \exists \gamma^{\prime} <
    \gamma^{\ast} \ \forall \tau_1^{\prime}, \tau_2^{\prime}, \eta^{\prime} \ \big (
    \gamma^{\prime} < \gamma_{\tau_2^{\prime}} \le \gamma^{\ast} \wedge
\]
\[
    \qquad \qquad \qquad \qquad \qquad \qquad
    \wedge A_4^{Mb}(\chi, \tau_1^{\prime}, \tau_2^{\prime}, \eta^{\prime}, X_1, X_2) \rightarrow
    \eta < \eta^{\prime} \big);
\]
в подобном случае, когда \ $A_{5.2}^{sc}$ \ выполняется, мы будем
говорить, что типы покрытия кардинала \ $\gamma^{\ast}$ \
неубывают до \ $\eta^{\ast}$ \ {\sl существенно}; таким образом
ординал  \ $\eta^{\ast}$ \ предельный;
\\
\quad \\
\emph{2.3.}\quad $A_{5.3}^{sc}(\chi, \gamma^{\ast}, \gamma^1,
\gamma, \eta^{\ast}, X_1, X_2)$:
\[
    \exists \tau_1, \tau \ \Big (
    \gamma_{\tau_1} = \gamma^1 \wedge
    \gamma_{\tau} = \gamma \wedge A_4^b(\chi, \tau_1,
    \tau, \eta^{\ast}, X_1, X_2) \wedge
\]
\[
    \forall \tau_1^{\prime}, \tau_2^{\prime}, \eta^{\prime} \ \big (
    \gamma^{\ast} < \gamma_{\tau_2^{\prime}} \le \gamma^1 \wedge
    A_4^{Mb}(\chi, \tau_1^{\prime}, \tau_2^{\prime}, \eta^{\prime}, X_1, X_2) \rightarrow
\]
\[
    \qquad \qquad \qquad \qquad \qquad \qquad \qquad \qquad \qquad \qquad
    \rightarrow \eta^{\prime} = \eta^{\ast} \big) \Big);
\]
теперь эти три условия следует собрать вместе в следующее
\\
\quad \\
\noindent \emph{2.4.} \ Условие подавляющего покрытия
\\
\quad \\
\hspace*{1.5em} $A_{5.4}^{sc}(\chi, \gamma, \eta^{\ast}, X_1,
X_2)$:
\[
    \exists \gamma^m, \gamma^{\ast}, \gamma^1 \ \Big (
    \gamma^m < \gamma^{\ast} < \gamma^1 < \gamma
    \wedge \eta^{\ast} < \chi^+ \wedge \quad
\]
\[
    \quad
    A_{5.1}^{sc}(\chi, \gamma^m, \gamma^1, X_1, X_2) \wedge
    A_{5.2}^{sc}(\chi, \gamma^m, \gamma^{\ast}, \eta^{\ast}, X_1, X_2)
    \wedge
\]
\[
    \qquad\qquad\qquad\qquad\qquad\qquad
    \wedge A_{5.3}^{sc}(\chi, \gamma^{\ast}, \gamma^1, \gamma,
    \eta^{\ast}, X_1, X_2) \Big);
\]
назовём покрытие \ $\gamma$, \ обладающее этим свойством, {\sl
подавляющим покрытием} для \ $\gamma$ \ типа \ $\eta^{\ast}$;
\\
таким образом, эти три условия \ $A_{5.1}^{sc}-A_{5.3}^{sc}$ \
вместе означают, что покрытие кардинала \ $\gamma$ \ разделяется
на три части: его типы неубывают до ординала \ $\eta^{\ast} <
\chi^+$ \ существенно слева от \ $\gamma^{\ast}$, \ затем они
стабилизируются от \ $\gamma^{\ast}$ \ до \ $\gamma^1$, \ то есть
  интервал \ $[\gamma^{\ast}, \gamma^1[$ \ покрыт
максимальными блоками постоянного типа \ $\eta^{\ast}$, \ также
существует блок \ $[\gamma^1, \gamma[$ \ того же типа\
$\eta^{\ast} < \chi^+$; очевидно, эти условия определяют ординалы
\ $\gamma^m < \gamma^{\ast} < \gamma^1 < \gamma$, $\eta^{\ast}$ \
единственным образом через \ $\gamma$ \ (если они существуют);
\\
\quad \\
\emph{2.5.}\quad $A_{5.5}^{sc}(\chi, \gamma, \eta^{\ast}, \alpha,
X_1, X_2)$:
\[
    \forall \gamma^{\prime} \ \big( \gamma \le \gamma^{\prime} <
    \alpha \rightarrow \exists \tau_1^{\prime},
    \tau_2^{\prime}, \eta^{\prime} \big(
    \gamma_{\tau_1^{\prime}}^{<\alpha} \le \gamma^{\prime} <
    \gamma_{\tau_2^{\prime}}^{<\alpha} \wedge \qquad
\]
\[
    \qquad\qquad\qquad\qquad
    \wedge A_4^{Mb \vartriangleleft \alpha} (\chi, \tau_1^{\prime},
    \tau_2^{\prime}, \eta^{\prime}, X_1, X_2) \wedge
    \eta^{\prime} \ge \eta^{\ast}  \big ) \big );
\]
здесь указывается, что интервал \ $[\gamma, \alpha[$ \ покрыт
максимальными блоками ниже
 \ $\alpha$ \ типов \
$\eta^{\prime} \ge \eta^{\ast}$.
\\

\noindent Теперь следует интегрировать  все эти условия в единое
\\
\quad \\
\noindent \emph{2.6.} \ Результирующее условие подавления
\\
\quad \\
\hspace*{1.5em} $A_{5}^{S,0}(\chi, a, \gamma, \alpha, \rho, S,
X_1^0, X_2^0, X_1, X_2)$: \vspace{-6pt}
\begin{multline*}
    \qquad a=0 \wedge SIN_n(\gamma) \wedge \rho < \chi^+ \wedge
    \sigma(\chi, \alpha, S) \wedge
\\
    \wedge \exists \eta^{\ast}, \tau < \gamma \Big(
    \gamma = \gamma_{\tau} \wedge
    A_{5.4}^{sc}(\chi, \gamma, \eta^{\ast}, X_1^0, X_2^0) \wedge
\\
    \wedge \forall \tau^{\prime} \big( \tau < \tau^{\prime}
    \wedge SIN_n(\gamma_{\tau^{\prime}}) \rightarrow
\\
    \rightarrow \exists \alpha^{\prime}, S^{\prime} \big[
    \gamma_{\tau^{\prime}} < \alpha^{\prime} <
    \gamma_{\tau^{\prime}+1} \wedge SIN_n^{<\alpha_{\chi}^{\prime \Downarrow}}
    (\gamma_{\tau^{\prime}}) \wedge \sigma(\chi,
    \alpha^{\prime}, S^{\prime}) \wedge
\\
    \wedge A_{5.5}^{sc}(\chi, \gamma, \eta^{\ast},
    \alpha_{\chi}^{\prime \Downarrow},
    X_1[\alpha_{\chi}^{\prime \Downarrow}],
    X_2[\alpha_{\chi}^{\prime \Downarrow}]) \big] \big) \Big);
\end{multline*}
это последнее условие накладывает на матрицу \ $S$ \ на её
носителе \ $\alpha$ \ тяжёлые требования, и если оно может быть
реализовано, то только в следующем очень особом случае:
\\
редуцированная матрица \ $S$ \ должна иметь {\sl нулевую}
характеристику на носителе \ $\alpha$, кардинал \ $\gamma$ \
должен быть \ $SIN_n$-{\sl субнедостижим}, база\ $\rho$ \ должна
быть {\sl строго меньше} кардинала  \ $\chi^+$, \ кардинал \
$\gamma$ \ должен быть {\sl покрыт подавляющем покрытием} типа \
$\eta^{\ast}$; \ более того, для всех \ $\gamma^{\prime} >
\gamma$, \mbox{$\gamma^{\prime} \in SIN_n$} \ существуют носители
 \ $\alpha^{\prime} > \gamma^{\prime}$ \ с
кардиналами предскачка \ $\alpha_{\chi}^{\prime \Downarrow}$,\
сохраняющими все \ $SIN_n$-сардиналы  \ $\le \gamma^{\prime}$, \ и
с интервалом \ $[\gamma, \alpha_{\chi}^{\prime \Downarrow}[$, \
покрытым максимальными блоками типов \ $\eta^{\prime} \ge
\eta^{\ast}$ \ ниже \ $\alpha_{\chi}^{\prime \Downarrow}$.

Далее такие случаи матрицы \ $S$ \ нулевой характеристики на её
носителе \ $\alpha$ \ будут систематически отвергаться в ходе
определения матричной функции и поэтому мы будем говорить, что эта
нулевая матрица \ $S$ \ здесь на \ $\alpha$ \  {\sl подавляется}
для \ $\gamma$.
\\
Соответственно, нулевая матрица \ $S$ \ на \ $\alpha$ \ с
диссеминатором \ $\delta$ \ и базой \ $\rho$ \ не подавляется для
 \ $\gamma$, \ если это условие нарушается; \ следовательно, всякая матрица \
$S$ \ на \ $\alpha$ \  {\sl не подавлена} для \ $\gamma$, \ если
она {\sl единичная}, или имеет базу \ $\rho \ge \chi^{+}$ \ на \
$\alpha$, \ или \ $\gamma$ \ не \ $SIN_n$-кардинал; так что надо
всегда иметь ввиду характеристику матрицы, её базу \ $\rho$ \ и
соответствующий кардинал \ $\gamma$.

\em

Приостановим в последний раз ненадолго это определение \ref{8.1.}
чтобы пояснить направление развития его финальной части; его суть
состоит в обычной диагональной конструкции, вызывающей финальное
противоречие.
\\
Чтобы  дальнейшая конструкция определения  \ref{8.2.} матричной
функции работала надлежащим образом, она должна руководствоваться
\ $\Pi_{n-2}$-формулой
\[
    U_{n-2}(\mathfrak{n}, x, \chi, a, \delta, \gamma, \alpha, \rho,
    S),
\]
универсальной для класса \ $\Pi_{n-2}$ \ формул с указанными
свободными переменными
\[
    x, \chi, a, \delta, \gamma, \alpha, \rho, S,
\]
и переменным гёделевым номером \ $\mathfrak{n}$ \ таких формул в
базовой модели  \ $\mathfrak{M}$ (\ см. Тарский~\cite{Tarski},
также Аддисон~\cite{Addison}).
\\
Когда этот номер \ $\mathfrak{n}$ \ и переменная \ $x$ \ примут
определённое специальное значение \ $\mathfrak{n}^{\alpha}$ \
одновременно:
\[
    \mathfrak{n} = x = \mathfrak{n}^{\alpha},
\]
тогда эта формула вместе с \ $\Sigma_n$-формулой \
$\mathbf{K}_n^{\forall}(\gamma, \alpha_{\chi}^{\Downarrow})$ \
станут утверждать, что \ $S$ \ -- это \ $\alpha$-матрица,
редуцированная к \ $\chi$ \ на носителе \ $\alpha$ \
характеристики \ $a$ \ с её диссеминатором \ $\delta$ \ и базой
 \ $\rho$, \ {\em допустимые} для \ $\gamma$ \ и подчиняющиеся
 определённым рекурсивным
условиям; напомним, что формула \
$\mathbf{K}_n^{\forall}(\gamma,\alpha)$ \
(определение~\ref{7.1.}\;) означает, что ординал  \ $\alpha$ \
сохраняет все \ $SIN_n$-кардиналы \ $\le \gamma$.
\\
Однако до тех пор, пока значение \ $\mathfrak{n}^{\alpha}$ \ не
будет придано переменным \ $\mathfrak{n}$, \ $x$, \ эта формула
будет действовать в определении~\ref{8.2.} при \ $\mathfrak{n} =
x$:
\[
    U_{n-2}(x, x, \chi, a, \delta, \gamma, \alpha, \rho, S).
\]
Также далее будут использоваться следующие функциональные
ограничения:
\[
    X|\tau^0 = \big\{ (\tau, \eta) \in X: \tau < \tau^0 \big \};
\]
\[
    X|^1 \alpha^0 = \big\{ ((\alpha, \tau), \eta) \in X: \alpha < \alpha^0 \big
    \}.
\]
Теперь вернёмся к определению \ref{8.1.} в последний раз. Условие
\emph{неподавления} \ $\neg A_5^{S,0}$ \ будет действовать в
следующей конъюнкции с формулой \ $U_{n-2}$, \ осуществляющей
``несущую конструкцию'' всего предстоящего определения:
\\

\em

\noindent \emph{III.} Несущее условие
\[
    3.1     \qquad U_{n-2}^{\ast}(\mathfrak{n}, x, \chi, a, \delta,
    \gamma, \alpha, \rho, S, X_1^0|\tau^{\prime}, X_2^0|\tau^{\prime},
    X_1|^1\alpha^0, X_2|^1\alpha^0): \qquad
\]
\[
    U_{n-2}(\mathfrak{n}, x, \chi, a, \delta, \gamma,
    \alpha, \rho, S) \wedge
    \qquad\qquad\qquad\qquad
\]
\[
    \quad \wedge \neg A_{5}^{S,0}(\chi, a, \gamma, \alpha, \rho, S,
    X_1^0|\tau^{\prime}, X_2^0|\tau^{\prime},
    X_1|^1\alpha^0, X_2|^1\alpha^0);
\]
это условие в конъюнкции с формулой \
$\mathbf{K}_n^{\forall}(\gamma, \alpha_{\chi}^{\Downarrow})$ \
после их \ $\vartriangleleft$-ограничения кардиналом \ $\alpha^0$
\ и для констант
\[
    x = \mathfrak{n}^{\alpha}, ~ \chi, ~ \delta, ~ \gamma, ~ \alpha, ~ \rho, ~ \tau^{\prime} < \alpha^0, ~
    S \vartriangleleft \rho
\]
будет описывать следующую ситуацию ниже \ $\alpha^0$: \ $S$ \ --
это матрица, редуцированная к \ $\chi$ \ на её носителе \ $\alpha$
\ характеристики\ $a$, \ допустимая для \ $\gamma$ \ вместе со
своим диссеминатором \ $\delta$ \ и с базой \ $\rho$, \ которая
\emph{неподавлена} для \ $\gamma$ \ ниже\ $\alpha^0$ -- и следует
подчеркнуть, что эта ситуация для каждой пары \ $(\alpha^0,
\tau^{\prime})$ \ будет определяться функциями
\[
    X_1^0|\tau^{\prime} = X_1[\alpha^0]|\tau^{\prime}, ~
    X_2^0|\tau^{\prime} = X_2[\alpha^0]|\tau^{\prime}, ~
    \mbox{\it \ и } ~
    X_1[\alpha_{\chi}^{\prime \Downarrow}], ~
    X_2[\alpha_{\chi}^{\prime \Downarrow}]
\]
для различных \ $\alpha_{\chi}^{\prime \Downarrow} < \alpha^0$ \ и
определённых на меньших парах, поэтому рекурсивная конструкция,
заданная этим условием, будет действовать корректно.
\\
\quad \\
$3.2    \qquad A^0(x, \chi, \tau)$:
\\
\[
    \exists \gamma \Big( \gamma = \gamma_{\tau} \wedge \neg
    \exists a, \delta, \alpha, \rho, S \big( \mathbf{K}_n^{\forall}(\gamma,
    \alpha_{\chi}^{\Downarrow}) \wedge \qquad
\]
\[
    \qquad\qquad\qquad\qquad U_{n-2}(x, x, a, \delta, \gamma, \alpha, \rho, S) \big) \Big);
\]
это условие для \ $x=\mathfrak{n}^{\alpha}$ \ будет означать, что
нет \ $\alpha$-матрицы\ $S$ \ на некотором носителе \ $\alpha$, \
допустимой для \ $\gamma_{\tau}$.
\\
$3.3    \qquad    A_2^0(x, \chi, \tau_1, \tau_2, \tau_3, X_1)$:
\\
\[
    A^0(x, \chi, \tau_1) \wedge A_2(\chi, \tau_1, \tau_2, \tau_3,
    X_1).
\]

\quad \\
\noindent \emph{IV.} Замыкающее условие
\quad \\
Это условие осуществит в дальнейшем замыкание диагонального
рассуждения, устанавливающего финальное противоречие.

\noindent \emph{4.1} \qquad $ \Big( a=0 \longrightarrow \forall
\tau_1^{\prime}, \tau_1^{\prime\prime},
    \tau_2^{\prime}, \tau_3^{\prime}, \eta^{\prime} <
    \alpha_{\chi}^{\Downarrow} \big[
    \gamma_{\tau_1^{\prime}}^{<\alpha_{\chi}^{\Downarrow}} \le
    \delta <
    \gamma_{\tau_3^{\prime}}^{<\alpha_{\chi}^{\Downarrow}} \wedge$
\[
    \wedge A_4^{M b \vartriangleleft \alpha_{\chi}^{\Downarrow} } \big(
    \chi, \tau_1^{\prime}, \tau_1^{\prime\prime}, \tau_2^{\prime},
    \tau_3^{\prime}, \eta^{\prime}, X_1[\alpha_{\chi}^{\Downarrow}],
    X_2[\alpha_{\chi}^{\Downarrow}] \big) \rightarrow
    \eta^{\prime} < \rho \vee \rho = \chi^{+} \big] \Big);
\]
эта формула имеет вот какое содержание для {\textit нулевой}
матрицы\ $S$ \ на носителе \ $\alpha$ \ с диссеминатором \
$\delta$ \ и базой \ $\rho$:

если этот диссеминатор попадает в максимальный блок \ $ \big[
\gamma_{\tau_1^{\prime}}^{<\alpha_{\chi}^{\Downarrow}},
\gamma_{\tau_3^{\prime}}^{<\alpha_{\chi}^{\Downarrow}} \big[$ \
ниже кардинала предскачка \ $\alpha_{\chi}^{\Downarrow}$, \ то
есть если
\[
    \gamma_{\tau_1^{\prime}}^{<\alpha_{\chi}^{\Downarrow}} \le
    \delta <
    \gamma_{\tau_3^{\prime}}^{<\alpha_{\chi}^{\Downarrow}},
\]
то эта база \ $\rho$ \ должна существенно возрасти и превзойти тип
\ $\eta^{\prime}$ \ этого самого блока, или даже принять
максимальное возможное значение:
\[
    \eta^{\prime} < \rho \vee \rho = \chi^{+},
\]
за неимением ничего лучшего;
\\
поэтому в подобных случаях интервал \ $[\tau_1^{\prime},
\tau_3^{\prime}[$ \ и соответствующий интервал
\[
    \big[ \gamma_{\tau_1^{\prime}}^{<\alpha_{\chi}^{\Downarrow}},
    \gamma_{\tau_3^{\prime}}^{<\alpha_{\chi}^{\Downarrow}} \big[
\]
будут считаться  {\sl ``обременительными''} для такой {\sl
нулевой} матрицы \ $S$ \ с таким диссеминатором на носителе \
$\alpha$ \ и будут {\sl препятствовать образованию матрицы} \ $S$
\ на этом носителе (с этим диссеминатором \ $\delta$).
\\
\quad \\

\noindent \emph{V.} Условие эквиинформативности
\\
\quad \\
\emph{5.1}\quad $A_{6}^{e}(\chi, \alpha^0)$:
\[
    \chi < \alpha^0 \wedge A_n^{\vartriangleleft \alpha^0}(\chi) =
    \| u_n^{\vartriangleleft \alpha^0}(\underline{l})\| \wedge
    SIN_{n-2}(\alpha^0) \wedge
\]
\[
    \qquad\qquad\qquad\qquad
    \wedge \forall \gamma < \alpha^0 \exists
    \gamma_1 \in [\gamma,\alpha^0[ \quad
    SIN_n^{< \alpha^0}(\gamma_1);
\]
кардинал \ $\alpha^0$ \ здесь с этим свойством называется,
напомним, {\sl эквиинформативным} с кардиналом \ $\chi$.
\\
\hspace*{\fill} $\dashv$
\\
\em Последнее понятие было использовано выше несколько раз (см.
также
\cite{Kiselev2,Kiselev3,Kiselev4,Kiselev5,Kiselev6,Kiselev7,Kiselev8},
\cite{Kiselev11}) и здесь оно акцентируется по причине его особой
важности: каждое \ $\Pi_n$-утверждение \ $\varphi(l)$ \
выполняется или нет в каждом генерическом расширении \
$\mathfrak{M}[l]$ \ ниже \ $\chi$ \ и также в этом же расширении
ниже \ $\alpha^0$ \ одновременно (см. комментарий после
(\ref{e7.1})); наилучший пример такого \ $\alpha^0$ -- это
кардинал предскачка \ $\alpha_{\chi}^{\Downarrow}$ \ после \
$\chi$ \ любого матричного носителя \ $\alpha > \chi$ \ (если этот
кардинал пределен для класса \ $SIN_n^{<
\alpha_{\chi}^{\Downarrow}}$).
\\

Теперь всё готово для того, чтобы собрать все введённые выше
фрагменты вместе в следующем \textit{интегрирующем}
определении~\ref{8.2.}, где рассматривается переменная матрица \
$S$ \ на её носителе \ $\alpha$.\
\\
Требования, которые накладываются в нём на матрицу \ $S$ \ на её
носителе \ $\alpha$ \ и на её диссеминатор \ $\delta$ \ с базой \
$\rho$, \ зависят от функций \ $X_i,~i=\overline{1,5}$, \ которые
уже будут рекурсивно  определены ниже кардинала предскачка \
$\alpha_{\chi}^{\Downarrow}$; \ они будут заданы на определённом
подмножестве множества
\begin{equation} \label{e8.1}
    \qquad\qquad
    \mathcal{A}_{\chi}^{\alpha_{\chi}^{\Downarrow}} = \Big\{
    (\alpha^0, \tau): \exists \gamma < \alpha^0 \big( \chi <
    \gamma=\gamma_{\tau}^{< \alpha^0} \wedge \qquad\qquad\qquad\qquad
\end{equation}
\[
    \qquad\qquad\qquad\qquad\qquad
    \wedge \alpha^0 \le \alpha_{\chi}^{\Downarrow} \wedge
    A_6^e(\chi,\alpha^0) \big) \Big\}
\]
и поэтому функции
\[
    X_i^0 = X_i[\alpha^0], ~ i=\overline{1,5}
\]
будут заданы на соответствующем подмножестве множества
\[
    \big\{ \tau: \gamma_{\tau} \in SIN_{n-1}^{<\alpha^0} \big\}
\]
для каждого кардинала \ $\alpha^0 \le \alpha_{\chi}^{\Downarrow}$
\ эквиинформативного с \ $\chi$. \ Это множество \
$\mathcal{A}_{\chi}^{\alpha_{\chi}^{\Downarrow}}$ \ полагается
канонически упорядоченным (с \ $\alpha^0$ \ как первой компонентой
в этом упорядочении и с \ $\tau$ \ как второй).
\\

Итак, переменная \ $X_{2}^0$ \ будет играть здесь роль
характеристической функции (то есть функции характеристик) \
$a_{f}^{<\alpha^0}$, \ определённой ниже кардинала \ $\alpha^0$; \
$X_{1}^0$ \ будет играть роль матричной функции \ $\alpha
S_{f}^{<\alpha^0} $; \ $X_{3}^0$ \ -- роль диссеминаторной функции
(то есть фунции диссеминаторов) \
$\widetilde{\delta}_{f}^{<\alpha^0}$; \ $X_{4}^0$ \ -- её базовой
функции (то есть функции баз данных) \ $\rho _{f}^{<\alpha^0}$; \
$X_{5}^0$ \ -- роль несущей функции (то есть функции носителей) \
$\alpha _{f}^{<\alpha^0}$; \ все они будут определены ниже \
$\alpha^0$.
\\

После того, как все эти функции будут будут определены для всех
таких кардиналов
\[
    \alpha^0 \le \alpha_{\chi}^{\Downarrow}
\]
-- тогда в заключение результирующее требование будет наложено на
саму матрицу \ $S$ \ на её носителе \ $\alpha$ \ вместе с её
диссеминатором \ $\delta$ \ и базой данных \ $\rho$ \ в
зависимости от расположения этого \ $\delta$, \ более точно -- в
зависимости от максимального блока
\[
    \big[ \gamma_{\tau_1^{\prime}}^{<\alpha_{\chi}^{\Downarrow}},
    ~ \gamma_{\tau_3^{\prime}}^{<\alpha_{\chi}^{\Downarrow}} \big[
\]
заключающего в себе этот \ $\delta$, \ который будет уже определён
ниже \ $\alpha_{\chi}^{\Downarrow}$.
\\
Здесь, напомним, действует замыкающее требование, наложенное на
матрицу \ $S$ \ на её носителе \ $\alpha$, \ упомянутое выше:
\\
если \ $S$ \ на \ $\alpha$ \ имеет {\textit нулевую}
характеристику и его допустимый диссеминатор \ $\delta$ \ попадает
в максимальный блок типа \ $\eta^{\prime}$ \ ниже \
$\alpha_{\chi}^{\Downarrow}$, \ то $\eta^{\prime} < \rho \wedge
\rho = \chi^{+}$; \ таким образом, в этом случае база данных \
$\rho$ \ должна значительно возрасти и мы увидим, что это
возможно, но всякий раз ведёт к противоречию.
\\

Ещё здесь потребуется формулировка \ $L j^{<\alpha}(\chi)$ \
понятия насыщенности кардинала \ $\chi$ \ (см. аргумент
перед~(\ref{e7.1}) или определение~6.9~4)~\cite{Kiselev11}\;);
напомним также, что  \ $\widehat{\rho}$ \ обозначает замыкание
 \ $\rho$ \ относительно функции пары.

Итак, рекурсивное определение матричной функции, заданное на
множестве \ $\mathcal{A}_{\chi}^{\alpha_{\chi}^{\Downarrow}}$, \
начинается \label{c10};  после определения 8.4 мы покажем, как это
рекурсивное определение действует:
\endnote{
\ стр. \pageref{c10}. \ Это определение использовалось ранее
(Киселев~\cite{Kiselev6,Kiselev7,Kiselev8}) в виде единого текста,
а теперь оно разделено на части чтобы прояснить его структуру.
\\
\quad \\
} %

\begin{definition}
\label{8.2.} \ \newline

\emph{1)}\quad Пусть
\[
    U_{n-2}(\mathfrak{n},x,\chi,a,\delta,\gamma,\alpha,\rho ,S)
\]
это \ $\Pi _{n-2}$-формула, универсальная для класса \ $\Pi _{n-2}
$, \ где \ $\mathfrak{n}$ \ это переменный гёделев номер \ $\Pi
_{n-2}$-формул со свободными переменными \ $x$, $\chi $, $a$,
$\delta $, $\gamma $, $ \alpha $, $\rho $, $S$, \ и пусть
\[
    U_{n-2}^{\ast}(x,\chi,a,\delta, \gamma,
    \alpha ,\rho ,S,
    X_1^0|\tau^{\prime}, X_2^0|\tau^{\prime},
    X_1|^1\alpha^0, X_2|^1\alpha^0).
\]
это формула
\begin{multline*}
    U_{n-2}(x,x,\chi,a,\delta,\gamma,\alpha,\rho,S) \wedge
\\
    \wedge \neg A_{5}^{S,0}(\chi, a, \gamma, \alpha, \rho, S,
    X_1^0|\tau^{\prime}, X_2^0|\tau^{\prime},
    X_1|^1\alpha^0, X_2|^1\alpha^0).
\end{multline*}

\quad \\
\emph{2)}\quad Пусть
\[
    A_7^{RC}(x, \chi, X_1, X_2, X_3, X_4, X_5, \alpha_{\chi}^{\Downarrow})
\]
это следующая \ $\Delta_1$-формула, составляющая требуемое {\sl
условие рекурсии}:

\begin{equation*}
\begin{array}{l}
\bigwedge _{1\leq i\leq 5} \Big( ( X_{i}\mbox{\it \ \ это
функция})\wedge X_{i}\vartriangleleft \alpha_{\chi}^{\Downarrow +}
\wedge
\quad  \\
\quad  \\
\hspace{0.1cm}\wedge dom(X_{i}) \subseteq \Big \{ (\alpha^0,
\tau): \exists \gamma < \alpha^0 \big( \chi < \gamma =
\gamma_{\tau}^{< \alpha^0} \wedge \qquad
\quad  \\
\quad  \\
\hspace{4.5cm} \wedge \alpha^0 \le \alpha_{\chi}^{\Downarrow}
\wedge A_6^{e}(\chi, \alpha^0) \big) \Big\} \Big) \wedge
\end{array}
\end{equation*}

\vspace{-6pt}

\begin{equation*}
\begin{array}{l}
\hspace{0.1cm} \wedge \forall \alpha^0 \Big( \big( \alpha^0 \le
\alpha_{\chi}^{\Downarrow} \wedge A_6^{e}(\chi,\alpha^0) \big) \longrightarrow
\quad  \\
\quad  \\
\hspace{0.5cm} \longrightarrow \exists X_1^0, X_2^0, X_3^0, X_4^0,
X_5^0, X_1^{1,0}, X_2^{1,0} ~ \Big[\bigwedge_{1 \le i \le 5} X_i^0
= X_i[\alpha^0] \wedge
\quad  \\
\quad  \\
\hspace{4.5cm} \wedge X_1^{1,0} = X_1 |^1 \alpha^0 \wedge X_2^{1,0}
= X_2 |^1 \alpha^0 \wedge
\quad  \\
\quad  \\
\hspace{0.1cm} \wedge \forall \tau ^{\prime },\gamma ^{\prime },
\gamma^{\prime\prime} < \alpha^0 \Big( \chi <\gamma ^{\prime
}\wedge \gamma ^{\prime }=\gamma _{\tau ^{\prime }}^{<\alpha^0}
\wedge \gamma^{\prime\prime} =
\gamma_{\tau^{\prime}+1}^{<\alpha^0} \longrightarrow
\end{array}
\end{equation*}

\vspace{-6pt}

\begin{equation*}
\begin{array}{l}
\forall a^{\prime} \Big( X_{2}^0(\tau ^{\prime })=a^{\prime
}\leftrightarrow
\quad  \\
\quad  \\
\hspace{0.1cm} \leftrightarrow a^{\prime }=\max_{\leq } \big\{
a^{\prime \prime }:\exists \delta ^{\prime \prime },\alpha
^{\prime \prime },\rho ^{\prime \prime }< \gamma^{\prime \prime }
\exists S^{\prime \prime } \vartriangleleft \chi^{+} \big(
\mathbf{ K}_{n}^{\forall <\alpha^0}( \gamma ^{\prime},\alpha
_{\chi }^{\prime \prime \Downarrow }) \wedge
\quad  \\
\quad  \\
\wedge U_{n-2}^{\ast \vartriangleleft \alpha^0}(x,\chi ,a^{\prime \prime
}, \delta ^{\prime \prime },\gamma^{\prime}, \alpha ^{\prime
\prime }, \rho^{\prime \prime}, S^{\prime\prime}, X_1^0|\tau^{\prime}, X_2^0|\tau^{\prime},
X_1^{1,0}, X_2^{1,0}) \big) \big\} \Big) \wedge  \\
\end{array}
\end{equation*}

\vspace{-6pt}

\begin{equation*}
\begin{array}{l}
\wedge \forall S^{\prime } \Big( X_{1}^0(\tau ^{\prime
})=S^{\prime }\longleftrightarrow \exists a^{\prime }~ \Big(
a^{\prime} = X_{2}^0(\tau ^{\prime })\wedge  \\
\quad  \\
\hspace{0.5cm}\wedge S^{\prime } =\min_{\underline{\lessdot }}
\big\{ S^{\prime \prime } \vartriangleleft \chi^{+} :\exists
\delta ^{\prime \prime },\alpha ^{\prime \prime },\rho ^{\prime
\prime }< \gamma^{\prime\prime} \big( \mathbf{K}_{n}^{\forall
<\alpha^0} (\gamma ^{\prime },\alpha _{\chi }^{\prime \prime
\Downarrow })\wedge
\quad  \\
\quad  \\
\wedge U_{n-2}^{\ast \vartriangleleft \alpha^0}(x,\chi ,a^{\prime
},\delta ^{\prime \prime}, \gamma^{\prime}, \alpha ^{\prime
\prime}, \rho ^{\prime \prime },S^{\prime \prime },
X_1^0|\tau^{\prime}, X_2^0|\tau^{\prime},
X_1^{1,0}, X_2^{1,0}) \big) \big\} \Big) \Big) \wedge  \\
\end{array}
\end{equation*}

\vspace{-6pt}

\begin{equation*}
\begin{array}{l}
\wedge \forall \delta ^{\prime } \Big( X_{3}^0(\tau ^{\prime
})=\delta ^{\prime }\longleftrightarrow \exists a^{\prime
},S^{\prime } \Big( a^{\prime }=X_1^0(\tau ^{\prime })\wedge
S^{\prime }=X_2^0(\tau ^{\prime })\wedge
\quad  \\
\quad  \\
\hspace{0.5cm}\wedge \delta ^{\prime }=\min_{\leq } \big\{
\delta^{\prime \prime} < \gamma^{\prime} : \exists \alpha ^{\prime
\prime },\rho ^{\prime \prime }< \gamma^{\prime\prime} \big(
\mathbf{K}_{n}^{\forall <\alpha^0}(\gamma ^{\prime },\alpha _{\chi
}^{\prime \prime \Downarrow })\wedge
\quad  \\
\quad  \\
\wedge U_{n-2}^{\ast \vartriangleleft \alpha^0}(x,\chi ,a^{\prime
},\delta ^{\prime \prime }, \gamma^{\prime}, \alpha ^{\prime
\prime },\rho ^{\prime \prime },S^{\prime }, X_1^0|\tau^{\prime},
X_2^0|\tau^{\prime}, X_1^{1,0}, X_2^{1,0}) \big) \big\} \Big) \Big) \wedge \\
\end{array}
\end{equation*}

\vspace{-6pt}

\begin{equation*}
\begin{array}{l}
\wedge \forall \rho ^{\prime }~ \Big( X_{4}^0(\tau ^{\prime
})=\rho ^{\prime }\longleftrightarrow \exists a^{\prime
},S^{\prime },\delta ^{\prime } \Big( a^{\prime }=X_1^0(\tau
^{\prime}) \wedge S^{\prime }=X_2^0(\tau ^{\prime})\wedge
\quad  \\
\quad  \\
\hspace{0.1cm}\wedge \delta ^{\prime }=X_{3}^0(\tau ^{\prime })
\wedge \rho ^{\prime }=\min_{\leq } \big \{ \rho ^{\prime \prime }
< \chi^{+} :\exists \alpha^{\prime\prime} < \gamma^{\prime\prime}
\big( \mathbf{K}_{n}^{\forall <\alpha^0}(\gamma ^{\prime },\alpha
_{\chi }^{\prime \prime \Downarrow })\wedge
\quad  \\
\quad  \\
\hspace{0.2cm}\wedge U_{n-2}^{\ast \vartriangleleft
\alpha^0}(x,\chi ,a^{\prime },\delta ^{\prime }, \gamma^{\prime},
\alpha ^{\prime\prime}, \rho ^{\prime \prime }, S^{\prime },
X_1^0|\tau^{\prime}, X_2^0|\tau^{\prime}, X_1^{1,0}, X_2^{1,0})
\big) \big\} \Big) \Big) \wedge
\end{array}
\end{equation*}

\vspace{-6pt}

\[
\wedge \forall \alpha ^{\prime } \Big( X_{5}^0(\tau ^{\prime
})=\alpha ^{\prime }\longleftrightarrow \exists a^{\prime
},S^{\prime },\delta ^{\prime },\rho ^{\prime } \Big(
a^{\prime}=X_1^0(\tau ^{\prime })\wedge S^{\prime }=X_2^0(\tau
^{\prime })\wedge \qquad \qquad \qquad
\]
\[
\hspace{0.1cm}\delta ^{\prime }=X_3^0(\tau ^{\prime })\wedge
\rho^{\prime }=X_4^0(\tau ^{\prime })\wedge \alpha ^{\prime
}=\min_{\leq} \big \{ \alpha^{\prime\prime}<
\gamma^{\prime\prime}: \mathbf{K}_{n}^{\forall <\alpha^0}(\gamma
^{\prime },\alpha _{\chi}^{\prime \prime \Downarrow })\wedge
\]
\[
\hspace{0.2cm}\wedge U_{n-2}^{\ast \vartriangleleft
\alpha^0}(x,\chi , a^{\prime }, \delta ^{\prime },
\gamma^{\prime}, \alpha^{\prime \prime}, \rho ^{\prime },S^{\prime
}, X_1^0|\tau^{\prime}, X_2^0|\tau^{\prime}, X_1^{1,0}, X_2^{1,0}
)  \big\} \Big) \Big) \Big) \Big]
\Big). \quad  \\
\]

\noindent \emph{3)}\quad Мы обозначаем через

\begin{equation*}
\alpha \mathbf{K}_{n+1}^{\exists }(x,\chi ,a,\delta ,\gamma ,\alpha ,\rho ,S)
\end{equation*}
\vspace{0pt}

\noindent  \ $\Pi _{n-2}$-формулу, которая эквивалентна следующей
формуле:

\vspace{6pt}
\begin{equation*}
\begin{array}{l}
(a=0\vee a=1)\wedge \sigma (\chi ,\alpha ,S)\wedge Lj^{<\alpha
}(\chi )\wedge \chi <\delta <\gamma <\alpha \wedge
\quad  \\
\quad  \\
\hspace{0.5cm}\wedge S\vartriangleleft \rho \leq \chi ^{+}\wedge
\rho = \widehat{\rho }\wedge SIN_{n}^{<\alpha _{\chi }^{\Downarrow
}}(\delta )\wedge SIN_{n+1}^{<\alpha _{\chi }^{\Downarrow }}\left[
<\rho \right] (\delta )\wedge
\end{array}
\end{equation*}

\vspace{-6pt}

\begin{equation*}
\begin{array}{l}
\hspace{1.0cm}\wedge \forall \gamma <\alpha _{\chi }^{\Downarrow
}~\exists \gamma ^{\prime }\in [ \gamma ,\alpha _{\chi
}^{\Downarrow }[ ~SIN_{n}^{<\alpha _{\chi }^{\Downarrow }}(\gamma
^{\prime })\wedge cf(\alpha _{\chi }^{\Downarrow })\geq \chi
^{+}\wedge
\quad  \\
\quad  \\
\wedge \exists X_{1},X_{2},X_{3},X_{4},X_{5}\Bigl\{ A_7^{RC}(x,
\chi, X_1, X_2, X_3, X_4, X_5, \alpha_{\chi}^{\Downarrow}) \wedge
\quad  \\
\quad  \\
\Big( a=0\longleftrightarrow \exists \tau_1^{\prime},
\tau_2^{\prime}, \tau_3^{\prime}<\alpha_{\chi}^{\Downarrow} ~
\big( A_{2}^{0\vartriangleleft \alpha _{\chi }^{\Downarrow }}
(\chi, \tau_1^{\prime}, \tau_2^{\prime}, \tau_3^{\prime},
X_1[\alpha_{\chi}^{\Downarrow}] ) \wedge
\quad  \\
\quad  \\
\hspace{0.2cm} \wedge \forall \tau^{\prime\prime} \big(
\tau_1^{\prime} < \tau^{\prime\prime} \le \tau_2^{\prime}
\rightarrow X_2[\alpha_{\chi}^{\Downarrow}](\tau^{\prime\prime})=1
\big) \wedge X_1[\alpha_{\chi}^{\Downarrow}] (\tau_{2}^{\prime})
=S \big) \Big)\wedge
\end{array}
\end{equation*}

\vspace{-6pt}

\begin{equation*}
\begin{array}{l}
\wedge \Big( a=0\longrightarrow \forall \tau_1^{\prime},
\tau_1^{\prime\prime}, \tau_2^{\prime}, \tau_3^{\prime},
\eta^{\prime} <\alpha _{\chi }^{\Downarrow }~\big[
\gamma_{\tau_{1}^{\prime}}^{<\alpha _{\chi }^{\Downarrow }} \leq
\delta < \gamma_{\tau_{3}^{\prime}}^{<\alpha _{\chi
}^{\Downarrow }}\wedge  \\
\quad  \\
\wedge A_{4}^{M b \vartriangleleft \alpha _{\chi }^{\Downarrow }}
(\chi ,\tau_1^{\prime}, \tau_1^{\prime\prime}, \tau_2^{\prime},
\tau_3^{\prime}, \eta^{\prime}, X_1[\alpha_{\chi}^{\Downarrow}],
X_2[\alpha_{\chi}^{\Downarrow}]) \rightarrow \\
\quad  \\
\qquad\qquad\qquad\qquad\qquad\qquad\qquad\qquad \rightarrow
\eta^{\prime} < \rho \vee \rho =\chi ^{+}\big] \Big) \Bigr\}.
\end{array}
\end{equation*}
\vspace{6pt}

Обозначим через \ $\mathbf{K}^{0}(\chi ,a,\delta ,\alpha ,\rho)$ \
последний конъюнктивный член в больших фигурных скобках $\{\;,\}$
в последней формуле, то есть замыкающее условие:
\[
    \Bigl(a=0\longrightarrow \forall \tau_1^{\prime},
    \tau_1^{\prime\prime}, \tau_2^{\prime}, \tau_3^{\prime}, \eta^{\prime}
    <\alpha _{\chi }^{\Downarrow } ~ \bigl[
    \gamma_{\tau_{1}^{\prime}}^{<\alpha _{\chi }^{\Downarrow }} \leq
    \delta < \gamma_{\tau_{3}^{\prime}}^{<\alpha _{\chi
    }^{\Downarrow }} \wedge  \qquad \qquad \\
    \quad  \\
\]
\[
    \wedge A_{4}^{M b\vartriangleleft \alpha _{\chi }^{\Downarrow
    }}(\chi ,\tau_1^{\prime}, \tau_1^{\prime\prime}, \tau_2^{\prime},
    \tau_3^{\prime}, \eta^{\prime}, X_1[\alpha_{\chi}^{\Downarrow}],
    X_2[\alpha_{\chi}^{\Downarrow}]) \rightarrow \qquad \qquad
\]
\hspace*{17em} $\rightarrow \eta^{\prime} < \rho \vee \rho =
\chi^{+}\bigr]\Bigr).$\label{c11}
\endnote{
\ стр. \pageref{c11}. \ Это замыкающее условие \ $\mathbf{K}^0$ \
действует здесь как условие \ $\mathbf{K}^0$, \ использованное
ранее (Киселев~\cite{Kiselev6,Kiselev7,Kiselev8}), но более
оперативным образом, так как оно действует теперь вполне успешно
без подформулы \
$(\gamma_{\tau_3^{\prime}}^{<\alpha_{\chi}^{\Downarrow}} = \gamma
\rightarrow \lim(\gamma))$, \ которая ранее вызывала значительное
усложнение доказательства.
\\
\quad \\
} %
\quad \\

\noindent Функции \ $X_1[\alpha_{\chi}^{\Downarrow}],
X_2[\alpha_{\chi}^{\Downarrow}]$ \ не упоминаются здесь в
обозначении \ $\mathbf{K}^{0}$ \ для краткости, так как они
определяются единственным образом в предшествующей части этой
формулы \ $\alpha \mathbf{K}_{n+1}^{\exists}$.

\emph{4)}\quad Формула \ $\alpha \mathbf{K}_{n+1}^{\exists }$ \ --
это \ $\Pi _{n-2}$-формула и поэтому она получает свой гёделев
номер \ $\mathfrak{n} ^{\alpha }$, \ то есть:
\begin{equation*}
\alpha \mathbf{K}_{n+1}^{\exists }(x,\chi ,a,\delta ,\gamma ,\alpha ,\rho
,S)\longleftrightarrow U_{n-2}(\mathfrak{n}^{\alpha },x,\chi ,a,\delta
,\gamma ,\alpha ,\rho ,S).
\end{equation*}
Присвоим значение \ $\mathfrak{n}^{\alpha }$ \ переменной \ $x$ \
здесь в этой эквивалентности и везде далее и соответственно этому
символы \ $\mathfrak{n}^{\alpha}$, $x$ \ будут опущены везде далее
в обозначениях.
\\
Мы обозначаем через \ $\alpha \mathbf{K}^{<\alpha _{1}}(\chi
,a,\delta ,\gamma ,\alpha ,\rho ,S)$ \ следующую \
$\Delta_1$-формулу:
\[
    \mathbf{K}_{n}^{\forall <\alpha _{1}}(\gamma ,\alpha _{\chi }^{\Downarrow
    })\wedge \alpha \mathbf{K}_{n+1}^{\exists \vartriangleleft \alpha _{1}}
    (\chi ,a,\delta ,\gamma ,\alpha ,\rho
    ,S)\wedge \alpha <\alpha _{1},
\]
и, соответственно, через\ $\alpha \mathbf{K}^{\ast <\alpha
_{1}}(\chi,a,\delta,\gamma, \alpha, \rho, S)$ -- формулу, которая
получается из формулы \ $\alpha \mathbf{K}^{<\alpha_1}$ \
присоединением к ней конъюнктивного условия {\slнеподавления}
матрицы \ $S$ \ на \ $\alpha$ \ для \ $\gamma$ \ (см. определение
\ref{8.1.}~2.6\;), но ниже \ $\alpha_1 < k$ \ (это делалось выше в
пунктах 2), 3) для  \ $\alpha^0$,\ \ $\alpha_{\chi}^{\Downarrow}$)
-- следующим образом:
\[
    \mathbf{K}_{n}^{\forall <\alpha _{1}}(\gamma ,\alpha _{\chi }^{\Downarrow
    })\wedge \alpha \mathbf{K}_{n+1}^{\exists \vartriangleleft  \alpha _{1}}
    (\chi,a,\delta,\gamma, \alpha, \rho, S) \wedge \alpha <\alpha_{1}
    \wedge
\]
\[
    \wedge \neg \Big( a=0 \wedge SIN_n^{<\alpha_1}(\gamma) \wedge
    \rho < \chi^+ \wedge \sigma(\chi, \alpha, S) \wedge \qquad\qquad\qquad
\]
\[
    \wedge \exists X_1, X_2, X_3, X_4, X_5 \Big(
    A_7^{RC}(\mathfrak{n}^{\alpha}, \chi, X_1, X_2, X_3, X_4, X_5,
    \alpha_1) \wedge
\]
\[
    \wedge \exists \eta^{\ast}, \tau < \gamma \Big( \gamma =
    \gamma_{\tau}^{<\alpha_1} \wedge A_{5.4}^{sc \vartriangleleft
    \alpha_1} (\chi, \gamma, \eta^{\ast}, X_1[\alpha_1]|\tau,
    X_2[\alpha_1]|\tau) \wedge
\]
\[
    \wedge \forall \tau^{\prime} \big( \tau < \tau^{\prime} \wedge
    SIN_n^{<\alpha_1}(\gamma_{\tau^{\prime}}^{<\alpha_1})
    \rightarrow
    \qquad\qquad\qquad \qquad\qquad\qquad
\]
\[
    \exists \alpha^{\prime}, S^{\prime} \big[
    \gamma_{\tau^{\prime}}^{<\alpha_1} < \alpha^{\prime} <
    \gamma_{\tau^{\prime}+1}^{<\alpha_1} \wedge
    SIN_n^{<\alpha_{\chi}^{\prime\Downarrow}}(\gamma_{\tau^{\prime}}^{<\alpha_1})
    \wedge \sigma(\chi, \alpha^{\prime}, S^{\prime}) \wedge
\]
\[
    \qquad \qquad \qquad
    \wedge A_{5.5}^{sc \vartriangleleft \alpha_1}(\chi, \gamma, \eta^{\ast},
    \alpha_{\chi}^{\prime \Downarrow},
    X_1[\alpha_{\chi}^{\prime \Downarrow}],
    X_2[\alpha_{\chi}^{\prime \Downarrow}]) \big] \big) \Big)
    \Big) \Big);
\]
здесь утверждается допустимость матрицы \ $S$ \ на \ $\alpha$, \ и
сверх того --- её неподавленность для \ $\gamma$ \ ниже \
$\alpha_1$. \ Таким образом, если \ $\alpha
\mathbf{K}^{<\alpha_1}$ \ выполняется, но \ $\alpha
\mathbf{K}^{\ast <\alpha_1}$ \ нарушается, то \ $S$ \ на \
$\alpha$ \ допустима, но подавлена (всё это ниже \
$\alpha_1$).\label{c12}
\endnote{
\ стр. \pageref{c12}. \ Эти понятия аналогичным образом можно
ввести и в нерялитивизованной форме для \ $\alpha_1 = k$, \ но в
этой форме они не используются в дальнейшем; кроме того, в этой
форме они нуждаются в привлечении более сложного (неэлементарного)
языка над \ $L_k$.
\\
\quad \\
} %

\emph{5)}\quad Если формула \ $\alpha \mathbf{K}^{<\alpha
_{1}}(\chi ,a,\delta ,\gamma ,\alpha ,\rho ,S)$ \ выполняется
константами \ $\chi $, $a$, $\delta $, $\alpha $, $\gamma $, $\rho
$, $S$, $\alpha _{1}$, \ то мы говорим, что \ $\chi $, $a$,
$\delta $, $\alpha $, $\rho $, $S$ \  {\sl очень сильно} допустимы
для \ $\gamma $ \ ниже \ $\alpha _{1}$.
\\
Если некоторые из них фиксированы или указываются в контексте, то
мы говорим, что остальные {\sl очень сильно} допустимы для них (и
для \ $\gamma $) ниже \ $\alpha _{1}$. \ Соответственно, через
\[
    \alpha \mathbf{K}^{<\alpha_1}(\chi, \gamma, \alpha, S)
\]
будет обозначаться формула
\[
    \exists a, \delta, \rho < \gamma ~ \alpha
    \mathbf{K}^{<\alpha_1}(\chi, a, \delta, \gamma, \alpha, \rho, S)
\]
означающая, что \ $S$ \ на \ $\alpha$ \ {\sl очень сильно}
допустима для \ $\gamma$ \ ниже \ $\alpha_1$.

\emph{6)} \ Матрица \ $S$ \ называется автоэкзорцизивной или,
коротко, \ \mbox{$\alpha $-матрицей} характеристики \ $a$, \ {\sl
очень сильно} допустимой на носителе \ $\alpha $ \ для \ $\gamma =
\gamma_{\tau}^{<\alpha_1}$ \ ниже \ $\alpha _{1}$, \ если она
обладает на \ $\alpha$ \ некоторым диссеминатором \
$\delta<\gamma$ \ с базой \ $\rho$ \ очень сильно допустимыми для
них также ниже \ $ \alpha _{1}$.

В каждом подобном случае \ $\alpha $-матрица обозначается общим
символом \ $\alpha S$ \ или \ $S$.

Если \ $a_1=k$, \ или \ $\alpha_1$ \ указывается в контексте, то
верхние индексы \ $< \alpha_1, \vartriangleleft \alpha_1$ \ и
другие упоминания о \ $\alpha_1$ \ опускаются.

Далее все понятия допустимости полагаются очеь сильными, поэтому
термин ``очень сильно'' будет в дальнейшем опускаться.
\hspace*{\fill} $\dashv$
\end{definition}

\noindent Таким образом, ограниченная формула
\[
    \alpha \mathbf{K}^{\ast < \alpha _{1}}
    (\chi,a,\delta,\gamma, \alpha, \rho, S)
\]
возникает из формулы \ $\alpha \mathbf{K}^{<\alpha_1}$ \
посредством присоединения к ней условия неподавления матрицы \ $S$
\ на \ $\alpha$ \ для \ $\gamma$ \ ниже \ $\alpha_1$, \ которое
получается из условия \ $\neg A_5^{S,0}$ \ описанным выше способом
--- \ $\vartriangleleft$-ограничением кардиналом \ $\alpha_1$ \ (то
есть \ $\vartriangleleft$-ограничением его индивидных переменных
этим кардиналом и заменой его подформул \ $SIN_n(\gamma)$ \ на \
$SIN_n^{<\alpha_1}(\gamma)$).
\\
Соответственно, если матрица \ $S$ \ на \ $\alpha$ \
\textit{подавлена} для кардинала \ $\gamma$ \ ниже \ $\alpha_1$, \
то \ $\gamma$ \ -- это\ $SIN_n^{< \alpha_1}$-кардинал и \ $S$ \
имеет \textit{нулевую} характеристику на \ $\alpha$ \ и базу \
$\rho < \chi^+$ \ ниже \ $\alpha_1$.

Везде далее \ $\chi =\chi ^{\ast }<\alpha _{1}$; \ мы будем часто
опускать обозначения функций \ $X_1=\alpha S_{f}^{<\alpha _{1}} $,
\ $X_2=a_{f}^{<\alpha _{1}}$ \ и символов \ $\chi ^{\ast }$, \
$\mathfrak{n} ^{\alpha }$ \ в написаниях формул
\begin{equation*}
A_{0}-A_5^{S,0}, \quad A^0, \quad A_2^0, \quad \alpha
\mathbf{K}_{n+1}^{\exists }, \quad \mathbf{K }^0, \quad \alpha
\mathbf{K}^{<\alpha _{1}}, \quad \alpha \mathbf{K}^{\ast <\alpha
_{1}}
\end{equation*}
(из определений \ref{8.1.}, \ref{8.2.}) и других обозначений для
некоторой краткости (если это не вызовет недоразумений); например,
всякий кардинал предскачка \ $\alpha_{\chi^{\ast}}^{\Downarrow}$ \
будет обозначаться через \ $\alpha^{\Downarrow}$ \ и так далее.

\noindent Касательно этих формул следует отметить, что определение
\ref{8.2.} было сформировано с целью получить ключевую формулу \
$\alpha \mathbf{K}_{n+1}^{\exists}$ \ из класса \ $\Pi_{n-2}$. \
Для этого все составляющие формулы были \
$\vartriangleleft$-ограничены кардиналами $\alpha^0$ \ или
$\alpha^{\Downarrow}$.
\\
Но в дальнейшем эти ограничения будут часто опускаться безо всякой
потери содержания этих формул, потому что их индивидные переменные
и константы будут в действительности ограничены указанными в
контексте кардиналами в ходе их применения.
\\

Очевидно, переменные \ $X_i,~ i=\overline{1,5}$ \ задаются в
определении \ref{8.2.} единственным образом через его параметры,
поэтому такие же функции можно рекурсивно определить той же
рекурсией на таком же множестве пар \ $(\alpha_1,\tau)$ \
(напомним множество (\ref{e8.1}))
\[
    \mathcal{A} = \bigl\{ (\alpha_1,\tau): \exists \gamma<\alpha^1
    (\chi^{\ast}<\gamma=\gamma_{\tau}^{<\alpha_1} \wedge
    A_6^e(\chi,\alpha_1) ) \bigr\}
\]
упорядоченном как и выше канонически  (с \ $\alpha_1$ \ как первой
компонентой в этом упорядочении и \ $\tau$ \ как второй).

\begin{definition}
\label{8.3.} \ \\
Пусть \ $\chi^{\ast}<\alpha_1$.
\\
\emph{1)}\quad Мы называем характеристической функцией уровня \
$n$ \ ниже \ $\alpha _{1}$ \ редуцированной к \ $\chi ^{\ast }$ \
следующую функцию
\[
    a_{f}^{<\alpha _{1}}=(a_{\tau }^{<\alpha _{1}})_{\tau },
\]
принимающую значения: \vspace{-6pt}
\begin{multline*}
    a_{\tau }^{<\alpha _{1}} =
\\
    = \max_{\leq } \bigl \{ a: \exists \delta ,\alpha,
    \rho < \gamma_{\tau+1}^{<\alpha_1} ~ \exists S \vartriangleleft \chi^{\ast +}
    ~ \alpha \mathbf{K}^{\ast <\alpha _{1}}(a,\delta,
    \gamma_{\tau }^{<\alpha _{1}}, \alpha ,\rho ,S)
    \bigr\};
\end{multline*}

\noindent \emph{2)}\quad мы будем называть матричной
автоэкзорцизивной (самоисключающейся в нарушении монотонности)
функцией или, короче, \ $\alpha $-функцией уровня \ $n$ \ ниже
$\alpha _{1}$ \ редуцированной к \ $\chi ^{\ast }$ \ следующую
функцию
\[
    \alpha S_{f}^{<\alpha _{1}}=(\alpha S_{\tau }^{<\alpha _{1}})_{\tau }
\]
принимающую значения: \vspace{-6pt}
\begin{multline*}
    \alpha S_{\tau }^{<\alpha _{1}}=
\\
    = \min_{\underline{\lessdot }}
    \bigl\{ S \vartriangleleft \chi^{\ast +}:
    \exists \delta ,\alpha ,\rho < \gamma_{\tau+1}^{<\alpha_1} ~  \alpha
    \mathbf{K}^{\ast <\alpha_{1}} (a_{\tau}^{<\alpha _{1}},\delta ,
    \gamma _{\tau }^{<\alpha _{1}}, \alpha, \rho,
    S) \bigr\};
\end{multline*}

\noindent \emph{3)} следующие сопрвождающие ординальные функции
определяются ниже~$\alpha _{1}$:

\hspace*{1.5em} диссеминаторная функция \hspace{\stretch{0.5000}}
$\widetilde{\delta }_{f}^{<\alpha _{1}}=(\widetilde{\delta }_{\tau
}^{<\alpha _{1}})_{\tau }$, \hspace*{1em}

\hspace*{1.5em} базовая функция \hspace{\stretch{0.5000}} $\rho
_{f}^{<\alpha _{1}}=(\rho _{\tau }^{<\alpha _{1}})_{\tau }$,
\hspace*{1em}

\hspace*{1.5em} несущая функция \hspace{\stretch{0.5000}} $\alpha
_{f}^{<\alpha _{1}}=(\alpha _{\tau }^{<\alpha _{1}})_{\tau }$,
\hspace*{1em}

\hspace*{1.5em} производящая диссеминаторная

\hspace*{1.5em} функция \hspace{\stretch{0.5} }
$\check{\delta}_{f}^{<\alpha _{1}}=(\check{\delta}_{\tau
}^{<\alpha _{1}})_{\tau }$, \hspace*{1em} \newline \quad \newline
принимающие значения для \ $a_{\tau} = a_{\tau }^{<\alpha _{1}}$,
$S_{\tau} = \alpha S_{\tau }^{<\alpha _{1}}$:
\begin{multline*}
    \widetilde{\delta }_{\tau }^{<\alpha _{1}} = \vspace{-12pt}
\\
    = \min_{\leq }
    \bigl \{ \delta < \gamma_{\tau}^{< \alpha_1}:
    \exists \alpha ,\rho <
    \gamma_{\tau+1}^{< \alpha_1} ~ \alpha
    \mathbf{K}^{\ast <\alpha_{1}}(a_{\tau},
    \delta, \gamma_{\tau }^{<\alpha_{1}},
    \alpha, \rho, S_{\tau}) \bigr \};
\end{multline*}

\vspace{-24pt}
\begin{multline*}
    \rho _{\tau }^{<\alpha _{1}} = \vspace{-12pt}
\\
    = \min_{\leq }
    \bigl \{ \rho \le \chi^{\ast +} :
    \exists \alpha < \gamma_{\tau+1}^{< \alpha_1}
    ~ \alpha \mathbf{K}^{\ast <\alpha_{1}}(a_{\tau},
    \widetilde{\delta }_{\tau}^{<\alpha _{1}},
    \gamma_{\tau}^{<\alpha_{1}},
    \alpha, \rho, S_{\tau})  \bigr \};
\end{multline*}

\vspace{-24pt}
\begin{multline*}
    \alpha _{\tau }^{<\alpha _{1}} = \vspace{-12pt}
\\
    = \min_{\leq }
    \bigl \{ \alpha < \gamma_{\tau+1}^{< \alpha_1} :
    \alpha \mathbf{K}^{\ast <\alpha _{1}}(a_{\tau},
    \widetilde{\delta }_{\tau}^{<\alpha _{1}},
    \gamma_{\tau}^{<\alpha_{1}},
    \alpha ,\rho_{\tau }^{<\alpha _{1}},
    S_{\tau})  \bigr \};
\end{multline*}

\noindent и для \ $\alpha^1 =
\alpha_{\tau}^{<\alpha_1^{\Downarrow}}$:

\vspace{-12pt}
\begin{multline*}
    \check{\delta}_{\tau }^{<\alpha _{1}}= \min_{\leq }
    \bigl \{ \delta < \gamma_{\tau}^{< \alpha_1} :
    SIN_{n}^{<\alpha ^{1}}(\delta )\wedge SIN_{n+1}^{<\alpha ^{1}}
    \left[ <\rho_{\tau }^{<\alpha _{1}}\right] (\delta )) \bigr \};
\end{multline*}
\vspace{0pt}

\noindent Значение \ $a_{\tau }^{<\alpha _{1}}$ \ называется,
напомним, характеристикой матрицы \ $\alpha S_{\tau }^{<\alpha
_{1}}$ \ на носителе \ $\alpha_{\tau}^{<\alpha_1}$ \ и самого
этого носителя ниже~$\alpha _{1}$.

\noindent Все функции
\[
    a_{f}^{<\alpha _{1}}, \quad
    \widetilde{\delta }_{f}^{<\alpha _{1}},\quad
    \check{\delta}_{f}^{<\alpha _{1}},\quad
    \rho_{f}^{<\alpha _{1}}
\]
называются, для некоторой краткости, атрибутами функций
\[
    \alpha_f^{<\alpha_1}, \quad \alpha S_f^{<\alpha_1},
\]
а их значения для индекса \ $\tau$ \ называются также атрибутами
значений
\[
    \alpha_{\tau}^{<\alpha_1}, \quad \alpha S_{\tau}^{<\alpha_1};
\]
аналогично функция \ $\alpha_f^{<\alpha_1}$ \ называется атрибутом
функции \ $\alpha S_f^{<\alpha_1}$, \ а её значение \
$\alpha_{\tau}^{<\alpha_1}$ \ -- атрибутом матрицы \ $\alpha
S_{\tau}^{<\alpha_1}$ \ ниже \ $\alpha_1$, \ и так далее.
\\
\hspace*{\fill} $\dashv$

\end{definition}

\noindent Понятие характеристики вводится в общем случае:

\begin{definition}
\label{8.4.} \hfill {} \newline \hspace*{1em} Мы называем
характеристикой матрицы \ $S$ \ на носителе \ $\alpha
>\chi ^{\ast }$ \ число \ $a(S,\alpha )=a$, \ определяемое следующим образом:
\[
    \hspace{0.4cm} \bigl(a=1\vee a=0 \bigr)\wedge \Big( a=0\longleftrightarrow
    \qquad\qquad\qquad\qquad\qquad\qquad\qquad\qquad
\]
\[
    \hspace{0.8cm} \longleftrightarrow \exists \tau_{1}^{\prime},
    \tau_{2}^{\prime},\tau _{3}^{\prime} < \alpha^{\Downarrow} \
    \big( A_{2}^{0 \vartriangleleft \alpha^{\Downarrow}}(\tau_{1}^{\prime},
    \tau_{2}^{\prime},\tau _{3}^{\prime},
    \alpha S_{f}^{<\alpha ^{\Downarrow }}) \wedge
\]
\[
    \qquad \qquad \qquad
    \wedge \forall \tau^{\prime\prime} \big( \tau_1^{\prime} <
    \tau^{\prime\prime} \le \tau_2^{\prime} \rightarrow
    a_{\tau^{\prime\prime}}^{<\alpha^{\Downarrow}}=1 \big) \wedge
    \alpha S_{\tau_2^{\prime}}^{<\alpha ^{\Downarrow }}=S \big) \Big).
\]

\noindent Матрица \ $S$ \ на её носителе \ $\alpha$ \ называется
единичной матрицей на \ $\alpha$, \ если она имеет единичную
характеристику на \ $\alpha$; \ иначе она называется нулевой
матрицей на \ $\alpha$. \hspace*{\fill} $\dashv$
\end{definition}

\noindent Таким образом, в ходе определения \ $\alpha $-функции \
$\alpha S_{f}^{<\alpha _{1}}$ \ приоритет отдаётся \
$\alpha$-матрицам, обладающим большей характеристикой.
\\
Это обстоятельство, хотя и позволяющее решить всю проблему больших
кардиналов, существенно усложняет теорию матричных функций в
целом, так как релятивизирующие рассуждения не работают теперь
вполне свободно: ситуации, связанные с нулевой характеристикой
могут не переноситься на нижние части универсума, например,
определённые единичной характеристикой, или по другим причинам,
связанным с подавлением.
\\
\quad \\

\begin{sloppypar}
Определение~\ref{8.3.}  \ $\alpha $-функции и сопровождающих её
ординальных функций следует рекурсивному определению~\ref{8.2.} и
так как функции \ \mbox{$X_i, \; i=\overline{1,5}$} \ определяются
этой рекурсией в формуле \ $\alpha \mathbf{K}_{n+1}^{\exists}$\ в
её подформуле \ $\ A_7^{RC}$\
 единственным
образом через её параметры, то нетрудно видеть, что функции   \
$X_i[\alpha^0]$, $i=\overline{1,5}$, \   в определении \ref{8.2.}
совпадают с соответствующими функциями
\end{sloppypar}
\begin{equation} \label{e8.2}
    a_{f }^{<\alpha^0},\quad
    \alpha S_{f }^{<\alpha^0},\quad
    \widetilde{\delta }_{f }^{<\alpha^0},\quad
    \rho_{f }^{<\alpha^0},\quad
    \alpha_{f }^{<\alpha^0},\quad
\end{equation}
для каждого кардинала \ $\alpha^0$ \ эквиинформативного с \
$\chi^{\ast}$. По этой причине мы будем использовать их
обозначения (\ref{e8.2}) вместо соответствующих обозначений этих
функций  \ $X_i[\alpha^0], ~ i=\overline{1,5}$ \ в формулах из
определения~\ref{8.1.}, то есть используя эти формулы, но для
функций \ \mbox{$X_i[\alpha^0], \; i=\overline{1,5}$}, \
заменённых на соответствующие функции (\ref{e8.2}) для \
$\alpha^0=\alpha_1$; \ мы даже будем опускать их для некоторой
краткости, когда это не будет вызывать недоразумений и когда
контекст будет очевидно на них указывать.

Например, формула \ $ A_0^{\vartriangleleft \alpha_1}(\tau_1,
\tau_2, \alpha S_f^{<\alpha_1})$ \ означает, что ниже \ $\alpha_1$
\ выполняется
\[
    \tau_1+1 < \tau_2 \wedge (\alpha S_f^{<\alpha_1} \mbox{\it \ это функция на
    } \left] \tau_1, \tau_2 \right[ ) \wedge
\]
\[
    \wedge \tau_1 = \min \bigl \{ \tau: \left] \tau, \tau_2 \right[ \subseteq
    dom(\alpha S_{f }^{<\alpha_1} )\wedge
\]
\[
    \wedge \chi^{\ast} \le  \gamma_{\tau_1}^{<\alpha_1} \wedge \gamma_{\tau_1}^{<\alpha_1}
    \in SIN_n^{<\alpha_1};
\]
формула \ $ A_1^{\vartriangleleft \alpha_1}(\tau_1, \tau_2, \alpha
S_f^{<\alpha_1})$ \ означает, что ниже \ $\alpha_1$ \ выполняется
\[
    A_0^{\vartriangleleft \alpha_1}(\tau_1, \tau_2, \alpha S_f^{< \alpha_1}) \wedge
    \gamma_{\tau_2}^{<\alpha_1} \in SIN_n^{<\alpha_1};
\]
формула
\[
    A_2^{0 \vartriangleleft \alpha^{\Downarrow}}(\tau_1, \tau_2, \tau_3,
    \alpha S_f^{<\alpha^{\Downarrow}}) \wedge
    \forall \tau \in \; ]\tau_1, \tau_2] ~
    a_{\tau}^{<\alpha^{\Downarrow}}=1
    \wedge \alpha S_{\tau_2}^{<\alpha^{\Downarrow}} = S
\]
означает, что здесь  \ $\alpha^{\Downarrow}$ \ -- это кардинал
предскачка носителя \ $\alpha $ \ после \ $\chi^{\ast}$, \ и что
не существует  \ $\alpha$-матрицы, допустимой для \
$\gamma_{\tau_1}^{<\alpha_1}$ \ ниже \ $\alpha^{\Downarrow}$, \ и
ниже того же \ $\alpha^{\Downarrow}$ \ выполняется \
\[
    A_1^{\vartriangleleft \alpha^{\Downarrow}} (\tau_1, \tau_3, \alpha
    S_f^{<\alpha^{\Downarrow}}),
\]
и \ $\tau_2 \in \left] \tau_1, \tau_3 \right[$ \ это первый
ординал, на котором нарушается монотонность на \ $]\tau_1,
\tau_3[$ \ матричной функции \ $\alpha S_f^{<\alpha^{\Downarrow}}$
\ (но уже ниже \ $\alpha^{\Downarrow}$ \ ) и, более того, \
$\alpha S_{\tau_2}^{<\alpha^{\Downarrow}} = S$ \ и все матрицы \
$\alpha S_{\tau}^{<\alpha^{\Downarrow}}$ \ обладают единичной
характеристикой на  \ $]\tau_1, \tau_2]$ \ -- \ и так далее.

\noindent Теперь необходимо сделать следующие два простых
замечания:

1. Все интервалы \ $[ \gamma_{\tau_1}^{<\alpha_1},
\gamma_{\tau_2}^{<\alpha_1} [$ \ определённости ниже \ $\alpha_1$,
\ рассмотренные в определении~\ref{8.1.} \ для функций
\[
    X_1 = \alpha S_f^{< \alpha_1}, \quad X_2 = a_f^{< \alpha_1},
\]
были различных видов и были определены различными условиями, но
все они включали условие максимальности интервала \
$[\gamma_{\tau_1}^{<\alpha_1}, \gamma_{\tau_2}^{<\alpha_1}[$ \
влево:
\[
    A_0^{\vartriangleleft \alpha_1}(\tau_1, \tau_2, \alpha
    S_f^{<\alpha_1}),
\]
которое утверждает, помимо прочего, что матричная функция \
$\alpha S_f^{<\alpha_1}$ \ ниже \ $\alpha_1$ \ определена на
интервале \ $\left] \tau_1, \tau_2 \right[$ \ и ординал \ $\tau_1$
\ является \textit{минимальным} с этим свойством и, более того, \
$\gamma_{\tau_1}^{<\alpha_1}$ \ является \
$SIN_n^{<\alpha_1}$-кардиналом. Благодаря этой минимальности
нетрудно видеть, что эта функция \ $\alpha S_f^{<\alpha_1}$ \  не
определена для самого этого ординала  \ $\tau_1$ \ !

2. Понятия допустимости, приоритетности и подавления следует
различать. Можно представить себе две матрицы \ $S^{\prime},
S^{\prime\prime}$ \ на их носителях \ $\alpha^{\prime},
\alpha^{\prime\prime}$ \ соответственно вместе с их
соответствующими атрибутами, обе допустимые для единого кардила \
$\gamma_{\tau}^{<\alpha_1}$; когда \ $S^{\prime}$ \ обладает
единичной характеристикой на \ $\alpha^{\prime}$ \ она всегда
неподавлена и имеет приоритет над \ $S^{\prime\prime}$ \ нулевой
характеристики на \ $\alpha^{\prime\prime}$. \ Но даже когда нет
такой матрицы \ $S^{\prime}$, \ всё-таки матрица \
$S^{\prime\prime}$ \ on \ $\alpha^{\prime\prime}$ \ может быть
подавлена, если выполняется условие подавления \ $A_5^{S,0}$ \
ниже \ $\alpha_1$; \ и в любом случае  каждая подавленная матрица
не может быть значением матричной функции \ $\alpha
S_f^{<\alpha_1}$.
\\
Таким образом, для интервала $[\gamma_{\tau_1}^{<\alpha_1},
\gamma_{\tau_2}^{<\alpha_1}[$, \ максимального влево ниже \
$\alpha_1$, \ не существует значения \ $\alpha
S_{\tau}^{<\alpha_1}$ \ для \ $\tau=\tau_1$, \ но тем не менее это
не исключает существования некоторой матрицы \textit{только
допустимой} (но подавленной) для
 \ $\gamma_{\tau_1}^{<\alpha_1}$ \ ниже \ $\alpha_1$.
\\
\quad \\

А теперь, имея эти замечания ввиду, рассмотрим, как
определение~\ref{8.2.} -- и, следовательно, определение~\ref{8.3.}
-- работает ниже \ $\alpha_1$ (мы рассматриваем, напомним, самый
важный случай, когда \ $\chi=\chi^{\ast}$, \ $\mathfrak{n} =
\mathfrak{n}^{\alpha}$).
\\

I. Итак, в его третьей части в начале формулы
\[
    \alpha \mathbf{K}_{n+1}^{\exists < \alpha_1}
    (a, \delta, \gamma, \alpha, \rho, S)
\]
утверждается, что \ $S$ \ -- это \ $\delta$-матрица на её носителе
 \ $\alpha > \chi^{\ast}$, $\alpha <
\alpha_1$, \ редуцированная к \ $\chi^{\ast}$ \ с диссеминатором
 \ $\delta<\gamma$ \ и базой \
$\rho$:
\[
    S \vartriangleleft \rho = \widehat{\rho} \le \chi^{\ast +};
\]
кардинал предскачка \ $\alpha^{\Downarrow} =
\alpha_{\chi^{\ast}}^{\Downarrow}$ \ является предельным для \
$SIN_n^{<\alpha^{\Downarrow}}$ \ и имеет конфинальность \ $\ge
\chi^{\ast +}$; \ диссеминатор \ $\delta$ \ обладает
субнедостижимостью ниже \ $\alpha^{\Downarrow}$ \ уровня \ $n$ \ и
даже уровня \ $n+1$ \ с базой \ $\rho$, \ то есть
\[
    \delta \in SIN_n^{<\alpha^{\Downarrow}} \cap
    SIN_{n+1}^{<\alpha^{\Downarrow}}[<\rho].
\]

II. Затем ниже \ $\alpha^{\Downarrow}$ \ определяются функции
 \ $X_i, \ i = \overline{1,5}$ \ на
парах \ $(\alpha^0, \tau^{\prime}) \in
\mathcal{A}_{\chi^{\ast}}^{\alpha^{\Downarrow}}$, \ где кардиналы
 \ $\alpha^0 \in \left] \chi^{\ast},
\alpha^{\Downarrow} \right]$ \ эквиинформативны с \ $\chi^{\ast}$
\ и существуют кардиналы  \ $\gamma_{\tau^{\prime}}^{<\alpha^0}$.
\\
Все эти функции рекурсивно определяются через определение функций
 \ $X_i^0, \ i = \overline{1,5}$, \ посредством
условия рекурсии \ $A_7^{RC}$ из второй части определения 8.2:
\[
    X_1[\alpha^0] = \alpha S_f^{<\alpha^0}, \quad
    X_2[\alpha^0] = a_f^{<\alpha^0}, \quad
    X_3[\alpha^0] = \widetilde{\delta}_f^{<\alpha^0},
\]
\[
    X_4[\alpha^0] = \rho_f^{<\alpha^0}, \quad
    X_5[\alpha^0] = \alpha_f^{<\alpha^0}.
\]
Цель этого определения -- получить результирующую матричную
функцию \ $\alpha S_f^{<\alpha^0}$, \ но первой определяется
именно характеристическая функция
\[
    X_2^0 = X_2[\alpha^0] = a_f^{<\alpha^0}.
\]
Эта функция принимает \textit{максимальные} возможные значения,
единичное или нулевое, которые являются характеристиками матриц,
допустимых ниже \ $\alpha^0$, \ \textit{но только не нулевую
характеристику подавленных нулевых матриц} \ $S^{\prime\prime}$ \
на их носителях \ $\alpha^{\prime\prime}$, \ которые удовлетворяют
условию своего подавления ниже \ $\alpha^0$:
\[
    A_5^{S,0 \vartriangleleft \alpha^0}(0, \gamma_{\tau^{\prime}}^{<\alpha^0}, \alpha^{\prime\prime},
    \rho^{\prime\prime}, S^{\prime\prime}, X_1^0|\tau^{\prime},
    X_2^0|\tau^{\prime}, X_1|^1\alpha^0, X_2|^1\alpha^0),
\]
где функции здесь
\[
    X_1^0|\tau^{\prime}, \quad X_2^0|\tau^{\prime}, \quad
    X_1|^1\alpha^0, \quad X_2|^1\alpha^0
\]
уже определены. И везде в дальнейшем такие подавленные нулевые
матрицы систематически отвергаются.
\\
После того, как характеристическая функция \ $X_2[\alpha^0] =
a_f^{<\alpha_0}$ \ определена, все оставшиеся функции
\[
    X_1[\alpha^0], \quad X_i[\alpha^0], \quad i=\overline{3,5}
\]
определяются последовательно одна за другой посредством
минимизации их допустимых и неподавленных значений.
\\
Таким образом, следующей по очереди определяется матричная функция
 \ $X_1[\alpha^0] =
\alpha S_f^{<\alpha^0}$, \ после этого соответствующая
диссеминаторная функция \ $X_3[\alpha^0] =
\widetilde{\delta}_f^{<\alpha^0}$, \ затем базовая функция \
$X_4[\alpha^0] = \rho_f^{<\alpha^0}$, \ и в последнюю очередь
определяется несущая функция \ $X_5[\alpha^0] =
\alpha_f^{<\alpha^0}$ .\
\\
При этом значения каждой из последующих из этих функций
существенно зависят от значений предыдущих функций.

III. После того, как эти функции сформированы для всякого
\[
    \alpha^0 \in ] \chi^{\ast}, \alpha^{\Downarrow}[\;,
\]
это определение переходит к кардиналу
\[
    \alpha^0 = \alpha^{\Downarrow}
\]
и здесь определяет характеристику самой матрицы  \ $S$ \ на её
носителе  \ $\alpha$:
\\
матрица $S$ \ на \ $\alpha$ \ получает \textit{нулевую
характеристику}, когда она участвует в следующем нарушении
монотонности матричной функции
\[
    X_1[\alpha^{\Downarrow}] = \alpha S_f^{<\alpha^{\Downarrow}}
\]
ниже \ $\alpha^{\Downarrow}$: \ когда ниже \ $\alpha^{\Downarrow}$
\ выполняется условие

\begin{multline*}
    \exists \tau_1^{\prime}, \tau_2^{\prime}, \tau_3^{\prime} <
    \alpha^{\Downarrow} \big(A_2^{0 \vartriangleleft \alpha^{\Downarrow}}
    (\tau_1^{\prime}, \tau_2^{\prime}, \tau_3^{\prime},
    \alpha S_f^{<\alpha^{\Downarrow}} ) \wedge
\\
    \wedge \forall \tau^{\prime\prime} \in \; ]\tau_1^{\prime}, \tau_2^{\prime}] ~
    a_{\tau^{\prime\prime}}^{<\alpha^{\Downarrow}}=1 \wedge
    \alpha S_{\tau_2^{\prime}}^{<\alpha^{\Downarrow}} = S \big);
\end{multline*}
иначе \ $S$ \ на \ $\alpha$ \ получает \textit{единичную
характеристику}.

IV. И в последнюю очередь это определение формирует
\textit{замыкающее условие} для \ $S$ \ на \ $\alpha$:
\\
Если \ $S$ \ -- это \textit{нулевая} матрица на \ $\alpha$ \ и её
допустимый диссеминатор \ $\delta$ \ попадает в некоторый
масимальный блок типа \ $\eta^{\prime}$ \ ниже \
$\alpha^{\Downarrow}$
\[
    [ \gamma_{\tau_1^{\prime}}^{<\alpha^{\Downarrow}},
    \gamma_{\tau_3^{\prime}}^{<\alpha^{\Downarrow}} [
\]
{\sl обременительный} для \ $S$ \ на \ $\alpha$, \ то есть если
выполняется
\[
    \gamma_{\tau_1^{\prime}}^{<\alpha^{\Downarrow}} \le \delta <
    \gamma_{\tau_3^{\prime}}^{<\alpha^{\Downarrow}} \wedge
    A_4^{M b \vartriangleleft \alpha^{\Downarrow}}
    (\tau_1^{\prime}, \tau_1^{\prime\prime}, \tau_2^{\prime},
    \tau_3^{\prime}, \eta^{\prime}, \alpha S_f^{<\alpha^{\Downarrow}},
    a_f^{<\alpha^{\Downarrow}} )
\]
ниже \ $\alpha^{\Downarrow}$, \ тогда требуется допустимая база
данных \ $\rho$ \ диссеминатора \ $\delta$ \ матрицы \ $S$ \ на \
$\alpha$, \ но только такая, что
\[
    \eta^{\prime} < \rho \vee \rho = \chi^{\ast +}.
\]
Таким образом, подобный случай \textit{затрудняет} существенно
использование такой матрицы \ $S$ \ на \ $\alpha$; \ кроме того, \
$S$ \ на \ $\alpha$ \ должна быть неподавлена; во всех других
случаях никаких требований на матрицу \ $S$ \ на \ $\alpha$ не
накладывается.
\\
Но напомним, что база \ $\rho = \chi^{\ast +}$ \ и каждая
\textit{единичная} матрица всегда допустимы и неподавлены; каждая
матрица неподавлена для \ $\gamma \notin SIN_n$ \ в любом случае.

После этого определение 8.2 в четвёртой части формирует конъюнкцию
 \ $\alpha
\mathbf{K}^{<\alpha_1}$:
\[
    \mathbf{K}_n^{\forall < \alpha_1} (\gamma,
    \alpha^{\Downarrow}) \wedge \alpha
    \mathbf{K}_{n+1}^{\exists \vartriangleleft \alpha_1}
    (a, \delta, \gamma, \alpha, \rho, S)
    \wedge \alpha < \alpha_1
\]
где дополнительно требуется, как обычно, что \
$\alpha^{\Downarrow}$ \ сохраняет все \
$SIN_n^{<\alpha_1}$-кардиналы \ $\le \gamma$ \ ниже \ $\alpha_1$;
\ и, наконец, возникает формула \ $\alpha \mathbf{K}^{\ast
<\alpha_1}$, \ получающаяся из формулы \ $\alpha
\mathbf{K}^{<\alpha_1}$ \ присоединением требования неподавления
нулевой матрицы \ $S$ \ на\ $\alpha$ \ ниже\ $\alpha_1$.
\\

Так как определение \ref{8.3.} матричной \ $\alpha$-функции и
сопрвождающих функций следует определению \ref{8.2.}, то
выполняется следующая очевидная лемма, которая в действительности
повторяет это определение 8.2. Здесь используется понятие
собственного порождающего диссеминатора \ $\check{\delta}^S$ \ для
произвольной матрицы \ $S$ \ на носителе \ $\alpha$, \ который,
напомним, является минимальным допустимым диссеминатором для \ $S$
\ на \ $\alpha$ \ с минимальной допустимой базой \ $\rho^S =
\widehat{\rho_1}$, $\rho_1 = Od(S)$ \ (см. \cite{Kiselev11},
\cite{Kiselev8}).

\begin{lemma}
\label{8.5.} \hfill {} \newline \hspace*{1em} Пусть \ $S$ \ --
произвольная \ $\alpha$-матрица, редуцированная к  \ $\chi ^{\ast
}$ \ и характеристики \ $a$ \ на носителе \ $\alpha < \alpha_1$,
{\sl допустимая} для \ $\gamma _{\tau }^{<\alpha _{1}}$ \ вместе с
её диссеминатором \ $\widetilde{ \delta }$, \ производящим
диссеминатором \ $\check{\delta}$ \ с базой \ $\rho $ и
порождающим собственным диссеминатором \ $\check{\delta}^S$ \ ниже
 \ $\alpha _{1}$, \ тогда для кардинала предскачка \ $\alpha ^{\Downarrow }$ \ после \
$\chi^{\ast}$ \ ниже \ $\alpha_1$ выполняется:
\newline

\noindent \medskip \emph{1)}\quad $\forall \gamma \leq \gamma
_{\tau }^{<\alpha _{1}}(SIN_{n}^{<\alpha _{1}}(\gamma
)\longrightarrow SIN_{n}^{<\alpha ^{\Downarrow }}(\gamma ))$~;
\newline

\noindent \medskip \emph{2)}\quad $\chi ^{\ast }<\widetilde{\delta }<\gamma
_{\tau }^{<\alpha _{1}}<\alpha ^{\Downarrow }\wedge S\vartriangleleft \rho
\leq \chi ^{\ast +}\wedge \rho =\widehat{\rho }$~; \newline

\noindent \medskip \emph{3)}\quad $\widetilde{\delta }\in
SIN_{n}^{<\alpha ^{\Downarrow }}\cap SIN_{n+1}^{<\alpha
^{\Downarrow }}\left[ <\rho \right] $; \ аналогично для
 \ $\check{\delta}$; \newline

\noindent \medskip \emph{4)}\quad $\sup SIN_{n}^{<\alpha
^{\Downarrow }}=\alpha ^{\Downarrow}\wedge cf(\alpha ^{\Downarrow
})\geq \chi ^{\ast +}$; \newline

\noindent \medskip \emph{5)}\quad $a=0\longleftrightarrow \exists
\tau_{1}^{\prime}, \tau_{2}^{\prime}, \tau_{3}^{\prime}
~\big(A_{2}^{0 \vartriangleleft
\alpha^{\Downarrow}}(\tau_{1}^{\prime}, \tau_{2}^{\prime},
\tau_{3}^{\prime}, \alpha S_f^{<\alpha^{\Downarrow}}) \wedge $
\newline

\noindent \medskip \hspace{10em} $ \wedge \forall
\tau^{\prime\prime} \in \; ]\tau_1^{\prime}, \tau_2^{\prime}] ~
a_{\tau^{\prime\prime}}^{<\alpha^{\Downarrow}}=1 \wedge \alpha
S_{\tau_{2}^{\prime}}^{<\alpha ^{\Downarrow }}=S\big)$;
\newline

\noindent \medskip \emph{6)}\quad $a=0\longrightarrow \forall
\tau_{1}^{\prime}, \tau_{1}^{\prime\prime}, \tau_{2}^{\prime},
\tau_{3}^{\prime}, \eta^{\prime} \bigl[ \gamma_{\tau_1^{\prime}}^{<\alpha
^{\Downarrow }}\leq \widetilde{\delta} <
\gamma_{\tau_3^{\prime}}^{<\alpha ^{\Downarrow }}\wedge $
\newline

\noindent \medskip \hspace{1em} $\wedge A_{4}^{M b
\vartriangleleft \alpha^{\Downarrow} } (\tau_{1}^{\prime},
\tau_{1}^{\prime\prime}, \tau_{2}^{\prime}, \tau_{3}^{\prime},
\eta^{\prime}, \alpha S_f^{<\alpha^{\Downarrow}},
a_f^{<\alpha^{\Downarrow}} ) \longrightarrow \eta^{\prime} <\rho \vee \rho
=\chi ^{\ast +} \bigr]$;
\newline

\noindent \medskip \emph{7)}\quad (i) \ $\check{\delta}^S \le
\check{\delta}\leq \widetilde{\delta }<\gamma _{\tau }^{<\alpha
_{1}}$; \newline

\noindent \medskip \qquad (ii) \ если \ $\widetilde{\delta }$ \
это {\sl минимальный} плавающий диссеминатор матрицы \ $S$ \ на \
$\alpha $ \ с минимальной базой \ $\rho$, \ допустимые для
 \ $\gamma _{\tau }^{<\alpha _{1}}$ \ вместе с \ $\rho $, \ тогда:
\[
    a=1\longrightarrow \widetilde{\delta } = \check{\delta}^S
    \wedge \rho = \rho^S = \widehat{\rho_1}, \wedge \rho_1 = Od(S),
\]
то есть когда \ $S$ \ это единичная матрица на \ $\alpha$, \ тогда
\ $\widetilde{\delta}$ \ это производящий собственный диссеминатор
\ $\check{\delta}^S$ \ матрицы \ $S$ \ на \ $\alpha$ \ с базой\
$\rho^S$;
\\

\noindent \medskip \emph{8)}\quad существует минимальный носитель
 \ $\alpha^{\prime} <
\gamma_{\tau+1}^{<\alpha_1}$ \ матрицы \ $S$ \ той же
характеристики \ $a$, \ допустимые для \
$\gamma_{\tau}^{<\alpha_1}$ \ вместе с теми же своими атрибутами
 \ $\widetilde{\delta}$, $\rho$ \
ниже \ $\alpha_1$:
\[
    \gamma _{\tau }^{<\alpha _{1}}<\alpha^{\prime} <
    \gamma _{\tau +1}^{<\alpha_{1}}~;
\]

\noindent аналогично для неподавленности \ $S$ \ для \
$\gamma_{\tau}$ \ вместе с теми же атрибутами.
\end{lemma}

\noindent \textit{Доказательство.} \ Осталось доказать последние
два утверждения; верхний индекс \ $< \alpha_1$ \ часто будет
опускаться.
\\
Итак, рассмотрим матрицу \ $S$ \ характеристики \ $a$ \ на её
носителе \ $\alpha < \alpha_1$, \ допустимую для \
$\gamma_{\tau}^{<\alpha_1}$ \ вместе с её диссеминатором \
$\widetilde{\delta}$ \ и базой \ $\rho$. \ Утверждение 7)~$(i)$
очевидно; переходя к  7)~$(ii)$ предположим \ $a=1$, \ тогда база

\[
    \rho = \rho^S = \widehat{\rho_1}, \ \rho_1 = Od(S)
\]
вместе с минимальным диссеминатором
\[
    \check{\delta}^{S} \in SIN_n^{<\alpha^{\Downarrow}} \cap
    SIN_{n+1}^{<\alpha^{\Downarrow}}[<\rho^S]
\]
очевидно выполняет все требования условия
\[
    \alpha \mathbf{K}(a, \check{\delta}^S, \gamma, \alpha, \rho^S, S)
\]
вплоть до последнего её конъюнктивного члена  \ $\mathbf{K}^0$.
\\
Но последнее выполняется тоже, так как для
 \ $a=1$ \ его посылка нарушается.
\\
Таким образом, вся формула \ $\alpha \mathbf{K}$ \ выполняется и \
$\widetilde{\delta} = \check{\delta}^S$, $\rho=\rho^S$.

\noindent Обращаясь к доказательству 8) легко применить лемму
 3.2~\cite{Kiselev11}  об ограничении, как это было сделано в доказательстве леммы
 5.17~2)~$(ii)$ . Однако подобное применение
составляет типичное рассуждение, которое будет использоваться
далее в различных важных случаях, поэтому следует представить его
себе в деталях.
\\
Во-первых, выше преполагается, что \ $\alpha_1$ \ -- предельный
кардинал для класса \ $SIN_{n-1}^{<\alpha_1}$ (напомним также
соглашение после (\ref{e7.1})), поэтому \
$\gamma_{\tau+1}^{<\alpha_1}$ \ всегда существует для каждого \
$\gamma_{\tau}^{<\alpha_1}$.
\\
Далее, предположим, что матрица \ $S$ \ с диссеминатором \
$\delta$ \ и базой \ $\rho$ \ на носителе
\[
    \alpha \in \; ]\gamma_\tau^{<\alpha_1}, \alpha_1 [
\]
допустимы для \ $\gamma_{\tau}^{<\alpha_1}$ \ ниже \ $\alpha_1$, \
тогда выполняется следующее утверждение \ $\varphi(\chi^{\ast},
\delta, \gamma_{\tau}^{<\alpha_1}, \rho, S)$:
\[
    \exists \alpha^{\prime} ~ ( \gamma_{\tau}^{<\alpha_1} < \alpha^{\prime}
    \wedge \alpha \mathbf{K}(\delta, \gamma_{\tau}^{<\alpha_1},
    \alpha^{\prime}, \rho, S))
\]
ниже \ $\alpha_1$, \ то есть после его  \
$\vartriangleleft$-ограничения кардиналом \ $\alpha_1$. \ Само это
утверждение \ $\varphi$ \ содержится в классе \ $\Sigma_n$, \ так
как оно включает в себя \ $\Sigma_n$-формулу\
$\mathbf{K}_n^{\forall}$. \ Но рассмотрим кардинал
\[
    \gamma_{\tau^n}=\sup\left\{\gamma \le \gamma_\tau^{<\alpha_1} :
    \gamma \in SIN_n^{<\alpha_1} \right\};
\]
по лемме~3.4~\cite{Kiselev11} \ $\gamma_{\tau^n}$ \ также
принадлежит классу \ $SIN_n^{<\alpha_1}$. \ Теперь заменим в
формуле \ $\alpha \mathbf{K}$ \ её подформулу \
$\mathbf{K}_n^{\forall}$ \ \ $\Delta_1$-формулой
\[
    SIN_n^{<\alpha^{\Downarrow}}(\gamma_{\tau^n}),
\]
тогда \ $\Sigma_n$-формула \ $\alpha \mathbf{K}$ \ преобразуется в
некоторую \ $\Pi_{n-2}$-формулу, которую обозначим через \ $\alpha
\mathbf{K}_{n-2}$. \ Соответственно, формула \ $\varphi$ \
преобразуется в некоторую\ $\Sigma_{n-1}$-формулу \
$\varphi_{n-2}(\chi^{\ast}, \delta, \gamma_{\tau^n},
\gamma_{\tau}^{<\alpha_1}, \rho, S)$:

\[
    \exists \alpha^{\prime} ~ (\gamma_{\tau}^{<\alpha_1} < \alpha^{\prime} \wedge
    \alpha \mathbf{K}_{n-2}(\delta, \gamma_{\tau^n}, \gamma_{\tau}^{<\alpha_1},
    \alpha^{\prime}, \rho, S))
\]
\\
в точности того же содержания ниже \ $\alpha_1$, \ и поэтому
выполняется
\[
    \varphi_{n-2}^{\vartriangleleft \alpha_1}(\chi^{\ast}, \delta,
    \gamma_{\tau^n}, \gamma_{\tau}^{<\alpha_1}, \rho, S).
\]
Последнее предложение содержит индивидные константы
\[
    \chi^{\ast}, \delta, \gamma_{\tau^n}, \gamma_{\tau}^{<\alpha_1}, \rho, S
\]
строго меньшие \ $SIN_{n-1}^{<\alpha_1}$-кардинала \
$\gamma_{\tau+1}^{<\alpha_1}$ \ и, следовательно, этот кардинал
ограничивает это предложение по лемме 3.2~\cite{Kiselev11} (где
 \ $n$ \ заменено на \ $n-1$), то есть
выполняется утверждение

\[
    \exists \alpha^{\prime} \in [\gamma_{\tau}^{<\alpha_1},
    \gamma_{\tau+1}^{<\alpha_1} [ \quad
    \alpha \mathbf{K}_{n-2}^{<\alpha_1}(\delta,
    \gamma_{\tau^n}, \gamma_{\tau}^{<\alpha_1},
    \alpha^{\prime}, \rho, S)
\]
\\
и \ $S$ \ получает свой носитель \ $\alpha^{\prime} \in
[\gamma_{\tau}^{<\alpha_1}, \gamma_{\tau+1}^{<\alpha_1}[$, \
допустимый для \ $\gamma_{\tau}^{<\alpha_1}$ \ вместе с теми же
диссеминатором и базой данных.
\\
Утверждение 8) о неподавленности не будет использоваться до \S 11
и там мы вернёмся к нему ещё раз.
\\
\hspace*{\fill} $\dashv$
\\
Нетрудно видеть, что введённые в определении \ref{8.3.} функции
обладают многими простыми свойствами \ $\delta $-функций и её
сопровождающих функций, поэтому доказательства следующих трёх лемм
вполне аналогичны доказательствам лемм \ref{7.3.}, \ref{7.4.} (или
лемм~5.16, 5.15 \cite{Kiselev11}) и леммы \ref{7.5.}:

\begin{lemma}
\label{8.6.} \hfill {} \newline \hspace*{1em} Для \ $\alpha
_{1}<k$ \ формулы \ $\alpha \mathbf{K}^{<\alpha _{1}}$, \ $\alpha
\mathbf{K}^{\ast <\alpha _{1}}$ \ принадлежат \ $\Delta _{1}$ \ и
поэтому все функции из определения~\ref{8.3.}:
\begin{equation*}
a_{f}^{<\alpha _{1}}, ~~\alpha S_{f}^{<\alpha _{1}},
~~\widetilde{\delta } _{f}^{<\alpha _{1}}, ~~\rho _{f}^{<\alpha
_{1}}, ~~\alpha _{f}^{<\alpha _{1}}, ~~\check{\delta}_{f}^{<\alpha
_{1}}
\end{equation*}
$\Delta _{1}$-определимы через \ $\chi ^{\ast },\alpha _{1}$. \
Для \ $\alpha _{1}=k$ \ формула \ $\alpha \mathbf{K}$ \
принадлежит \ $\Sigma _{n+1}$.
\\
\hspace*{\fill} $\dashv$
\end{lemma}

\begin{lemma}
\label{8.7.} \emph{(Об абсолютности \ $\alpha $-функций)}
\newline \hspace*{1em} Пусть \ $\chi ^{\ast }<\gamma _{\tau
+1}^{<\alpha _{1}}<\alpha _{2}<\alpha _{1}\leq k$, \quad $\alpha
_{2}\in SIN_{n-2}^{<\alpha _{1}}$ \ и
\begin{equation*}
(\gamma _{\tau }^{<\alpha _{1}}+1)\cap SIN_{n}^{<\alpha
_{2}}=(\gamma _{\tau }^{<\alpha _{1}}+1)\cap SIN_{n}^{<\alpha
_{1}},
\end{equation*}

\noindent \medskip \emph{1)}\quad тогда на множестве
\[
    T = \{\tau^{\prime }:~\chi ^{\ast }\leq
    \gamma _{\tau ^{\prime }}^{<\alpha_{2}}\leq
    \gamma _{\tau }^{<\alpha _{1}}\}
\]
допустимость ниже \ $\alpha_2$ \ равносильна допустимости ниже \
$\alpha_1$: \ для всякого \ $\tau^{\prime} \in T$ \ и матрицы \
$S^{\prime}$ \ на её носителе \ \mbox{$\alpha^{\prime} \in \; ]
\gamma_{\tau^{\prime}}^{<\alpha_2},
\gamma_{\tau^{\prime}+1}^{<\alpha_1}[$}
\[
    \alpha \mathbf{K}^{<\alpha_2}(
    \gamma_{\tau^{\prime}}^{<\alpha_2}, \alpha^{\prime},
    S^{\prime}) \leftrightarrow
    \alpha \mathbf{K}^{<\alpha_1}(
    \gamma_{\tau^{\prime}}^{<\alpha_2}, \alpha^{\prime},
    S^{\prime});
\]
\emph{2)}\quad на множестве
\[
    \big \{ \tau^{\prime}: \chi^{\ast} \le \gamma_{\tau^{\prime}}^{<\alpha_2}
    \le \gamma_{\tau}^{<\alpha_1} \wedge ( a_{\tau^{\prime}}^{<\alpha_2} = 1 \vee
    \neg SIN_n^{<\alpha_2}(\gamma_{\tau^{\prime}}^{<\alpha_2}) ) \big \}
\]
функции (\ref{e8.2}) ниже \ $\alpha^0 = \alpha_{2}$ \ тождественно
совпадают соответственно с функциями (\ref{e8.2}) ниже \ $\alpha^0
= \alpha _{1}$.
\\
\hspace*{\fill} $\dashv$
\end{lemma}

\begin{lemma}
\label{8.8.} \emph{(О диссеминаторе)} \newline \emph{1)}\quad
Пусть

\quad \medskip (i) \quad $] \tau_{1}, \tau_{2}[\; \subseteq dom
\bigl( \alpha S_{f}^{< \alpha_{1}} \bigr ) $, \quad
$\gamma_{\tau_{2}} \in SIN_{n}^{< \alpha_{1}}$;

\quad \medskip (ii) \quad $\tau_{3} \in dom \bigl ( \alpha
S_{f}^{< \alpha_{1}} \bigr ),\quad \tau_{2} \leq \tau_{3}$;

\quad \medskip (iii) \quad $\widetilde{\delta}_{\tau_{3}}^{<
\alpha_{1}} < \gamma_{\tau_{2}}^{< \alpha_{1}} $ \ и \
$a_{\tau_{3}}^{< \alpha_{1}}=0$.
\newline
Тогда
\begin{equation*}
\widetilde{\delta}_{\tau_{3}}^{ < \alpha_{1}} \leq
\gamma_{\tau_{1}}^{ < \alpha_{1}}.
\end{equation*}
Аналогично для \ $\check{\delta}_{\tau_{3}}^{ < \alpha_{1}}$.
\\
\quad \\
\emph{2)}\quad Пусть \ $\alpha$-матрица \ $S$ \ характеристики \
$a$ \ на носителе \ $\alpha$ \ допустима для \
$\gamma_{\tau}^{<\alpha_1}$ \ вместе с её диссеминатором \
$\widetilde{\delta}$ \ и базой \ $\rho$ \ ниже \ $\alpha_1$, \
тогда

\[
    \{ \tau^\prime : \widetilde{\delta}  <
    \gamma_{\tau^\prime}^{ < \alpha_{1}} < \gamma_{\tau}^{ <
    \alpha_{1}} \} \subseteq dom \bigl ( \alpha S_{f}^{< \alpha_{1}}
    \bigr ).
\]
\end{lemma}

\noindent \textit{Доказательство} \ 1) Предстоящее рассуждение
аналогично доказательству леммы 7.5~1), но теперь некоторые
специальные свойства диссеминаторов матриц единичной или нулевой
характеристики вызывают особые обстоятельства. Поэтому здесь
следует использовать следующий аргумент, который будет
использоваться в дальнейшем в различных типичных ситуациях;
верхние индексы \ $< \alpha_{1} $, $\vartriangleleft \alpha_{1}$ \
будут опускаться для краткости.
\\
Предположим, что 1) нарушается; рассмотрим матрицу \ $S^3 = \alpha
S_{\tau_{3}} $ \ характеристики \ $a^3 = a_{\tau_{3}}=0 $ \ на
носителе \ $\alpha_{\tau_3}$ \ с кардиналом предскачка \ $\alpha^3
= \alpha_{\tau_3}^{\Downarrow}$, \ обладающую диссеминаторами \
$\check{\delta}^3 = \check{\delta}_{\tau_{3}}$, \
$\widetilde{\delta}^3 = \widetilde{\delta}_{\tau_{3}}$ \ с базой
 \ $\rho^3 = \rho_{\tau_{3}}$, \ и предположим что
\\
\begin{equation} \label{e8.3}
\gamma_{\tau_1} < \widetilde{\delta}^3 < \gamma_{\tau_2}.
\end{equation}
Здесь следует иметь ввиду минимальный ординал \ $\tau_1$, \
выполняющий условие $(i)$.
\\
По определению \ref{8.3.} \ выполняется утвеждение \ $\alpha
\mathbf{K} $ \ и поэтому выполняется утвеждение \
$\mathbf{K}^{0}$:

\begin{equation*}
\begin{array}{l}
a^3 = 0 \longrightarrow \forall \tau_1^{\prime},
\tau_1^{\prime\prime}, \tau_2^{\prime}, \tau_3^{\prime},
\eta^{\prime} < \alpha^3 ~ \bigl[ \gamma
_{\tau_1^\prime}^{<\alpha^3} \leq
\widetilde{\delta}^3 < \gamma_{\tau_3^\prime}^{<\alpha^3} \wedge \\
\quad \\
~ \wedge A_{4}^{M b \vartriangleleft \alpha^3} ( \tau_1^{\prime},
\tau_1^{\prime\prime}, \tau_2^{\prime}, \tau_3^{\prime},
\eta^{\prime}, \alpha S_f^{<\alpha^3}, a_f^{<\alpha^3})
\longrightarrow \eta^{\prime} < \rho^3 \vee \rho^3 = \chi^{\ast +}
\bigr].
\end{array}
\end{equation*}

\noindent Предположим, что существуют ординалы \ $\tau_1^{\prime},
\tau_1^{\prime\prime}, \tau_2^{\prime}, \tau_3^{\prime},
\eta^{\prime} < \alpha^3$,\
 выполняющие посылку
этого утверждения: \vspace{6pt}
\begin{equation} \label{e8.4}
\gamma_{\tau_{1}^{\prime}}^{< \alpha^{3}} \leq
\widetilde{\delta}^3 < \gamma_{\tau_{3}^{\prime}}^{< \alpha^{3}}
\wedge A_{4}^{M b \vartriangleleft \alpha^3} ( \tau_1^{\prime},
\tau_1^{\prime\prime}, \tau_2^{\prime}, \tau_3^{\prime},
\eta^{\prime}, \alpha S_{f}^{< \alpha^{3}}, a_{f}^{< \alpha^{3}}
)~.
\end{equation}
\vspace{0pt}

\noindent Следует отметить снова, что благодаря \ $A_{4}^{M b
\vartriangleleft \alpha^3} $ \ эти ординалы определяются через\
$\widetilde{\delta}^3, \alpha^{3} $ \ \textit{единственным
образом}. Так как \ $\gamma_{\tau_{2}} \in SIN_{n} $ \ и \
$\gamma_{\tau_{1}}$ \ минимален, то можно видеть, что из
предположения (\ref{e8.3}) следует
\begin{equation} \label{e8.5}
    \gamma_{\tau_{1}^{\prime}} < \check{\delta}^3 =
    \widetilde{\delta}^3 < \gamma_{\tau_{2}}
\end{equation}
как результат минимизации диссеминатора \
$\widetilde{\delta}_{\tau_{3}} $ \ в интервале \ $[
\gamma_{\tau_{1}^{\prime}}, \gamma_{\tau_{2}} [ $ \ согласно
определению~8.3. Теперь мы приходим к ситуации из доказательства
леммы~7.5~1) и остаётся повторить его аргументы, то есть
использовать \ $\underline{\lessdot}$-минимальную матрицу \ $S^m
\lessdot S^3$ \ на некотором носителе \ $\alpha^m \in ]
\gamma_{\tau_3}, \alpha^3 [$ \ характеристики \ $a^m$, \
допустимую и неподавленную для \ $\gamma_{\tau_3}$ \ вместе со
своим минимальным диссеминатором \ $\widetilde{\delta}^m <
\gamma_{\tau_2}$ \ и базой \ $\rho^m < Od(S^3)$, \ потому что
подавленность матрицы \ $S^m$ \ для \ $\gamma_{\tau_3}$ \ влечёт
подавленность самой матрицы \ $S^3$ \ для \ $\gamma_{\tau_3}$, \
хотя она неподавлена по определению (ниже \ $\alpha_1$).
\\
Это вызывает противоречие: так как \ $S^m \lessdot S^3$ \ и \
$a^3=0$, \ то по определению \ref{8.3.} матрица \ $S^3$ \ не может
быть минимальным значением \ $\alpha S_{\tau_3}$.

Если же таких ординалов \ $\tau_1^{\prime}, \tau_1^{\prime\prime},
\tau_2^{\prime}, \tau_3^{\prime}, \eta^{\prime}$ \ нет, то
утверждение \ $\mathbf{K}^{0} $ \ очевидно сохраняется при
минимизации диссеминатора \ $\widetilde{\delta}_{\tau_{3}} $ \ в \
$[ \gamma_{\tau_{1}}, \gamma_{\tau_{2}} [ $ \ и тогда\
$\widetilde{\delta}_{\tau_{3}} \leq \gamma_{\tau_{1}} $, \ иначе
снова выполняется \ $\gamma_{\tau_{1}} < \check{\delta}_{\tau_{3}}
= \widetilde{\delta}_{\tau_{3}} < \gamma_{\tau_{2}} $ \ и прежнее
рассуждение вызывает прежнее противоречие.

Обращаясь к утверждению 2 следует просто заметить, что это
утверждение повторяет предшествующие леммы~5.17~2)
\cite{Kiselev11}, 7.5~2) в следующей форме:
\\
матрица \ $S$, \ допустимая для \ $\gamma_{\tau}$ \ на её носителе
\ $\alpha_\tau$, \ по лемме \ref{8.5.}~8) и определению \ref{8.2.}
остаётся допустимой и неподавленной для всякого \
$\gamma_{\tau^\prime} < \gamma_{\tau}$, \ такого, что \
$\widetilde{\delta} < \gamma_{\tau^\prime}$, \ вместе с теми же
сопровождающими ординалами \ $a$, $\widetilde{\delta}$, $\rho$,
$\alpha$. \ Для единичной характеристики \ $a=1$ \ это очевидно;
для \ $a=0$ \ эта лемма будет использоваться только в \S 11 и там
мы вернёмся к её доказательству, изложенному более подробно.
\hspace*{\fill} $\dashv$
\\
\quad \\

Следующие леммы подтверждают дальнейшее распространение теории \
\mbox{$\delta$-функций} на \ $\alpha$-функции и аналогичны
леммам~7.6, 7.7 об определённости \ $\delta$-фукций на
заключительном интервале недостижимого кардинала \ $k$.
\\
Итак, следующая лемма показывает, что существует кардинал \
$\delta <k$ \ такой, что
\begin{equation*}
\{ \tau: \delta < \gamma _{\tau } < k \} \subseteq  dom ( \alpha
S_{f} );
\end{equation*}
более точно:
\\
\begin{lemma}
\label{8.9.} \emph{(Об определённости \ $\alpha $-функции)}
\newline \hspace*{1em}\\Существуют ординалы \ $\delta < \gamma <
k$ \ такие, что для каждого \ \mbox{$SIN_n$-}кардинала $\alpha_1
> \gamma$, $\alpha_1 < k$, \  предельного для \ $SIN_n \cap \alpha_1$,\ функция \ $\alpha
S_f^{<\alpha_1}$ \ определена на непустом множестве
\[
    T^{\alpha_1} = \{\tau: \delta < \gamma_{\tau}^{<\alpha_1} < \alpha_1\}.
\]
Минимальный из таких ординалов \ $\delta$ \ обозначается через \
$\alpha \delta ^{\ast }$, \ следующий за ним в \ $SIN_n$  \
кардинал -- через \ $\alpha \delta ^{\ast 1}$; \ также вводятся
следующие соответствующие ординалы:
\begin{align*}
    & \alpha \tau_1^{\ast} = \tau(\alpha \delta^{\ast}),
    \ \alpha \tau^{\ast 1} = \tau(\alpha \delta^{\ast 1}),
\\
    & \mbox{\it так что \ \ } \alpha \delta^{\ast} = \gamma_{\alpha
    \tau_1^{\ast}}, \ \alpha \delta^{\ast 1} = \gamma_{\alpha
    \tau^{\ast 1}},
\\
    & \mbox{\it и \ \ } \alpha^{\ast 1} = \alpha_{\alpha \tau^{\ast
    1}}^{< \alpha_1\Downarrow}, \ \alpha \rho^{\ast 1} = \rho_{\alpha \tau^{\ast
    1}}^{< \alpha_1}.
\end{align*}
\end{lemma}

\noindent \textit{Доказательство} состоит в применении леммы~6.14
\cite{Kiselev11}, как это было сделано в доказательстве леммы~7.6,
но для \textit{большего} редуцирующего кардинала
\begin{equation*}
\chi = ( \chi^{\ast})^{+ \omega_{0}} \mbox{\it \quad и \quad}
\alpha_{1} = k, \quad m=n+1.
\end{equation*}
Получающуюся функцию \ $\mathfrak{A}$, \ определённую на некотором
непустом множестве
\begin{equation*}
T = \{ \tau : \gamma_{\tau_{0}} \leq \gamma_{\tau} < k \},
\end{equation*}
нужно рассмотреть следующим образом:
\\
Рассмотрим согласно лемме~6.14 \cite{Kiselev11} матрицу \
$S_{\tau}^{1} = \mathfrak{A} (\tau) $, \ редуцированную к
кардиналу \ $\chi=(\chi^{\ast})^{+ \omega_0}$ \ на носителе \
\mbox{$\alpha_{\tau}^{1}> \gamma_{\tau} $;} \ она имеет допустимый
производящий собственный диссеминатор \ $
\check{\delta}_{\tau}^{1} < \gamma_{\tau} $ \ с базой \
$\rho_{\tau}^{1} \vartriangleright S_{\tau}^{1} $. \ Можно видеть,
что \ $\rho_{\tau}^{1} > \chi^{\ast+} $ \ и поэтому \
$\check{\delta}_{\tau}^{1} $ \ может быть рассмотрен как
допустимый диссеминатор для \ $S_{\tau}^1$ \ на \
$\alpha_{\tau}^1$ \ с базой \ $\chi^{*+}$.
\\
Теперь обратимся к кардиналу предскачка
\[
    \alpha^{1} = \alpha_{\tau}^{1 \Downarrow};
\]
по той же лемме \ $cf(\alpha^{1}) \geq \chi^{\ast +} $ \ и можно
ввести матрицу \ $S_{\tau}^{2} $ \ \textit{редуцированную} к \ $
\chi^{\ast} $ \ и обладающую тем же кардиналом предскачка скачка \
$\alpha^{1} $ \ и поэтому тем же диссеминатором  \
$\check{\delta}_{\tau}^{1} $ \  с той же базой \ $ \chi^{\ast +} $
,\ используя лемму 5.12~\cite{Kiselev11} следующим образом:
\\
Если выполняется утверждение
\begin{equation} \label{e8.6}
\exists \alpha \in [ \alpha^{1}, \alpha_{\tau}^{1}  [ \ \ \sigma
(\chi^{\ast}, \alpha),
\end{equation}
тогда пусть \ $S_{\tau}^{2} $ \ будет матрица, редуцированная к \
$\chi^{\ast} $ \ на минимальном носителе \ $\alpha_{\tau}^{2}$ \ и
порождённая кардиналом \ $\alpha^{1}$, \ так что  \ $\alpha^1 =
\alpha_\tau ^{2 \Downarrow} $ \ (как это было сделано в
доказательстве леммы~6.12 \cite{Kiselev11} с кардиналом \
$\alpha_0$, \ играющем здесь роль \ $\alpha^1$ \ ).
\\
В противном случае, когда (\ref{e8.6}) нарушается, можно видеть,
что так как  утверждение леммы~5.12 \cite{Kiselev11} сохраняется
ниже \ $\alpha^{1} $,\ то матрица \ $S_{\tau}^{1} $ \ сохраняет
кардинал предскачка \ $\alpha_{\tau}^{\downarrow} $ \ (и, значит,
сохраняет \ $ \alpha^{1}$), \ то есть он сохраняется (остаётся
кардиналом предскачка) при редуцировании матрицы \ $S_{\tau}^{1} $
\ на носителе \ $\alpha_{\tau}^{1} $ \ к кардиналу \ $\chi^{\ast}
$; \ поэтому мы можем определить следующую матрицу (см.
определения 4.1, 5.1, 5.5~\cite{Kiselev11})
\begin{equation*}
S_{\tau}^{2} \Rightarrow \widetilde{\mathbf{S}}_{n}^{\sin
\triangleleft \alpha_{\tau}^{1}} \overline{\overline{\lceil }}
\chi^{\ast} \mbox{\it \quad на носителе \quad} \alpha_{\tau}^{2} =
\sup dom \bigl ( \widetilde{\mathbf{S}}_{n}^{\sin \triangleleft
\alpha_{\tau}^{1}} \overline{ \overline{\lceil }}  \chi^{\ast}
\bigr ) .
\end{equation*}
Эта матрица -- сингулярная на носителе \ $\alpha_\tau^2$: условия
1), 3) определения сингулярности 5.7~\cite{Kiselev11} очевидно
выполняются, а условие 2) можно установить с помощью
\textit{расщепляющего метода}, повторяя дословно аргумент из
даказательства леммы 5.12~\cite{Kiselev11} (где \ $\alpha_1$, \
$\chi$ \ заменяются на \ $\alpha_\tau^2$, \ $\chi^\ast$ \
соответственно).
\\
В любом случае \ $\alpha^{1} = \alpha_{\tau}^{2 \Downarrow} $ \ и
\ $ S_{\tau}^{2} $ \ оказывается допустимой на
 \ $\alpha_{\tau}^{2} $ \ для \
$\gamma_{\tau} $ \ вместе с тем же диссеминатором \
$\check{\delta}_{\tau}^{1} $ \ и его базой \ $\chi^{\ast +} $, \
так как все условия утверждения
 \ $\mathbf{K}^0$ \ из определения
\ref{8.2.} очевидно выполняются когда \ $\rho=\chi^{\ast +}$. \
Точно также такая матрица \ $S_{\tau}^2$ \ на носителе \
$\alpha_{\tau}^2$ \ неподавлена благодаря базе \ $\rho=\chi^{\ast
+}$. \ Она может быть единичной или нулевой, но в любом случае
существует некоторая \ $\alpha$-матрица, редуцированная к \
$\chi^{\ast}$, \ допустимая и неподавленная для \ $\gamma_{\tau}$
\ вместе со своими атрибутами.
\\
Теперь нужно взять любой достаточно большой кардинал  \ $\gamma$ \
такой, что для всякого \ $\gamma_{\tau} > \gamma$ \ существует
некоторая матрица \ $S_{\tau}^2$ \ с базой \ $\rho = \chi^{\ast
+}$; \ она допустима и неподавлена для \ $\gamma_{\tau}$ \ ниже \
$\alpha_1$ \ для каждого \ $\alpha_1 \in SIN_n$, \ $\alpha_1
> \gamma$, $\alpha_1 < k$, \ по
определению.
\\
Таким образом, после минимизации таких получающихся матриц и их
атрибутов следуя определению~\ref{8.3.} возникает функция \
$\alpha S_{f}^{<\alpha_1} $ \ и её сопровождающие ординальные
функции, определённые на  \ $T^{\alpha_1}$ \ для всякого \
$\alpha_1 \in SIN_n$, $\alpha_1 > \gamma$.
\\
\hspace*{\fill} $ \dashv$
\\

\noindent В заключение этого раздела применяя метод доказательства
леммы~\ref{7.7.} устанавливается

\begin{lemma}
\label{8.10.}
\begin{equation*}
\alpha \delta ^{\ast } \in SIN_{n} \cap SIN_{n+1}^{<\alpha^{\ast
1}}[< \alpha \rho^{\ast 1}].
\end{equation*}
\end{lemma}
\noindent \textit{Доказательство.} \ Будем использовать
обозначения из предыдущей леммы~\ref{8.9.}.
\\
Сначала начинает действовать рассуждение из доказательства
леммы~\ref{7.7.}, рассматривая диссеминатор \
$\widetilde{\delta}_{\alpha \tau^{\ast 1}}$ \ с базой \ $\alpha
\rho^{\ast 1}$ \ матрицы \ $\alpha S_{\alpha \tau^{\ast 1}}$ \ на
носителе \ $\alpha_{\alpha \tau^{\ast 1}}$ \ с кардиналом
предскачка \ $\alpha^{\ast 1} = \alpha_{\alpha \tau^{\ast
1}}^{\Downarrow}$. \ Так как
\[
    \alpha \delta^{\ast 1} \in SIN_n, \
    \widetilde{\delta}_{\alpha \tau^{\ast 1}} <
    \alpha \delta^{\ast 1}
\]
и
\[
    \widetilde{\delta}_{\alpha \tau^{\ast 1}} \in
    SIN_n^{<\alpha^{\ast 1}} \cap  SIN_{n+1}^{<\alpha^{\ast 1}}
    [< \alpha \rho^{\ast 1}],
\]
то лемма~3.8 влечёт\ $\widetilde{\delta}_{\alpha \tau^{\ast 1}}
\in SIN_n$.

\noindent Теперь предположим, что эта лемма~\ref{8.10.} неверна:
\[
    \alpha \delta^{\ast} \notin SIN_n,
\]
тогда
\[
    \widetilde{\delta}_{\alpha \tau^{\ast 1}} <
    \alpha \delta^{\ast} = \gamma_{\alpha \tau_1^{\ast}}.
\]
Благодаря лемме 3.2~\cite{Kiselev11} \ можно ограничить
\mbox{$\Sigma_{n-1}$-}утверждение о существовании носителя матрицы
\ $\alpha S_{\alpha \tau^{\ast 1}}$, \ допустимом вместе с теми же
\ $\widetilde{\delta}_{\alpha \tau^{\ast 1}}$, \ $\alpha
\rho^{\ast 1}$, \ имеющимся \ $SIN_{n-1}$-кардиналом \
$\gamma_{\alpha \tau_1^{\ast}+1}$, \ как это было сделано в
доказательстве леммы~\ref{8.5.}~8).
\\
Тогда матрица \ $\alpha S_{\alpha \tau^{\ast 1}}$ \ получает снова
некоторый свой носитель
\[
    \alpha^{\prime} \in \; ] \gamma_{\alpha \tau^{\ast 1}},
    \gamma_{\alpha \tau^{\ast 1} + 1} [,
\]
допустимый для \ $\gamma_{\alpha \tau^{\ast 1}}$ \ вместе с теми
же диссеминатором и базой данных.
\\
Но благодаря минимальности \ $\alpha \tau_1^{\ast}$ \ выполняется
\[
    \alpha \tau_1^{\ast} \notin dom(\alpha S_f).
\]
Это может быть только когда матрица \ $\alpha S_{\tau^{\ast 1}}$ \
на \ $\alpha^{\prime}$ \ допустима, но подавлена для \
$\gamma_{\alpha \tau_1^{\ast}}$; \ в свою очередь если это может
быть, то только когда
\[
    \alpha \delta^{\ast} = \gamma_{\alpha \tau_1^{\ast}} \in
    SIN_n,
\]
вопреки предположению.

Что касается оставшейся части леммы:
\[
    \alpha \delta^{\ast} \in SIN_n^{<\alpha^{\ast 1}} [ < \alpha
    \rho^{\ast 1} ],
\]
то она не используется в дальнейшем вплоть до \S 11 и поэтому мы
мы вернёмся к ней в \S 11 . \hspace*{\fill} $\dashv$

\newpage

\section{Анализ монотонности \ $\protect\alpha$\,-функций}
\setcounter{equation}{0}

\hspace*{1em} В этом разделе первая составляющая требуемого
противоречия -- монотонность \ $\alpha$-функций -- исследуется в
различных важных случаях.

Как мы увидим, это свойство довольно сильное; в частности всякий
интервал \ $[ \tau_1, \tau_2 [$ \ монотонности такой функции не
может быть ``слишком длинным'', -- соответствующий интервал \ $]
\gamma_{\tau_1}, \gamma_{\tau_2} [$ \ не может содержать никаких
 \ \mbox{$SIN_n$-кардиналов}, и если
 \ $\gamma_{\tau_2} \in SIN_n$, \ то такая функция получает
некоторые постоянные характеристики и \textit{стабилизируется} на
таком \ $[ \tau_1, \tau_2 [$.
\\
Мы начинаем с последней ситуации стабилизации:

\begin{definition}
\label{9.1.} \hfill {}

Функция \ $\alpha S_{f}^{<\alpha_{1}}$ \ называется монотонной на
интервале \ $\left[ \tau _{1},\tau _{2}\right[ $ \ и на
соответствующем интервале \ $[ \gamma_{\tau_1}^{<\alpha_1},
\gamma_{\tau_2}^{<\alpha_1} [$ \ ниже \ $ \alpha _{1}$, \ если \
$\tau _{1}+1<\tau _{2}$, \ $\left] \tau _{1},\tau _{2} \right[
\subseteq dom ( \alpha S_{f}^{<\alpha _{1}} ) $ \ и

\begin{equation*}
\forall \tau ^{\prime },\tau ^{\prime \prime }  ( \tau _{1} <
\tau^{\prime }<\tau ^{\prime \prime } <  \tau _{2}\longrightarrow
\alpha S_{\tau ^{\prime }}^{<\alpha _{1}} \underline{\lessdot
}\alpha S_{\tau ^{\prime \prime }}^{<\alpha _{1}} )~.
\end{equation*}
\vspace{-12pt} \hspace*{\fill} $\dashv$
\end{definition}

\noindent Чтобы оперировать этим понятием, удобно использовать
следующие \ \mbox{$\Delta_1$-формулы}, которые будут играть
основную роль в этом разделе:
\\

\noindent \quad $A_0^{1 \vartriangleleft \alpha_1}(\chi, \tau_1,
\tau_2, \alpha S_f^{<\alpha_1})$: \vspace{-6pt}
\begin{multline*}
    A_0^{\vartriangleleft \alpha_1}(\chi, \tau_1, \tau_2, \alpha S_f^{<\alpha_1}) \wedge
    \forall \tau^{\prime}, \tau^{\prime\prime} \big( \tau_1 < \tau^{\prime}
    < \tau^{\prime\prime} < \tau_2 \rightarrow
\\
    \rightarrow \alpha S_f^{<\alpha_1}(\tau^{\prime})
    \underline{\lessdot} \alpha
    S_f^{<\alpha_1}(\tau^{\prime\prime}) \big);
\end{multline*}
таким образом, здесь утверждается, что функция\ $\alpha
S_f^{<\alpha_1}$ \ определена на интервале\ $] \tau_1,\tau_2 [$, \
обладающим свойством \ $A_0$ \ (напомним определение
 \ref{8.1.}~1.0\;), и более того --  она монотонна
на интервале \ $[ \tau_1,\tau_2[\;$; \ поэтому мы будем называть
его и соответствующий интервал \ $[
\gamma_{\tau_1},\gamma_{\tau_2} [$ \ \emph{интервалами
монотонности} функции \ $\alpha S_f^{<\alpha_1}$; \
\\

\noindent \quad $A_1^{1 \vartriangleleft \alpha_1}(\chi, \tau_1,
\tau_2, \alpha S_f^{<\alpha_1})$:
\\
\[
    \qquad \qquad A_0^{1 \vartriangleleft \alpha_1}(\chi, \tau_1, \tau_2, \alpha S_f^{<\alpha_1})
    \wedge \exists \gamma^2 \big(
    \gamma^2 = \gamma_{\tau_2} \wedge SIN_n^{<\alpha_1}(\gamma^2) \big);
\]
\quad %

\noindent далее символы функции \ $\alpha S_f^{<\alpha_1}$ \ будут
опускаться в таких обозначениях (если она будет подразумеваться в
контексте ).

Теперь ещё не все готово для доказательства (тотальной)
монотонности функции \ $\alpha S_f$ \ -- второго компонента
заключительного противоречия -- но некоторые её фрагменты очевидны
аналогично леммам 5.17~1)~\cite{Kiselev11}, \ref{7.9.}. Например,
из леммы 3.2~\cite{Kiselev11} непосредственно следует

\begin{lemma}
\label{9.2.} \emph{(О монотонности $\alpha $-функции)}
\newline \hspace*{1em} Пусть
\begin{equation*}
\tau _{1}<\tau _{2}, ~~ a_{\tau _{2}}^{<\alpha _{1}}= 1  ~~, ~~
\widetilde{\delta }_{\tau _{2}}^{<\alpha _{1}} <  \gamma_{\tau
_{1}},
\end{equation*}
\noindent тогда
\begin{equation*}
\alpha S_{\tau _{1}}^{<\alpha _{1}}\underline{\lessdot }  \ \
\alpha S_{\tau _{2}}^{<\alpha_{1}}  ,\quad \quad a_{\tau_{1}}^{<
\alpha_{1}} = 1 .
\end{equation*}

\noindent Аналогично для нулевой характеристики \ $a_{\tau
_{1}}^{<\alpha _{1}} = a_{\tau _{2}}^{<\alpha _{1}} = 0 $.
\\
\hspace*{\fill} $\dashv$
\end{lemma}

\begin{lemma}
\label{9.3.} \emph{(О стабилизации \ $\alpha $-функции)}
\newline \hspace*{1em} Пусть
\\
(i) \ $\alpha S_{f}^{<\alpha _{1}}$ \ монотонна на \ $[ \tau
_{1},\tau _{2} [$ \ ниже \ $\alpha_1$:
\[
    A_1^{1 \vartriangleleft \alpha_1}(\tau_1, \tau_2);
\]
(ii) \ $\gamma _{\tau _{2}}^{<\alpha_{1}}$ \ это наследник в \
$SIN_{n}^{<\alpha _{1}}$. \newline Тогда функция \ $\alpha
S_{f}^{<\alpha _{1}}$ \ стабилизируется на \ $[ \tau _{1},\tau
_{2} [ $, \ то есть существуют  \ $ S^{0}$ \ и \ $\tau _{0}\in \;
] \tau_{1},\tau _{2} [ \;$ \ такие, что

\vspace{6pt}
\begin{equation*}
\forall \tau \in \left[ \tau _{0},\tau _{2}\right[ \quad \alpha
S_{\tau }^{<\alpha_{1}}=S^{0}.
\end{equation*}
\vspace{0pt}

\noindent Наименьший из таких ординалов \ $\tau _{0}$ \ называется
{\sl ординалом стабилизации} функции \ $\alpha S_{f}^{<\alpha
_{1}}$ \ для \ $\tau_2$ \ ниже \ $\alpha_1$ \ и обозначается через
\ $\tau_2^{s<\alpha_1} $.
\end{lemma}

\noindent \textit{Доказательство} снова представляет собой
типичное применение леммы 3.2~\cite{Kiselev11}; мы будем опускать
верхние индексы \ $<\alpha _{1}$, $\vartriangleleft \alpha _{1}$.
\ Предположим, что эта лемма неверна; рассмотрим ординал
\[
    \rho ^{0}= \sup \{ Od ( \alpha S_{\tau } ) :  \tau _{1}<\tau
    <\tau_{2}\}.
\]
Мы будем применять метод рассуждений, использованный выше в
доказательстве леммы \ref{8.5.}~8); введём для этого кардиналы
\[
    \gamma_{\tau_2^n} = \sup \{ \gamma < \gamma_{\tau_2}: \gamma
    \in SIN_n \};
\]
\[
    \gamma_{\tau_{1,2}^n} = \max \{ \gamma_{\tau_1}, \gamma_{\tau_2^n}
    \}.
\]

\noindent Затем следует повторить определение~\ref{8.3.} матричной
функции \ $\alpha S_f$ \ и сопровождающих его ординальных функций
ниже \ $\alpha_1$ \ на множестве

\begin{equation} \label{e9.1n}
    T_{\tau_{1,2}^n}^{\alpha_1} = \{ \tau: \gamma_{\tau_{1,2}^n}
    < \gamma_{\tau} < \alpha_1 \},
\end{equation}
но сохраняя \textit{только $SIN_n$-кардиналы} \ $\le
\gamma_{\tau_2^n}$; \ это можно сделать следующим образом:
\\
Определение~\ref{8.3.} базируется на формуле
\begin{equation} \label{e9.2n}
\alpha \mathbf{K}^{\ast}(a, \delta, \gamma_{\tau}, \alpha, \rho,
S)
\end{equation}
ниже \ $\alpha_1$ \ (см. определение \ref{8.2.}~4)\;), которая
означает, напомним, что \ $S$ \ это \ $\alpha$-матрица на её
носителе \ $\alpha$ \ характеристики \ $a$ \ с диссеминатором \
$\delta$ \ и базой \ $\rho$, \ допустимая для \ $\gamma_{\tau}$ \
и, более того, неподавленная на этом  \ $\alpha$ \ для \
$\gamma_{\tau}$ \ ниже \ $\alpha_1$; \ но так как для каждого \
$\tau \in T_{\tau_{1.2}^n}^{<\alpha_1}$ \ выполняется \
$\gamma_{\tau} \notin SIN_n$, \ то условие неподавления \ $\neg
A_5^{S,0}$ \ в \ $\alpha \mathbf{K}^{\ast}$ \ тривиально
выполняется и может быть опущено, а тогда \ $\alpha
\mathbf{K}^{\ast}$ \ трансформируется в формулу $\alpha
\mathbf{K}$.
\\
Эта последняя формула относится к классу \ $\Sigma_n$, \ потому
что она включает в себя \ $\Sigma_n$-формулу \
$\mathbf{K}_n^{\forall}$. \ Но рассмотрим кардинал \
$\gamma_{\tau_2^n}$ \ и заменим в формуле~(\ref{e9.2n}) её
поформулу \ $\mathbf{K}_n^{\forall}$ \  \ $\Delta_1$-формулой
\[
    SIN_n^{<\alpha^{\Downarrow}}(\gamma_{\tau_2^n}),
\]
тогда \ $\Sigma_n$-формула (\ref{e9.2n}) преобразуется в некоторую
\ $\Pi_{n-2}$-формулу, которую будем обозначать через
\[
    \alpha \mathbf{K}_{n-2}^{\ast 1}(a, \delta, \gamma_{\tau},
    \alpha, \rho, S).
\]
Таким образом матричная функция, определённая на множестве \
$T_{\tau_{1,2}^n}^{\alpha_1}$ \ (\ref{e9.1n}) как в
определении~\ref{8.3.}, но посредством формулы (\ref{e9.2n})
заменённой на \ $\alpha\mathbf{K}_{n-2}^{\ast 1}$, \ очевидно
совпадает с функцией \ $\alpha S_f$ \ на интервале \
$]\tau_{1,2}^n, \tau_2 [\;$; \ мы будем обозначать эту последнюю
функцию через \ $\alpha S_f^1$.
\\
Далее, так как эта \ $\alpha S_f^1$ \ монотонна на on \
$]\tau_{1,2}^n,~\tau_2[$, \ но не стабилизируется на этом
интервале, то ординал \ $\rho_0$ \ -- \textit{предельный} и
поэтому выполняется следующее предложение ниже \
$\gamma_{\tau_2}$:

\vspace{6pt}
\begin{equation*}
\forall \tau ~ ( \tau_{1,2}^n < \tau \longrightarrow \exists S~ (
S=\alpha S_{\tau }^{1}\wedge S\vartriangleleft  \rho ^{0} ) )~.
\end{equation*}
\vspace{0pt}

\noindent Оно может быть сформулировано в \ $\Pi_n$-форме:

\begin{equation*}
\forall \tau, \gamma^{\prime}, \gamma^{\prime\prime} \Bigl[
\gamma_{\tau_{1,2}^n} < \gamma^{\prime} = \gamma_{\tau} <
\gamma^{\prime\prime} = \gamma_{\tau+1} \rightarrow
\qquad\qquad\qquad\qquad\qquad\qquad
\end{equation*}
\begin{equation} \label{e9.3n}
~ \rightarrow \Big( \exists \delta, \alpha, \rho <
\gamma^{\prime\prime} \exists S \vartriangleleft \rho ~ \big(
\alpha \mathbf{K}_{n-2}^{\ast 1} (1, \delta, \gamma^{\prime},
\alpha, \rho, S) \wedge S \vartriangleleft \rho^0 \big) \vee
\end{equation}
\begin{equation*}
~ \vee \big( \exists \delta, \alpha, \rho < \gamma^{\prime\prime}
\exists S \vartriangleleft \rho ~ \alpha \mathbf{K}_{n-2}^{\ast 1}
(0, \delta, \gamma^{\prime}, \alpha, \rho,  S) \wedge S
\vartriangleleft \rho^0 \wedge \quad \quad
\end{equation*}
\begin{equation*}
\quad \wedge \forall \delta^{\prime}, \alpha^{\prime},
\rho^{\prime} < \gamma^{\prime\prime} \forall
S^{\prime}\vartriangleleft \rho^{\prime} ~ \neg \alpha
\mathbf{K}_{n-2}^{\ast 1} (1, \delta^{\prime}, \gamma^{\prime},
\alpha^{\prime}, \rho^{\prime}, S^{\prime}) \big) \Big) \Bigr ] ~.
\end{equation*}
\vspace{0pt}

\noindent Теперь возникает противоречие:
\\
С одной стороны, \ $SIN_n$-кардинал \ $\gamma _{\tau _{2}}$ \
продолжает утверждение (\ref{e9.3n}) до \ $\alpha _{1}$ \ и в
результате возникает\emph{ минимальная} матрица \ $ \alpha S_{\tau
_{2}}^{1}\vartriangleleft \rho ^{0}$ .\
\\
Но сдругой стороны \ $\rho ^{0} $ \ -- это предельный ординал и
поэтому существует
\[
    \tau _{1,3}^{n} \in ] \tau _{1,2}^{n}, \tau _{2} [
    \mbox{\it \ \ такой, что \ }
    \alpha S_{\tau _{1,3}^{n}}^{1} \gtrdot \alpha S_{\tau _{2}}^{1}.
\]
Поэтому ниже \ $\gamma _{\tau _{2}}$ \ выполняется следующее
предложение:

\vspace{6pt}
\begin{equation*}
\forall \tau ~ \big( \tau _{1,3}^{n }<\tau \longrightarrow \forall
S ( S=\alpha S_{\tau }^{1}\longrightarrow S \gtrdot \alpha
S_{\tau_{2}}^{1} ) \big)~.
\end{equation*}
\vspace{0pt}

\noindent Оно тоже может быть сформулировано в \ $\Pi_{n}$-форме:

\begin{equation*}
\forall \tau, \gamma^{\prime}, \gamma^{\prime\prime} \Bigl[
\gamma_{\tau_{1,3}^n} < \gamma^{\prime} = \gamma_{\tau} <
\gamma^{\prime\prime} = \gamma_{\tau+1} \rightarrow \quad \quad
\quad \quad \quad \quad \quad \quad \quad \quad \quad \quad \quad
\end{equation*}
\begin{equation*}
\rightarrow \Big( \forall \delta, \alpha, \rho <
\gamma^{\prime\prime} ~ \forall S \vartriangleleft \rho \big(
\alpha \mathbf{K}_{n-2}^{\ast 1} (1, \delta, \gamma^{\prime},
\alpha, \rho, S) \rightarrow \alpha S_{\tau_2}^{1} \lessdot S
\big) \wedge
\end{equation*}
\begin{equation} \label{e9.4}
\quad \wedge \forall \delta, \alpha, \rho < \gamma^{\prime\prime}
~ \forall S \vartriangleleft \rho  \big( \alpha
\mathbf{K}_{n-2}^{\ast 1} (0, \delta, \gamma^{\prime}, \alpha,
\rho, S) \wedge
\end{equation}
\begin{equation*}
\wedge \forall \delta^{\prime}, \alpha^{\prime},  \rho^{\prime} <
\gamma^{\prime\prime} ~ \forall S^{\prime} \vartriangleleft
\rho^{\prime}  \neg \alpha \mathbf{K}_{n-2}^{\ast 1} (1,
\delta^{\prime}, \gamma^{\prime}, \alpha^{\prime}, \rho^{\prime},
S^{\prime}) \rightarrow
\end{equation*}
\begin{equation*}
\qquad \qquad \qquad \qquad \qquad \qquad \qquad \qquad \qquad
\qquad \rightarrow \alpha S_{\tau_2}^{1} \lessdot S \big) \Big)
\Bigr ] ,
\end{equation*}

\vspace{6pt}

\noindent которое \ $\gamma _{\tau _{2}}$ \ продолжает до \
$\alpha _{1} $ \ и поэтому
\[
    \alpha S_{\tau _{2}}^{1} \lessdot \alpha S^1_{\tau _{2}}.
\]
\hspace*{\fill} $\dashv$

\quad \\

Напомним, что символы \ $\mathfrak{n}^{\alpha}$, \ $\chi^{\ast}$,
\ $ \alpha S_{f}^{<\alpha _{1}}$, \ $a_{f}^{<\alpha _{1}}$ \ в
написаниях формул могут часто опускаться для сокращения записей.
Кроме того, условие эквиинформативности \ $A_6^e(\alpha_1)$
\begin{multline*}
    \chi^{\ast} < \alpha_1 \wedge A_n^{<\alpha_1} (\chi^\ast) =  \|
    u_n^{\vartriangleleft \alpha_1} ( \underline{l}) \| \wedge SIN_{n-2}
    (\alpha_1) \wedge
\\
    \wedge \forall \gamma < \alpha_1 \exists \gamma_1 \in [\gamma,
    \alpha_1[ ~ SIN_n^{<\alpha_1}(\gamma_1)
\end{multline*}
\vspace{0pt}

\noindent всегда накладывается на ограничивающие кардиналы \
$\alpha_1$.
\newline

Стабилизационное свойство чрезвычайно существенно для дальнейшего;
более того, оказывается, что аналогичное свойство возникает у
характеристической функции, играющее далее ключевую роль --
\\
усложняя определённым образом доказательство леммы~\ref{9.3.}
можно доказать аналогичное характеристическое свойство:

\begin{lemma}
\label{9.4.} \emph{(О стабилизации характеристики)}
\newline \hspace*{1em} Пусть \medskip

(i) \quad $A_{1}^{1 \vartriangleleft \alpha_1} ( \tau _{1},
\tau_{2} ) ;$
\medskip

(ii) \quad $\forall \tau <\tau _{2}\quad \exists \tau^{\prime} \in \left[
\tau ,\tau_{2}\right[ \quad a_{\tau ^{\prime }}^{<\alpha _{1}}=1$; \medskip
\newline
Тогда
\begin{equation*}
\forall \tau ^{\prime }\in \left] \tau _{1},  \tau_{2}\right[
\quad a_{\tau^{\prime} }^{<\alpha _{1}}=1 .
\end{equation*}
В таком случае мы будем говорить, что единичная характеристика
стабилизируется на \ $[ \tau _{1},\tau _{2} [ $ \ ниже \ $\alpha
_{1}$.
\\
Аналогично для нулевой характеристики.
\end{lemma}

\noindent \textit{Доказательство} проводится индукцией по \ $(
\alpha_{1}, \tau_{2} )$ \ (напомним, что множество таких пар
рассматривается как канонически упорядоченное как и ранее, с \
$\alpha_1$ \ как первым компонентом и \ $\tau_2$ \ вторым).
Предположим, что эта пара -- минимальная нарушающая эту лемму.
Нетрудно видеть, что \ $\gamma_{\tau_{2}}^{< \alpha_{1}} $ \
является наследником в классе \ $SIN_{n}^{< \alpha_{1}} $; \
именно этот случай используется в дальнейшем. Напомним, что для
каждой матрицы единичной характеристики на её носителе и каждого
$\gamma_{\tau} \notin SIN_n$ \    условие неподавления  \ $\neg
A_5^{S,0}$ \ тривиально выполняется и может быть опущено; таким
образом для единичной характеристики \ $a=1$ \ формула \ $\alpha
\mathbf{K}^{\ast <\alpha_1}$ \ эквивалентна $\alpha
\mathbf{K}^{<\alpha_1}$; \ верхние индексы \ $< \alpha_1$,
$\vartriangleleft \alpha_1$ \ будут опускаться для краткости.
\\
Согласно предыдущей лемме существует стабилизационный ординал \
$\tau_{2}^s $ \ функции \ $\alpha S_{f} $ \ на \ $[ \tau_{1},
\tau_{2} [ $ \ и матрица \ $S^0$ \ такие, что

\begin{equation*}
\forall \tau \in [ \tau_{2}^s, \tau_{2}[ \quad  \alpha S_{\tau} =
S^{0} .
\end{equation*}

\noindent По  условию существует минимальный ординал \ $\tau^{1}
\in [ \tau_{2}^s, \tau_{2}[ $ \ такой, что \ $ a_{\tau^{1}} = 1 $.
\ Дальнейшее рассуждение разделяется на две части:
\\
1. \ Сначала докажем, что \ $\forall \tau \in [ \tau^{1},
\tau_{2}[ \ \ a_{\tau} = 1 $. \newline Допустим это неверно, тогда
существует минимальный ординал \ $\tau^{0} \in ] \tau^{1},
\tau_{2}[ $ \ для которого \ $a_{\tau^{0}} = 0 $; \ так что \
$a_{\tau} \equiv 1 $ \ on \ $[ \tau^{1}, \tau^{0}[ $ \ (можно
заметить кстати, что здесь матрица \ $S^0$ \ на разных носителях
может обладать разными характеристиками). Рассмотрим следующие
подслучаи для
\begin{equation*}
S^0 = \alpha S_{\tau^0}, \quad \check{\delta}^{0} =
\check{\delta}_{\tau^{0}} , \quad \alpha^{0} =
\alpha_{\tau^{0}}^{\Downarrow} ~:
\end{equation*}
\newline
1a. \quad $\check{\delta}^{0} \notin SIN_{n} $. \ Так как \
$\gamma_{\tau_{1}} \in SIN_{n} $, \ то лемма~3.8~1)
\cite{Kiselev11} влечёт
\begin{equation*}
\gamma_{\tau_1} < \check{\delta}^{0}, \quad \check{\delta}^{0} \in
\left ( SIN_{n}^{< \alpha^{0}} - SIN_{n} \right )
\end{equation*}
и можно использовать кардинал
\begin{equation*}
\gamma_{\tau^{2}} = min \left (SIN_{n}^{< \alpha^{0}} -  SIN_{n} \right ) .
\end{equation*}
Благодаря лемме 3.8 \cite{Kiselev11} нетрудно видеть, что \
$\gamma_{\tau^{2}} $ \ это\textit{наследник} \ $ SIN_{n}^{<
\alpha^{0}} $ \ некоторого кардинала
\begin{equation*}
\gamma_{\tau^{3}} \geq \gamma_{\tau_{1}}, \quad  \gamma_{\tau^{3}}
\in SIN_{n} \mbox{\it \ \ ниже\ \ } \alpha_1,
\end{equation*}
и функция \ $\alpha S_{f}^{< \alpha^{0}} $ \ монотонна на
интервале \ $[ \tau^3, \tau^2 [\;$. \ Так как \ $a_{\tau} \equiv 1
$ \ на \ $[ \tau^{1}, \tau^{0} [ $, \ то лемма
3.2~\cite{Kiselev11} влечёт, что интервал \ $ ]\gamma_{\tau_{1}},
\gamma_{\tau^{2}} [ $ \ содержит допустимые носители матриц
единичной характеристики, расположенные \emph{конфинально}
кардиналу \ $\gamma_{\tau^{2}} $ ,\ потому что \
$SIN_n^{<\alpha^0}$-кардинал \ $\gamma_{\tau^2}$ \ ограничивает \
$\Sigma_n$-утверждение о существовании таких носителей. После
этого кардинал \ $\gamma_{\tau^{2}} $ \ \textit{продолжает}
единичную характеристику до \ $\alpha_0$, \ и таким образом \
$S^0$ \ на \ $\alpha_{\tau_0}$ \ становится единичной матрицей
вопреки предположению. \noindent
\\Метод этого аргумента состоит в ограничениях и продолжениях,
применяемых по очереди, поэтому мы будем называть его методом
\textit{ограничения-и-продолжения}. Он будет часто использоваться
в дальнейшем в разнообразных типичных ситуациях, поэтому сейчас
следует остановиться на нём подробнее:
\\
Рассмотрим произвольный кардинал \ $\gamma < \gamma_{\tau_2}$;\
имеется единичная матрица \ $S^0$ \ на некотором носителе \
$\alpha> \gamma$ \ и она остаётся единичной ниже \ $\alpha^0$ \
благодаря лемме~\ref{8.7.} об абсолютности. Затем рассуждение
переходит к ситуации ниже \ $\alpha^0$; \  очевидно, выполняется
утверждение ниже \ $\alpha^0$:
\[
    \exists \delta, \gamma_{\tau}, \alpha, \rho ~ \bigl( \gamma <
    \gamma_{\tau} \wedge \alpha \mathbf{K}(1, \delta, \gamma_{\tau}, \alpha,
    \rho, S^0) \bigr),
\]
оно относится \ $\Sigma_n$ \ и содержит только константы
\[
    \chi^{\ast}, ~ \gamma < \gamma_{\tau^2}, ~ S^0 \vartriangleleft \chi^{\ast
    +} < \gamma_{\tau^2}.
\]
Поэтому \ $\Pi_n$-кардинал \ $\gamma_{\tau_2}$ \ ниже \ $\alpha_0$
\ ограничивает его по лемме 3.2~\cite{Kiselev11}, то есть оно
выполняется после его ограничения кардиналом $\gamma_{\tau^2}$:
\[
    \exists \delta, \gamma_{\tau}^{<\gamma_{\tau^2}}, \alpha,
    \rho < \gamma_{\tau^2} ~ \bigl( \gamma <
    \gamma_{\tau}^{<\gamma_{\tau^2}} \wedge \alpha
    \mathbf{K}^{\vartriangleleft
    \gamma_{\tau^2}}(1, \delta, \gamma_{\tau}^{<\gamma_{\tau^2}},
    \alpha, \rho, S^0) \bigr).
\]
Но здесь верхние индексы \ $<\gamma_{\tau^2}$, \ $\vartriangleleft
\gamma_{\tau^2}$ \ можно опустить благодаря \
$\Pi_n^{<\alpha^0}$-субнедостижимости кардинала \
$\gamma_{\tau^2}$ \ и в результате здесь появляются допустимые
носители
\[
    \alpha \in \; ] \gamma, \gamma_{\tau_2} [
\]
матрицы \ $S^0$ \  \textit{единичной} характеристики ниже \
$\alpha^0$ \ для произвольного \ $\gamma < \gamma_{\tau_2}$.

\noindent Тогда по индуктивной гипотезе \ $a_{\tau} \equiv 1 $ \
на \ $] \tau_{1}, \tau^{2} [ $, \ и ниже \ $\gamma_{\tau^{2}} $ \
выполняется утверждение
\begin{equation*}
\forall \tau \ \ (\gamma_{\tau^{3}} <  \gamma_{\tau}\longrightarrow a_{\tau}
= 1 )
\end{equation*}
которое может быть сформулировано в \ $\Pi_{n} $-форме:
\begin{multline*}
\forall \gamma \Bigl( \gamma_{\tau^{3}} < \gamma  \wedge
SIN_{n-1}(\gamma) \longrightarrow
\\
\longrightarrow \exists \delta, \alpha, \rho, S \Bigl(
SIN_{n}^{<\alpha^{\Downarrow}} (\gamma_{\tau^{3}}) \wedge \alpha
\mathbf{K}_{n+1}^{\ast \exists} (1, \delta, \gamma, \alpha,  \rho,
S) \Bigr) \Bigr) \ .
\end{multline*}
\vspace{0pt}

\noindent Кардинал \ $\gamma_{\tau^{2}} \in SIN_{n}^{< \alpha^{0}}
$ \ продолжает это предложение до \ $\alpha^{0} $ \ и ниже \ $
\alpha^{0} $ \ появляется матрица единичной характеристики на
некотором носителе \ $\in \; ]\gamma_{\tau^{0}}, \alpha^{0} [ $, \
допустимая вместе со своим диссеминатором \ $<\gamma_{\tau_{0}} $
\ и его базой для \ $ \gamma_{\tau_{0}} $.
\\
Таким образом, \ $a_{\tau_{0}}=1 $ \ вопреки предположению, и мы
переходим к следующему подслучаю:
\\

\noindent 1b. \quad $\check{\delta}^{0} \in SIN_{n} $. \ Так как
 \ $a_{\tau} \equiv 1 $ \ на\ $[ \tau^{1}, \tau^{0} [ $,
\ то существует матрица
\[
    S^0 = \alpha S_{\tau^{1,0}} \mbox{\it \ на носителе \ }
    \alpha_{\tau^{1,0}} \in [ \check{\delta}^0, \gamma_{\tau^0}[
\]
\textit{единичной характеристики} \ $a_{\tau^{1,0}}=1$ \ и можно
рассмотреть ситуацию ниже \ $\alpha^{1,0} =
\alpha_{\tau^{1,0}}^{\Downarrow}$ \ следующим образом.
\\
Предстоящее рассуждение применяется далее не один раз, поэтому
нужно остановиться на нём внимательнее.
\\
Начнём с матрицы \ $S^0$ \ на on \ $\alpha_{\tau_0}$. \ По
лемме~\ref{8.5.}~5) нулевая характеристика матрицы \ $S^0$ \ на \
$\alpha_{\tau^0}$ \ значает, что выполняется

\begin{equation}  \label{e9.5}
    \qquad
    \exists \tau_1^{\prime}, \tau_2^{\prime}, \tau_3^{\prime} <
    \alpha^0 ~ \big( A_2^{0 \vartriangleleft \alpha^0} (\tau_1^{\prime},
    \tau_2^{\prime}, \tau_3^{\prime}, \alpha S_f^{<\alpha^0} )
    \wedge \qquad \qquad \qquad
\end{equation}
\[
    \qquad \qquad \qquad
    \wedge \forall \tau^{\prime\prime} \in \; ]\tau_1^{\prime},
    \tau_2^{\prime}] ~ a_{\tau^{\prime\prime}}^{<\alpha^0} = 1
    \wedge \alpha S_{\tau_2^{\prime}}^{<\alpha^0} = S^0 \big).
\]
\vspace{0pt}

\noindent Поэтому могут быть использованы некоторые ординалы
\[
    \tau_1^{\prime} < \tau_2^{\prime} < \tau_3^{\prime} <
    \alpha^0
\]
такие, что выполняется

\begin{equation}  \label{e9.6}
    A_2^0( \tau_{1}^{\prime}, \tau_{2}^{\prime}, \tau_{3}^{\prime},
    \alpha S_{f})
    \wedge \forall \tau^{\prime\prime} \in \; ]\tau_1^{\prime},
    \tau_2^{\prime}] ~ a_{\tau^{\prime\prime}}^{<\alpha^0} = 1
    \wedge \alpha S_{\tau_{2}^{\prime}} = S^{0}
\end{equation}

\noindent ниже \ $\alpha^0$, \ то есть после \
$\vartriangleleft$-ограничения кардиналом \ $\alpha^0$.

\noindent Ключевую роль будет играть так называемый
\textit{медиатор}: это некоторый \ $SIN_n^{<\alpha^0}$-кардинал \
$\gamma^0$ \ такой, что

\begin{equation}  \label{e9.7}
    \gamma_{\tau_1^{\prime}}^{<\alpha^0} <
    \gamma_{\tau_2^{\prime}}^{<\alpha^0} <
    \gamma_{\tau_3^{\prime}}^{<\alpha^0} <
    \gamma^0 < \alpha^0
\end{equation}

\noindent который существует благодаря лемме~\ref{8.5.}~4). По
лемме~\ref{8.7.} об абсолютности допустимости и единичной
характеристики матричной \ $\alpha$-функции эти значения и их
атрибуты ниже \ $\alpha^0$ \ и ниже \ $\gamma^0$ \ сопадают на
множестве
\[
    \bigl \{ \tau: \gamma_{\tau}^{<\alpha^0} < \gamma^0
    \wedge a_{\tau}^{<\alpha^0} = 1
    \bigr \}
\]
и поэтому (\ref{e9.7}), (\ref{e9.6}) влекут следующее \
$\Sigma_{n+1}$-утверждение ниже \ $\alpha^0$:
\[
    \exists \gamma^0 \exists \tau_1^{\prime}, \tau_2^{\prime},
    \tau_3^{\prime} < \gamma^0 \big( SIN_n(\gamma^0) \wedge
    \gamma_{\tau_1^{\prime}}^{<\gamma^0} <
    \gamma_{\tau_2^{\prime}}^{<\gamma^0} <
    \gamma_{\tau_3^{\prime}}^{<\gamma^0} < \gamma^0 \wedge
    \qquad
\]
\begin{equation} \label{e9.8}
    \wedge  A_2^{0 \vartriangleleft \gamma^0} (\tau_1^{\prime},
    \tau_2^{\prime}, \tau_3^{\prime}, \alpha S_f^{<\gamma^0})
    \wedge \forall \tau^{\prime\prime} \in \; ]\tau_1^{\prime},
    \tau_2^{\prime}] ~ a_{\tau^{\prime\prime}}^{<\gamma^0} = 1
    \wedge
    \quad
\end{equation}
\[
    \qquad \qquad \qquad \qquad \qquad \qquad \qquad \qquad
    \wedge \alpha S_{\tau_2^{\prime}}^{<\gamma^0} = S^0 \big).
\]

\noindent По лемме~\ref{8.5.}~3)мы имеем
\[
    \check{\delta}^0 \in SIN_{n+1}^{<\alpha^0}[< \rho_{\tau^0}]
\]
и тогда по лемме 3.2~\cite{Kiselev11} существует некоторое \
$\gamma^0$, \ для которого выполняется (\ref{e9.8}), но
\textit{уже ниже} \ $\check{\delta}^0$. \ Отсюда и из леммы
3.8~\cite{Kiselev11} следует, что \ $SIN_n$-субнедостижимость \
$\check{\delta}^0$ \ влечёт существование \ $SIN_n$-кардинала \
$\gamma^{0 \prime} < \check{\delta}^0$ \ с тем же свойством
(\ref{e9.8}); отметим, что \ $\gamma^{0 \prime}$ \ обладает той же
\ $SIN_n$-субнедостижимостью, что и \ $\check{\delta}^0$. \ Таким
образом, для некоторых кардиналов
\begin{equation} \label{e9.9}
    \gamma_{\tau_1^{\prime\prime}}^{<\gamma^{0 \prime}} <
    \gamma_{\tau_2^{\prime\prime}}^{<\gamma^{0 \prime}} <
    \gamma_{\tau_3^{\prime\prime}}^{<\gamma^{0 \prime}} <
    \gamma^{0 \prime}
\end{equation}
выполняется
\begin{equation} \label{e9.10}
    \qquad \qquad
    A_2^{0 \vartriangleleft \gamma^{0 \prime}} (\tau_1^{\prime\prime},
    \tau_2^{\prime\prime}, \tau_3^{\prime\prime}, \alpha S_f^{<\gamma^{0 \prime}})
    \wedge
    \qquad\qquad\qquad\qquad\qquad\qquad\qquad
\end{equation}
\[
    \qquad\qquad
    \wedge \forall \tau^{\prime\prime\prime} \in \; ]\tau_1^{\prime\prime},
    \tau_2^{\prime\prime}] ~ a_{\tau^{\prime\prime\prime}}^{<\gamma^{0\prime}} = 1
    \wedge \alpha S_{\tau_2^{\prime\prime}}^{<\gamma^{0 \prime}} =
    S^0.
\]

\noindent Так как \ $\gamma^{0 \prime}$ \ это \ $SIN_n$-кардинал,
то везде в (\ref{e9.9}), (\ref{e9.10}) \ $\vartriangleleft$- и \
$<$-ограничения этим кардиналом \ $\gamma^{0 \prime}$ могут быть
опущены по тем же леммам 3.8, \ref{8.7.}.

Начиная с этого места следует повторить проведённые выше
рассуждения, но в обратном порядке, и не для \ $\alpha^0$, \
\textit{но для} \ $\alpha^{1,0}$. \ Тогда (\ref{e9.9}),
(\ref{e9.10}) влекут (\ref{e9.5}), где \ $\alpha^0$ \ заменён на
 \ $\alpha^{1,0}$ \ и поэтому матрица
\ $S^0$ \ на\ $\alpha_{\tau^{1,0}}$ \ получает \textit{нулевую
характеристику} в противоречии с предположением.
\\\noindent 2. Итак, утвеждение 1. здесь доказано; осталось проверить ординал
\[
    \tau^{1,2} = \min \bigl \{\tau \in [ \tau_1, \tau_2 [ : \forall
    \tau^{\prime} \in ] \tau, \tau_2 [ ~ a_{\tau^{\prime}} = 1 \bigr \}
\]
и доказать, что он совпадает с \ $\tau_1$.
\\
Допустим, это не так и \ $\tau_1 < \tau^{1,2}$, \ тогда нужно
рассмотреть две единичные матрицы
\[
    S^1 = \alpha S_{\tau^{1,2}}, ~ S^2 = \alpha S_{\tau^{1,2}+1},
\]
используя матрицу \ $S^2$ \ на её носителе \
$\alpha_{\tau^{1,2}+1}$ \ с её производящим диссеминатором \
$\check{\delta}^2 = \check{\delta}_{\tau^{1,2}+1} $ \  следующим
образом. Согласно леммам~\ref{8.5.}~7)~$(ii)$, \ \ref{8.8.}~2)
получается
\[
    \gamma_{\tau_1} \le \check{\delta}^2 =
    \widetilde{\delta}_{\tau^{1,2}+1}
\]
и поэтому возникают только три подслучая:

\noindent 2a. \quad $\gamma_{\tau_1} = \check{\delta}^2$. \ Тогда
по определению~\ref{8.3.}
\[
    \forall \tau \in ] \tau_1, \tau^{1,2} [ ~ a_{\tau} = 1
\]
в противоречии с предположением.

\noindent 2b. \quad $\gamma_{\tau_1} < \check{\delta}^2$, \
$\check{\delta}^2 \notin SIN_n$. \ Тогда действует метод
ограничения-и-продолжения, буквально как выше в случае 1а. этого
доказательства, но для
\[
    S^2, ~ \check{\delta}^2 ~ \mbox{\it \ вместо}
    ~ \alpha S_{\tau_0}, ~ \check{\delta}^0
\]
и снова получается \ $a_{\tau} \equiv 1$ \ on \ $] \tau_1,
\tau^{1,2} ]$.

\noindent 2c. \quad $\gamma_{\tau_1} < \check{\delta}^2$, \
$\check{\delta}^2 \in SIN_n$. \ Здесь снова используется метод
ограничения-и-продолжения, но в несколько другой манере. Сначала
по лемме 3.2~\cite{Kiselev11} матрица \ $S^1$ \ получает свои
допустимые носители, расположенные конфинально кардиналу
 \ $\check{\delta}^2$, \ так что по индуктивной гипотезе
\[
    a_{\tau} \equiv 1 \mbox{\it \ на множестве \ }
    \{ \tau : \gamma_{\tau_{1}} < \gamma_{\tau} <
    \check{\delta}^2 \}.
\]
Тогда ниже диссеминатора \ $\check{\delta}^2$ \ выполняется
следующее \ $ \Pi_{n+1} $-утверждение

\begin{equation*}
\forall \gamma \big( \gamma_{\tau_1} < \gamma  \wedge
SIN_{n-1}(\gamma) ~ \rightarrow ~ \exists \delta,  \alpha, \rho, S
\quad \alpha \mathbf{K} (1, \delta,  \alpha, \gamma, \rho, S)
\big)
\end{equation*}

\noindent которое распространяется этим диссеминатором до \ $
\alpha_{\tau^{1,2}+1}^{\Downarrow} $ \ согласно лемме
6.6~\cite{Kiselev11} (для \ $ m=n+1 $, \ $\delta =
\check{\delta}^2$, $\alpha_{1} =
\alpha_{\tau^{1,2}+1}^{\Downarrow}$) \ -- и снова получается \ $
a_{\tau} \equiv 1 $ \ на том же множестве \ $]\tau_1,
\tau^{1,2}]$.
\\
И в любом случае мы приходим к \ $\tau_1 = \tau^{1,2}$.
\hspace*{\fill} $\dashv$
\\

Следующая важная лемма будет доказана снова методом
\textit{ограничения-и-продолжения}, но в несколько более
синтезированной форме.

Однако предварительно хорошо ввести следующие довольно удобные
понятия, использующие понятия редуцированных спектров и матриц
(напомним определения~4.1, 5.1~\cite{Kiselev11}).

В дальнейшем основной технический приём в рассуждениях будет
состоять в рассмотрении некоторой матрицы \ $S$ \ на её
\textit{различных носителях по очереди}. Такую трансформацию
редуцированной матрицы \ $S$ \ от одного её носителя \ $\alpha$ \
к другому её носителю \ $\alpha^1$ \ мы будем называть
\textit{переносом} матрицы \ $S$ \ от \ $\alpha$ \ к \ $\alpha^1$.
\\
Такая техника будет часто применяться в дальнейшем и уже была
использована в доказательствах лемм \ref{7.5.}, \ref{8.8.},
\ref{9.4.}.

В ходе такого переноса редуцированной матрицы \ $S$ \ от \
$\alpha$ \ к \ $\alpha^1$ \ некоторые свойства универсума,
ограниченного кардиналами скачка или предскачка матрицы \ $S$ \ на
\ $\alpha$, \ сохраняются и поэтому мы будем называть их
\textit{внутренними свойствами} матрицы \ $S$; \ другие свойства
при этом нарушаются и поэтому они будут называться
\textit{внешними свойствами}.

Более точно: свойство или признак матрицы \ $S$, \ редуцированной
к \ $\chi^{\ast}$ \ на её носителе \ $\alpha$, \ будет называться
\textit{внешним} свойством или признаком этой \ $S$ \ (на \
$\alpha$), \ если оно определимо ниже некоторого кардинала скачка
или предскачка спектра
\[
    dom  \bigl( \widetilde{\mathbf{S}}_{n}^{\sin \vartriangleleft
    \alpha}\overline{\overline{\lceil}}\chi^{\ast} \bigr )
\]
через её некоторые другие кардиналы скачка или предскачка;
аналогично для других объектов из \ $L_{\alpha}$; \ во всех других
случаях они будут называться \textit{внешними свойствами или
признаками} матрицы \ $S$.
\\
Эти понятия вводятся в действие леммой~5.11~\cite{Kiselev11} о
матричной информативности, которая означает, напомним, что такие
\textit{внутренние} свойства сохраняются во время переноса
матрицы\ $S$ \ от одного её носителя к другому.

И вот очень важный пример \textit{внешнего} свойства -- свойство
\textit{характеристики}; оно включает в себя \textit{всю} матрицу\
$S$ \ на её носителе \ $\alpha$, \ а не только некоторые её
кардиналы скачка.
\\
Действительно, возьмём произвольную матрицу \ $S$ \ на её носителе
\ $\alpha$ \ \textit{нулевой} характеристики (если он существует),
тогда по лемме~\ref{8.5.}~5) выполняется \vspace{-6pt}
\begin{multline*}
    \exists \tau_1^{\prime}, \tau_2^{\prime}, \tau_3^{\prime} <
    \alpha^{\Downarrow} \bigl(
    A_2^{0 \vartriangleleft \alpha^{\Downarrow}}
    ( \tau_1^{\prime}, \tau_2^{\prime}, \tau_3^{\prime},
    \alpha S_f^{< \alpha^{\Downarrow}}) \wedge
\\
    \wedge \forall \tau^{\prime\prime} \in \; ]\tau_1^{\prime},
    \tau_2^{\prime}] ~ a_{\tau^{\prime\prime}}^{<\alpha^{\Downarrow}} = 1
    \wedge \alpha S_{\tau_2^{\prime}}^{< \alpha^{\Downarrow}} = S \bigr),
\end{multline*}
где \ $S$ \ получает \textit{меньший} носитель \
$\alpha_{\tau_2^{\prime}}^{< \alpha^{\Downarrow}}$, \ но уже
 \textit{единичной } характеристики, благодаря
условию \ $\forall \tau^{\prime\prime} \in \; ]\tau_1^{\prime},
\tau_2^{\prime}] ~ a_{\tau^{\prime\prime}}^{<\alpha^{\Downarrow}}
= 1$.
\\
Но в дальнейшем все другие матричный свойства --
{\emph{внутренние}}, и одно из них релизует метод
ограничения-и-продолжения в доказательстве нижеследующей леммы
9.5.

\noindent Эта лемма использует удобную функцию, которая уже была
использована в доказательстве леммы~9.3:
\[
    Od \alpha S_{f}^{<\alpha _{1}} ( \tau _{1},\tau _{2} )  = \sup \{
    Od ( \alpha S_{\tau }^{<\alpha _{1}} ) :  \tau_{1} < \tau <\tau
    _{2}\};
\]
она необходима для формирования так называемых \textit{лестниц} --
семейств интервалов, которые будут служить главным техническим
средством в доказательстве основной теоремы. Для этого необходимы
следующие формулы ниже\ $\alpha_1$:
\\

\noindent 1.\quad $A_{1.1}^{m \vartriangleleft \alpha_1}(\tau_1,
\tau_2, \alpha S_f^{<\alpha_1})$:
\[
    A_1^{1 \vartriangleleft \alpha_1}
    (\tau_1, \tau_2, \alpha S_f^{<\alpha_1})
    \wedge \tau_2 = \sup \big \{\tau: A_1^{1 \vartriangleleft \alpha_1}
    (\tau_1, \tau, \alpha S_f^{<\alpha_1}) \big \};
\]
\quad \\
здесь интервал \ $\left[ \tau_1,\tau_2 \right[$ \ -- это
максимальный интервал монотонности с левым \
$SIN_n^{<\alpha_1}$-концом \ $\gamma_{\tau_1}^{<\alpha_1}$ \ и с
правым \ $SIN_n^{<\alpha_1}$-концом  \
$\gamma_{\tau_2}^{<\alpha_1}$ \ \ (максимальный в том смысле, что
он не включается ни в какой другой такой интервал), поэтому мы
будем называть его и соответствующий интервал \ $[
\gamma_{\tau_1}^{<\alpha_1}, \gamma_{\tau_2}^{<\alpha_1} [$ \
\textit{максимальным интервалом монотонности} функции \ $\alpha
S_f^{<\alpha_1}$ \
 ниже\ $\alpha_1$.
\\
\quad %

\noindent 2.\quad $A_{1.1}^{m 1 \vartriangleleft \alpha_1}(\tau_1,
\tau_2, \alpha S_f^{<\alpha_1}, a_f^{<\alpha_1})$:
\\
\[
    \
    A^{0 \vartriangleleft \alpha_1}(\tau_1) \wedge
    A_{1.1}^{m \vartriangleleft \alpha_1}(\tau_1, \tau_2, \alpha S_f^{<\alpha_1})
    \wedge \forall \tau \big( \tau_1 < \tau < \tau_2
    \rightarrow a_{\tau}^{<\alpha_1} = 1 \big);
\]
\quad \\
в добавление  к \ $A_{1.1}^{m \vartriangleleft \alpha_1}$ \ здесь
утверждается, что не существует \ $\alpha$-матриц, допустимых для
 \ $\gamma_{\tau_1}^{<\alpha_1}$ \ ниже \ $\alpha_1$ \ и функция \ $\alpha S_f^{<\alpha_1}$ \ имеет на \
$]\tau_1,\tau_2[$ \ значения \ $\alpha S_{\tau}^{<\alpha_1}$ \
только единичной характеристики \ $a_{\tau}^{<\alpha_1}=1$; \ в
этом случае единичная характеристика \ $a=1$ \ стабилизируется на
интервале \ $[\tau_1,\tau_2[$ \ и на соответствующем интервале \
$[ \gamma_{\tau_1}^{<\alpha_1}, \gamma_{\tau_2}^{<\alpha_1} [$ \
ниже \ $\alpha_1$ \ по лемме~\ref{9.4.}.
\\
\quad \medskip %

\noindent 3.\quad $A_{1.1}^{st \vartriangleleft \alpha_1}(\tau_1,
\tau_{\ast}, \tau_2,\alpha
S_f^{<\alpha_1},a_f^{<\alpha_1})$:\newline
\[
    A_{1.1}^{m 1 \vartriangleleft \alpha_1}(\tau_1, \tau_{\ast},
    \alpha S_f^{<\alpha_1}, a_f^{<\alpha_1})
    \wedge \tau_1 < \tau_{\ast} \le \tau_2 \wedge
    A_{1}^{\vartriangleleft \alpha_1}(\tau_1, \tau_2, \alpha S_f^{<\alpha_1});
\]
\quad \\
здесь говорится, что функция \ $\alpha S_f^{<\alpha_1}$ \
определена на интервале \ $]\tau_1,\tau_2[$, \ но на его начальном
максимальном подинтервале монотонности \ $]\tau_1,\tau_{\ast}[$ \
с \ $\gamma_{\tau_{\ast}} \in SIN_n^{<\alpha_1}$ \ она имеет даже
\textit{единичную} характеристику, стабилизирующуюся на нём;
поэтому этот интервал и соответствующий интервал \ $[
\gamma_{\tau_1}^{<\alpha_1}, \gamma_{\tau_2}^{<\alpha_1} [$ \
будут называться далее (единичными) ступенями ниже \ $\alpha_1$
(как мы увидим, этот термин оправдывается ростом гёделевской
функции \emph{Od} на таких ступенях); в этом случае ординал
\[
    Od \alpha S_f^{<\alpha_1}(\tau_1, \tau_{\ast})
\]
будет называться \textit{высотой} такой ступени.
\\
\quad \medskip %

\noindent 4.\quad $A_{1.1}^{Mst \vartriangleleft \alpha_1}(\tau_1,
\tau_{\ast}, \tau_2,\alpha S_f^{<\alpha_1},a_f^{<\alpha_1})$:
\[
    \qquad
    A_{1.1}^{st \vartriangleleft \alpha_1}(\tau_1, \tau_{\ast}, \tau_2,
    \alpha S_f^{<\alpha_1}, a_f^{<\alpha_1})    \wedge
    A_{1.1}^{M \vartriangleleft \alpha_1}(\tau_1, \tau_2, \alpha S_f^{<\alpha_1});
\]
\quad \\
в дополнение здесь указывается, что интервал \ $] \tau_1, \tau_2[$
\ -- максимальный с\ \mbox{$\gamma_{\tau_2}^{<\alpha_1} \in
SIN_n^{<\alpha_1}$}, \ поэтому мы будем называть интервал \ $[
\tau_1, \tau_2 [$ \ и соответствующий интервал \ $[
\gamma_{\tau_1}^{<\alpha_1}, \gamma_{\tau_2}^{<\alpha_1} [$ \
\textit{максимальными (единичными) ступенями} ниже \ $\alpha_1$.

Эти наблюдения приводят к следующему понятию \emph{лестницы}:
\\
\quad \\
5.\quad $A_{8}^{\mathcal{S} t \vartriangleleft \alpha_1}
(\mathcal{S} t, \chi, \alpha S_f^{<\alpha_1}, a_f^{<\alpha_1})$:
\[
    (\mathcal{S} t \mbox{\it \ \ это функция на} \chi^{\ast +})
    \wedge \qquad\qquad\qquad\qquad\qquad\qquad\qquad\qquad
\]
\[
    \wedge \forall \beta < \chi^{\ast +} ~ \exists \tau_1, \tau_{\ast},
    \tau_2 \big( \mathcal{S} t(\beta) = (\tau_1,\tau_{\ast},\tau_2)
    \wedge
    \qquad\qquad\qquad\qquad
\]
\[
    \qquad\qquad\qquad
    \wedge A_{1.1}^{M st \vartriangleleft \alpha_1}(\tau_1,\tau_{\ast},\tau_2,
    \alpha S_f^{<\alpha_1}, a_f^{<\alpha_1}) \wedge
\]
\[
    \wedge \forall \tau_1, \tau_{\ast}, \tau_2 \bigl(
    A_{1.1}^{M st \vartriangleleft \alpha_1}(\tau_1,\tau_{\ast},\tau_2,
    \alpha S_f^{<\alpha_1}, a_f^{<\alpha_1})
    \longrightarrow
    \qquad\qquad\qquad
\]
\[
    \qquad\qquad\qquad\qquad\qquad
    \longrightarrow
    \exists \beta < \chi^{\ast +} ~ \mathcal{S}t(\beta) =
    (\tau_1, \tau_{\ast}, \tau_2) \bigr) \wedge
\]
\[
    \wedge \forall \beta_1, \beta_2 < \chi^{\ast +}
    ~ \forall \tau_1^{\prime}, \tau_{\ast}^{\prime}, \tau_2^{\prime}
    ~ \forall \tau_1^{\prime\prime}, \tau_{\ast}^{\prime\prime},
    \tau_2^{\prime\prime} \big( \beta_1 < \beta_2
    \wedge \qquad\qquad\qquad
\]
\[
    \wedge \mathcal{S} t(\beta_1) = (\tau_1^{\prime},
    \tau_{\ast}^{\prime}, \tau_2^{\prime})
    \wedge \mathcal{S} t(\beta_2) = (\tau_1^{\prime\prime},
    \tau_{\ast}^{\prime\prime}, \tau_2^{\prime\prime} ) \rightarrow
    \tau_2^{\prime} < \tau_1^{\prime\prime} \wedge
\]
\[
    \qquad\qquad\qquad
    \wedge Od~\alpha S_f^{<\alpha_1}(\tau_1^{\prime}, \tau_{\ast}^{\prime})
    < Od~\alpha S_f^{<\alpha_1}(\tau_1^{\prime\prime}, \tau_{\ast}^{\prime\prime})
    \big) \wedge
\]
\[
    \wedge \sup \big\{ Od~\alpha S_f^{<\alpha_1}(\tau_1, \tau_{\ast}) : \exists
    \beta, \tau_2 ~ \mathcal{S} t(\beta) = (\tau_1,\tau_{\ast},\tau_2) \big\} =
    \chi^{\ast +};
\]
здесь указывается, что \ $\mathcal{S} t$ \ это функция на \
$\chi^{\ast +}$ \ и что её значения -- это все тройки \
$(\tau_1,\tau_{\ast},\tau_2)$ \ такие, что интервалы \ $[ \tau_1,
\tau_2 [$ \ являются максимальными единичными ступенями,
расположенными последовательно, одна за другой. Поэтому такая
функция \ $\mathcal{S} t$ \ будет называться \textit{лестницей},\
а интервалы \ $[ \tau_1, \tau_2 [$ \ и соответствующие интервалы
 \ $[ \gamma_{\tau_1}^{<\alpha_1},
\gamma_{\tau_2}^{<\alpha_1} [$ будут называться её ступенями ниже
\ $\alpha_1$.
\\
Это понятие оправдывается строгим возрастанием высот таких
ступеней; мы будем также говорить, что эта лестница \ $\mathcal{S}
t$ \ состоит из этих ступеней, или содержит эти ступени.
\\
Соответственно этому, ординал
\[
    h(\mathcal{S} t) = \sup \big \{ Od~\alpha S_f^{<\alpha_1}(\tau_1,\tau_{\ast}):
    \exists \beta, \tau_2 ~ \mathcal{S} t(\beta) =
    (\tau_1, \tau_{\ast}, \tau_2) \big \}
\]
будет называться \textit{высотой} всей лестницы \ $\mathcal{S} t$.
\ Таким образом, здесь требуется, чтобы высота всей лестницы \
$\mathcal{S} t$ \ возрастала до \ $\chi^{\ast +}$, \ то есть чтобы

\[
    h(\mathcal{S} t) = \chi^{\ast +}.
\]
Также кардинал
\[
    \upsilon = \sup \left\{ \gamma_{\tau_2}: \exists \beta, \tau_1,
    \tau_{\ast} ~ \mathcal{S} t(\beta) = (\tau_1, \tau_{\ast}, \tau_2)
    \right\}
\]
будет называться \emph{завершающим} кардиналом лестницы \
$\mathcal{S} t$ \ и будет обозначаться через
\[
    \upsilon(\mathcal{S} t);
\]
так что мы будем говорить, что лесница  \ $\mathcal{S} t$ \
завершается в этом кардинале \ $\upsilon(\mathcal{S} t)$.

\noindent Если такая лестница \ $\mathcal{S} t$ \ существует ниже
кардинала \ $\alpha_1$, \ то мы будем говорить, что этот \
$\alpha_1$ \  \textit{обладает} этой лестницей \ $\mathcal{S} t$.
\\
Когда кардинал \ $\alpha > \chi^{\ast}$ \ это носитель матрицы \
$S$ \ и его кардинал предскачка \ $\alpha_1 =
\alpha_{\chi{\ast}}^{\Downarrow}$ \ после \ $\chi^{\ast}$ \
обладает некоторой лестницей \ $\mathcal{S} t$, \ то мы будем
говорить, что эта \ $S$ \ на \ $\alpha$ \ обладает этой лестницей.

И вот очень важный пример \emph{внутреннего свойства} матрицы:
\textit{внутреннее} свойство обладания матрицы \ $S$ \ некоторой
лестницей.
\\
Это свойство для матрицы \ $S$ \ на её носителе \ $\alpha$ \
определяется формулой
\[
    \exists \mathcal{S}t \vartriangleleft \alpha^{\Downarrow +} ~
    A_{8}^{\mathcal{S} t \vartriangleleft \alpha^{\Downarrow} }
    (\mathcal{S} t, \alpha S_f^{< \alpha^{\Downarrow}},
    a_f^{< \alpha^{\Downarrow}}),
\]
которая может быть ограничена кардиналом скачка \
$\alpha^{\downarrow}$ \ её носителя \ $\alpha$ \ после \
$\chi^{\ast}$. \ Поэтому согласно лемме 5.11~\cite{Kiselev11} тем
же свойством матрица \ $S$ \ обладает на любом другом своём
носителе \ $\alpha^1 > \chi^{\ast}$:
\[
    \exists \mathcal{S}t^1 \vartriangleleft \alpha^{1 \Downarrow +} ~
    A_{8}^{\mathcal{S} t \vartriangleleft \alpha^{ 1 \Downarrow} }
    (\mathcal{S} t^1, \alpha S_f^{< \alpha^{1 \Downarrow}},
    a_f^{< \alpha^{1 \Downarrow}} ),
\]
которое теперь может быть ограничено кардиналом скачка \
$\alpha^{1 \downarrow}$ \ носителя  \ $\alpha^{1}$ \ после \
$\chi^{\ast}$, \ и таким образом \ $S$ \ на \ $\alpha^1$ \ снова
обладает некоторой лестницей \ ${\mathcal{S}t}^1$ \ как и раньше.

\begin{lemma}
\label{9.5.} \emph{(О срезании лестницы сверху)}  \\
\hspace*{1em} Пусть
\medskip

(i) \quad $A_{1}^{1 \vartriangleleft \alpha_1}(\tau _{1},\tau
_{2})$;
\newline

(ii) \quad $\tau_2 \le \tau_3$ \ и \ $S^3$ \ это матрица
характеристики \ $a^3$ \  на носителе
\[
    \alpha_3 \in \; ]\gamma_{\tau_3}^{<\alpha_1}, \alpha_1 [
\]
с диссеминатором \ $\widetilde{\delta}^3$ \ и базой данных \
$\rho^3$, \ допустимые для \ $\gamma_{\tau_3}^{<\alpha_1}$ \ ниже
\ $\alpha_1$, \ с производящим собственным диссеминатором \
$\check{\delta}^{S^3}$ \ на \ $\alpha^3$;
\newline

(iii) \quad $\forall \tau <\tau _{2}\quad \exists \tau ^{\prime
}\in \left[ \tau ,\tau _{2}\right[ \quad a_{\tau ^{\prime
}}^{<\alpha _{1}} = a^3$. \medskip

\noindent Тогда \medskip

1. \ $Od\alpha S_{f}^{<\alpha _{1}}(\tau _{1},\tau _{2}) <
Od(S^3)$;
\\

2a. \ таким образом, если \ $a^3 = 1$, \ то не существует лестницы
ниже \ $\alpha_1$, \ завершающейся в некотором \
$SIN_n^{<\alpha_1}$-кардинале \ $\upsilon <
\alpha_3^{\Downarrow}$;
\\

2b. \ следовательно, если существует некоторая единичная мактрица
\ $S^0$ \ на её носителях ниже кардинала \ $\alpha_1$, \
расположенных конфинально этому кардиналу \ $\alpha_1$:
\vspace{-6pt}
\begin{multline*}
    \forall \gamma < \alpha_1 ~ \exists \gamma^1 \in \; ]\gamma,
    \alpha_1[ ~ \exists \delta, \alpha, \rho < \alpha_1 ~ \big(
    SIN_{n-1}^{<\alpha_1}(\gamma^1) \wedge
\\
    \wedge \alpha \mathbf{K}^{<\alpha_1}(1, \delta, \gamma^1, \alpha, \rho,
    S^0) \big),
\end{multline*}
то этот кардинал \ $\alpha_1$ \ не обладает никакой лестницей;
\\

3. \ если \ $S^3$ \ это матрица, \
$\underline{\lessdot}$-минимальная изо всех матриц той же
характеристики \ $a^3$ \ на носителях\ $\in \; ]
\gamma_{\tau_3}^{<\alpha_1}, \alpha_1[$ \ допустимых для \
$\gamma_{\tau_3}^{<\alpha_1}$, \ то
\[
    \gamma_{\tau _{2}}^{<\alpha _{1}}<\check{\delta}^{S^3} \leq
    \widetilde{\delta }^3 < \gamma _{\tau_{3}}.
\]
\end{lemma}

\noindent \textit{Доказательство.} \ Проведём это рассуждение для
\ $a^3=1$, \ именно этот случай используется в дальнейшем; в этом
важном случае\ $\widetilde{\delta}^3 = \check{\delta}^{S^3}$ \ и
условие $(iii) $ можно ослабить до \ $a^3=1$ \ по лемме
3.2~\cite{Kiselev11}. В этом случае условие неподавления \ $\neg
A_5^{S,0}$ \ для единичной матрицы \ $S$ \ на её носителях может
быть опущено, потому что такая матрица \ $S$ \ всегда неподавлена
и формулы \ $\alpha \mathbf{K}^{\ast < \alpha_1}$, $\alpha
\mathbf{K}_{n+1}^{\ast < \alpha_1}$ \ эквивалентны формулам \
$\alpha \mathbf{K}^{< \alpha_1}$, $\alpha \mathbf{K}_{n+1}^{<
\alpha_1}$; \ верхние индексы \ $< \alpha_{1} $, $\vartriangleleft
\alpha_{1}$ \ будут опускаться, как всегда.
\\
По этому условию этой леммы матрица \ $S^3$ получает единичную
характеристику на её допустимых носителях, расположенных
конфинально кардиналу \ $ \gamma_{\tau_{2}}$, \ как это было в
доказательстве части 1а. леммы~\ref{9.4.},  где  \
$\gamma_{\tau^2}$, \ $S^0$ \ на \ $\alpha_{\tau^1}$ \ следует
заменить на \ $\gamma_{\tau_2}$, \ $S^3$ \ на $\alpha_3$. \ Таким
образом,  $(i)$ и лемма~\ref{9.4.} влекут
\[
    Od \alpha S_{f} (\tau_{1}, \tau_{2}) \leq Od (S^3)
    \mbox{\it \ \  и \ \ } a_{\tau} \equiv 1 \mbox{\it \ \ на \ \ }
    ]\tau_{1}, \tau_{2} [ ~.
\]
Теперь предположим, что функция \ $\alpha S_f$ \ стабилизируется
на \ $[ \tau_1, \tau_2[$ \ и пусть \ $\tau_{2}^s $ \ -- это
стабилизационный ординал для \ $\alpha S_{f} $ \ на \ $[\tau_{1},
\tau_{2} [ $, \ так что существует матрица  \ $S^0 $ \ такая, что
\[
    \alpha S_{\tau} \equiv \alpha S_{\tau_2^s} = S^{0}
    \mbox{\it \ \ на \ \ } [ \tau_2^s, \tau_2[ ~.
\]
Мы применим сейчас  метод ограничения-и-продолжения, который уже
был использован выше несколько раз. Для этого вернёмся к матрице \
$S^{0} $ \ на носителе \ $\alpha_{\tau_2^s+1} $ \ с кардиналом
предскачка \ $\alpha^{1} = \alpha_{\tau_2^s+1}^{\Downarrow} $ \ и
диссеминатором \ $\check{\delta}^{1} =
\check{\delta}_{\tau_2^s+1}$. \ Та же самая матрица \ $S^{0} $ \
на носителе \ $\alpha_{\tau_2^s} $ \ единичной характеристики по
лемме 3.2~\cite{Kiselev11} об ограничении получает единичные
характеристики на её допустимых носителях, расположенных
конфинально кардиналу  \ $\check{\delta}^{1} $ \ и поэтому ниже \
$\check{\delta}^{1} $ \ выполняется следующее \ $\Pi_{n+1}
$-утверждение для \ $S=S^{0}$:
\begin{equation} \label{e9.11}
    \qquad\qquad
    \forall \gamma \exists \gamma^{\prime} ~ \bigl( \gamma <
    \gamma^{\prime} \wedge SIN_{n-1}(\gamma^{\prime}) \wedge
    \qquad\qquad\qquad\qquad\qquad\qquad\qquad
\end{equation}
\[
    \qquad\qquad\qquad
    \wedge \exists \delta, \alpha, \rho ~~ \alpha\mathbf{K} (1,
    \delta, \gamma^{\prime}, \alpha, \rho, S) \bigr ) ~.
\]
\vspace{-6pt}

\noindent Диссеминатор \ $\check{\delta}^{1} $ \ продолжает его до
\ $ \alpha^{1} $ \ и поэтому матрица \ $S^{0} $ \ получает
единичную характеристику на её допустимых носителях, расположенных
конфинально кардиналу \ $\alpha^{1} $, \ то есть (\ref{e9.11})
выполняется  матрицей \ $S = S^{0} $ \ при ограничении \ $
\vartriangleleft \alpha^{1} $. \ После минимизации таких матриц \
$S $ \ мы получаем матрицу \ $S=S^{1} $ \ со свойством
(\ref{e9.11}) ниже \ $ \alpha^{1} $ \ и по лемме
4.6~\cite{Kiselev11} \ $S^{1} \lessdot S^{0} $. \ Нужно отметить,
что утверждение (\ref{e9.11}), \ $\vartriangleleft$-ограниченное
кардиналом\ $\alpha^1$ \ с \ $S = S^1$, \ это \textit{внутреннее
свойство} матрицы \ $S^0$.
\\
Если теперь

\[
    Od \alpha S_{f} (\tau_{1},\tau_{2}) = Od ( S^3 ), \mbox{\it \
    \ то есть \ \ } S^{0} = S^3,
\]
то матрица \ $S^{1} $ \ по лемме 5.11~\cite{Kiselev11} об
информативности получает свои допустимые носители той же единичной
характеристики, расположенными конфинально кардиналу предскачка \
$\alpha^{3}=\alpha_3^{\Downarrow} $, \ так как \ $S^{0} $ \ на
носителе \ $\alpha_{\tau_2^s+1} $ \ обладает тем же свойством.
\\
После этого снова по лемме 3.2~\cite{Kiselev11} такие носители
появляются, расположенные конфинально \ $\gamma _{\tau_{2}} $. \
Таким бразом, наконец, благодаря $(i)$ появляется противоречие:

\begin{equation}  \label{e9.12}
Od \alpha S_{f} (\tau_{1},\tau_{2})\leq Od (S^{1})< Od (S^{0}).
\end{equation}
\vspace{0pt}

\noindent Если же функция \ $\alpha S_f$ \ не стабилизируется на \
$[\tau_1, \tau_2[$, то ординал
\[
    \rho=Od\alpha S_f(\tau_1, \tau_2)
\]
-- предельный. Теперь для того, чтобы закончить доказательства
утверждения 1., можно заметить, что гёделевская конструктивная
функция \ $F$ \ принимает значения \ $F (\alpha ) = F | \alpha$ \
для предельных ординалов \ $\alpha$ \ (см. Гёдель \cite{Godel}).
Нетрудно видеть, что \ $Od(S)$ \ не может быть предельным и
поэтому \ $\rho \leq Od(S) $ \ влечёт \ $\rho < Od(S)$.
\\

\noindent Обращаясь к утверждению 3., допустим, что оно неверно и
\[
    \check{\delta}^{3} = \check{\delta}^{S^3} \leq \gamma_{\tau_{2}}
\]
и, стоя на кардинале \ $\alpha^3 = \alpha_3^{\Downarrow}$, \
рассмотрим получающуюся ситуацию ниже \ $ \alpha^{3} $ .\ Здесь
следует рассмотреть два случая:
\\
Случай 1. \quad $\gamma_{\tau_{1}} < \check{\delta}^{3} \leq
\gamma_{\tau_{2}} $.
\\
Так как \ $\gamma_{\tau_{2}} \in SIN_{n} $ \ и \
$\check{\delta}^{3} \in SIN_{n}^{< \alpha^{3}} $, \ то из лемм
~3.8~\cite{Kiselev11}, 8.5~1), 8.7, 9.4 следует, что \ $\check{
\delta}^{3} \in SIN_{n} $ \ и
\[
    \alpha S_{\tau}^{< \alpha^{3}} \equiv \alpha S_{\tau}
    \mbox{\it \quad на \quad}  \{
    \tau : \gamma_{\tau_{1}} < \gamma_{\tau} < \check{\delta}^{3} \}
\]
и затем диссеминатор \ $\check{\delta}^{3} $ \ продолжает до \
$\alpha^{3} $ \ следующее \ $\Pi_{n+1} $-утверждение об
определённости функции \ $ \alpha S_{f} $ \ \textit{единичной}
характеристики со значениями \ $\lessdot S^3$ \ благодаря части
1.:
\begin{equation}  \label{e9.13}
    \qquad
    \forall \gamma^{\prime} ~ \Big( \gamma_{\tau_1} < \gamma^{\prime} \wedge
    SIN_{n-1}(\gamma^{\prime}) ~ \longrightarrow
    \qquad \qquad \qquad \qquad \qquad \qquad \qquad
\end{equation}
\[
    \qquad
    \longrightarrow \exists \delta, \alpha, \rho, S ~~
    \big ( S \lessdot S^3 \wedge
    \alpha \mathbf{K}(1, \delta, \gamma^{\prime}, \alpha, \rho, S) \big) \Big)
\]
\vspace{0pt}

\noindent и поэтому существует  допустимая матрица \ $\alpha
S_{\tau_{3}}^{< \alpha^{3}} $ \ на носителе

\[
    \alpha_{\tau_3}^{<\alpha^3} \in \; ] \gamma_{\tau_3}, \alpha^3 [
\]
единичной характеристики и со значением \ $\alpha S_{\tau_{3}}^{<
\alpha^{3}} \lessdot S^3 $ \ вопреки \ $\underline{\lessdot}
$-минимальности самой матрицы \ $ S^3$ \ на \ $\alpha_3$. \
Осталось рассмотреть следующий
\\

\noindent Случай 2. \quad $\check{\delta}^{3} \leq
\gamma_{\tau_{1}} $.
\\
Нужно отметиь, что условие \ $\underline{\lessdot} $-минимальности
матрицы \ $S^3$ \ не используется в этом случае. Здесь матрица \
$S^0 = \alpha S_{\tau_1 + 2}$ \ единичной характеристики должна
рассматриваться на носителе \ $\alpha_{\tau_1 + 2}$ \ с кардиналам
предскачка \ $\alpha^1 = \alpha_{\tau_1 + 2}^{\Downarrow}$ \ и
диссеминатором \ $\check{\delta}^1 = \check{\delta}_{\tau_1 + 2}$,
\ как это было сделано выше для \ $S^0 = \alpha S_{\tau_2^s + 1}$,
\ $\alpha^1 = \alpha_{\tau_2^s + 1}^{\Downarrow}$, \
$\check{\delta}^1 = \check{\delta} _{\tau_2^s + 1}$ \ в
доказательстве утверждения 1. (договоримся сохранить прежние
обозначения для удобства). И снова матрица \ $\alpha S_{\tau_{1} +
1} $ \ получает единичные характеристики на своих допустимых
носителях расположенных конфинально диссеминатору\
$\check{\delta}^1 = \gamma_{\tau_1}$ \ и он продолжает утверждение
(\ref{e9.11}) для \ $S = \alpha S_{\tau_1 + 1}$ \ до \ $\alpha^1$;
\ таким образом, это вызывает появление минимальной матрицы \ $S^1
\lessdot S^0$ \ с прежними свойствами: она получает единичную
характеристику на её допустимых носителях, расположенных
конфинально кардиналу \ $\alpha^1$.
\\
По лемме 3.2~\cite{Kiselev11} появляются носители матрицы \ $S^{1}
$ \ единичной характеристики, расположенные конфинально
диссеминатору \ $\check{\delta}^3$, \ то есть (\ref{e9.11})
выполняется для \ $S=S^{1} $ \ при \
$\vartriangleleft$-ограничении этим диссеминатором \
$\check{\delta }^3$; \ поэтому \ $\check{\delta}^{3} $ \
продолжает это утверждение до \ $\alpha^{3} $. \ После этого
кардинал \ $ \gamma_{\tau_{2}} \in SIN_{n}$ \ ограничивает это
утверждение с \ $\gamma $ \ заменённым на произвольный постоянный
кардинал \ $\gamma^1 < \gamma_{\tau_{2}} $. \ В результате матрица
\ $S^{1} $ \ получает единичную характеристику на допустимых
носителях, расположенных конфинально кардиналу \
$\gamma_{\tau_{2}} $ \ и мы снова приходим к противоречию
(\ref{e9.12}).

Обращаясь к характеристике \ $\alpha^3=0$ \ следует повторить это
рассуждение, но для \textit{нулевых} матриц \ $S$ \ на их
носителях \ $\alpha$, \ допустимых для рассматриваемых кардиналов
\ $\gamma_{\tau}$, \ но только для \ $\gamma_{\tau} \notin SIN_n$.
\ Во всех подобных случаях такие матрицы \ $S$ \ неподавляются по
определению и снова поэтому условие неподавления \ $\neg
A_5^{S,0}$ \ может быть опущено и снова формулы \ $\alpha
\mathbf{K}^{\ast < \alpha_1}$, $\alpha \mathbf{K}_{n+1}^{\ast <
\alpha_1}$ \ можно заменить на формулы \ $\alpha \mathbf{K}^{<
\alpha_1}$, $\alpha \mathbf{K}_{n+1}^{< \alpha_1}$. \ Именно такие
матрицы \ $S$ \ на их носителях \ $\alpha$ \ следует использовать
в методе ограничения-и-подавления, применённом выше, что
обеспечивает доказательство для характеристики \ $a_3=0$.

И, наконец, утверждение 2а. следует из 1. почти очевидно.
Предположим, оно не выполняется, то есть существует некоторая
лестница \ $\mathcal{S}t$ \ ниже \ $\alpha_1$, \ завершающаяся в \
$SIN_n$-кардинале \ $\upsilon(\mathcal{S}t) <
\alpha_3^{\Downarrow}$; \ отсюда следует, что
\[
    \upsilon(\mathcal{S}t) < \gamma_{\tau_3}^{<\alpha_1}.
\]
Дальнейшее рассуждение можно представить себе как \emph{``срезание
этой лестницы сверху''} матрицей \ $S^3$ \ и тем самым исключение
этой лестницы:
\\
По определению эта лестница состоит из \textit{единичных} ступеней
ниже \ $\alpha_1$
\[
    \mathcal{S}t(\tau^{\prime}) = (\tau_1^{\prime},
    \tau_2^{\prime\ast}, \tau_2^{\prime})
\]
и каждая из них обладает свойством
\[
    A_1^{1 \vartriangleleft \alpha_1} (\tau_1^{\prime},
    \tau_2^{\prime\ast})
\]
с единичной характеристикой, стабилизирущейся на \ $[
\tau_1^{\prime}, \tau_2^{\prime\ast} [$ \ (см. определение
лестницы перед леммой \ref{9.5.}). Поэтому доказанное утверждение
1. влечёт, что высота всех её ступеней, то есть высота  \
$h(\mathcal{S}t)$ \ всей лестницы \ $\mathcal{S}t$, \ ограничена
сверху  ординалом
\[
    Od(S^3) < \chi^{\ast +}
\]
(``срезается'' этим ординалом), хотя по определению  лестницы \
$\mathcal{S}t$ \ высота её ступеней возрастает до \ $\chi^{\ast
+}$, \ то есть \ $h(\mathcal{S}t) = \chi^{\ast +}$. \
\\
Это же рассуждение составляет и доказательство утверждения 2b.,
если матрицу \ $S^0$ \ использовать вместо матрицы \ $S^3$.
\\
\hspace*{\fill} $\dashv$

 Отметим, что для характеристики \ $a^3 =
1$ \ условие минимальности матрицы \ $S^3$ \ в пункте 3. это леммы
~\ref{9.5.} можно опустить, используя немного изменённые аргументы
из обсуждения случая 1.
\\

Следующее очевидное следствие показывает, что высоты таких
ступеней строго возрастают:

\begin{corollary}
\label{9.6.} \hfill {} \newline \hspace*{1em} Пусть \medskip

(i) \quad $A_{1}^{1 \vartriangleleft \alpha_1} ( \tau _{1},\tau
_{2})$, \ $A_{1}^{1 \vartriangleleft \alpha_1} ( \tau _{3},\tau
_{4}) $, \ $\tau _{2}<\tau _{4}$; \newline

(ii)\quad $\forall \tau <\tau _{4} \quad \exists \tau ^{\prime }\in \left[
\tau ;\tau _{4}\right[ \quad a_{\tau ^{\prime }}^{<\alpha _{1}}=1$. \medskip

\noindent Тогда \medskip

1) \ $\forall \tau \in \; ]\tau_1, \tau_2[ \; \cup \; ]\tau_3, \tau_4[ \quad
a_\tau^{<\alpha_1} = 1$; \newline

2) \ $Od\alpha S_{f}^{<\alpha _{1}} ( \tau _{1},\tau _{2} ) <Od\alpha
S_{f}^{<\alpha _{1}} ( \tau _{3},\tau _{4} )$; \newline

3) \ $\forall \tau \in \; ]\tau_3, \tau_4[ \quad \gamma _{\tau
_{2}}^{<\alpha_{1}} < \check{\delta}_{\tau}^{S} = \widetilde{
\delta}_{\tau}^{ < \alpha _{1}} < \gamma^{ < \alpha _{1}}_{\tau}$,
\\
где \ $\check{\delta}_{\tau}^{S}$ \  -- это производящий
собственный диссеминатор матрицы \ $\alpha S_{\tau}^{< \alpha_1}$
\ на \ $\alpha_{\tau}^{< \alpha_1}$; \ поэтому
\[
    \gamma_{\tau_2}^{<\alpha_1} < \gamma_{\tau_3}^{<\alpha_1}.
\]
\end{corollary}

\noindent \textit{Доказательство.} \ Из условий $(i)$, $(ii)$ и
лемм 3.2~\cite{Kiselev11}, \ref{9.4.} следует, что \
$a_{\tau}^{<\alpha_1} \equiv 1$ \ на интервалах \ $]\tau_1,
\tau_2[$, \ $]\tau_3, \tau_4[$\;. \ Поэтому лемма~\ref{9.5.} (где
\ $\tau_3$ \ играет роль любого \ $\tau \in \; ]\tau_3,
\tau_4[$\;) влечёт утверждения~2),~3). Для \ $\tau = \tau_3 + 1$ \
тогда получаается \ $\gamma_{\tau_2}^{<\alpha_1} <
\check{\delta}_{\tau}^{<\alpha_1}$ \ и в то же время согласно
леммам~\ref{8.5.}~7)~$(ii)$, \ref{8.8.}~2) \ -- \
$\check{\delta}_{\tau}^{S} = \gamma_{\tau_3}^{<\alpha_1}$; \ таким
образом \ $\tau_2 < \tau_3$.
\\
\hspace*{\fill} $\dashv$

\begin{corollary}
\label{9.7.} \hfill {} \newline \hspace*{1em} Пусть \medskip

(i) \quad $A_{1}^{1 \vartriangleleft \alpha_1}(\tau _{1},\tau
_{2})$; \newline

(ii) \quad $\tau _{3}\in dom(\alpha S_{f}^{<\alpha _{1}})$, \
$\tau _{3}\geq \tau _{2}$;
\\

(iii) \quad матрица \ $\alpha S_{\tau_3}^{< \alpha_1}$ \ на \
$\alpha_{\tau_3}^{< \alpha_1}$ \ имеет производящий собственный
диссеминатор
\[
    \check{\delta}_{\tau _{3}}^{S}\leq
    \gamma _{\tau _{2}}^{<\alpha _{1}}
\]
ниже \ $\alpha_1$.
\\
Тогда \medskip

1) \ $a_{\tau }^{<\alpha _{1}}\equiv 1$ \ on \ $\left] \tau
_{1},\tau _{2} \right[ $, \quad $a_{\tau _{3}}^{<\alpha _{1}}=0$;
\newline

2) \ $\check{\delta}_{\tau _{3}}^S \leq \gamma _{\tau
_{1}}^{<\alpha _{1}}$ \quad и \quad \newline

3) \ $Od\alpha S_{f}^{<\alpha _{1}}(\tau _{1},\tau _{2})>Od(\alpha
S_{\tau _{3}}^{<\alpha _{1}})$.
\\

\noindent Аналогично для производящего диссеминатора \
$\check{\delta}_{\tau _{3}}^{<\alpha_1} $ \ матрицы \ $\alpha
S_{\tau _{3}}^{<\alpha_1}$ \ на \ $\alpha_{\tau
_{3}}^{<\alpha_1}$.
\end{corollary}

\noindent \textit{Доказательство.} \ Мы будем опускать верхние
индексы \ $< \alpha_{1} $, $\vartriangleleft \alpha_{1}$. \ По
лемме~\ref{9.5.} для \ $S^3 = \alpha S_{\tau_3}^{< \alpha_1}$ \
условие $(iii)$ влечёт, что для некоторого \ $\tau < \tau_{2}$
\begin{equation*}
\forall \tau^{\prime} \in [\tau, \tau_{2} [ \quad  a_{\tau^{\prime}} \neq
a_{\tau_{3}} ;
\end{equation*}
благодаря лемме 3.2~\cite{Kiselev11} это возможно только когда
\begin{equation*}
\forall \tau^{\prime} \in [\tau, \tau_{2} [ \quad
a_{\tau^{\prime}} =1, \quad a_{\tau_{3}} = 0,
\end{equation*}
а тогда по лемме 9.4 \ $a_{\tau} \equiv 1 $ \ on \ $]\tau_{1},
\tau_{2} [ $\;.
\\
Если \ $\check{\delta}_{\tau_{3}}^S \in \; ] \gamma_{\tau_{1}},
\gamma_{\tau_{2}}[\;$, \ то можно получить \ $a_{\tau_{3}} =1 $, \
используя снова рассуждение из доказательства леммы~\ref{9.5.}, и
продолжая утверждение~(\ref{e9.13}) без его подформулы \ $S
\lessdot S^3$ \ диссеминатором \ $\check{\delta}_{\tau_{3}}^S$ \
до \ $\alpha_{\tau_{3}}^{\Downarrow}$, \ а это влечёт \
$a_{\tau_3}^{< \alpha_1} = 1$.
\\
После этого достаточно провести рассуждение из конца этого
доказательства (случай 2.), буквально повторяя его посредством
метода ограничения-и-повторения.
\\
\hspace*{\fill} $\dashv$
\\
\quad \\
Немедленное следствие этой леммы для \ $\tau_{2}= \tau_{3}$ \
составляет следующая
\newline
\quad \\
\quad \\
\quad \\
\quad \\
 \noindent \textbf{Теорема 1}\quad \newline
\emph{\hspace*{1em} Пусть
\medskip }

\emph{(i) \quad $\alpha S_f^{<\alpha_1}$ \ монотонна на \
$[\tau_1, \tau_2[ $ \ ниже \ $\alpha_1$;
\newline }

\emph{(ii) \quad $\tau_1 = \min \{ \tau: \; ]\tau, \tau_2[ \; \subseteq dom
(\alpha S_f^{<\alpha_1}) \} $. \medskip }

\emph{\noindent Тогда
\begin{equation*}
]\gamma_{\tau_1}^{< \alpha_1}, \gamma_{\tau_2}^{< \alpha_1} [  ~
\cap ~ SIN_n^{<\alpha_1} = \varnothing.
\end{equation*}
\vspace{0pt} }

\noindent \textit{Доказательство.} \ Предположим, что наоборот --
существует \ $SIN_n^{<\alpha_1}$-кардинал \
$\gamma_{\tau_2^{\prime}}^{< \alpha_1} \in ]\gamma_{\tau_1}^{<
\alpha_1}, \gamma_{\tau_2}^{< \alpha_1} [ $~. \newline Тогда
кардинал \ $\gamma_{\tau_1}^{< \alpha_1}$ \ также принадлежит
классу \ $SIN_n^{< \alpha_1}$; \ это можно увидеть, повторяя
доказательство леммы~\ref{8.10.}. Таким образом, выполняется
утверждение \ $A_1^{1 \vartriangleleft \alpha_1}(\tau_1,
\tau_2^\prime)$; \ остаётся применить следствие~\ref{9.7.},
используя \ $\tau_2^{\prime}$ \ как \ $\tau_2 = \tau_3$, \ так как
\ $\check{\delta}_{\tau_2^\prime}^{S} <
\gamma_{\tau_2^\prime}^{<\alpha_1}$ \ по определению.
\\
\hspace*{\fill} $\dashv$

\newpage

\section{Анализ немонотонности \ $\protect\alpha $\,-функции}

\setcounter{equation}{0}

\hspace*{1em} Итак, всякий интервал монотонности \
$\alpha$-функции не может быть ``слишком длинным'' по теореме 1.
\\
Однако такая функция может быть определена на ``довольно длинных''
интервалах; например, функция \ $\alpha S_f^{<\alpha_1}$ \ может
быть определена на заключительном промежутке \ $T^{\alpha_1}$ \
для всякого достаточно большого \ $SIN_n$-кардинала \ $\alpha_1 <
k$ \ (лемма 8.9\;). Следовательно, её монотонность нарушается на
некоторых ординалах из такого промежутка.
\newline Как это случается? В этом разделе анализируются все
существенные случаи подобных нарушений. Для этого нужно вспомнить
формулу \ $A_2^{\vartriangleleft \alpha_1}(\tau_1, \tau_2,
\tau_3)$ \ (см. определение~\ref{8.1.}~1.4 для \ $X_1 = \alpha
S_f^{<\alpha_1}$): \vspace{-6pt}
\begin{multline*}
    A_1^{\vartriangleleft \alpha_1}(\tau_1, \tau_3, \alpha
    S_f^{<\alpha_1}) \wedge \tau_1 + 1 < \tau_2 <  \tau_3 \wedge
\\
    \wedge \tau_2 = \sup \bigl \{ \tau < \tau_3 :  \forall
    \tau^\prime, \tau^{\prime\prime} (\tau_1 < \tau^\prime <
    \tau^{\prime\prime} < \tau \rightarrow  \alpha
    S_{\tau^\prime}^{<\alpha_1} \underline{\lessdot} \alpha
    S_{\tau^{\prime\prime}}^{<\alpha_1} \bigr ) \}.
\end{multline*}

\noindent Здесь указывается, что \ $\tau_2$ \ -- это
\textit{минимальный} ординал, на котором нарушается монотонность
функции \ $\alpha S_f^{<\alpha_1}$ \ на интервале \ $ [ \tau_1,
\tau_3 [$\;. \ Таким образом, во всех рассуждениях этого параграфа
рассматривается некоторая немонотонность \ $A_2^{\vartriangleleft
\alpha_1}(\tau_1, \tau_2, \tau_3)$ \ на интервалах \ $[\tau_1,
\tau_3[$ \ в различных ситуациях (но условие \
$SIN_n^{\vartriangleleft \alpha_1}(\gamma_{\tau_3})$ \ может быть
опущено везде, кроме последней леммы~\ref{10.5.}\;).

\begin{lemma}
\label{10.1.} \hfill {}
\\
Пусть\medskip

(i) \ $A_{2}^{\vartriangleleft \alpha_1}(\tau _{1},\tau _{2},\tau
_{3})$;
\newline

(ii) \ $SIN_{n}^{<\alpha _{\tau _{2}}^{\Downarrow }}\cap \gamma _{\tau
_{2}}^{<\alpha _{1}}\subseteq SIN_{n}^{<\alpha _{1}}$. \medskip

\noindent Тогда \medskip

1) \ $a_{\tau }^{<\alpha _{1}}\equiv 1$ \ на \ $] \tau _{1},
\tau_{2} [ $, \ $a_{\tau _{2}}^{<\alpha _{1}}=0$ \quad и
\newline

2) \ $\widetilde{\delta }_{\tau _{2}}^{<\alpha _{1}} \le \gamma
_{\tau _{1}}^{<\alpha _{1}}$. \label{c13}
\endnote{
\ стр. \pageref{c13}. \ Можно доказать, что здесь \
$\widetilde{\delta}_{\tau_2}^{<\alpha_1} =
\gamma_{\tau_1}^{<\alpha_1}$.
\\
\quad \\
} %
\end{lemma}

\noindent \textit{Доказательство.} \ Верхние индексы \ $<
\alpha_{1} $, $\vartriangleleft \alpha_{1} $ \ будут опускаться.
Так как функция\ $\alpha S_f$ \ монотонна на \ $]\tau_1,\tau_2[$,
\ то из теоремы~1 следует, что
\begin{equation} \label{e10.1}
    ]\gamma_{\tau_1}, \gamma_{\tau_2}[ \;  \cap \; SIN_n =
    \varnothing.
\end{equation}
Стоя на кардинале \ $\alpha^2 = \alpha_{\tau_2}^{\Downarrow}$ \
рассмотрим ниже \ $\alpha^{2}$ \ функцию \ $\alpha
S_{f}^{<\alpha^{2}}$. \ По условию $(ii)$ и лемме~\ref{8.7.} об
абсолютности  она совпадает с \ $\alpha S_{f}$ \ на \ $] \tau
_{1}, \tau _{2} [ $ \ и монотонна на этом интервале.
\\
Поэтому \ $\widetilde{\delta}_{\tau_2} \le \gamma_{\tau_1}$, \
иначе \ $\widetilde{\delta}_{\tau_2} \in \; ] \gamma_{\tau_1},
\gamma_{\tau_2}[$ \ вопреки $(ii)$, (\ref{e10.1}).
\\
Если теперь \ $a_{\tau _{2}}=1$, \ то по лемме~\ref{9.2.}
\begin{equation*}
Od\alpha S_{f} ( \tau _{1},\tau _{2} ) \leq  Od ( \alpha S_{\tau _{2}} )
\end{equation*}
несмотря на (\textit{i}) и поэтому \ $a_{\tau_2} = 0$. \ То же
самое получается, если
\[
    \forall \tau <\tau _{2} \quad \exists \tau ^{\prime}
    \in \left[ \tau ;\tau _{2}\right[ \quad a_{\tau ^{\prime}}=0,
\]
так как в этом случае благодаря условию $(i)$ можно рассмотреть
 \mbox{$\tau_1^2 \in \; ]\tau_1, \tau_2[$} \
такой, что для \ $S^2 = \alpha S_{\tau_2}$
\begin{equation} \label{e10.2}
    a_{\tau_1^2} = 0, \quad \alpha S_{\tau_1^2} \gtrdot S^2.
\end{equation}
Благодаря заключению 2) и лемме 3.2~\cite{Kiselev11} об
ограничении
 \textit{нулевая} матрица \ $S^2$ \ получает
некоторые допустимые носители \ $\alpha \in \; ]\gamma_{\tau_1^2},
\gamma_{\tau_1^2+1}[$ \ как результат ограничения \
$SIN_{n-1}$-кардиналом \ $\gamma_{\tau_1^2+1}$ \ следующего \
$\Sigma_{n-1}$-утверждения
\begin{equation} \label{e10.3}
    \quad
    \exists \alpha \Big( \gamma_{\tau_1^2} < \alpha \wedge
    \exists \delta, \alpha, \rho ~ \bigl(
    \delta \le \gamma_{\tau_1} \wedge
    SIN_n^{<\alpha^{\Downarrow}}(\gamma_{\tau_1}) \wedge
    \qquad \qquad \qquad \qquad
\end{equation}
\begin{equation*}
    \qquad\qquad\qquad\qquad\qquad
    \wedge \alpha \mathbf{K}_{n+1}^{\exists}(0, \delta, \gamma_{\tau_1^2}, \alpha, \rho, S^2)
    \big) \Big),
\end{equation*}
которое выполняется ниже \ $\gamma_{\tau_1^2 + 1}$, \ так как оно
выполняется для \ $\alpha=\alpha_{\tau_2}$ \ ниже \ $\alpha_1$.
\\
Следовательно, из (\ref{e10.2}), (\ref{e10.3}) следует, что
матрица \ $S^2$ \ запрещается для определения матричного значения
\ $\alpha S_{\tau_1^2}$ \ по определению~\ref{8.3.}~2).
\\
Но это может случиться, только если \ $S^2$ \ на \ $\alpha$ \
\textit{подавлена} для \ $\gamma_{\tau_1^2}$, \ что влечёт
 \ $SIN_n(\gamma_{\tau_1^2})$ \ вопреки
(\ref{e10.1}).
\\
Поэтому
\[
    \exists \tau < \tau_2 ~ \forall \tau^{\prime} \in [\tau,\tau_2[
    ~~ a_{\tau^{\prime}}=1,
\]
и лемма \ref{9.4.} влечёт \ $a_{\tau} \equiv 1$ \ на \ $]\tau_1,
\tau_2[\;$.
\\
\hspace*{\fill} $\dashv$
\\

\noindent Отсюда и из теоремы~1 сразу же вытекает

\begin{corollary}
\label{10.2.} \hfill {}
\\
Пусть \medskip

(i) \ $A_{2}^{\vartriangleleft \alpha_1}(\tau _{1},\tau _{2},\tau
_{3})$;
\newline

(ii) \ $\left] \gamma _{\tau _{1}}^{<\alpha _{1}}, \gamma _{\tau
_{2}}^{<\alpha _{1}}\right] \cap SIN_{n}^{<\alpha _{1}} \neq
\varnothing $.
\medskip

\noindent Тогда \medskip

1) \ $\gamma _{\tau _{2}}^{<\alpha _{1}}$ \ -- это наследник
кардинала \ $\gamma _{\tau _{1}}^{<\alpha _{1}}$ \ в \
$SIN_{n}^{<\alpha _{1}}$; \medskip

2) \ $a_{\tau }^{<\alpha _{1}}\equiv 1$ на \ $\left] \tau
_{1},\tau _{2} \right[ $, \ $a_{\tau _{2}}^{<\alpha _{1}}=0$ \quad
и \medskip

3) \ $\widetilde{\delta }_{\tau _{2}}^{<\alpha _{1}} \le \gamma
_{\tau _{1}}^{<\alpha _{1}}$.  \label{c14}
\endnote{
\ стр. \pageref{c14}. \ Снова в действительности здесь \
$\widetilde{\delta}_{\tau_2}^{<\alpha_1} =
\gamma_{\tau_1}^{<\alpha_1}$.
\\
\quad \\
} %
\\
\hspace*{\fill} $\dashv$
\end{corollary}

\begin{lemma}
\label{10.3.} \hfill {}
\\
Пусть \medskip

(i) \ $A_{2}^{\vartriangleleft \alpha_1}(\tau _{1},\tau _{2},\tau
_{3})$; \newline

(ii) \ $a_{\tau _{2}}^{<\alpha _{1}}=1$. \medskip

\noindent Тогда для кардинала предскачка \ $\alpha ^{2}=\alpha
_{\tau _{2}}^{<\alpha _{1}\Downarrow }$ \ существует ординал
\begin{equation*}
\tau _{\ast }=\min \{\tau \in \; ]\tau_1, \tau_2[ \; : \gamma
_{\tau }^{<\alpha _{1}}\in SIN_{n}^{<\alpha ^{2}}\}
\end{equation*}
такой, что: \medskip

1) \ $\gamma _{\tau _{\ast }}^{<\alpha _{1}}<\widetilde{\delta }_{\tau
_{2}}^{<\alpha _{1}}$, \ $\gamma _{\tau _{\ast }}^{<\alpha _{1}}\notin
SIN_{n}^{<\alpha _{1}}$; \medskip

2) \ $a_{\tau }^{<\alpha ^{2}} \equiv a_{\tau }^{<\alpha_1} \equiv
1$ \ on \ $\left] \tau _{1},\tau _{\ast }\right[ $;
\medskip

3) \ $\alpha S_{f}^{<\alpha ^{2}}$ \ монотонна на \ $[\tau
_{1},\tau _{\ast }[$ \quad и \medskip

4) \ $Od\alpha S_{f}^{<\alpha _{1}}(\tau _{1},\tau _{\ast })>Od(\alpha
S_{\tau _{2}}^{<\alpha _{1}})$. \medskip
\end{lemma}

\noindent \textit{Доказательство.} \ Верхние индексы \ $<
\alpha_{1} $, $\vartriangleleft \alpha_{1} $ \ будут опускаться.
Во-первых, нужно отметить, что \ $\gamma_{\tau_{1}} <
\widetilde{\delta}_{\tau_{2}} $, \ иначе $(ii)$ и лемма~\ref{9.2.}
нарушают условие $(i)$.

\noindent Затем, из леммы~\ref{10.1.} и условий (\textit{i}),
(\textit{ii}) следует, что имеется  следующий ординал ниже
кардинала предскачка \ $\alpha^2 = \alpha_{\tau_2}^{\Downarrow}$:

\vspace{0pt}
\begin{equation*}
    \tau_{\ast} = \min \bigl \{ \tau >\tau _{1}:  \gamma _{\tau }\in
    \bigl ( SIN_{n}^{<\alpha^2}  - SIN_{n} \bigr ) \bigr \} .
\end{equation*}
\vspace{0pt}

\noindent На интервале \ $[ \tau _{1},\tau_{\ast} [ $ \ функция
\quad $\alpha S_{f}$ \ монотонна и по теореме~1 (для \ $\alpha^2 $
\ вместо \ $\alpha_{1} $) кардинал \ $\gamma _{\tau_{\ast}}$ \
является наследником  \ $\gamma _{\tau _{1}}$ \ в классе \
$SIN_{n}^{<\alpha^2}$. \ Благодаря лемме~\ref{9.3.} функция \
$\alpha S_{f} $ \ стабилизируется на интервале \ $[ \tau
_{1},\tau_{\ast} [ $, \ поэтому для некоторых \ $\tau _{0}\in
\left] \tau _{1} ,\tau_{\ast} \right[$ \ и \quad $S^{0}$ \
выполняется утверждение \ $\forall \tau \geq \tau _{0} \ \alpha
S_{\tau }=S^{0}$ \ ниже \ $ \gamma _{\tau_{\ast}}$. \
Следовательно, ниже \ $\gamma_{\tau_{\ast}}$ \ справедливо более
слабое утверждение:
\begin{equation*}
    \forall \tau \Bigl( \tau_0 < \tau ~ \rightarrow \exists S  ( S = \alpha
    S_{\tau} \wedge S \underline{\gtrdot} S^0 )\Bigr).
\end{equation*}
Оно может быть сформулировано в \ $\Pi_n$-форме, как это было
сделано в доказательстве леммы~\ref{9.3.} посредством утверждения
 (\ref{e9.4}), где\ $\tau_{1,3}^n$,
\ $\alpha S_{\tau_2}^1$ \ следует заменить на \ $\tau_0$, \ $S^0$
\ соответственно.
\newline
Затем кардинал \ $\gamma _{\tau_{\ast}}$ \ продолжает это
утверждение до
 \ $ \alpha^2$ \ и поэтому по
условию $(i)$
\begin{equation*}
\rho =Od\alpha S_{f} ( \tau _{1},\tau_{\ast} ) >  Od ( \alpha S_{\tau _{2}}
) .
\end{equation*}
\vspace{0pt}

\noindent Далее, обсудим единичную характеристику. Лемма
3.2~\cite{Kiselev11} и условия (\textit{i}), (\textit{ii}) влекут,
что существуют некоторые допустимые носители \ $\alpha $-матриц
единичной характеристики, расположенные конфинально кардиналу \
$\gamma _{\tau _{1}}$, \ как это происходило несколько раз выше.
Поэтому такие носители должны быть также в интервале \ $]
\gamma_{\tau_{\ast}}, \alpha^2 [\;$, \ иначе кардинал \ $\gamma
_{\tau _{\ast}}$ \ был бы определим ниже \ $\alpha^2$ \ вместе с
ординалом  \ $\rho $, \ а тогда по лемме 4.6~\cite{Kiselev11}
получается противоречие:
\[
    \rho <Od (\alpha S_{\tau _{2}} ).
\]
Осталось применить лемму 3.2~\cite{Kiselev11}, \ref{9.4.} (где
 \ $\tau_{\ast}$, \ $ \alpha^2$ \ играют роль \ $\tau_2$, \
$\alpha_1$ \ соответственно) что завершает доказательство, так как
\ $\alpha S_{f}^{<\alpha^1}$, \ $a_{f}^{<\alpha^1}$ \ совпадают с
 \ $ \alpha S_{f}^{<\alpha^2}$, \ $a_{f}^{<\alpha^2}$
\ на \ $[\tau_1, \tau_\ast [ $ \ согласно лемме~\ref{8.7.} об
абсолютности.
\\
\hspace*{\fill} $\dashv$
\\
\quad \\
С помощью рассуждений, вполне аналогичных доказательствам
лемм~\ref{10.1.}-\ref {10.3.}, легко устанавливается

\begin{lemma}
\label{10.4.} \hfill {}
\\
Пусть \medskip

(i) \ $A_{2}^{\vartriangleleft \alpha_1}(\tau _{1},\tau _{2},\tau
_{3})$;
\medskip

(ii) \ $\forall \gamma <\gamma _{\tau _{2}}^{<\alpha _{1}}\exists
\tau ~(\gamma <\gamma _{\tau }^{<\alpha ^{2}}\wedge a_{\tau
}^{<\alpha ^{2}}=1)$ \ для \ $\alpha ^{2}=\alpha _{\tau
_{2}}^{<\alpha _{1}\Downarrow }$ ;
\medskip

(iii) \ $a_{\tau _{2}}^{<\alpha _{1}}=0$. \medskip

\noindent Тогда \medskip

1) \ $\widetilde{\delta }_{\tau _{2}}^{<\alpha _{1}} \le \gamma
_{\tau _{1}}^{<\alpha _{1}}$ \quad и \medskip

2) \ $\exists \tau \in \left] \tau _{1},\tau _{2}\right[ \quad
(a_{\tau }^{<\alpha _{1}}=1\wedge \alpha S_{\tau }^{<\alpha
_{1}}\gtrdot \alpha S_{\tau _{2}}^{<\alpha _{1}})$. \label{c15}
\endnote{
\ стр. \pageref{c15}. \ В действительности десь \
$a_{\tau}^{<\alpha_1} \equiv 1$ \ на \ $]\tau_1, \tau_2[$\ \ и
снова \ $\widetilde{\delta}_{\tau_2}^{<\alpha_1} =
\gamma_{\tau_1}^{<\alpha_1}$.
\\
\quad \\
} %
\end{lemma}

\noindent \textit{Доказательство.} \ Верхние индексы \ $<
\alpha_{1} $, $\vartriangleleft \alpha_{1}$ \ будут опускаться.
Как обычно, рассмотрим ситуацию ниже, стоя на \ $\alpha^{2} =
\alpha_{\tau_{2}}^{\Downarrow} $. \ Предположим, что \
$\widetilde{\delta}_{\tau_{2}} \in ]\gamma_{\tau_{1}},
\gamma_{\tau_{2}} [ $; \ здесь достаточно рассмотреть следующие
два случая:
\\
Случай 1. \quad $[\gamma_{\tau_1}, \gamma_{\tau_2} [ \; \cap \;
SIN_n^{<\alpha^2} \subseteq SIN_n$, \ тогда снова (буквально как
это было в доказательстве теоремы 1\;) функция \ $ \alpha S_{f} $
\ монотонна на интервале \ $[\tau_1, \tau_2[$, \ в то время как
интервал \ $]\gamma_{\tau_1}, \gamma_{\tau_2}[$ \ содержит
 \ $ SIN_{n} $-кардинал  \ $\widetilde{\delta}_{\tau_{2}}
\ $, \ в противоречии с теоремой 1.
\\
Случай 2. \quad $[\gamma_{\tau_1}, \gamma_{\tau_2} [ \; \cap \;
SIN_n^{<\alpha^2} \nsubseteq SIN_n$. \ В этом случае нужно снова
применить технику ограничения-и-продолжения, буквально как это
было сделано в доказательстве части~1a. леммы~\ref{9.4.}. Сначала
повторим аргумент из доказательства леммы 10.3 касательно функции
\ $\alpha S_{f} $, \ определённой на интервале \ $]\tau_{1},
\tau_\ast [\; $, \ где
\begin{equation*}
\gamma_{\tau_\ast} = \min \left ( SIN_{n}^{ < \alpha^{2}} -
SIN_{n} \right ) \mbox{\it \quad и \quad}  a_{\tau} \equiv 1
\mbox{\it \ \ on \ \ } ]\tau_{1},\tau_\ast [
\end{equation*}
благодаря условию $(ii)$. По теореме 1. (для \ $\tau_\ast$,
$\alpha^2$ как $\tau_2$, $\alpha_1$) и условию $(i)$,\ кардинал \
$\gamma_{\tau_\ast}$ \ является наследником  кардинала \
$\gamma_{\tau_1}$ \ в \ $SIN_n^{<\alpha^2}$. Следовательно,
утверждение
\begin{equation*}
\forall \tau > \tau_1 \quad a_{\tau} = 1
\end{equation*}
выполняется ниже \ $\gamma_{\tau_\ast}$; \ в этом нетрудно
убедиться с помощью условия $(ii)$ и лемм~3.2~\cite{Kiselev11},
\ref{9.4.}. Это условие можно сформулировать в \ $\Pi_n$-форме
буквально как это было сделано в доказательстве части~1a.
леммы~\ref{9.4.} с \ $\tau^3$ \ как \ $\tau_1$:

\vspace{0pt}
\begin{equation*}
\forall \gamma \Bigl ( \gamma_{\tau_1} < \gamma  \; \wedge \; SIN_{n-1}
(\gamma) \longrightarrow \qquad \qquad \qquad  \qquad \qquad \qquad \qquad
\qquad \qquad
\end{equation*}
\begin{equation*}
\qquad \qquad \longrightarrow \exists \delta, \alpha,  \rho, S ~
\bigl (  SIN_n^{<\alpha^{\Downarrow}}(\gamma_{\tau_1}) \wedge
\alpha \mathbf{K}_{n+1}^{\exists}(1, \delta, \gamma, \alpha,\rho,
S)  \bigr ) \Bigr ).
\end{equation*}
\vspace{0pt}

\noindent После этого \ $SIN_n^{<\alpha^2} $-кардинал \
$\gamma_{\tau_\ast} $ \ продолжает это предложение до \
$\alpha^{2} $ \ и поэтому \ $a_{\tau_{2}}=1 $ \ вопреки условию
 $(iii)$.
\\
Итак, \ $\widetilde{\delta}_{\tau_{2}} \leq \gamma_{\tau_{1}}$; \
для завершения доказательства следует применить лемму~9.2.
Предположим, что
\begin{equation*}
\exists \tau \in \; ]\tau_{1},\tau_{2} [ ~~ \forall \tau^{\prime}
\in [\tau ,\tau_{2} [ \ \ \ a_{\tau^{\prime}} = 0,
\end{equation*}
тогда по этой лемме монотонность функции  \ $\alpha S_{f} $ \ на
 \ $]\tau_{1},\tau_{2}[$ \ влечёт
\[
    Od \alpha S_{f} (\tau_{1},\tau_{2}) \leq Od (\alpha
    S_{\tau_2}),
\]
нарушая условие $(i)$. Это противоречие вместе с $(i)$
устанавливает 2) и завершает доказательство этой леммы.
\\
\hspace*{\fill} $\dashv$

Следующая лемма будет использована в конце доказательства основной
теоремы, снова  полагаясь на формулу \ $A^{0 \vartriangleleft
\alpha_1}(\tau)$ \ (напомним определение \ref{8.1.}~3.2\;):
\vspace{-6pt}
\begin{multline*}
    \exists \gamma < \alpha_1 \Big( \gamma =
    \gamma_{\tau}^{<\alpha_1} \wedge \neg \exists a, \delta,
    \alpha, \rho < \alpha_1 \exists S \vartriangleleft \rho ~ \big(
    \mathbf{K}_n^{\forall <\alpha_1}(\gamma,\alpha_{\chi}^{\Downarrow})
    \wedge
\\
    \wedge \alpha \mathbf{K}_{n+1}^{\exists \vartriangleleft
    \alpha_1}(a, \delta, \gamma, \alpha, \rho, S) \big) \Big),
\end{multline*}
означающую, что не существует \ $\alpha$-матрицы на каком-то
носителе \ $\alpha
> \gamma_{\tau}^{<\alpha_1}$, \ допустимой для \
$\gamma_{\tau}^{<\alpha_1}$ \ ниже \ $\alpha_1$.
\\
Соответственно, через \ $A_1^{0 \vartriangleleft \alpha_1}(\tau_1,
\tau_2, \alpha S_f^{<\alpha_1})$ \ обозначается формула (напомним
определение  \ref{8.1.}~1.1 для \ $X_1 = \alpha S_f^{<\alpha_1}$):
\[
    A^{0 \vartriangleleft \alpha_1}(\tau_1) \wedge
    A_1^{\vartriangleleft \alpha_1}(\tau_1, \tau_2, \alpha
    S_f^{<\alpha_1});
\]
также будет использоваться формула $A_2^{0 \vartriangleleft
\alpha_1}(\tau_1, \tau_2^{\prime}, \tau_3, \alpha
S_f^{\vartriangleleft \alpha_1})$ \ (напомним определение
\ref{8.2.}~3.3\;):
\[
    A^{0 \vartriangleleft \alpha_1}(\tau_1) \wedge
    A_2^{\vartriangleleft \alpha_1}(\tau_1, \tau_2^{\prime}, \tau_3, \alpha
    S_f^{<\alpha_1}).
\]

\begin{lemma}
\label{10.5.} \hfill {}
\\
Пусть \medskip

(i)\quad $A_{1}^{0 \vartriangleleft \alpha_1}(\tau _{1},\tau _{2},
\alpha S_f^{<\alpha_1})$;
\medskip

(ii)\quad $\tau_2 \le \tau_3$ \ и \ $S^3$ \ это матрица
характеристики \ $a^3$ \ на носителе
\[
    \alpha_3 \in \; ] \gamma_{\tau_3}^{<\alpha_1}, \alpha_1 [
\]
сохраняющая \ $SIN_n^{<\alpha_1}$-кардиналы \ $\le
\gamma_{\tau_2}^{<\alpha_1}$ \ ниже \ $\alpha_1$ \ и обладающая
производящим собственным диссеминатором \ $\check{\delta}^{S^3}$;
\medskip

(iii)\quad $\check{\delta}^{S^3} \leq \gamma _{\tau _{1}}^{<\alpha
_{1}}$. \medskip

\noindent Тогда \quad $a^3=0$.
\\
Аналогично для всякого диссеминатора \ $\widetilde{\delta}$ \
матрицы \ $S^3$ \ на \ $\alpha^3$ \ с любой базой \ $\rho \ge
\rho^{S^3}$.
\end{lemma}

\noindent \textit{Доказательство.} \ Как обычно, мы будем
исследовать ситуацию ниже, стоя на кардинале предскачка \ $a^3 =
\alpha_3^{\Downarrow}$ \ и рассматривая диссеминатор \ $\check{
\delta}^{S^3} $ \ с базой данных \ $\rho^3 = \rho^{S^3} =
\widehat{\rho_1}$, $\rho_1 = Od(S^3)$; \ верхние индексы \
$<\alpha _{1}$, $\vartriangleleft \alpha_1$ \ будут опускаться для
некоторого удобства.
\\
Допустим, что эта лемма неверна и \ $a^3 = 1$, \ тогда
$\check{\delta}^{S^3}$ \ допустимый и неподавленый диссеминатор
матрицы \ $S^3$ \ на \ $\alpha_3$ \ для каждого \ $\gamma_{\tau}
\in \; ]\tau_1, \tau_3[$ \ и по лемме~\ref{9.2.} \ $a_{\tau }
\equiv 1$ \ на \ $ \left] \tau _{1}, \tau _{2}\right[ $ \ и
\begin{equation*}
Od \alpha S_{f}(\tau_{1},\tau_{2})\leq Od(\alpha  S_{\tau_{3}}).
\end{equation*}
Из леммы~\ref{9.5.} следует, что благодаря условию $(iii)$ функция
\ $\alpha S_{f}$ \ немонотонна на \ $[ \tau _{1},\tau _{2} [ $ \
(напомним случай~2. в доказательстве пункта 3. леммы~\ref{9.5.}\;)
и поэтому существует \ $\tau _{2}^{\prime }$, \ для которого
выполняется
\begin{equation}  \label{e10.4}
A_2^0 ( \tau _{1},\tau _{2}^{\prime }, \tau_{2} ).
\end{equation}
Теперь следует повторить рассуждение из  случая~1b. доказательства
леммы~\ref{9.4.}. Ниже \ $\alpha^3$ \ справедливо следующее \
$\Sigma _{n+1}$-утверждение по лемме~\ref{8.7.} об абсолютности
 (напомним (\ref{e9.8})\;):
\[
    \exists \gamma^0 \exists \tau _{1}^{\prime },  \tau _{2}^{\prime
    \prime },\tau _{3}^{\prime } < \gamma^0 ~  \Bigl( SIN_{n} ( \gamma^0 ) \wedge
    \gamma_{\tau_1^{\prime}}^{<\gamma^0} <
    \gamma_{\tau_2^{\prime\prime}}^{<\gamma^0} <
    \gamma_{\tau_3^{\prime}}^{<\gamma^0} < \gamma^0 \wedge
    \qquad \qquad
\]
\begin{equation}  \label{e10.5}
    \qquad
    \wedge A_2^{0 \vartriangleleft \gamma^0}
    ( \tau _{1}^{\prime },\tau _{2}^{\prime \prime },
    \tau_{3}^{\prime },\alpha S_{f}^{<\gamma^0 } )
    \wedge \forall \tau^{\prime\prime\prime} \in \;
    ]\tau_1^{\prime}, \tau_2^{\prime\prime}] ~
    a_{\tau^{\prime\prime\prime}}^{<\gamma^0}  = 1 \wedge
\end{equation}
\[
    \qquad \qquad \qquad \qquad \qquad \qquad \qquad \qquad \wedge
    \alpha S_{\tau _{2}^{\prime \prime }}^{<\gamma^0 } =  \alpha
    S_{\tau_{2}^{\prime }} \Bigr).
\]
\vspace{0pt}

\noindent Оно содержит индивидные константы \ $< \rho^3 $ \ и \
$\alpha S_{\tau _{2}^{\prime }} \vartriangleleft \rho^3 $, \
поэтому диссеминатор \ $\check{\delta}^{S^3}$ \ ограничивает это
предложение и оно выполняется уже ниже \ $\check{ \delta}^{S^3}$.

Теперь перейдём к ситуации ниже кардинала предскачка \ $\alpha^2 =
\alpha_{\tau_2^\prime}^\Downarrow$\ единичной матрицы \ $\alpha
S_{\tau _{2}^{\prime }} $\ на носителе \
$\alpha_{\tau_2^\prime}$.\ По условию $(i)$ \ $\gamma_{\tau_1} \in
SIN_n$, \ поэтому по лемме ~\ref{8.5.}~1) получается \
$\gamma_{\tau_1} \in SIN_n^{<\alpha^3}$. \ Так как \
$\check{\delta}^{S^3} \le \gamma_{\tau_1}$ \ и \
$\check{\delta}^{S^3} \in SIN_n^{<\alpha^3}$, \ то лемма
3.8~\cite{Kiselev11} (для \ $\alpha^3$, \ $\gamma_{\tau_1}$ \ как
\ $\alpha_1$, \ $\alpha_2$) \ влечёт \ $\check{\delta}^{S^3} =
\gamma_{\tau_1}$ \ или \ $\check{\delta}^{S^3} \in
SIN_n^{<\gamma_{\tau_1}}$; \ затем из той же леммы (для \
$\gamma_{\tau_1}$ \ как \ $\alpha_2$) вытекает \
$\check{\delta}^{S^3} \in SIN_n$.
\\
Отсюда и из леммы~\ref{8.5.}~1) следует \ $\check{\delta}^{S^3}
\in SIN_n^{<\alpha^2}$; \ значит, в утверждении (\ref{e10.5})
можно заменить \ $\gamma^0 $ \ на \ $\alpha^{2}$ \ по
лемме~\ref{8.7.} об абсолютности, а тогда по лемме \ref{8.5.}~5) (
для \ $\alpha_{\tau_2^\prime}$ \ как \ $\alpha$) \ мы получаем
противоречие: \ $a_{\tau _{2}^{\prime }}=0$ \ вопреки тому, что \
$a_{\tau } \equiv 1$ \ на всём интервале \ $ \left] \tau _{1},
\tau _{2}\right[ $ .\
\\
\hspace*{\fill} $\dashv$
\\

Теперь специальная теория матричных функций разработана достаточно
для того, чтобы приступить к доказательству основной теоремы.

\newpage

\chapter{Приложения специальной теории}

\setcounter{section}{10}

\section{Доказательство основной теоремы}

\setcounter{equation}{0}

\hspace*{1em} Противоречие, доказывающее основную теорему,
достигается следующим диагональным рассуждением:
\\
С одной стороны, по лемме 8.9 функция \ $\alpha S_f^{<\alpha_1}$ \
определена на непустом множестве
\begin{equation*}
T^{\alpha_1} = \{\tau:\alpha\delta^{\ast}<\gamma_{\tau}\ <
\alpha_1 \}
\end{equation*}
для каждого достаточно большого кардинала \ $\alpha_1 \in SIN_n$.
\\
Её монотонность на этом множестве исключена по теореме 1.
\\
Но с другой стороны, эта монотонность устанавливается следующей
теоремой для каждого \ $SIN_n$-кардинала \
$\alpha_1>\alpha\delta^{\ast}$ \ всякой достаточно большой
конфинальности. Напомним, что ограничивающие кардиналы \
$\alpha_1$ \ всегда предполагаются эквиинформативными с \
$\chi^{\ast}$, \ то есть что всегда выполняется \
$A_6^e(\chi^{\ast}, \alpha_1)$ \ (см. определение \ref{8.1.}~5.1
для \ $\chi=\chi^{\ast}$, $\alpha^0=\alpha_1$).
\\

\noindent \textbf{Теорема2.}\quad \newline \emph{\hspace*{1em}
Пусть функция \ $\alpha S_f^{<\alpha_1}$ \ определена на непустом
множестве
\begin{equation*}
T^{\alpha_1} =\{\tau:\gamma_{\tau_1}^{<\alpha_1} <
\gamma_{\tau}^{<\alpha_1} < \alpha_1 \}
\end{equation*}
таком, что  \ $\alpha_1 < k $ \ и: \medskip }

\emph{(i) \ $\tau_1=\min\{\tau:\forall\tau^{\prime}
(\gamma_{\tau}^{< \alpha_1} <
\gamma_{\tau^{\prime}}^{<\alpha_1}\longrightarrow\tau^{\prime}\in
dom(\alpha S_f^{<\alpha_1}))\};$ \newline }

\emph{(ii) \ $\sup SIN_n^{<\alpha_1}=\alpha_1;$ \newline }

\emph{(iii) \ $cf(\alpha_1)\geq\chi^{\ast +}.$ \newline
}

\emph{\noindent Тогда эта функция \ $\alpha S_f^{<\alpha_1}$ \
монотонна на этом множестве:
} %

\[
    \forall \tau_1, \tau_2 \in T^{\alpha_1} \bigl( \tau_1 <
    \tau_2 \rightarrow \alpha S_{\tau_1}^{<\alpha_1}
    \underline{\lessdot} \alpha S_{\tau_2}^{<\alpha_1} \bigr).
\]

\quad \\
\noindent \textit{Доказательство.} \ План этого доказательства
состоит в следующем: \\ Рассуждение будет проводиться индукцией по
кардиналу \ $\alpha_1$.
\\
Предположим, что эта теорема неверна и кардинал \ $\alpha_1^\ast$
\ -- \textit{минимальный}, нарушающий эту теорему, то есть функция
\ $\alpha S^{<\alpha_1^\ast}_f$ \ \emph{немонотонна} на \
множестве
\begin{equation*}
T^{\alpha_1^\ast} = \bigl \{ \tau:
\gamma_{\tau_1^\ast}^{<\alpha_1^\ast} <
\gamma_{\tau}^{<\alpha_1^\ast} < \alpha_1^{\ast} \bigr \}
\end{equation*}
с указанными свойствами $(i)$--$(iii)$ для некоторого \
$\tau_1^{\ast}$; \ таким образом, принимается \textit{первая
индуктивная гипотеза}:
\\
для каждого \ $\alpha_1 < \alpha_1^{\ast}$ \ функция \ $\alpha
S_f^{<\alpha_1}$ \ \emph{монотонна} на множестве \ $T^{\alpha_1}$
\ со свойствами $(i)$--$(iii)$.
\\
Из теоремы 1 сразу же следует, что этот \ $\alpha_1^{\ast}$ \ --
это просто \textit{минимальный} кардинал \ $\alpha_1$, \ для
которого такое множество \ $T^{\alpha_1}$ \ существует, потому что
для каждого \ $\alpha_1 < \alpha_1^{\ast}$ \ это множество не
существует, так как функция \ $\alpha S_f^{<\alpha_1}$ \ на таком
\ $T^{\alpha_1}$ \ \textit{немонотонна} по теореме~1 и в то же
время \textit{монотонна} по первой индуктивной гипотезе о
минимальности \ $\alpha_1^{\ast}$.

Всё рассуждение будет проводиться ниже \ $\alpha_1^\ast$ \ (и все
переменные будут ограничиваться этим \ $\alpha_1^{\ast}$), \ или
ниже ограничивающих кардиналов \ $\alpha_1 \le \alpha_1^{\ast}$, \
поэтому верхние индексы \ $< \alpha_1^\ast$, \ $\vartriangleleft
\alpha_1^\ast$ \ будут опускаться для некоторого удобства вплоть
до конца доказательства теоремы 2.

Сначала отметим, что в условиях теоремы 2 выполняется
\begin{equation*}
\gamma_{\tau_1}^{<\alpha_1}\in SIN_n^{<\alpha_1};
\end{equation*}
чтобы убедиться в этом достаточно ещё раз повторить аргумент, уже
применённый выше несколько раз (сначала в доказательствах лемм
7.7, 8.10\;). Благодаря этому нетрудно видеть, что для
\textit{каждого} достаточно большого ординала \ $ \tau_3^\ast \in
T^{\alpha_1^\ast}$ \ интервал \ $[\gamma_{\tau_1^\ast},\
\gamma_{\tau_3^\ast}[$ \ может рассматриваться \textit{как блок},
то есть существуют некоторые ординалы \ $\tau_1^{\ast\prime},\
\tau_2^\ast,\ \eta^{\ast 3}$ \ которые выполняют утверждение
(напомним определение \ref{8.1.}~1.6 для \ $X_1 = \alpha S_f$,
$X_2 = a_f$):

\begin{equation*}
A_4^b (\tau_1^\ast,\ \tau_1^{\ast\prime},\ \tau_2^\ast,\
\tau_3^\ast,\ \eta^{\ast 3}, \alpha S_f, a_f).
\end{equation*}

\noindent Здесь (согласно этому определению~\ref{8.1.}\;) \
$\tau^{\ast\prime}_1$ \ это индекс матрицы \ $\alpha
S_{\tau_1^{\ast \prime}}$ \  \textit{единичной характеристики} \
$a_{\tau^{\ast\prime}_1}=1$ \ на её носителе \ $
\alpha_{\tau^{\ast\prime}_1}$ \ и \ $\eta^{\ast 3}$ \ --- \ тип
этого интервала.
\\
Далее, благодаря условию \ ($iii$) \ этой теоремы 2 мы можем
использовать индекс \ $ \tau_3^\ast \in T^{\alpha_1^\ast}$ \
такой, что интервал \ $[\gamma_{\tau_1^\ast}$, \
$\gamma_{\tau_3^\ast}[$ \  имеет именно  тип

\begin{equation*}
\eta^{\ast 3} > Od(\alpha S_{\tau_1^{\ast\prime}}), \quad
\eta^{\ast 3} < \chi^{\ast +}.
\end{equation*}

\noindent Теперь формула \ $\mathbf{K}^0$ \ начинает действовать и
\textit{замыкает диагональное рассуждение}:
\\
Рассмотрим матрицу\ $\alpha S_{\tau_3^\ast}$ \ на носителе \
$\alpha_{\tau_3^\ast}$ \ вместе с её диссеминатором \
$\widetilde{\delta} ^{\ast 3}=\widetilde{\delta}_{\tau_3^\ast}$ \
и базой данных \ $\rho^{\ast 3}=\rho_{\tau_3^\ast}$; \ мы увидим,
что по лемме~\ref{10.5.} она имеет \textit{нулевую} характеристику
на этом носителе.
\\
Стоя на кардинале предскачка \ $\alpha^{\ast
3}=\alpha_{\tau_3^\ast}^{\Downarrow}$ \ следует рассмотреть
ситуацию ниже этого \ $\alpha^{\ast 3}$:
\\
Почти очевидно, что по лемме \ref{8.8.} диссеминатор \
$\widetilde{\delta}^{\ast 3}$ \ попадает в некоторый
\emph{максимальный} блок \ $[\gamma_{\tau_1^\ast},\ \gamma^{\ast
3}[$ \ ниже \ $\alpha^{\ast 3}$ \ типа \ $\eta^{\ast 3 \prime}<
\chi^{\ast +}$, \ где \ $\gamma^{\ast 3}$ \ это некоторый \ $
\gamma_{\tau^{\ast\prime}_3}^{<\alpha^{\ast 3}}$. \ Нетрудно
видеть, что
\begin{equation*}
\gamma_{\tau_3^\ast} \leq \gamma^{\ast 3}\ \ \wedge\ \  \eta^{\ast
3} \leq \eta^{\ast 3 \prime},
\end{equation*}
поэтому выполняется
\\
\begin{equation*}
A^{M b \vartriangleleft \alpha^{\ast 3} } _4(\tau_1^\ast,
\tau_1^{\ast \prime}, \tau_2^\ast, \tau_3^{\ast\prime}, \eta^{\ast
3 \prime}, \alpha S_f^{<\alpha^{\ast 3}}, a_f^{<\alpha^{\ast 3}}).
\end{equation*}
\newline
Вс эти факты вместе составляют посылку леммы 8.5~6):
\\
\[
    a_{\tau_3^\ast}=0\wedge\gamma_{\tau_1^\ast}^{<\alpha^{\ast 3}}
    \leq \widetilde{\delta}^{\ast 3}<
    \gamma_{\tau^{\ast\prime}_3}^{<\alpha^{\ast3}} \wedge
    \qquad\qquad\qquad\qquad
\]
\[
    \qquad\qquad\qquad\qquad
    \wedge A_4^{M b \vartriangleleft \alpha^{\ast 3}} (\tau_1^\ast,
    \tau_1^{\ast \prime}, \tau_2^\ast, \tau_3^{\ast\prime}, \eta^{\ast
    3 \prime}, \alpha S_f^{<\alpha^{\ast 3}}, a_f^{<\alpha^{\ast 3}}).
\]

\noindent Следовательно, эта лемма влечёт

\begin{equation*}
\eta^{\ast 3 \prime}<\rho^{\ast 3}\vee\rho^{\ast 3} =\chi^{\ast +}
\end{equation*}

\noindent и во всяком случае \newline
\begin{equation*}
Od(\alpha S_{\tau_1^{\ast\prime}} ) < \eta^{\ast 3} \leq
\eta^{\ast 3 \prime} < \rho^{\ast 3}.
\end{equation*}

\noindent Но мы увидим скоро, что это невозможно, так как по
лемме~\ref{9.5.} о срезании лестницы сверху и по лемме \ref{11.3.}
(установленной ниже) выполняется:
\begin{equation*}
\rho^{\ast 3} \leq Od(\alpha S_{\tau_1^{\ast\prime}} ).
\end{equation*}
\vspace{0pt}

\noindent Это противоречие завершит доказательство теоремы~2.
\quad \\

  Чтобы осуществить изложенный план нужна некоторая
допонительная информация.
\\
Рассуждение, описанное выше, полагается на следущие несложные
вспомогательные леммы 11.1, 11.3, которые представляют собой его
``несущую конструкцию'' и описывают некоторые существенные
свойства поведения нулевых матриц; они не были представлены ранее
по причине их довольно специального характера. С этой целью
следует вспомнить формулу (см. определение \ref{8.1.}~1.1 для \
$X_1 = \alpha S_f^{<\alpha_1}$)
\quad \\
\quad \\
\ $A_1^{\vartriangleleft \alpha_1} (\tau_1, \tau_2, \alpha
    S_f^{<\alpha_1})$:

\[
    \tau_1+1 < \tau_2 \wedge \tau_1 = \min \bigl\{ \tau:
    \; ]\tau, \tau_2[ \; \subseteq dom(\alpha S_f^{<\alpha_1} )\bigr\}
\]
\[
    \qquad \qquad \qquad \qquad \qquad \wedge
    \gamma_{\tau_1}^{<\alpha_1} \in SIN_n^{<\alpha_1} \wedge
    \gamma_{\tau_2}^{<\alpha_1} \in SIN_n^{<\alpha_1};
\]
\vspace{0pt}Напомним также, что мы часто опускаем обозначения
функций \ $\alpha S_f^{<\alpha_1}$, $a_f^{<\alpha_1}$ \ в
написаниях формул ниже \ $\alpha_1$; \ напомним также, что тип
интервала \ $[ \gamma_{\tau_1}^{<\alpha_1},
\gamma_{\tau_2}^{<\alpha_1}[$ \ ниже \ $\alpha_1$ \ -- это тип
множества (см. определение \ref{8.1.}~1.3\;):
\[
    \bigl \{ \gamma: \gamma_{\tau_1}^{<\alpha_1} < \gamma <
    \gamma_{\tau_2}^{<\alpha_1} \wedge SIN_n^{<\alpha_1}(\gamma)
    \bigr \}.
\]

\noindent Предварительно следует остановиться на следующих
дополнительных аргументах, удобных для сокращения последующих
рассуждений; с этой целью нужно ввести такие понятия:

 Интервал \ $[\tau_1, \tau_2[$ \ и соответствующий интервал \
$[\gamma_{\tau_1}^{<\alpha_1}, \gamma_{\tau_2}^{<\alpha_1}[$ \
будут называться интервалами матричной допустимости, или просто
\textit{интервалами допустимости}, ниже \ $\alpha_1$, \ если для
каждого \ $\tau^{\prime} \in \; ]\tau_1,\tau_2[$ \ существует
некоторая \ $\alpha$-матрица \ $S$ \ на некотором носителе \ $>
\gamma_{\tau^{\prime}}^{<\alpha_1}$, \ \textit{допустимом} для \
$\gamma_{\tau^{\prime}}^{<\alpha_1}$ \ ниже \ $\alpha_1$:
\[
    \forall \tau^{\prime} \in \; ]\tau_1, \tau_2[ ~ \exists
    a^{\prime}, \delta^{\prime}, \alpha^{\prime}, \rho^{\prime},
    S^{\prime} ~ \alpha \mathbf{K}^{<\alpha_1} (a^{\prime},
    \delta^{\prime}, \gamma_{\tau^{\prime}}^{<\alpha_1},
    \alpha^{\prime}, \rho^{\prime}, S^{\prime}),
\]
и \ \mbox{$\gamma_{\tau_1}^{<\alpha_1} \in SIN_n$},
\mbox{$\gamma_{\tau_2}^{<\alpha_1} \in SIN_n$} \ и \ $\tau_1$ \ --
это минимальный ординал с этими свойствами.
\\
Затем, нужно обратиться к следующим свойствам произвольной нулевой
матрицы \ $S$ на её носителе \ $\alpha$,\ допустимыми для \
$\gamma_{\tau}^{<\alpha_1}$ \ вместе со своими \textit{минимальным
допустимым} диссеминатором \ $\widetilde{\delta}$ \ с базой \
$\rho$ \ ниже \ $\alpha_1$ \ для \ $\alpha_1 \le \alpha_1^{\ast}$
:\

(1a.) если \ $\gamma_{\tau_1}^{<\alpha_1} <
\gamma_{\tau_2}^{<\alpha_1} \le  \gamma_{\tau}^{<\alpha_1}$ \ и \
$\widetilde{\delta}$ \ попадает в интервал допустимости \
$[\gamma_{\tau_1}^{<\alpha_1}, \gamma_{\tau_2}^{<\alpha_1}[\;$, \
то есть \ $ \gamma_{\tau_1}^{<\alpha_1} \le \widetilde{\delta} <
\gamma_{\tau_2}^{<\alpha_1}$, \ то \ $\gamma_{\tau_1}^{<\alpha_1}
= \widetilde{\delta}$;

(1b.) если существует некоторая нулевая матрица \ $S^1$ \ на
некотором \textit{другом} носителе \ $\alpha^1 \neq \alpha$, \
тоже допустимая для того же \ $\gamma_{\tau}^{<\alpha_1}$ \ вместе
с своими минимальным допустимым диссеминатором \
$\widetilde{\delta}^1$ \ с базой \ $\rho^1$, \ то \ $S$ \ на \
$\alpha$ \ \textit{неподавлена} для \ $\gamma_{\tau}^{<\alpha_1}$
\ вместе со своими \ $\widetilde{\delta}$, $\rho$ \ ниже \
$\alpha_1$.

\noindent Проверка этих свойств будет проводиться индукцией по
тройкам \ $(\alpha_1, \alpha, \tau)$, \ упорядоченным канонически
как обычно (с \ $\alpha_1$ \ как первым компонентом, \ $\alpha$ \
-- вторым и с \ $\tau$ \ -- третьим).
\\
Предположим, тройка $(\alpha_1^0, \alpha^0, \tau^0)$ \ --
минимальная нарушающая (1a.) или (1b.); таким образом принимается
\textit{вторая индуктивная гипотеза}:
\\
для всякой меньшей тройки $(\alpha_1, \alpha, \tau)$ выполняются
(1a.) и (1b.).
\\
Мы увидим, что это вызывает противоречия; предстоящее
доказательство этого будет проводиться ниже \ $\alpha_1^0$, \
поэтому верхние индексы \ $<\alpha_1^0$, $\vartriangleleft
\alpha_1^0$ \ как обычно будут опускаться (вплоть до специального
замечания, если контекст не укажет явно на другой случай).

1. \ Начнём со свойства (1a.); \ допустим, оно нарушается, то есть
существует нулевая матрица \ $S^0$ \ на её носителе \ $\alpha^0
> \gamma_{\tau^0}$ \ с её минимальным диссеминатором \
$\widetilde{\delta^0}$ \ с базой \ $\rho^0$, \ все допустимые для
\ $\gamma_{\tau^0}$, \ и \ $\widetilde{\delta^0}$ \ попадает в
\textit{интервал допустимости} \ $[\gamma_{\tau_1^0},
\gamma_{\tau_2^0}[$, \ но
\begin{equation} \label{e11.1}
    \gamma_{\tau_1^0} < \widetilde{\delta}^0 < \gamma_{\tau_2^0}
    \le \gamma_{\tau^0}\;,
    \mbox{\it \ \ то есть \ } \widetilde{\delta}^0 = \gamma_{\tau_3^0}\;,
    ~ \tau_1^0 < \tau_3^0 < \tau_2^0.
\end{equation}
Отсюда и из леммы~3.8 \cite{Kiselev11} немедленно следует, что
\[
    \widetilde{\delta}^0 \in SIN_n,
\]
так как \ $\widetilde{\delta}^0 < \gamma_{\tau_2^0}$, \
$\widetilde{\delta}^0 \in SIN_n^{<\alpha^{0 \Downarrow}}$, \
$\gamma_{\tau_2^0} \in SIN_n$.
\\
По определению интервала допустимости существует некоторая матрица
\ $\alpha S_{\tau_3^0}$ \ на её носителе \ $\alpha_{\tau_3^0}$, \
допустимые для \ $\gamma_{\tau_3^0}$ \ вместе с её минимальным
диссеминатором \ $\widetilde{\delta}_{\tau_3^0}$ \ с базой \
$\rho_{\tau_3^0}$ \ (всё это ниже \ $\alpha_1^0$).
\\
Из второй индуктивной гипозы следует, что можно рассматривать \
$\widetilde{\delta}_{\tau_3^0} = \gamma_{\tau_1^0}$; \ поэтому из
леммы 3.2~\cite{Kiselev11} следует, что для каждого \
\mbox{$\gamma_{\tau} \in \; ]\gamma_{\tau_1^0},
\gamma_{\tau_3^0}[$} \ матрица \ $\alpha S_{\tau_3^0}$ \ обладает
многими носителями \ $\alpha \in \; ]\gamma_{\tau},
\gamma_{\tau+1}[$, \ допустимыми для \ $\gamma_{\tau}$, \ которые
неподавлены для этого \ $\gamma_{\tau}$ \ благодаря той же
индуктивной гипотезе, поэтому выполняется утверждение \
$A_1^0(\tau_1^0, \tau_3^0, \alpha S_f^{<\alpha_1^0})$:
\[
    A^0(\tau_1^0) \wedge
    A_1(\tau_1^0, \tau_3^0, \alpha S_f^{<\alpha_1^0})
\]
ниже \ $\alpha_1^0$. \ Те же аргументы действуют ниже кардинала
предскачка \ $\alpha^{0 \Downarrow}$, \ поэтому ниже \ $\alpha^{0
\Downarrow}$ \ также выполняется:
\[
    A_1^{0 < \alpha^{0 \Downarrow}}(\tau_1^0, \tau_3^0, \alpha S_f^{<\alpha^{0 \Downarrow}}).
\]
Для производящего диссеминатора \ $\check{\delta}^0$ \ матрицы \
$S^0$ \ на \ $\alpha^0$ \ с той же базой \ $\rho^0$ это влечёт:
\begin{equation} \label{e11.2}
    \check{\delta}^0 \le \gamma_{\tau_1^0},
\end{equation}
так как в противном случае \ $\check{\delta}^0$ \ попадает строго
в интервал допустимости \ $]\gamma_{\tau_1^0},
\gamma_{\tau_2^0}]$:
\begin{equation} \label{e11.3}
    \gamma_{\tau_1^0} < \check{\delta}^0 \le
    \gamma_{\tau_3^0},
\end{equation}
а тогда \ $\check{\delta^0}$ \ продолжает до \ $\alpha^{0
\Downarrow}$ \  \ $\Pi_{n+1}$-предложение о допустимости некоторых
матриц для каждого   \ $\gamma_{\tau}^{<\alpha^{0 \Downarrow}} >
\gamma_{\tau_1^0}$  \ например так, как это было сделано в
доказательстве леммы 7.5 1), где в формуле (7.2) нужно теперь
связать \ $\rho^4$, \ $S^4$ \ кванторами существования. Такие
матрицы становятся даже неподавленными для всех таких \
$\gamma_{\tau}^{<\alpha^{0 \Downarrow}}$ \ по второй индуктивной
гипотезе (всё это ниже $\alpha^{0 \Downarrow}$); \ следовательно,
возникает множество \ $T^{\alpha^{0 \Downarrow}}$ \ со свойствами
$(i)$--$(iii)$, указанными в теореме~2, в противоречии с первой
индуктивной гипотезой и теоремой~1, то есть с минимальностью \
$\alpha_1^{\ast}$.

\noindent Начиная с этого места рассуждение переходит к ситуации
ниже \ $\alpha^{0 \Downarrow}$ \ и верхние индексы \ $<\alpha^{0
\Downarrow}$, $\vartriangleleft \alpha^{0 \Downarrow}$ \ будут
опускаться.
\\
Ниже \ $\alpha^{0 \Downarrow}$ \ функция \ $\alpha S_f$ \
определена на интервале \ $]\tau_3^0, \tau_3^1[$, \ где \
$\gamma_{\tau_3^1}$ \ наследник кардинала \ $\widetilde{\delta}^0$
\ в\ $SIN_n$, \ по лемме~\ref{8.7.} об абсолютности. Отсюда и из
(\ref{e11.2}) следует
\begin{equation} \label{e11.4}
\tau_3^0 \notin dom(\alpha S_f),
\end{equation}
иначе снова возникает (\ref{e11.3}), или \ $\widetilde{\delta}^0 =
\gamma_{\tau_1^0}$ \ как результат минимизации \
$\widetilde{\delta}^0$ \ внутри \ $[\gamma_{\tau_1^0},
\gamma_{\tau_3^0}[$ \ в противоречии с предположением (всё это
ниже \ $\alpha^{0 \Downarrow}$).
\\
Но (\ref{e11.4}) возможно только если допустимая матрица \ $\alpha
S_{\tau_3^0}$ \  \textit{подавлена} для \ $\gamma_{\tau_3^0}$, \
то есть если выполняется условие подавления \ $A_5^{S,0}$ \ для \
$\alpha S_{\tau_3^0}$ \ на \ $\alpha_{\tau_3^0}$ \ характеристики
\ $a_{\tau_3^0}$ \ с базой \ $\rho_{\tau_3^0}$ (см. определение
\ref{8.1.}~2.6\;) ниже \ $\alpha^{0 \Downarrow}$, --- и теперь все
ограничения нужно подробно указать:
\[
    a_{\tau_3^0} = 0 \wedge SIN_n^{<\alpha^{0 \Downarrow}}(\gamma_{\tau_3^0}) \wedge
    \rho_{\tau_3^0} < \chi^{\ast +} \wedge \sigma(\chi^{\ast},
    \alpha_{\tau_3^0}, S_{\tau_3^0}) \wedge \qquad
\]
\[
    \wedge \exists \eta^{\ast} < \gamma_{\tau_3^0} \Big(
    A_{5.4}^{sc \vartriangleleft \alpha^{0 \Downarrow}}(\gamma_{\tau_3^0},
    \eta^{\ast}, \alpha S_f^{< \alpha^{0 \Downarrow}}|\tau_3^0,
    a_f^{< \alpha^{0 \Downarrow}}|\tau_3^0) \wedge \qquad
\]
\begin{equation} \label{e11.5}
    \wedge \forall \tau^{\prime} \big( \gamma_{\tau_3^0} <
    \gamma_{\tau^{\prime}}^{< \alpha^{0 \Downarrow}} \wedge
    SIN_n^{<\alpha^{0 \Downarrow}} (\gamma_{\tau^{\prime}}^{< \alpha^{0 \Downarrow}})
    \rightarrow \qquad \qquad
\end{equation}
\[
    \rightarrow \exists \alpha^{\prime}, S^{\prime} \big[
    \gamma_{\tau^{\prime}}^{< \alpha^{0 \Downarrow}} < \alpha^{\prime} <
    \gamma_{\tau^{\prime}+1}^{< \alpha^{0 \Downarrow}} \wedge
    SIN_n^{<\alpha^{\prime \Downarrow}} (\gamma_{\tau^{\prime}}^{< \alpha^{0 \Downarrow}} )
    \wedge \sigma(\chi^{\ast}, \alpha^{\prime}, S^{\prime}) \wedge
    \qquad \qquad
\]
\[
    \qquad \qquad \qquad \qquad \wedge
    A_{5.5}^{sc \vartriangleleft \alpha^{0 \Downarrow}} (\gamma_{\tau_3^0}, \eta^{\ast},
    \alpha^{\prime \Downarrow},
    \alpha S_f^{< \alpha^{\prime \Downarrow}},
    a_f^{< \alpha^{\prime \Downarrow}}) \big] \big) \Big).
\]

\noindent Следовательно, существуют кардиналы
\[
    \gamma^m < \gamma^{\ast} \le \gamma_{\tau_1^0} <
    \gamma_{\tau_3^0} \mbox{\it \ \ и предельный тип \ }
    \eta^{\ast},
\]
которые выполняют все составляющие его условия\ $A_{5.1}^{sc} -
A_{5.5}^{sc}$ \ ниже \ $\alpha^{0 \Downarrow}$ \ (см. определение
\ref{8.1.}~2.1--2.5\;); в частности, интервал \
$[\gamma_{\tau_1^0}, \gamma_{\tau_3^0}[$ \ является блоком типа \
$\eta^{\ast}$ \ благодаря условию
\[
    A_4^{b < \alpha^{0 \Downarrow}}(\tau_1^0, \tau_3^0,
    \eta^{\ast}, \alpha S_f^{< \alpha^{0 \Downarrow}} | \tau_3^0,
    a_f^{< \alpha^{0 \Downarrow}} | \tau_3^0)
\]
из условия \ $A_{5.4}^{sc}$ \ (см. определение
\ref{8.1.}~2.4,~2.3\;). Более того, существует \emph{последующий
за ним} максимальный блок
\begin{equation} \label{e11.6}
    [\gamma_{\tau_3^0}, \gamma_{\tau_2^2}^{< \alpha^{0
    \Downarrow}}[ \mbox{\it \   \emph{типа}\ } \ge \eta^{\ast}
    \mbox{\it \ \ ниже \ \ } \alpha^{0 \Downarrow}.
\end{equation}
Действительно, рассмотрим кардинал
\[
    \gamma^{\prime} = \gamma_{\tau_2^{\prime}}^{< \alpha^{0
    \Downarrow}} \in SIN_n^{< \alpha^{0 \Downarrow}},
    \tau_2^{\prime} > \tau_2
\]
такой, что ниже \ $\alpha^{0 \Downarrow}$
\begin{equation} \label{e11.7}
    \tau_2^{\prime} \not \in  dom(\alpha S_f^{< \alpha^{0
    \Downarrow}}).
\end{equation}
Тогда согласно  (\ref{e11.5}) существует некоторая сингулярная
матрица \ $S^{\prime}$ \ на её носителе \ $\alpha^{\prime} >
\gamma^{\prime}$ \ с кардиналом предскачка \ $\alpha^{\prime
\Downarrow}$, \ сохраняющим все \ $SIN_n^{< \alpha^{0
\Downarrow}}$-кардиналы \ $\le \gamma^{\prime}$ \ и выполняющим
ниже \ $\alpha^{\prime \Downarrow}$ условие (напомним определение
\ref{8.1.}~2.5\;):
\[
    A_{5.5}^{sc}(\gamma_{\tau_3^0}, \eta^{\ast}, \alpha^{\prime
    \Downarrow}, \alpha S_f^{< \alpha^{\prime \Downarrow}},
    a_f^{< \alpha^{\prime \Downarrow}}) ;
\]
оно означает, что весь интервал \ $[\gamma_{\tau_3^0},
\alpha^{\prime \Downarrow}[$ \ покрыт блоками типов \ $\ge
\eta^{\ast}$ ниже \ $\alpha^{\prime \Downarrow}$ \ . \ Среди них
имеется \textit{последующий} блок
\[
    [\gamma_{\tau_3^0}, \gamma_{\tau_2^3}^{< \alpha^{\prime
    \Downarrow}}[ \mbox{\it \ \ типа \ \ } \ge \eta^{\ast},
\]
поэтому можно рассмотреть его подблок \ $[\gamma_{\tau_3^0},
\gamma_{\tau_2^4}^{< \alpha^{\prime \Downarrow}}[$ \  в точности
типа \ $\eta^{\ast}$.
\\
Напомним, тип \ $\eta^{\ast}$ \ -- предельный, поэтому для каждого
\ $\gamma_{\tau}$ \ из этого подблока существует много разных
матричных носителей, допустимых для таких \ $\gamma_{\tau}$ \ по
лемме 3.2~\cite{Kiselev11} об ограничении; благодаря второй
индуктивной гипотезе все они неподавлены для всех соответствующих
 \ $\gamma_{\tau}$ \ -- и
всё это ниже\ $\alpha^{\prime \Downarrow}$.
\\
Тот же аргумент действует ниже \ $\alpha^{0 \Downarrow}$ \ и мы
возвращаемся к ситуации ниже этого кардинала. Из (\ref{e11.7})
следует
\[
    \gamma_{\tau_2^4}^{< \alpha^{\prime \Downarrow}} =
    \gamma_{\tau_2^4}^{< \alpha^{0 \Downarrow}}, \quad
    \gamma_{\tau_2^4}^{< \alpha^{0 \Downarrow}}
    \in SIN_n^{< \alpha^{0 \Downarrow}},
\]
поэтому интервал \ $[\gamma_{\tau_3^0}, \gamma_{\tau_2^4}^{<
\alpha^{0 \Downarrow}}[$ \ действительно является блоком типа \
$\eta^{\ast}$, \ но уже ниже \ $\alpha^{0 \Downarrow}$, \ который
содержит допустимый диссеминатор \ $\widetilde{\delta}^0$ \
матрицы \ $S^0$ \ на \ $\alpha^0$.

\noindent Но это составляет противоречие. С одной стороны, матрица
\ $S^0$ \ допустима для \ $\gamma_{\tau^0}$ \ и тогда согласно
замыкающему условию \ $\mathbf{K}^0$ \ она имеет диссеминатор \
$\widetilde{\delta}^0$ \ с базой \ $\rho^0 > \eta^{\ast}$. \ Но с
другой стороны, предшествующий блок \ $[\gamma_{\tau_1^0},
\gamma_{\tau_3^0}[$ \  ниже \ $\alpha^{0 \Downarrow}$ \ имеет тот
же тип \ $\eta^{\ast}$ \ и согласно (\ref{e11.2}) его левый конец
$\gamma_{\tau_1^0}$ может служить допустимым диссеминатором для \
$S^0$ \ на \ $\alpha^0$ \ с той же базой \ $\rho^0$, \ и поэтому \
$\widetilde{\delta^0} \le \gamma_{\tau_1^0}$ \ благодаря
минимальности \ $\widetilde{\delta}^0$, \ в противоречии с
предположением (\ref{e11.1}).

2. \ Итак, предложение (1a.) выполняется для \ $(\alpha_1^0,
\alpha^0, \tau^0)$ \ и осталось предположить, что (1b.) нарушается
для этой тройки, и мы возвращаемся к ситуации ниже \ $\alpha_1^0$;
\  это означает:
\\
имеется некоторая матрица \ $S^{0 1}$ \ на носителе \ $\alpha^{0
1}\neq \alpha^0$\ нулевой характеристики, допустимая для \
$\gamma^0 = \gamma_{\tau^0}$ \ вместе со своим минимальным
диссеминатором \ $\widetilde{\delta}^{0 1}$ \ и производящим
диссеминатором \ $\check{\delta}^{0 1}$ \ с базой \ $\rho^{0 1}$,
\\
но всё-таки \ $S^0$ \ на \ $\alpha^0$ \ подавлена для \ $\gamma^0
= \gamma_{\tau^0}$ \ (ниже \ $\alpha_1^0$); \ мы рассмотрим
\textit{минимальный} носитель \ $\alpha^{0 1}$ \ такой \ $S^{0
1}$. \
\\
Так как матрица  \ $S^0$ \ на \ $\alpha^0$ \ допустима для \
$\gamma^0$, \ то это подавление означает, что выполняется условие
подавления (\ref{e11.5}) ниже \ $\alpha_1^0$, \ то есть для\
$\alpha^{0 \Downarrow}$, \ $\gamma_{\tau_3^0}$, \ заменённых на
 \ $\alpha_1^0$, \ $\gamma^0$ \ соответственно везде в (\ref{e11.5}).

\noindent Отсюда следует, что
\[
    \alpha^0 < \alpha^{0 1},
\]
потому что если \ $\alpha^0 > \alpha^{0 1}$, \ то вторая
индуктивная гипотеза утверждает, что \ $S^{0 1}$ \ на \ $\alpha^{0
1}$ неподавлена для \ $\gamma^0$ \ ниже \ $\alpha_1^0$, \ и в то
же время она подавлена тем же условием подавления . Кроме того, \
$S^0$ \ на \ $\alpha^0$ \ --  это единственная матрица, допустимая
для \ $\gamma^0$ \ на носителе \ $\alpha^0 \in \; ]\gamma^0,
\alpha^{0 1}[$ \ согласно минимальности \ $\alpha^{0 1}$.
\\
Теперь это условие (\ref{e11.5}), с \ $\alpha_1^0$, $\gamma^0$ \
вместо \ $\alpha^{0 \Downarrow}$, $\gamma_{\tau_3^0}$ \
соответственно, утверждает существование кардиналов (мы сохраняем
предыдущие обозначения, чтобы подчеркнуть аналогию с рассуждениями
в части~1.):
\[
    \gamma^m < \gamma^{\ast} \le \gamma_{\tau_1^0} < \gamma^0
    \mbox{\it \ \ и предельного типа \ } \eta^{\ast},
\]
выполняющих все составляющие его условия \ $A_{5.1}^{sc} -
A_{5.5}^{sc}$; \ в частности, интервал \ $[\gamma^m,
\gamma^{\ast}[$ \ покрыт максимальными блоками типов, существенно
неубывающих до предельного ординала \ $\eta^{\ast}$; \
$[\gamma^{\ast}, \gamma_{\tau_1^0}[$ \ покрыт максимальными
блоками в точности типа \ $\eta^{\ast}$; \ $[\gamma_{\tau_1^0},
\gamma^0[$ \ тоже блок того же типа \ $\eta^{\ast}$ \ -- и так
далее.

Эти условия определяют \ $\gamma^m$, $\gamma^{\ast}$,
$\gamma_{\tau_1^0}$, \ $\eta^{\ast}$ \ единственным образом через
 \ $\gamma^0$ \ ниже \ $\alpha_1^0$ \ и
задают особый тип подавляющий этого покрытия; чтобы обсуждать его
удобны следующие вспомогательные \ $\Sigma_n$-формулы,
использующие только понятие допустимости (напомним определение
\ref{8.2.}~5\;):
\\
\quad \\
$
    \alpha \mathbf{K}^1(\gamma): \quad \exists \alpha^{\prime},
    S^{\prime} ~ \alpha \mathbf{K}(\gamma, \alpha^{\prime},
    S^{\prime});
$
\quad \\
\quad \\
$
    \alpha \mathbf{K}^2(\gamma): \quad \exists \alpha^{\prime},
    S^{\prime} ~ \exists \alpha^{\prime\prime}, S^{\prime\prime}
    \big( \alpha^{\prime} \neq \alpha^{\prime\prime} \wedge
$
\[
    \qquad\qquad\qquad
    \wedge \alpha \mathbf{K}(\gamma, \alpha^{\prime},
    S^{\prime}) \wedge \alpha \mathbf{K}(\gamma, \alpha^{\prime\prime},
    S^{\prime\prime}) \big).
\]
Первое из них означает, что существует по крайней мере
\textit{один} матричный носитель \ $\alpha^{\prime}$, \ допустимый
для \ $\gamma$; \ второе -- что существует \textit{более одного}
такого носителя \ $\alpha^{\prime} \neq \alpha^{\prime\prime}$; \
таким образом \ $\neg \alpha \mathbf{K}^1(\gamma)$ \ означает, что
таких носителей не существует вовсе.

Так как тип \ $\eta^{\ast}$ \ -- предельный, то каждый
максимальный блок \ $[\gamma_{\tau_1}, \gamma_{\tau_2}[$ \ из
покрытия интервала \ $[\gamma^{\ast}, \gamma_{\tau_1^0}[$ \
обладает следующими двумя свойствами:
\\
(i) \ если \ $\gamma_{\tau}$ \ -- внутренний кардинал в \
$[\gamma_{\tau_1}, \gamma_{\tau_2}[$, $\tau_1 < \tau < \tau_2$, \
то выполняется \ $\alpha \mathbf{K}^2(\gamma_{\tau})$; \ это
следует из второй индуктивной гипотезы и леммы
3.2~\cite{Kiselev11} об ограничении;
\\
(ii) \ если \ $\gamma_{\tau}$ \ -- это левый или правый конец
этого блока, то \ $\alpha \mathbf{K}^1(\gamma_{\tau})$ \
нарушается.
\\
В этом можно убедиться следующим образом. Предположим, что
 \ $\gamma_{\tau}$ \ это правый конец,
\ $\gamma_{\tau} = \gamma_{\tau_2}$, \ тогда существование
некоторой \ $S^{\prime}$ \ н \ $\alpha^{\prime}$, \ допустимой для
\ $\gamma_{\tau_2}$, \ влечёт объединение этого блока и следующего
за ним блока\ $[\gamma_{\tau_2}, \gamma_{\tau_3}[$ \ в единый
интервал допустимости \ $[\gamma_{\tau_1}, \gamma_{\tau_3}[$ \
типа \ $2 \eta^{\ast}$. \ И снова по второй индуктивной гипотезе и
лемме 3.2 [27] существует несколько матричных носителей \
$\alpha^{\prime}$, \ допустимых для \ $\gamma_{\tau_2}$, \ которые
становятся неподавленными для \ $\gamma_{\tau_2}$ \ и поэтому
функция \ $\alpha S_f$ \ становится определённой на всём интервале
\ $[\tau_1, \tau_3[$, \ хотя \ $[\gamma_{\tau_1},
\gamma_{\tau_2}[$ \ это \textit{максимальный} блок (всё это ниже \
$\alpha_1^0$). \ Левый конец \ $\gamma_{\tau} = \gamma_{\tau_1}$ \
нужно рассмотреть аналогичным образом.

\noindent Следовательно, для каждого \ $\gamma_{\tau} \in
[\gamma^{\ast}, \gamma^0[$ \ выполняется \ $\Delta_{n+1}$-формула:

\begin{equation} \label{e11.8}
    \alpha \mathbf{K}^2(\gamma_{\tau}) \vee \neg \alpha
    \mathbf{K^1}(\gamma_{\tau});
\end{equation}
нетрудно видеть, что та же ситуация выполняется ниже \ $\alpha^{0
1 \Downarrow}$ \ по тем же причинам.
\\
Далее производящий диссеминатор \ $\check{\delta}^{0 1}$ \ матрицы
\ $S^{0 1}$ \ на \ $\alpha^{0 1}$ \ осуществляет метод
ограничения-и-продолжения.
\\
Во-первых,
\[
    \widetilde{\delta}^{0 1} \le \gamma_{\tau_1^0};
\]
в противном случае
\[
    \gamma_{\tau_1^0} < \check{\delta}^{0 1} =
    \widetilde{\delta}^{0 1} < \gamma^0
\]
и \ $\check{\delta}^{0 1}$ \ продолжает до \ $\alpha^{0 1
\Downarrow}$ \  \ $\Sigma_{n+1}$-утверждение
\[
    \forall \gamma_{\tau} > \gamma_{\tau_1^0} \quad \alpha
    \mathbf{K}^2(\gamma_{\tau}).
\]
Этот факт вместе со второй индуктивной гпотезой влечёт
определённость функции \ $\alpha S_f^{<\alpha^{0 1 \Downarrow}}$ \
на некотором непустом множестве \ $T^{\alpha^{0 1 \Downarrow}}$ \
со свойствами $(i)$--$(iii)$, указанными в теореме~2, вопреки
минимальности \ $\alpha_1^{\ast}$.
\\
Далее, из \ $\widetilde{\delta}^{0 1} \le \gamma_{\tau_1^0}$ \
следует
\begin{equation} \label{e11.9}
\gamma^{\ast} < \check{\delta}^{0 1} \le \widetilde{\delta}^{0 1}
\le \gamma_{\tau_1^0}.
\end{equation}
Действительно, блок \ $[\gamma_{\tau_1^0}, \gamma^0[$ \ очевидно
обеспечивает справедливость следующего \
$\Sigma_{n+1}$-утверждения \ $\varphi(\tau_1^0, \tau^0,
\eta^{\ast})$ \ ниже \ $\alpha^{0 1}$:
\[
    \exists \gamma \big( \gamma_{\tau_1^0} < \gamma_{\tau^0} \le
    \gamma \wedge SIN_n(\gamma) \wedge A_{1.2}^{\vartriangleleft
    \gamma} (\tau_1^0, \tau^0, \eta^{\ast}) \wedge
\]
\[
    \wedge \forall \tau \in \; ]\tau_1, \tau^0[ \quad \alpha
    \mathbf{K}^{2 \vartriangleleft \gamma} (\gamma_{\tau}) \big);
\]
напомним, здесь \ $A_{1.2}(\tau_1, \tau^0, \eta^{\ast})$ \
означает, что интервал \ $[\gamma_{\tau_1^0}, \gamma_{\tau^0}[$ \
имеет тип \ $\eta^{\ast}$.
\\
Диссеминатор \ $\widetilde{\delta}^{0 1}$ \ попадает в интервал
 \ $[\gamma_{\tau_1^0}, \gamma^0[$ \ и поэтому
\ $\widetilde{\delta}^{0 1} = \gamma_{\tau_1^0}$, \ иначе
 \ $\widetilde{\delta}^{0 1} < \gamma_{\tau_1^0}$ \ и по
лемме 3.2~\cite{Kiselev11} появляется много носителей матрицы \
$S^{0 1}$, \ допустимых для этого \ $\gamma_{\tau_1^0}$; \ затем
по второй индуктивной гипотезе все они неподавлены для \
$\gamma_{\tau_1^0}$; \ поэтому матричная функция \ $\alpha S_f$ \
становится определённой для \ $\tau_1^0$ \ в противоречии с
минимальностью левого конца \ $\gamma_{\tau_1^0}$ \ по определению
понятия блока.

\noindent Так как \ $\widetilde{\delta}^{0 1} =
\gamma_{\tau_1^0}$, \ то замыкающее условие \ $\mathbf{K}^0$ \ для
\ $S^{0 1}$ \ на \ $\alpha^{0 1}$ \ влечёт \ $\eta^{\ast} <
\rho^{0 1}$ \ для базы \ $\rho^{0 1}$ \ диссеминатора \
$\widetilde{\delta}^{0 1}$.
\\
Но тогда производящий диссеминатор \ $\check{\delta}^{0 1}$ \ с
этой базой ограничивает \ $\Sigma_{n+1}$-утверждение
\[
    \exists \tau_1^{\prime}, \tau^{\prime} ~ \varphi(\tau_1^{\prime},
    \tau^{\prime}, \eta^{\ast}),
\]
потому что оно содержит только индивидные константы, ограниченные
этой базой \ $\rho^{0 1}$.
\\
Поэтому ниже \ $\check{\delta}^{0 1}$ \ появляются блоки типов \
$\ge \eta^{\ast}$ \ (снова благодаря второй индуктивной гипотезе).
\\
Теперь если \ $\check{\delta}^{0 1} \le \gamma^{\ast}$, \ то
нарушается условие \ $A_{5.2}^{sc}$ \ существенного неубывания
типов покрытия интервала  \ $[\gamma^m, \gamma^{\ast}[$ \ до \
$\eta^{\ast}$ \  (см. определение 8.1 2.2).

Таким образом, (\ref{e11.9}) установлено. Благодаря (\ref{e11.8})
ниже \ $\check{\delta}^{0 1}$ \ выполняется \
$\Pi_{n+1}$-утверждение
\[
    \forall \tau \big( \gamma^{\ast} < \gamma_{\tau} \rightarrow
    (\alpha \mathbf{K}^2(\gamma_{\tau}) \vee \neg \alpha
    \mathbf{K}^1(\gamma_{\tau} ) ) \big)
\]
и диссеминатор \ $\check{\delta}^{0 1}$ \ продолжает его до \
$\alpha^{0 1 \Downarrow}$ \ по лемме 6.6~\cite{Kiselev11} (для \
$m=n+1$, $\delta=\check{\delta}^{0 1}$, $\alpha_0 =
\gamma^{\ast}$, $\alpha_1 = \alpha^{0 1 \Downarrow}$). \ Но это
вызывает противоречие: (\ref{e11.8}) выполняется для \
$\gamma_{\tau} = \gamma^0$, \ хотя существует \emph{в точности
одна матрица} \ $S^0$ \ на \ $\alpha^0$ \ ниже \ $\alpha^{0 1
\Downarrow}$, допустимая для \  \ $\gamma^0$.

Итак, свойства (1a.) и (1b.) установлены. Теперь можно обратиться
к леммам \ref{8.5.}~8), \ref{8.8.}, \ref{8.10.} (для \ $\alpha_1
\le \alpha_1^{\ast}$):

(2) \ Сначала можно остановиться на лемме~\ref{8.8.}~1) \ и на
утверждении (\ref{e8.5}). Для этого сравним два интервала

\[
    ]\gamma_{\tau_1}, \gamma_{\tau_2}[\;, \quad
    ]\gamma_{\tau_1^{\prime}}, \gamma_{\tau_2}[\;.
\]
Благодаря (\ref{e8.3}), (\ref{e8.4}) \ $\widetilde{\delta}^3$ \
содержится в каждом из них, что вызывает
\[
    \gamma_{\tau_1^{\prime}} \le \gamma_{\tau_1},
\]
иначе \ $\gamma_{\tau_1} < \gamma_{\tau_1^{\prime}}$ \ и тогда
условие (i)  и (\ref{e8.3}) влекут существование некоторой матрицы
\ $S^1$ \ на её носителе \ $\alpha^1$ \ с диссеминатором \
$\widetilde{\delta}^1$ \ и базой \ $\rho^1$, \ допустимой для
 \ $\gamma_{\tau_1^{\prime}}$. \ По лемме
3.2~\cite{Kiselev11} матрица \ $S^1$ \ получает свои носители,
допустимые для каждого \ $\gamma_{\tau} \in \;
]\widetilde{\delta}^1, \gamma_{\tau_1^{\prime}}]$ \ с теми же \
$\widetilde{\delta}^1$, $\rho^1$. \ Таким образом возникает
некоторый интервал допустимости \
$[\gamma_{\tau_1^{\prime\prime}}, \gamma_{\tau_3}[$ \ с \
$\gamma_{\tau_1^{\prime\prime}} < \gamma_{\tau_1^{\prime}}$ \ и по
уже доказанному (1a.) минимальный допустимый  \
$\widetilde{\delta}^3 $ \ совпадает с \
$\gamma_{\tau_1^{\prime\prime}}$, \ \ $\widetilde{\delta}^3 =
\gamma_{\tau_1^{\prime\prime}}$,\ вопреки (\ref{e8.4}). Поэтому
выполняется \ $\gamma_{\tau_1^{\prime}} \le \gamma_{\tau_1}$ \ и
вместе с (\ref{e8.3}), (\ref{e8.4}) это влечёт (\ref{e8.5}), что
обеспечивает оставшуюся часть доказательства леммы \ref{8.8.}~1).
\\
Обращаясь к лемме \ref{8.8.}~2) рассмотрим \ $\alpha$-матрицу \
$S$ \ характеристики \ $a$ \ на носителе \ $\alpha$, \ допустимые
для
 \ $\gamma_{\tau}^{<\alpha_1}$ \ вместе с
диссеминатором \ $\widetilde{\delta}$ \ и базой \ $\rho$ \ ниже
 \ $\alpha_1$; \ можно установить
\[
    \{ \tau^{\prime}: \widetilde{\delta} <
    \gamma_{\tau^{\prime}}^{<\alpha_1} < \gamma_{\tau}^{<\alpha_1}
    \} \subseteq dom(\alpha S_f^{<\alpha_1})
\]
рассуждениями, уже использованными выше:
\\
для каждого \ $\gamma_{\tau^{\prime}}^{<\alpha_1} \in \;
]\widetilde{\delta}, \gamma_{\tau}^{<\alpha_1}[$ \ существует
много допустимых носителей матрицы \ $S$ \ по лемме
3.2~\cite{Kiselev11}, поэтому все они неподавлены благодаря (1b.)
и поэтому \ $\tau^{\prime} \in dom(\alpha S_f^{<\alpha_1})$.
\\
Такие же рассуждения устанавливают, что в лемме \ref{8.10.}
функция \ $\alpha S_f^{<\alpha_1}$ \ определена на всём интервале
\ $]\alpha \tau_1^{\ast}, \alpha \tau^{\ast 1}[$ \ для каждого \
$\alpha_1 > \alpha \delta^{\ast 1}$, $\alpha_1 \in SIN_n$, \ и \
$\alpha \delta^{\ast} = \gamma_{\alpha \tau_1^{\ast}}$ \ -- это
диссеминатор матрицы \ $\alpha S_{\alpha \tau^{\ast
1}}^{<\alpha_1}$ \ на её носителе \ $\alpha_{\alpha \tau^{\ast
1}}^{<\alpha_1}$ \ с базой \ $\alpha \rho^{\ast 1} = \rho_{\alpha
\tau^{\ast 1}}^{<\alpha_1}$.
\\
Такие же рассуждения следует использовать и в доказательстве леммы
\ref{8.5.}~8). Чтобы закончить  это доказательство для
неподавленности достаточно заметить, что если \ $S$ \ вместе с  \
$\delta$, $\rho$ \ имеет носитель \ $\alpha$, \ допустимый и
неподавленный для \ $\gamma_{\tau}^{<\alpha_1}$ \ только в \
$[\gamma_{\tau+1}^{<\alpha_1}, \alpha_1[\;$, \ то \
$\gamma_{\tau+1}^{<\alpha_1}$ \ ограничивает \
$SIN_{n-1}^{<\alpha_1}$-утверждение
\[
    \exists \alpha^{\prime} \big( \gamma < \alpha^{\prime} \wedge
    \alpha \mathbf{K}_{n-2}(\delta,
    \gamma_{\tau^n}, \gamma_{\tau}^{<\alpha_1},
    \alpha^{\prime}, \rho, S) \big)
\]
(см. доказательство леммы 8.5. 8)) для каждого \ $\gamma \in \;
]\gamma_{\tau}^{<\alpha_1}, \gamma_{\tau+1}^{<\alpha_1}[\;$, \
поэтому \ $S$ \ получает много своих носителей в \
$]\gamma_{\tau}^{<\alpha_1}, \gamma_{\tau+1}^{<\alpha_1}[$, \
допустимых для \ $\gamma_{\tau}^{<\alpha_1}$, \ которые также
становятся неподавленными для \ $\gamma_{\tau}^{<\alpha_1}$ \ ниже
\ $\alpha_1$ \ благодаря (1b.).

\begin{sloppypar}
Следующая лемма показывает, что интервалы \
$[\gamma_{\tau_1}^{<\alpha_1}, \gamma_{\tau_3}^{<\alpha_1}[$ \
определённости матричной функции \ $\alpha S_f^{<\alpha_1}$ \ с
минимальными левыми концами \ $\gamma_{\tau_1}^{<\alpha_1} \in
SIN_n^{<\alpha_1}$ \ имеют особое строение: для каждого \
$SIN_n^{<\alpha_1}$-кардинала \ \mbox{$\gamma_{\tau}^{<\alpha_1}
\in \; ]\gamma_{\tau_1}^{<\alpha_1},
\gamma_{\tau_3}^{<\alpha_1}[$} \ матрица \ $\alpha
S_{\tau}^{<\alpha_1}$ \ имеет \emph{нулевую} характеристику и
диссеминаторы \ $\check{\delta}_{\tau} \le
\widetilde{\delta}_{\tau} = \gamma_{\tau_1}^{<\alpha_1}$ \ ниже \
$\alpha_1$:
\end{sloppypar}

\begin{lemma}
\label{11.1.} \hfill {}
\\
Пусть \medskip

(i) \ $A_{1}^{< \alpha_1}(\tau _{1},\tau _{2})$;
\newline

(ii) \ $S^2$ \ это\ $\alpha$-матрица характеристики\ $a^2$ \ на
носителе
\[
    \alpha^2 \in \; ] \gamma_{\tau_2}^{<\alpha_1}, \alpha_1 [,
\]
допустимая для \ $\gamma_{\tau_2}^{<\alpha_1}$ \ ниже \ $\alpha_1$
\ вместе со своими минимальным допустимым диссеминатором \
$\widetilde{\delta}^2$ \ с базой \ $\rho^2$ \ и производящим
собственным  диссеминатором
 \ $\check{\delta}^{S^2}$;
\\
Тогда
\[
    \check{\delta}^{S^2} \le \widetilde{\delta}^2 = \gamma_{\tau_1}^{<\alpha_1}
    \mbox{\it \ и \ } a^3 = 0.
\]
\end{lemma}

\noindent \textit{Доказательство.} \ Верхние индексы \ $<
\alpha_1, \vartriangleleft \alpha_1$ \ будут опускаться как
обычно.
\\
Рассмотрим \ $\alpha$-матрицу \ $S^2$ \ на её носителе \
$\alpha^2$, \ допустимую для \ $\gamma_{\tau_2}$ \ вместе с её
минимальным диссеминатором \ $\widetilde{\delta}^2$ \ с базой \
$\rho^2$ \ и с производящим собственным диссеминатором \
$\check{\delta}^2 = \check{\delta}^{S^2}$, \ и рассмотрим ситуацию
ниже кардинала предскачка \ $\alpha^{2 \Downarrow}$.
\\
1. \ Предположим, что наоборот -- эта лемма неверна и \
$\check{\delta^2} \nleq \gamma_{\tau_1}$, \ так что
\[
    \gamma_{\tau_1} < \check{\delta}^2 < \gamma_{\tau_2}.
\]
По условию (i) \ $\gamma_{\tau_2} \in SIN_n$ \ и, следовательно,
 \ $\gamma_{\tau_2} \in SIN_n^{<\alpha^{2
\Downarrow}}$. \ Благодаря этому и лемме \ref{8.7.} допустимость
ниже \ $\alpha_1$ \ равносильна допустимости ниже \ $\alpha^{2
\Downarrow}$ \ для всякого \ $\gamma_{\tau} \in \; ]\chi^{\ast},
\gamma_{\tau_2}[$.
\\
Тогда производящий собственный диссеминатор \ $\check{\delta}^2$ \
продолжает до \ $\alpha^{2 \Downarrow}$ \  \
$\Pi_{n+1}$-предложение о существовании допустимых \
$\alpha$-матриц \ $\alpha S_f^{<\alpha^{2 \Downarrow}}$, \ как это
было несколько раз ранее, например, в форме:
\begin{equation*}
    \forall \gamma^{\prime}
    \Bigl( \gamma_{\tau_1} < \gamma^{\prime} \wedge
    SIN_{n-1}(\gamma^{\prime}) \rightarrow
    \exists \alpha, S ~
    \alpha \mathbf{K} (\gamma^{\prime}, \alpha, S) \Bigr).
\end{equation*}
В результате появляется функция \ $\alpha S_f^{<\alpha^{2
\Downarrow}}$, \ определённая на множестве со свойствами
$(i)$--$(iii)$ из теоремы~2:
\[
    T^{\alpha^{2 \Downarrow}} = \bigl \{ \tau: \gamma_{\tau_1} <
    \gamma_{\tau}^{<\alpha^{2 \Downarrow}} < \alpha^{2 \Downarrow}
    \bigl \},
\]
потому что для каждого \ $\tau \in T^{\alpha^{2 \Downarrow}}$ \
появляются некоторые \ $\alpha$-матрицы на многих носителях,
допустимые для \ $\gamma_{\tau}^{< \alpha^{2 \Downarrow}}$, \
которые неподавлены для \ $\gamma_{\tau}^{< \alpha^{2
\Downarrow}}$ \ благодаря (1b.), (1a.); но это противоречит
минимальности \ $\alpha_1^{\ast}$.
\\
2. \ Итак, \ $\check{\delta}^2 \le \gamma_{\tau_1}$; \ более того
-- выполняется \ $A^0(\tau_1)$. \ Допустим, это не так и
существует некоторая \ $\alpha$-матрица \ $S^1$ \ на её носителе
 \ $\alpha^1$, \ допустимая для \
$\gamma_{\tau_1}$ \ вместе с её минимальным диссеминатором \
$\widetilde{\delta}^1 = \gamma_{\tau_1^1}$ \ с базой \ $\rho^1$; \
поэтому \ $\tau_1^1 < \tau_1$. \ Согласно лемме \ref{8.8.} 2)  \
$]\tau_1^1, \tau_1[\; \subseteq dom(\alpha S_f)$ \ и поэтому для
каждого \ $\tau \in \; ]\tau_1^1, \tau_2[$ \ существует некоторая
\ $\alpha$-матрица на носителе, допустимые для \ $\gamma_{\tau}$.
\ Следовательно, возникает некоторый интервал допустимости \
$]\tau_1^{1 \prime}, \tau_2[$ \ с \ $\tau_1^{1 \prime} \le
\tau_1^1 < \tau_2$. \ Благодаря (1a.) \ получается
$\widetilde{\delta}^2 = \gamma_{\tau_1^{1 \prime}}$ \ и снова по
лемме  \ref{8.8.} 2) получается \ $]\tau_1^{1 \prime}, \tau_2[ \;
\subseteq dom(\alpha S_f)$ \ вопреки минимальности \ $\tau_1$, \
требуемой здесь в условии $(i)$.
\\
3. \ Таким образом, выполняется \ $A^0(\tau_1)$ \ и,
следовательно, выполняется \ $A_1^0(\tau_1, \tau_2)$; \ поэтому
лемма~\ref{10.5.} влечёт \ $a^2=0$. \ И, наконец, снова из (1a.)
следует \ $\widetilde{\delta}^2 = \gamma_{\tau_1}$.
\\
\hspace*{\fill} $\dashv$
\\

Теперь обратимся к следующему удобному понятию, которое уже было
несколько раз использовано выше, но в дальнейшем оно будет играть
ключевую роль и поэтому на него обращает особое внимание

\begin{definition} \label{11.2.}
\hfill {} \newline \hspace*{1em} Пусть \ $S$ \ это матрица на
некотором носителе \ $\alpha$ \ вместе со своим диссеминатором
 $\widetilde{\delta} < \gamma_{\tau
}^{<\alpha _{1}}$ \ с базой \ $\rho$.
\\

{\em 1)}\quad Мы будем говорить, что матрица  \ $S$ \ опирается на
этот диссеминатор \ $\widetilde{\delta}$ \ на её носителе \
$\alpha $ \ ниже \ $\alpha_1$, \ если  \ $\widetilde{\delta}$ \
попадает в некоторый блок \ $[\gamma_{\tau_1}^{<\alpha_1},
\gamma_{\tau_3}^{<\alpha^{\Downarrow}}[ $ \ типа \ $\eta$, \ то
есть если существуют ординалы \ $\tau _{1}$, $\tau _{1}^{\prime
}$, $\tau _{2}$, $\tau _{3}, \eta$ \ такие, что
\[
    \gamma _{\tau _{1}}^{<\alpha _{1}} \leq \widetilde{\delta} <
    \gamma _{\tau _{3}}^{<\alpha^{\Downarrow }} \wedge
    A_{4}^{b \vartriangleleft \alpha^{\Downarrow}}
    (\tau _{1},\tau _{1}^{\prime },\tau_{2}, \tau_{3}, \eta).
\]
Если вдобавок этот блок -- максимальный ниже \
$\alpha^{\Downarrow}$ \ и имеет тип \ $\eta < \rho $, \ то мы
будем говорить, что матрица \ $S$ \ опирается на \
$\widetilde{\delta}$ \ очень сильно.
\\

{\em 2)}\quad Пусть существует покрытие кардинала \ $\upsilon \in
\; ]\chi^{\ast}, \alpha^{\Downarrow}]$ \ блоками; это покрытие
будет называться \ $\eta$-ограниченным (ниже \
$\alpha^{\Downarrow}$) \ если типы всех его блоков в \
$]\chi^{\ast}, \alpha^{\Downarrow}]$ ограничены некоторым
постоянным ординалом \ $\eta < \chi^{\ast +}$:
\\
\quad \\
\hspace*{1em} $\forall \gamma^{\prime} < \upsilon ~ \forall
\tau_1^{\prime}, \tau_2^{\prime}, \eta^{\prime} \big( \chi^{\ast}
< \gamma_{\tau_1^{\prime}}^{<\alpha^{\Downarrow}} <
\gamma_{\tau_2^{\prime}}^{<\alpha^{\Downarrow}} < \upsilon \wedge$
\\
\quad \\ \hspace*{12em} $\wedge A_4^{M b\vartriangleleft
\alpha^{\Downarrow}}(\tau_1^{\prime}, \tau_2^{\prime},
\eta^{\prime}) \rightarrow \eta^{\prime} \le \eta \big).$ \
\\
\hspace*{\fill} $\dashv$
\end{definition}

Часть (I) следующей леммы играет роль ``несущей конструкци''  в
дальнейших рассуждениях; часть (II) будет использована в самом
конце доказательста теоремы 2 как решающий аргумент. Здесь нужно
вспомнить понятие лестницы и его разнообразных атрибутов, которые
были введены выше перед леммой \ref{9.5.} посредством формул
1.--8.; такая лестница, определённая ниже кардинала предскачка \
\mbox{$\alpha_1 = \alpha^{\Downarrow}$} \ носителя \ $\alpha$ \
матрицы \ $S$ \  формулой \ $A_8^{\mathcal{S}t \vartriangleleft
\alpha^{\Downarrow}}(\mathcal{S}t, \alpha
S_f^{<\alpha^{\Downarrow}}, a_f^{<\alpha^{\Downarrow}})$, \ должна
быть использована как функция на \ $\chi^{\ast +}$
\[
    \mathcal{S}t = \bigl( (\tau_1^{\beta}, \tau_{\ast}^{\beta},
    \tau_2^{\beta}) \bigr) _{\beta < \chi^{\ast +}}
\]
такая, что для всяких \ $\beta$, $\beta_1$, $\beta_2$: \\
\quad \\
(i) \ $ \beta < \chi^{\ast +} \rightarrow \tau_{1}^{\beta} <
\tau_{\ast }^{\beta } \leq \tau_{2}^{\beta} \wedge A_{1.1}^{Mst
\vartriangleleft \alpha^{\Downarrow}} (\tau_1^{\beta},
\tau_{\ast}^{\beta}, \tau_2^{\beta}, \alpha
S_f^{<\alpha^{\Downarrow}}, a_f^{<\alpha^{\Downarrow}}),$
\\
\quad \\
то есть \ $[ \gamma_{\tau_{1}^{\beta}}^{< \alpha^{\Downarrow}},
\gamma_{\tau_{2}^{\beta}}^{< \alpha^{ \Downarrow}} [$ \ это
максимальная единичная ступень ниже \ $\alpha^{\Downarrow}$:
\[
    A_{1.1}^{st \vartriangleleft
    \alpha^{\Downarrow}} (\tau_{1}^{\beta}, \tau_{\ast}^{\beta},
    \tau_{2}^{\beta}, \alpha S_f^{<\alpha^{\Downarrow}},
    a_f^{<\alpha^{\Downarrow}}) \wedge
    A_{1.1}^{M \vartriangleleft \alpha^{\Downarrow}}
    (\tau_{1}^{\beta}, \tau_{2}^{\beta}, \alpha S_f^{<\alpha^{\Downarrow}});
\]
(ii) \ $\beta _{1}<\beta _{2}<\chi ^{\ast +}\longrightarrow$
\[
    \longrightarrow \tau _{2}^{\beta _{1}}<\tau _{1}^{\beta
    _{2}}\wedge Od\alpha S_{f}^{<\alpha ^{\Downarrow }} ( \tau
    _{1}^{\beta _{1}},\tau _{\ast }^{\beta _{1}} ) <Od\alpha
    S_{f}^{<\alpha ^{\Downarrow }} ( \tau _{1}^{\beta
    _{2}},\tau _{\ast }^{\beta _{2}} ) ,
\]
то есть такие ступени располагаются последовательно одна после
другой и при этом их высоты строго возрастают;
\\
\quad \\
(iii) \ $\sup_{\beta } Od\alpha S_{f}^{<\alpha ^{\Downarrow }} (
\tau _{1}^{\beta },\tau _{\ast }^{\beta } ) =\chi ^{\ast +}$,
\\
\quad \\
то есть \ $h(\mathcal{S}t) = \chi^{\ast +}$ \ и высоты этих
ступеней строго возрастают до \ $\chi^{\ast +}$;
\\
\quad \\
(iv) \ для каждой максимальной единичной ступени \ $[
\gamma_{\tau_{1}}^{< \alpha^{\Downarrow}}, \gamma_{\tau_{2}}^{<
\alpha^{ \Downarrow}} [$ \ ниже \ $\alpha^{\Downarrow}$ \
соответствующая тройка ординалов \ $(\tau_1, \tau_{\ast}, \tau_2)$
\ является значением этой функции. \label{c16}
\endnote{
\ стр. \pageref{c16}. \ Последнее условие не является необходимым,
но всё-таки принимается чтобы обеспечить единственность такой
лестницы для некоторого удобства.
\\
\quad \\
} %

\quad \\
Соответственно, эта лестница \ $\mathcal{S}t$ \ завершается в
кардинале \ $\upsilon(\mathcal{S}t) = \alpha^{\Downarrow}$, \ если
её ступени располагаются конфинально этому кардиналу \
$\alpha^{\Downarrow}$, \ то есть если выполняется условие \
$H(\alpha^{\Downarrow})$: \vspace{-6pt}
\begin{multline*}
    \forall \gamma < \alpha^{\Downarrow} \exists \beta < \chi^{\ast
    +} \exists \tau_1^{\beta}, \tau_{\ast}^{\beta}, \tau_1^{\beta},
    \bigl(\gamma < \gamma_{\tau_1^{\beta}}^{<\alpha^{\Downarrow}} <
    \gamma_{\tau_2^{\beta}}^{<\alpha^{\Downarrow}} < \alpha^{\Downarrow}
    \wedge
\\
    \wedge \mathcal{S}t(\beta) = (\tau_1^{\beta},
    \tau_{\ast}^{\beta}, \tau_2^{\beta}) \bigr).
\end{multline*}

\begin{lemma}
\label{11.3.} \hfill {} \newline \hspace*{1em} Для каждой матрицы
 \ $S$ \ нулевой характеристики на носителе \ $\alpha
> \chi^{\ast}$:
\\
\quad \\
(I) Матрица \ $S$ \ на носителе \ $\alpha$ \ обладает некоторой
лестницей \ $\mathcal{S}t$.
\\
\quad \\
(II) Эта лестница \ $\mathcal{S}t$ \ завершается в кардинале
 \ $\alpha^{\Downarrow}$, \ то есть
\[
    \upsilon(\mathcal{S}t) =
    \alpha^{\Downarrow} = \sup \big\{ \gamma_{\tau_2}^{<\alpha^{\Downarrow}} :
    \exists \beta, \tau_1, \tau_{\ast} ~~ \mathcal{S} t(\beta) = (\tau_1, \tau_{\ast}, \tau_2)
    \big\}.
\]
\end{lemma}

\noindent \textit{Доказательство} I. \ Рассмотрим любой носитель
 \ $\alpha_{0}>\chi ^{\ast }$ \ матрицы
 \ $S^0$ \ нулевой характеристики на носителе \
$\alpha_{0} $ \ и его кардинал предскачка \
$\alpha^{0}=\alpha_{0}^{\Downarrow} $.
\\
По лемме~\ref{8.5.}~5) существуют \ $\tau _{1}^{\prime }$, $\tau
_{2}^{\prime }$, $\tau _{3}^{\prime }$ \ такие, что ниже \
$\alpha^0$
\[
    A_{2}^{0\vartriangleleft \alpha^0}
    ( \tau _{1}^{\prime }, \tau _{2}^{\prime }, \tau _{3}^{\prime},
    \alpha S_f^{< \alpha^0} ) \wedge
    \forall \tau^{\prime\prime} \in \; ]\tau_1^{\prime},
    \tau_2^{\prime}] ~ a_{\tau^{\prime\prime}}^{<\alpha^0} = 1
    \wedge \alpha S_{\tau _{2}^{\prime }}^{< \alpha
    ^{0}} = S^0
\]
и, следовательно, \ $a_{\tau_{2}^{\prime }}^{<\alpha ^{0}}=1$. \
Теперь перейдём к носителю \ $\alpha _{2}=\alpha _{\tau
_{2}^{\prime }}^{<\alpha ^{0}}$ \ и рассмотрим матрицу \ $S^0$ \
на этом носителе \ $\alpha _{2}=\alpha _{\tau _{2}^{\prime
}}^{<\alpha ^{0}}$ \ и его кардинал предскачка \
$\alpha^{2}=\alpha_{2}^{\Downarrow}$. \ По лемме~\ref{10.3.}
существуют ординалы \ $\tau _{\ast }^{\prime }$, $\tau
_{2}^{\prime \prime }$ \ такие, что \ $\tau _{1}^{\prime }<\tau
_{\ast }^{\prime }\leq \tau _{2}^{\prime \prime }$ \ и

\begin{equation} \label{e11.10}
    A_{1.1}^{st \vartriangleleft \alpha^2} (\tau_{1}^{\prime },
    \tau_{\ast}^{\prime}, \tau_{2}^{\prime \prime }, \alpha S_f^{< \alpha^2},
    a_f^{< \alpha^2})
    \wedge A_{1.1}^{M \vartriangleleft \alpha^2}
    (\tau _{1}^{\prime},\tau _{2}^{\prime \prime }, \alpha S_f^{< \alpha^2} );
\end{equation}

\begin{equation} \label{e11.11}
    Od\alpha S_{f}^{<\alpha ^{2}} ( \tau _{1}^{\prime},
    \tau _{\ast }^{\prime } ) > Od ( S^0 ).
\end{equation}
\vspace{0pt}

\noindent Занумеруем без пропусков все тройки ординалов \ $( \tau
_{1}^{\prime },\tau _{\ast }^{\prime },\tau _{2}^{\prime \prime} )
$, \ обладающих свойством (\ref{e11.10}) в порядке возрастания их
первых компонент, то есть определим функцию
\[
    \mathcal{S}t = \bigl( ( \tau _{1}^{\beta },\tau _{\ast }^{\beta },
    \tau_{2}^{\beta } ) \bigr) _{\beta },
\]
обладающую свойством (i) лестницы, представленным выше для \
\mbox{$\alpha^{\Downarrow} = \alpha^2$}; \ свойство (ii)
получается затем из следствия~\ref{9.6.} для \ $\alpha_1 =
\alpha^2$.
\\
Отсюда и из (\ref{e11.11}) следует, что ординал \ $Od\alpha
S_{f}^{<\alpha ^{2}} ( \tau _{1}^{\beta },\tau _{\ast }^{\beta } )
$ \  возрастает вместе с \ $\beta $ \ до \ $\chi ^{\ast +}$; \ в
противном случае можно определить ниже \ $\alpha ^{2}$ \ верхнюю
границу множества таких ординалов
\[
    \rho \in \bigl[ Od\alpha S_{f}^{<\alpha ^{2}} ( \tau _{1}^{\prime
    },\tau _{\ast }^{\prime } ) ;\chi ^{\ast +} \bigr[,
\]
а тогда по лемме 4.6~\cite{Kiselev11} о спектральном типе
\[
    \rho <Od ( S^0 )
\]
вопреки (\ref{e11.11}). Поэтому \ $dom(\mathcal{S}t) = \chi^{\ast
+}$; \ свойство (iv) очевидно из конструкции функции \
$\mathcal{S} t$.
\\
Таким образом, свойства \ (i)--(iv) установлены для носителя \
$\alpha _{2}$ \ матрицы \ $S^0$ \ и на \ $\alpha_2$ \ она обладает
этой лестницей. Тогда по лемме 5.11~\cite{Kiselev11} об
информативности матрица \ $S^0$ \ точно так же обладает некоторой
лестницей \ $\mathcal{S}t^0$ \ на её носителе \ $\alpha_0$, \
потому что это свойство является \textit{внутренним} свойством
матрицы \ $S^0$ (см. комментарий перед леммой~\ref{9.5.}\;).

  II. Обращаясь к утвеждению (II) предположим, что оно ложное и эта лестница  \ $\mathcal{S}t^0$ \ завершается в некотором
кардинале \ $\upsilon^{0} < \alpha^0 = \alpha_0^{ \Downarrow}$:
\[
    \upsilon^0 = \sup \big\{ \gamma_{\tau_2}^{<\alpha^0} :
    \exists \beta, \tau_1, \tau_{\ast} ~~ \mathcal{S} t^0(\beta) = (\tau_1, \tau_{\ast}, \tau_2)
    \big\};
\]
очевидно, \ $\upsilon^{0}$ \ содержится в  \ $SIN_n^{<\alpha^0}$ \
и имеет конфинальность \ $\chi^{\ast +}$.

\noindent Оставшаяся часть доказательства этой леммы основывается
на методе, который может быть назван \emph{зашивающм методом}; вот
его общее описание (ниже \ $\alpha^0$):
\\
Рассматривая некоторый кардинал \ $\upsilon$ \ можно столкнуться с
ситуацией, когда существуют кардиналы \ $\gamma_{\tau} < \upsilon$
, расположенные ``близко'' к этому \  \ $ \upsilon$, \ но такие,
что функция  \ $\alpha S_f$ \ предполагается не определённой для
соответствующих \ $\tau$; \ поэтому такие кардиналы \
$\gamma_{\tau} < \upsilon$ \ могут быть названы ``проколами'' во
множестве
\[
    \upsilon \cap \big \{ \gamma_{\tau}: \tau \in dom (\alpha S_f) \big \}.
\]
С целью преодолеть эту ситуацию и, тем не менее, установить
определённоть функции \ $\alpha S_f$ \ для таких проколов, нужно
осуществить следующие два действия:
\\
Нужно найти некоторую \ $\alpha$-матрицу \ $S$ \ на носителе \ $
\alpha \geq \upsilon$ \ некоторой характеристики \ $a$ \ вместе с
производящим диссеминатором \ $\check{\delta}^{\rho} < \upsilon$ \
и базой \ $\rho$ \ такими, что интервал \ $]\check{\delta}^{\rho},
\upsilon[$ \ содержит такие проколы.
\\
Одновременно с этим нужно обнаружить некоторый кардинал
\begin{equation*}
    \gamma^{\delta} \in [ \check{\delta}^{\rho}, \upsilon [ \;  \cap \; SIN_n
\end{equation*}
который вместе с \ $S$, \ $\rho$ \ \textit{нарушает} посылку
замыкающего условия
\begin{equation*}
    \mathbf{K}^0 (a, \gamma^\delta, \alpha, \rho)
\end{equation*}
или выполняет его заключение и тем самым выполняет его в целом и
по этой причине  по лемме 6.8~1)~\cite{Kiselev11} \
$\gamma^\delta$ \ становится также допустимым диссеминатором
матрицы \ $ S$ \ на
 \ $\alpha$ \ с той же базой. Более того, можно видеть, что этот
диссеминатор  \textit{допустим  и неподавлен для каждого}
\begin{equation*}
    \gamma_{\tau} \in \{ \gamma_{\tau}: \gamma^\delta <  \gamma_{\tau} <
    \upsilon \},
\end{equation*}
потому что это условие тривиально выполняется для \
$\gamma_{\tau}$ \ и поэтому всё условие допустимости
\begin{equation*}
    \alpha \mathbf{K} (a, \gamma^{\delta}, \gamma_\tau, \alpha,  \rho, S)
\end{equation*}
выполняется для многих носителей\ $\alpha
> \gamma_{\tau}$ \ матрицы \ $S$ \ тоже. Поэтому благодаря (1b.) функция \ $\alpha S_f$ \ оказывется
определённой на всём множестве
\begin{equation*}
    \{ \tau: \gamma^\delta < \gamma_\tau < \upsilon \},
\end{equation*}
и таким образом производится  ``зашивание'' интервала \ $[
\gamma^\delta, \upsilon [$ \ -- это означает, что это множество
включается в \ $dom(\alpha S_f)$ \ и этот интервал не содержит
никаких проколов вопреки предположению.
\\
Противоречие этого типа поможет продвинуть далее доказательство
леммы~\ref{11.3.}, и тем самым доказательство теоремы 2 на каждой
его решающей стадии.
\\

Итак, рассмотрим в качестве подобного \ $\upsilon$ \ кардинал \
$\upsilon^1 \in SIN_n^{<\alpha^0}$, \ являющийся \ $\chi^{\ast +}$
\--им по порядку в \ $SIN_n^{<\alpha^0}$, \ то есть множество
\[
    \upsilon^1 \cap SIN_n^{<\alpha^0}
\]
имеет порядковый тип \ $\chi^{\ast +}$; \ этот кардинал \
$\upsilon^1 \le \upsilon^0$ \ действительно существует благодаря
 \ $\upsilon^0 < \alpha^0$, $cf(\upsilon^0) =
\chi^{\ast +}$.

Так как \ $\upsilon^1 \in SIN_n^{<\alpha^0}$ \ и \ $cf(\upsilon^1)
= \chi^{\ast +}$, \ то существует \ $\delta $-матрица \ $S^1$ \
характеристики \ $a^1$, \ редуцированная к \ $\chi^{\ast}$ \ и
порождённая кардиналом \ $\upsilon ^{1}$ \ на носителе \ $\alpha^1
< \alpha^0$ \ с кардиналом предскачка \ $\alpha^{1 \Downarrow} =
\upsilon^1$ \ с производящим собственным диссеминатором \
$\check{\delta}^1 = \check{\delta}^{S^1} < \upsilon^{1}$ \ с базой
\ $\rho^1 = \rho^{S^1}$ \ по лемме~6.13~\cite{Kiselev11},
использованной здесь для \ $m=n+1$, $\alpha_{0}=\upsilon^1$,
$\alpha_1 = \alpha^0$ \ и функции
\[
    f(\beta)= OT(\beta \cap SIN_n^{<\upsilon^1}) \ ;
\]
мы рассмотрим минимальный из таких \ $\alpha^1$ \ для некоторой
определённости.

\noindent Мы увидим, что это вызывает противоречие: возникает
определённое множество
\[
    T^{\upsilon^1} = \{ \tau: \gamma < \gamma_{\tau}^{\upsilon^1} <
    \upsilon^1 \} \subseteq dom(\alpha S_f^{<\upsilon^1})
\]
выполняющее все условия теоремы~2 (для \ $\upsilon^1$ \ вместо
 \ $\alpha_1$), \ в противоречии с с минимальностью
 \ $\alpha_1^{\ast}$; \ это
противоречие устанавливает, что на самом деле лестница\
$\mathcal{S}t^0$ \ завершается в кардинале \ $\upsilon^0 =
\alpha^0 = \alpha_0^{\Downarrow}$. \ Этот результат будет
достигнут методом зашивания, применённым к \ $\upsilon^1$.
\\
Во-первых, возникает покрытие интнрвала \ $[\check{\delta}^1,
\upsilon^1[$ \ максимальными блоками \ (ниже \ $\upsilon^1$). \
Допустим, что это не так, тогда существует некоторый кардинал
\[
    \gamma^1 \in SIN_n^{<\upsilon^1} \cap \; ]\check{\delta}^1, \upsilon^1[
\]
который не принадлежит ни одному блоку (ниже \ $\upsilon^1$). \
Тогда этот \ $\gamma^1$ \ может служить диссеминатором
\[
    \widetilde{\delta}^1 = \gamma^1
\]
с той же базой \ $\rho^1$ \ для той же матрицы \ $S^1$ \
характеристики \ $a^1$ \  на носителе \ $\alpha^1$\   по лемме
6.8~\cite{Kiselev11} (для \ $m=n+1$), \ допустимым для каждого \
$\gamma_{\tau} \in \; ]\widetilde{\delta}^1, \upsilon^1[\;$, \ так
как выполняется замыкающее \ $\Delta_1$-условие \
$\mathbf{K}^0(\alpha^1, \widetilde{\delta}^1, \alpha^1, \rho^1)$

\begin{equation} \label{e11.12}
    \big( a^1 = 0 \rightarrow \forall \tau_1^{\prime},
    \tau_1^{\prime\prime}, \tau_2^{\prime}, \tau_3^{\prime},
    \eta{\prime} < \alpha^{1 \Downarrow} \big[
    \gamma_{\tau_1^{\prime}}^{<\alpha^{1 \Downarrow}} \le
    \widetilde{\delta}^1 < \gamma_{\tau_3^{\prime}}^{<\alpha^{1
    \Downarrow}} \wedge
    \qquad\qquad\qquad\qquad
\end{equation}
\[
    \wedge A_4^{M \vartriangleleft \alpha^{1 \Downarrow}}
    (\tau_1^{\prime}, \tau_1^{\prime\prime},
    \tau_2^{\prime}, \tau_3^{\prime}, \eta{\prime}, \alpha
    S_f^{<\alpha^{1 \Downarrow}}, a_f^{<\alpha^{1 \Downarrow}})
    \rightarrow \eta^{\prime} < \rho^1 \vee \rho^1 = \chi^{\ast +}
    \big] \big)
\]
\vspace{0pt}

\noindent благодаря ложности его посылки \ $A_4^{M
\vartriangleleft \alpha^{1 \Downarrow}}$. \ Теперь действует
зашивающий метод: для каждого \ $\gamma_{\tau} \in \;
]\widetilde{\delta}^1, \upsilon^1[$ \ и для
\[
    \gamma_{\tau^n} = \sup\{\gamma \le \gamma_{\tau}: \gamma \in
    SIN_n^{<\upsilon^1} \}
\]
выполняется \ $\Pi_{n-2}$-условие \ $\varphi(a^1,
\widetilde{\delta}^1, \gamma_{\tau^n}, \gamma_{\tau}, \alpha^1,
\rho^1, S^1)$:
\[
    \gamma_{\tau} < \alpha^1 \wedge SIN_n^{<\alpha^{1
    \Downarrow}}(\gamma_{\tau^n}) \wedge
    \alpha \mathbf{K}_{n+1}^{\exists}(a^1, \widetilde{\delta}^1,
    \gamma_{\tau}, \alpha^1, \rho^1, S^1)
\]
утверждающее, что \ $S^1$ \ на \ $\alpha^1$ \ допустима для  \
$\gamma_{\tau}$ \ вместе с теми же \ $a^1$,
$\widetilde{\delta}^1$, $\rho^1$. \ Тогда \ $SIN_{n-1}$-кардинал \
$\gamma_{\tau+1}$ \ ограничивает \ $\Sigma_{n-1}$-утверждение
\begin{equation} \label{e11.13}
    \exists \alpha \big( \gamma_{\tau} < \alpha \wedge \varphi(a^1,
    \widetilde{\delta}^1, \gamma_{\tau^n}, \gamma_{\tau},
    \alpha, \rho^1, S^1) \big)
\end{equation}
и поэтому в \ $]\gamma_{\tau}, \gamma_{\tau+1}[$ \ появляется
много допустимых для \ $\gamma_{\tau}$ \ носителей \ $\alpha$ \ с
тем же свойством (\ref{e11.13}) и это вызывает противоречие ниже \
$\upsilon^1$:
\\
все они неподавлены для \ $\gamma_{\tau}$ \ благодаря (1b.) и
функция  \ $\alpha S_f^{< \upsilon^1}$ \ становится определенной
для \ $\gamma_{\tau}$ \ и, таким образом, становится определённой
для всего интервала \ $[\check{\delta}^1, \upsilon^1[$ \ (то есть
совершилось зашивание этого интервала); но это и составляет
противоречие -- появляется некоторое множество \ $T^{\upsilon^1}$
\ определённости этой функции \ $\alpha S_f^{< \upsilon^1}$ \ со
свойствами $(i)$-$(iii)$ из теоремы 2 (для \
$\alpha_1=\upsilon^1$), \ в противоречии с минимальностью \
$\alpha_1^{\ast}$.

Итак, интервал \ $[\check{\delta}^1, \upsilon^1[$ \ покрывается
максимальными блоками ниже \ $\upsilon^1$ \ и выполняется условие
\[
    A_{5.1}^{sc \vartriangleleft \upsilon^1} (\gamma^m, \alpha
    S_f^{< \upsilon^1}, a_f^{< \upsilon^1})
\]
утверждающее покрытие интервала \ $[\gamma^m, \upsilon^1[$ \
максимальными блоками ниже \ $\upsilon^1$ \ и минимальность
кардинала \ $\gamma^m$ \ с этим свойством (напомним определение
 \ref{8.1.}~2.1a.,~2.1b.).

\noindent Оставшаяся часть доказательства этой леммы~\ref{11.3.}
проводится ниже \ $\upsilon^1$ \ и верхние индексы \ $<
\upsilon^1$, \ $\vartriangleleft \upsilon^1$ \ и обозначения
функций \ $\alpha S_f^{< \upsilon^1}$, $a_f^{< \upsilon^1}$ \
будут как обычно опускаться (когда контекст будет очевидно на них
указывать).
\\
Мы переходим к заключительному противоречию этого доказательства:
\\
это покрытие не может быть \ $\eta$-ограниченным, и в то же время
оно должно быть \ $\eta$-ограниченным ниже \ $\upsilon^1$ \ (см.
определение \ref{11.2.}~2) для \ $\alpha^{\Downarrow} =
\upsilon^1$).

Действительно, это покрытие не может быть\ $\eta$-ограниченным,
так как в противном случае существует некоторый постоянный тип его
максимальных блоков, располагающихся конфинально \ $\upsilon^1$. \
Минимальный тип \ $\eta^1$ \ из таких типов очевидно определяется
ниже \ $\upsilon^1 = \alpha^{1 \Downarrow}$ \ и по лемме
4.6~\cite{Kiselev11} о спектральном типе это влечёт
\[
    \eta^1 < Od(S^1) < \rho^1.
\]
Но тогда действует метод зашивания. Пусть \ $[\gamma_{\tau_1},
\gamma_{\tau_2}[$ \ это максимальный блок в \
$[\widetilde{\delta}^1, \upsilon^1[$ \ типа \ $\eta^1$ \ с
минимальным левым концом \ $\gamma_{\tau_1}$, \ тогда \
$SIN_n$-кардинал \ $\gamma_{\tau_1}$ \ может служить
диссеминатором \ $\widetilde{\delta}^{1 \prime} = \gamma_{\tau_1}$
\ для \ $S^1$ \ на on \ $\alpha^1$ \ с той же базой \ $\rho^1$.
\\
И снова выполняется утверждение~(\ref{e11.12}) (где \
$\widetilde{\delta}^1$ \ следует заменить на
$\widetilde{\delta}^{1 \prime}$), \ но теперь потому, что \ $S^1$
\ на \ $\alpha^1$ \ опирается на \ $\widetilde{\delta}^{1 \prime}$
\ очень сильно: существует единственный максимальный блок \
$[\gamma_{\tau_1^{\prime}}^{<\alpha^{1 \Downarrow}},
\gamma_{\tau_3^{\prime}}^{<\alpha^{1 \Downarrow}}[\;$, \
совпадающий с \ $[\gamma_{\tau_1}, \gamma_{\tau_2}[$ \ типа \
$\eta^1$, \ который содержит \ $\widetilde{\delta}^{1 \prime}$ \ и
который выполняет его заключение \ $\eta^1 < \rho^1$. \ И снова
выполняется (\ref{e11.13}) для каждого \ $\gamma_{\tau} \in \;
]\widetilde{\delta}^{1 \prime}, \upsilon^1[$ \ и поэтому таким же
образом  существует множество \ $T^{\upsilon^1}$, \ вопреки
минимальности \ $\alpha_1^{\ast}$.
\\
Следовательно, покрытие интервала \ $[\gamma^m, \upsilon^1[$ \ не
должно быть \ $\eta$-ограниченным, а тогда типы его максимальных
блоков должны существенно неубывать до \ $\chi^{\ast +}$ \
 (ниже \ $\upsilon^1$), то есть:
\begin{equation} \label{e11.14}
    \forall \eta < \chi^{\ast +} ~ \exists \gamma^{\prime} <
    \upsilon^1 ~ \forall \tau_1^{\prime}, \tau_2^{\prime},
    \eta^{\prime} ~ \big( \gamma^{\prime} < \gamma_{\tau_2^{\prime}} \le
    \upsilon^1 \wedge \qquad \qquad
\end{equation}
\[
    \qquad \qquad \wedge A_4^{M b}(\tau_1^{\prime}, \tau_2^{\prime},
    \eta^{\prime}, \alpha S_f^{<\upsilon^1}, a_f^{<\upsilon^1})
    \rightarrow \eta < \eta^{\prime} \big );
\]
в противном случае снова существует некоторый постоянный тип его
максимальных блоков, располагающихся конфинально \ $\upsilon^1$ и
приводящий к прежнему противоречию.
\\ Но и это приводит к противоречию: бесконечно много значений матричной функции \
$\alpha S_f^{<\upsilon^1}$ \ становятся подавленными (ниже \
$\upsilon^1$), хотя они неподавлены по определению~\ref{8.3.} этой
матричной функции.
\\
\noindent Чтобы убедиться в этом, нужно применить следующий метод
рассуждения, который может быть назван ``отсечение блоков справа''
и который состоит в ``укорачивании слишком длинных блоков'' с
правого  конца посредством ``отсечения'' их конечных подинтервалов
справа.
\\
Этот метод будет действовать здесь вполне успешно, потому что это
покрытие не является \ $\eta$-ограниченным и поэтому действует на
такие блоки как подавляющее покрытие, выполняя условие (см.
определение \ref{8.1.}~2.4 для \ $X_1 = \alpha
S_f^{<\upsilon^1}|\tau$, $X_2 = a_f^{<\upsilon^1}|\tau$):
\[
    A_{5.4}^{sc}(\gamma_{\tau}, \eta^{\ast}, \alpha S_f^{<\upsilon^1}|\tau,
    a_f^{<\upsilon^1}|\tau)
\]
\
для бесконечно многих кардиналов \ $\gamma_{\tau}$, \
расположенных конфинально \ $\upsilon^1$, \ и для некоторых
соответствующих \ $\gamma^m$, $\gamma^{\ast}$, $\gamma^1$,
$\eta^{\ast}$ .

\noindent Здесь кардинал \ $\gamma^m$ \ уже определён выше как
минимальный из левых концов блоков покрытия кардинала \
$\upsilon^1$.
\\
Далее, ординалы \ $\gamma^{\ast}$, $\eta^{\ast}$ \ могут быть
определены здесь разными способами, например, согласно
(\ref{e11.14}) как предельные значения следующих
последовательностей  ниже \ $\upsilon^1$ \ (обозначения \ $\alpha
S_f^{<\upsilon^1}$, $a_f^{<\upsilon^1}$ \ будут опущены):
\[
    \gamma_0 = \gamma^m; \qquad\qquad\qquad\qquad\qquad\qquad\qquad\qquad\qquad\qquad
\]
\[
    \eta_i = \sup \big \{ \eta: ~ \exists \tau_1, \tau_2 \big (
    \gamma_{\tau_1} < \gamma_{\tau_2} < \gamma_i \wedge A_4^{M
    b}(\tau_1, \tau_2, \eta) \big) \big \};
\]
\[
    \gamma_{i+1} = \min \big \{ \gamma: ~ \exists \tau_1, \tau_2, \eta \big (
    \gamma_i < \gamma_{\tau_1} < \gamma_{\tau_2} = \gamma \wedge
    \eta_i < \eta \wedge
\]
\begin{equation} \label{e11.15}
    \qquad \qquad \wedge A_4^{M b}(\tau_1, \tau_2, \eta) \wedge
    \forall \tau_1^{\prime}, \tau_2^{\prime}, \eta^{\prime} \big(
    \gamma \le \gamma_{\tau_1^{\prime}} < \gamma_{\tau_2^{\prime}}
    \wedge
\end{equation}
\[
    \qquad \qquad \qquad \wedge A_4^{M b}(\tau_1^{\prime}, \tau_2^{\prime}, \eta^{\prime})
    \rightarrow \eta \le \eta^{\prime} \big ) \big) \big \};
\]
\[
    \qquad\qquad\qquad\qquad
    \eta^{\ast} = \sup_{i \in \omega_0} \eta_i; \quad \gamma^{\ast} =
    \sup_{i \in \omega_0} \gamma_i.
\]
Так как  типы максимальных блоков этого покрытия существенно
неубывают до \ $\chi^{\ast +}$, \ то существует максимальный блок
в \ $[\gamma^{\ast}, \upsilon^1[$
\[
    [\gamma_{\tau_1^{\ast}}, \gamma_{\tau_2^{\ast}}[ \mbox{\it \ \ большего типа \ }
    \eta^{\ast 1} > \eta^{\ast}
\]
и следует рассмотреть такой блок с \textit{минимальным} левым
концом \ $\gamma_{\tau_1^{\ast}}
> \gamma^{\ast}$. \ Очевидно, этот блок включает в себя свой начальный подинтервал
\[
    [\gamma_{\tau_1^{\ast}}, \gamma^{\ast 1}[ \mbox{\it \ \  типа в точности \ }
    \eta^{\ast},
\]
который тоже является блоком (не максимальным) с правым концом \
$\gamma^{\ast 1} \in SIN_n$, \mbox{$\gamma^{\ast 1} =
\gamma_{\tau_2^{\ast 1}}$}.
\\
Следовательно, существует матрица
\[
    S^{\ast 1} = \alpha S_{\tau_2^{\ast 1}}
\]
на её носителе \ $\alpha^{\ast 1} = \alpha_{\tau_2^{\ast 1}}$ \
характеристики \ $a^{\ast 1} = a_{\tau_2^{\ast 1}}$, \ которая
допустима и \textit{неподавлена} для \ $\gamma^{\ast 1}$ \ вместе
со своим диссеминатором \ $\widetilde{\delta}^{\ast 1} =
\widetilde{\delta}_{\tau_2^{\ast 1}}$ \ с базой \ $\rho^{\ast 1} =
\rho_{\tau_2^{\ast 1}}$ \ по определению (всё это ниже \
$\upsilon^1$).
\\
Но в то же время эта матрица \ $S^{\ast 1}$ \ и все её атрибуты,
наоборот, подавлены для \ $\gamma^{\ast 1}$, \  так как для них
выполняется условие подавления \ $A_5^{S,0}$ \  ниже \
$\upsilon^1$ \ (напомним определение \ref{8.1.}~2.6, а также
(\ref{e11.5})\;), которое для этой ситуации можно сформулировать
так:
\[
    a^{\ast 1} = 0 \wedge SIN_n^{<\upsilon^1}(\gamma^{\ast 1})
    \wedge \rho^{\ast 1} < \chi^{\ast +} \wedge \sigma(\chi^{\ast},
    \alpha^{\ast 1}, S^{\ast 1}) \wedge
\]
\[
    \wedge \exists \eta^{\ast}, \tau < \gamma^{\ast 1} \Big(
    \gamma^{\ast 1} = \gamma_{\tau}^{<\upsilon^1} \wedge
    A_{5.4}^{sc \vartriangleleft \upsilon^1} (\gamma^{\ast 1},
    \eta^{\ast}, \alpha S_f^{<\upsilon^1}|\tau, a_f^{<\upsilon^1}|\tau)
    \wedge
\]
\begin{equation} \label{e11.16}
    \wedge \forall \tau^{\prime} \Big( \tau < \tau^{\prime} \wedge
    SIN_n^{< \upsilon^1} (\gamma_{\tau^{\prime}}^{< \upsilon^1}) \rightarrow
    \qquad \qquad \qquad
\end{equation}
\[
    \rightarrow \exists \alpha^{\prime}, S^{\prime} \big [
    \gamma_{\tau^{\prime}}^{<\upsilon^1} < \alpha^{\prime} <
    \gamma_{\tau^{\prime} + 1}^{<\upsilon^1} \wedge
    SIN_n^{<\alpha^{\prime \Downarrow}} (\gamma_{\tau^{\prime}}^{< \upsilon^1})
    \wedge \sigma(\chi^{\ast}, \alpha^{\prime}, S^{\prime}) \wedge
\]
\[
    \qquad\qquad\qquad\qquad\qquad\qquad \wedge A_{5.5}^{sc \vartriangleleft \upsilon^1}(\gamma^{\ast 1}, \eta^{\ast},
    \alpha^{\prime \Downarrow}, \alpha S_f^{<\alpha^{\prime \Downarrow}},
    a_f^{<\alpha^{\prime \Downarrow}} ) \big] \Big) \Big).
\]
Здесь действительно \ $a^{\ast 1} = 0$ \ по лемме \ref{11.1.}; \
$SIN_n^{<\upsilon^1}({\gamma^{\ast 1}})$ \ выполняется по
построению; \ $\rho^{\ast 1} < \chi^{\ast +}$ \ так как \
$\upsilon^1$ \ это \ $\chi^{\ast +}$--\ кардинал по порядку в
$SIN_n^{<\upsilon^1}$; \ $\sigma(\chi^{\ast}, \alpha^{\ast 1},
S^{\ast 1})$ \ выполняется благодаря допустимости \ $S^{\ast 1}$ \
на \ $\alpha^{\ast 1}$ \ для \ $\gamma^{\ast 1}$; \ $A_{5.4}^{sc
\vartriangleleft \upsilon^1}$ \ выполняется потому, что типы
покрытия \ $\gamma^{\ast}$ \ существенно неубывают до \
$\eta^{\ast}$ \ по (\ref{e11.15}); и максимальные блоки из
интервала \ $[\gamma^{\ast}, \gamma_{\tau_1^{\ast}}[$ \ имеют
постоянный тип \ $\eta^{\ast}$ \ благодаря минимальности \
$\gamma_{\tau_1^{\ast}}$ -- и осталось проверить условие \
$A_{5.5}^{sc}$ \ из (\ref{e11.16}). Для этого следует применить
обычный ограничивающий аргумент:
\\
Каждый максимальный блок \ $[\gamma_{\tau_1}, \gamma_{\tau_2}[$ \
в \ $[\gamma^{\ast}, \upsilon^1[$ \ имеет тип \ $\eta \ge
\eta^{\ast}$ \ согласно (\ref{e11.15}) и поэтому выполняется
следующее \ $\Pi_{n-2}$-утверждение \ $\psi(\gamma^{\ast 1},
\eta^{\ast}, \alpha^1, S^1, \alpha S_f^{<\alpha^{1 \Downarrow}},
a_f^{<\alpha^{1 \Downarrow}})$:
\[
    \sigma(\chi^{\ast}, \alpha^1, S^1) \wedge A_{5.5}^{sc}
    (\gamma^{\ast 1}, \eta^{\ast}, \alpha^{1 \Downarrow}, \alpha
    S_f^{<\alpha^{1 \Downarrow}}, a_f^{<\alpha^{1 \Downarrow}}),
\]
где \ $A_{5.5}^{sc}$ \ -- это \ $\Delta_1$-формула (см.
определение
 \ref{8.1.}~2.5\;):
\begin{multline*}
    \forall \gamma^{\prime} \Big( \gamma^{\ast 1} \le
    \gamma^{\prime} < \alpha^{1 \Downarrow} \rightarrow
    \exists \tau_1^{\prime}, \tau_2^{\prime}, \eta^{\prime} \big(
    \gamma_{\tau_1^{\prime}}^{<\alpha^{1 \Downarrow}} \le
    \gamma^{\prime} < \gamma_{\tau_2^{\prime}}^{<\alpha^{1
    \Downarrow}}\wedge
\\
    \wedge A_4^{M \vartriangleleft \alpha^{1
    \Downarrow}} (\tau_1^{\prime}, \tau_2^{\prime}, \eta^{\prime},
    \alpha S_f^{<\alpha^{1 \Downarrow}}, a_f^{<\alpha^{1 \Downarrow}}
    ) \wedge \eta^{\prime} \ge \eta^{\ast} \big ) \Big).
\end{multline*}
\vspace{0pt}

\noindent Теперь рассмотрим любой кардинал \
$\gamma_{\tau^{\prime}}
> \gamma^{\ast 1}$, $\gamma_{\tau^{\prime}} \in
SIN_n^{<\upsilon^1}$; \ согласно лемме 3.2~\cite{Kiselev11} об
ограничении \ $SIN_{n-1}$-кардинал \ $\gamma_{\tau^{\prime}+1}$ \
ограничивает \ $\Sigma_{n-1}$-утверждение \ $\exists
\alpha^{\prime} ~ \psi_1(\gamma^{\ast 1}, \eta^{\ast},
\alpha^{\prime}, \gamma_{\tau^{\prime}})$, где \ $\psi_1$ это
формула: \vspace{-6pt}
\begin{multline*}
     \exists S^{\prime} \big [
    \gamma_{\tau^{\prime}} < \alpha^{\prime} \wedge SIN_n^{<\alpha^{\prime
    \Downarrow}}(\gamma_{\tau^{\prime}}) \wedge
\\
    \wedge \psi(\gamma^{\ast 1},
    \eta^{\ast}, \alpha^{\prime}, S^{\prime}, \alpha S_f^{<\alpha^{\prime \Downarrow}},
    a_f^{<\alpha^{\prime \Downarrow}}) \big ],
\end{multline*}
\vspace{0pt}  и поэтому некоторые носители  \ $\alpha^{\prime}$ \
матрицы
 \ $S^{\prime}$ \ с этим свойством появляются в \
$]\gamma_{\tau^{\prime}}, \gamma_{\tau^{\prime}+1}[$. \
Следовательно, ниже \ $\upsilon^1$ \ выполняется утверждение:

\[
    \forall \tau^{\prime} \big ( \tau_2^{\ast 1} < \tau^{\prime}
    \wedge SIN_n^{<\upsilon^1}(\gamma_{\tau^{\prime}}) \
    \rightarrow \exists \alpha^{\prime} < \gamma_{\tau^{\prime}+1}
    ~ \psi_1^{\vartriangleleft \upsilon^1} (\gamma^{\ast 1},
    \eta^{\ast}, \alpha^{\prime}, \gamma_{\tau^{\prime}}) \big) .
\]
\vspace{-6pt}

\noindent В результате условие подавления (\ref{e11.16})
выполняется в целом для матрицы \ $S^{\ast 1}$ \ на её носителе \
$\alpha^{\ast 1}$ \ и она не может быть значением матричной
функции \ $\alpha S_f^{<\upsilon^1}$ \ ниже \ $\upsilon^1$ \
вопреки предположению.
\\
\hspace*{\fill} $\dashv$
\\
\quad \\

Теперь заключительная часть доказательства теоремы~2 приходит к
своему завершению.  Снова все рассуждения будут релятивизированы к
 \ $\alpha _{1}^{\ast
}$ \ и поэтому верхние индексы \ $<\alpha _{1}^{\ast }$,
$\vartriangleleft \alpha_{1}^{\ast }$ \ и обозначения функций \
$\alpha S_f^{<\alpha_1^{\ast}}$, $a_f^{<\alpha_1^{\ast}}$ \ в
записях формул будут опускаться.
\\
По предположению эта теорема нарушается для минимального кардинала
 \ $\alpha_{1}^{\ast }$, \ поэтому существуют \ $\tau _{2}^{\ast }$,
$\tau _{3}^{\ast }$ \ \ $<\alpha _{1}^{\ast }$ такие, что
выполняется
\[
    A_{2} (\tau_{1}^{\ast }, \tau_{2}^{\ast }, \tau _{3}^{\ast }
    )  ,
\]
где, напомним, \ $\tau_1^{\ast}$ \ это минимальный ординал во
множестве \ $T^{\alpha_1^{\ast}}$ \ и где \ $\tau _{2}^{\ast } $ \
это минимальный ординал на котором нарушается монотонность функции
\ $\alpha S_{f}$ \ на \ $T^{\alpha_1^{\ast}}$; \ рассмотрим
произвольно большой ординал \ $\tau_3^{\ast}$ \ из множества
\[
    Z^{\ast }=\{ \tau : \gamma _{2}^{\ast } < \gamma _{\tau } <
    \alpha_{1}^{\ast } \wedge \gamma _{\tau } \in SIN_{n} \},
\]
и рассмотрим соответствующие кардиналы
\[
    \gamma _{i}^{\ast } = \gamma_{\tau _{i}^{\ast }}, ~
    i=\overline{1,3}\;,
\]
а также матрицу \ $S^{\ast 2}=\alpha S_{\tau _{2}^{\ast }}$ \ на
носителе \ $\alpha_{\tau_{2}^{\ast }}$ \ характеристики \ $a^{\ast
2} = a_{\tau_{2}^{\ast}}$ \ с кардиналом предскачка \
$\alpha^{\ast 2} = \alpha_{\tau_2^{\ast}}^{\Downarrow}$ \ и её
производящий собственный диссеминатор\ $\check{\delta}^{\ast 2} =
\check{\delta}_{\tau_{2}^{\ast }}^{S}$.
\\
Но главную роль будет играть далее матрица
\[
    S^{\ast 3} = \alpha S_{\tau_3^{\ast}} \mbox{\it \ \ на носителе \ \ }
    \alpha_{\tau_3^{\ast}}
\]
для этого \ $\tau_3^{\ast} \in Z^{\ast}$ \ с кардиналом предскачка
\ $\alpha^{\ast 3} = a_{\tau_3^{\ast}}^{\Downarrow}$ \ и с
производящим и плавающим диссеминаторами
\[
    \check{\delta}^{\ast 3} = \check{\delta}_{\tau_3^{\ast}},~
    \widetilde{\delta}^{\ast 3} =
    \widetilde{\delta}_{\tau_3^{\ast}}
    \mbox{\it \ \ с базой \ \ } \rho^{\ast 3} =
    \rho_{\tau_3^{\ast}}.
\]

\noindent Из леммы \ref{11.1.} (для \ $\tau_1^{\ast}, \tau, \alpha
S_{\tau}, \alpha_{\tau}, \alpha_1^{\ast}$ \ как \ $\tau_1, \tau_2,
S^2, \alpha^2, \alpha_1$) \ следует:
\begin{equation}
\label{e11.6o} \forall \tau \in Z^{\ast }\  ( a_{\tau }=0\wedge
\widetilde{\delta }_{\tau }=\gamma _{1}^{\ast } ).
\end{equation}

\noindent Теперь возникают следующие случаи:

Случай 1.\quad $a^{\ast 2}=1$.\quad  Тогда по лемме~\ref{10.3.}
существует \ $\tau _{1}^{\ast \prime }$ \ такой, что
\begin{equation}
\label{e11.7o} A_{3} (\tau _{1}^{\ast }, \tau_{1}^{\ast \prime },
\tau _{2}^{\ast }, \tau_{3}^{\ast } )
\end{equation}
(см. определение 81. 1.5) где, напомним, матрица \ $\alpha
S_{\tau_1^{\ast \prime}}$ \ имеет единичную характеристику на её
носителе \ $\alpha_{\tau_1^{\ast \prime}}$.
\\
Начиная с этого места следует использовать только такой ординал \
$\tau_3^{\ast}$, \ что соотвествующий интервал
\[
    [ \gamma_1^{\ast}, \gamma_{\tau_3^{\ast}}[
\]
 имеет тип
\[
    \eta^{\ast 3} > Od(\alpha S_{\tau_{1}^{\ast \prime}});
\]
существование таких ординалов \ $\tau_3^{\ast}$ \ следует из
условий $(i)$, $(iii)$ этой теоремы~2.

\noindent Рассмотрим следующие подслучаи: \\ Подслучай 1a. \quad
Предположим, что
\[
    \alpha S_{\tau_1^{\ast \prime}} \vartriangleleft \rho^{\ast 3}.
\]
Но это исключается рассуждением ограничения-и-продолжения,
обеспечивающим здесь следующий аргумент, который может быть назван
``срезание лестницы сверху'' и который состоит в срезании высот
ступеней лестницы:
\\
Матрица \ $\alpha S_{\tau_1^{\ast \prime}}$ \ единичной
характеристики имеет допустимый носитель
\[
    \alpha_{\tau_1^{\ast \prime}} \in \; ]
    \widetilde{\delta}_3^{\ast}, \gamma_3^{\ast}[ \; ,
\]
так как \ $\widetilde{\delta}_3^{\ast} = \gamma_1^{\ast}$. \ По
лемме 3.2~\cite{Kiselev11} об ограничении эта матрица получает
свои допустимые носители, расположенные конфинально диссеминатору
\ $\check{\delta}^{\ast 3} \le \widetilde{\delta}^{\ast 3}$, \ и
поэтому выполняется следующее \ $\Pi_{n+1}$-предложение ниже
 \ $\check{\delta}^{\ast 3}$:
\[
    \forall \gamma ~ \exists \gamma^{1} > \gamma ~
    \exists \delta, \alpha, \rho ~ \big( SIN_{n-1}(\gamma^1) \wedge \alpha
    \mathbf{K}(1, \delta, \gamma^{1}, \alpha, \rho,
    \alpha S_{\tau_1^{\ast \prime}}) \big)
\]
и благодаря  \ $ S_{\tau_1^{\ast \prime}} \vartriangleleft
\rho^{\ast 3} $\ и лемме 6.6~\cite{Kiselev11} (для \ $m=n+1$) \
этот диссеминатор продолжает это утверждение до кардинала
предскачка \ $\alpha^{\ast 3}$ \ и поэтому \ $\alpha
S_{\tau_1^{\ast \prime}}$ \ получает свои допустимые носители ниже
\ $\alpha^{\ast 3}$, \ расположенные конфинально \  \
$\alpha^{\ast 3}$.
\\
По лемме \ref{11.3.} существует лестница \ $\mathcal{S}t$ \ ниже \
$\alpha^{\ast 3}$, \ завершающаяся в \ $\alpha^{\ast 3}$, \ но по
лемме~\ref{9.5.}~2b. (о срезании лестницы сверху, где \
$\alpha_1$, $S^0$ \ нужно заменить на \ $\alpha^{\ast 3}$, \
$\alpha S_{\tau_1^{\ast \prime}}$) \ эта лестница \ $\mathcal{S}
t$ \ невозможна, так как высоты всех ступеней \
$\mathcal{S}t(\beta)$ \ этой лестницы ограничиваются ординалом
\[
    Od( \alpha S_{\tau_1^{\ast \prime}}) < \chi^{\ast
    +}
\]
(срезаются этим ординалом), хотя они возрастают до \ $\chi^{\ast
+}$ \ по определению.
\\
Подслучай 1b. \quad Таким образом
\[
    \rho^{\ast 3} \le Od( \alpha S_{\tau_1^{\ast \prime}} ).
\]
Но напомним, что здесь используется тип
\[
    \eta^{\ast 3} > Od(\alpha S_{\tau_1^{\ast\prime}})
\]
блока \ $[\gamma_1^{\ast}, \gamma_3^{\ast}[\;$.
\\
Согласно (\ref{e11.7o}) существуют ординалы \ $\tau_{3}^{\prime
}$, $\eta _{3}^{\prime }$ \ такие, что для \ $\alpha^{\ast 3} =
\alpha_{\tau_{3}^{\ast}}^{\Downarrow}$ \ выполняется
\begin{equation} \label{e11.8o}
    A_{4}^{M b \vartriangleleft \alpha^{\ast 3}}
    (\tau_{1}^{\ast}, \tau_{1}^{\ast \prime},
    \tau_{2}^{\ast}, \tau_{3}^{\prime },
    \eta_{3}^{\prime} ).
\end{equation}
Эти ординалы \ $\tau_{1}^{\ast}, \tau_{1}^{\ast \prime},
\tau_{2}^{\ast}, \tau_{3}^{\prime}, \eta_{3}^{\prime} $ \
однозначно определяются через \ $\chi^{\ast}$, $\gamma_{1}^{\ast}
= \widetilde{\delta }^{\ast 3}$ \ ниже \ $\alpha^{\ast 3} $ \ и
нетрудно видеть, что
\[
    \tau_3^{\ast} \le \tau_3^{\prime}, ~~ \eta^{\ast 3} \le
    \eta_3^{\prime}.
\]
Из допустимости \ $S^{\ast 3}$ \ на \ $\alpha_{\tau_3^{\ast}}$ \ и
леммы \ref{8.5.}~6) следует замыкающее условие \ $\mathbf{K}^{0}
(a^{\ast 3}, \widetilde{\delta}^{\ast 3},
\alpha_{\tau_{3}^{\ast}}, \rho^{\ast 3})$:
\begin{eqnarray*}
    \qquad a^{\ast 3} =0 \longrightarrow \forall \tau _{1}^{\prime },
    \tau_{1}^{\prime \prime },\tau _{2}^{\prime },
    \tau _{3}^{\prime },\eta ^{\prime }
    \left[ \gamma _{\tau _{1}^{\prime }}^{<\alpha^{ \ast 3}}
    \leq \widetilde{\delta}^{\ast 3} <
    \gamma _{\tau _{3}^{\prime}}^{<\alpha ^{\ast 3}}
    \wedge \right. \qquad \qquad
    \\
    \left. \wedge A_{4}^{M b \vartriangleleft \alpha^{\ast 3}}
    (\tau _{1}^{\prime }, \tau_{1}^{\prime \prime },
    \tau _{2}^{\prime },\tau _{3}^{\prime },
    \eta^{\prime }, \alpha S_f^{<\alpha^{\ast 3}}, a_f^{<\alpha^{\ast 3}})
    \rightarrow \eta^{\prime} < \rho^{\ast 3} \vee \rho^{\ast 3} = \chi^{\ast +}
    \right],
\end{eqnarray*}
которое и заканчивает это рассуждение следующим образом: Так как
блок \ $[ \gamma_{\tau_1^{\prime}}^{<\alpha^{\ast 3}},
\gamma_{\tau_3^{\prime}}^{<\alpha^{\ast 3}}[$ \ однозначно
определяется через
 \ $\widetilde{\delta}^{\ast 3} =
\gamma_1^{\ast}$ \ ниже  \ $\alpha^{\ast 3} $, \ то
\[
    \widetilde{\delta}^{\ast 3} = \gamma_1^{\ast} =
    \gamma_{\tau_1^{\prime}}^{<\alpha^{\ast 3}}, \quad
    \eta^{\ast 3} \le \eta_3^{\prime} = \eta^{\prime}.
\]
Отсюда, из (\ref{e11.6o}) и (\ref{e11.8o}) следует, что \ $S^{\ast
3}$ \ на \ $\alpha_{\tau_3^{\ast}}$ \ опирается на \
$\widetilde{\delta}^{\ast 3}$ \ очень сильно, то есть
\[
    \alpha S_{\tau _{1}^{\ast \prime }} \vartriangleleft
    \rho^{\ast 3},
\]
в противоречии с  условием этого подслучая.
\\
Случай 2. \quad $a^{\ast2}=0$. \ В этом случае ниже кардинала
предскачка \ $\alpha^{\ast 2} =
\alpha_{\tau_2^{\ast}}^{\Downarrow}$ \ выполняется
\begin{equation*}
    \forall \gamma < \gamma _{2}^{\ast } ~ \exists \tau ~
    ( \gamma < \gamma_{\tau}^{<\alpha^{\ast 2}} \wedge
    a_{\tau }^{<\alpha ^{\ast 2}}=1 )  .
\end{equation*}

\noindent Это утверждение очевидно следует из леммы \ref{11.3.},
потому что существует некоторая лестница \ $\mathcal{S}t$ \
единичных ступеней, завершающаяся в \ $\alpha^{\ast 2}$. \ Отсюда
и из леммы \ref{10.4.} следует существование ординала \ $\tau
_{1}^{\ast \prime }$, \ для которого (\ref{e11.7o}) снова
выполняется (мы сохраняем здесь обозначения из случая~1. для
некоторого удобства). Остаётся буквально повторить рассуждение,
следующее за (\ref{e11.7o}), и получить прежнее противоречие.
Доказательство теоремы~2 закончено.
\\
\hspace*{\fill} $\dashv$
\\

\noindent Подведём итоги.
\\
Все рассуждения проводились в теории
\[
    ZF + \exists k ~ (k \mbox{\it \ слабо недостижимый кардинал});
\]
в ней  рассматривалась счётная стандартная модель
\[
    \mathfrak{M} = (L_{\chi^0}, \in, =)
\]
теории
\[
    ZF + V = L + \exists k ~ (k \mbox{\it \ слабо недостижимый кардинал}),
\]
где всякий слабо недостижимый кардинал является сильно
недостижимым.
\\
В этой модели рассматривались матричные функции; такая функция \
$\alpha S_f^{<\alpha_1}$ \ определена на всяком непустом множестве
\ $T^{\alpha_1}$, \ которое существует для всякого достаточно
большого кардинала \ $\alpha_1 < k$, $\alpha_1 \in SIN_n$ \
согласно лемме~\ref{8.9.}.
\\
Это влечёт заключительное противоречие: рассмотрим любой \
$SIN_n$-кардинал \ $\alpha_1 > \alpha \delta^{\ast}$ \ предельный
для \ $SIN_n \cap \alpha_1$ \ и конфинальности \ $cf(\alpha_1) \ge
\chi^{\ast +}$, \ доставляющий такое непустое множество \
$T^{\alpha_1}$ \ со свойствами $(i)$--$(iii)$ из теоремы~2, тогда
функция \ $\alpha S_f^{<\alpha_1}$ \ немонотонна на этом \
$T^{\alpha_1}$ \ по теореме~1 и в то же время монотонна на этом
множестве по теореме~2.
\\
Это противоречие завершает доказательство основной теоремы.
\\
\hspace*{\fill} $\dashv$

\newpage

\section{Некоторые следствия}

\setcounter{equation}{0}

Вернёмся к началу во Введение в~\cite{Kiselev11}, где приводились
различные хорошо известные связи между Гипотезами больших
кардиналов, Аксиомой детерминированности, регулярными свойствами
континуальных множеств и так далее (см. Дрейк~\cite{Drake},
Канамори~\cite{Kanamori}). Здесь мы укажем некоторые простые
следствия из таких результатов и основной теоремы.

\paragraph{I. Гипотезы больших кардиналов}

\hfill {} \newline \hspace*{1em} Рассмотрим здесь Гипотезы больших
кардиналов, утверждающие существование больших кардиналов
какого-либо вида. Иерархия больших кардиналов располагает их по
``степени недостижимости'' и базируется на (слабо) недостижимых
кардиналах. Существование некоторых из них (кардиналов Мало, слабо
компактных и т.д.) прямо исключается основной теоремой.
Следовательно, не существует кардиналов, обладающих более сильными
партиционными свойствами, например, неописываемых кардиналов,
кардиналов Рамсея, Эрдёша и других; измеримые кардиналы также не
существуют, так как они являются кардиналами  Рамсея. В некоторых
случаях в доказательстве несовместности Гипотез больших кардиналов
может быть использована \ $AC$, \ но можно обойтись и без этого,
извлекая из таких гипотез существование модели \ $ZFC+\exists ~
$недостижимый кардинал
 (см., например, Сильвер~\cite{Silver}).

Мы опускаем переформулировку таких результатов в терминах
фильтров, деревьев, инфинитарных языков и т.д.

\paragraph{II. Сингулярные кардиналы. Шарпы.}

\hfill {} \newline \hspace*{1em} Так как недостижимые кардиналы не
существуют, то всякий несчётный предельный кардинал сингулярен.
Известно, что в \ $ZFC$ \ каждый наследный кардинал регулярен.
Следовательно, каждый несчётный кардинал сингулярен в точности
тогда, когда он пределен.

Хорошо известен замечательный результат Йенсена: отрицание
гипотезы недостижимых кардиналов влечёт Гипотезу сингулярных
кардиналов (см. также Стерн~\cite{Stern}). А несуществование
внутренней модели с измеримым кардиналом влечёт Покрывающую лемму
для стержневой модели \ $K $: \ для каждого несчётного \ $X
\subseteq On $ \ существует \ $Y \in K $ \ такое, что \ $X
\subseteq Y $ \ и \ $|X | = |Y |$.

\noindent Отсюда следует Гипотеза сингулярных кардиналов (Додд,
Йенсен~ \cite{Dodd,Dodd1}). Таким образом, несуществование
недостижимых кардиналов устанавливает Покрывающую лемму и Гипотезу
сингулярных кардиналов.
\newline Эта ситуация проливает новый свет на проблему шарпов.
\newline Хорошо известно, что существование \ $0^{\sharp} $ \
влечёт существование недостижимого в \ $L $ \ кардинала (Гитик,
Магидор, Вудин~\cite {Gitik}). Следовательно, \ $0^{\sharp} $ \ не
существует; этот результат влечёт Покрывающую лемму для \ $L $ \
по выдающейся теореме Йенсена (см.  Девлин, Йенсен~\cite{Devlin}).
\newline Отсюда и из знаменитого результата Кюнена, устанавливающего эквивалентность существования элементарного вложения \ $L \prec L
$ \ и существования \ $0^{\sharp} $, \ следует, что не существует
элементарных вложений \ $ L \prec L $ \ и, далее, не существует
элементарных вложений \ $L_{\alpha} \prec L_{\beta} $ \ с
критической точкой меньшей \ $|\alpha | $.

\paragraph{III. Аксиома детерминированности}

\hfill {} \newline \hspace*{1em} Хорошо известно, что \ $AD$ \
влечёт некоторые Гипотезы больших кардиналов. Например, Соловай
установил, что \ $AD$ \ влечёт измеримость кардинала \ $\omega
_{1}$; \ кардиналы \ $\omega _{2}$, \ $\omega _{\omega +1}$, \
$\omega _{\omega +2} $ \ тоже измеримы (см. также
Клейнберг~\cite{Kleinberg}, Mignone~\cite{Mignone}). Кроме того, \
$AD $ \ устанавливает, что кардиналы \ $\omega _{1}$, \ $\omega
_{2} $ \ являются \ $\delta $-суперкомпактными для некоторого
недостижимого кардинала \ $\delta $ \ (Беккер~\cite {Becker}).
Мыциельский~\cite{Mycielsky} получил выдающийся результат:
совместность
\begin{equation*}
ZF+AD
\end{equation*}
влечёт совместность
\begin{equation*}
ZFC + \exists \mbox{\it\;недостижимый кардинал.}
\end{equation*}
Следовательно, аксиома \ $AD$ \ несовместна, но можно более точно
выделить её несовместную часть:
\\
Именно, аксиома \ $AD ( \Sigma _{2}^{1} ) $ \ равнонепротиворечива
с Гипотезой измеримых кардиналов (Louveau~\cite{Louveau}). Поэтому
существуют недетерминированные \ \mbox{$\Sigma_{2}^{1}$-игры}.

\noindent Также, используя \ $AC_{\omega}({}^{\omega} \omega) $, \
можно доказать \ $\neg Det (\Pi_{1}^{1}) $ \ (см.
Канамори~\cite{Kanamori}). Поэтому в
\begin{equation*}
ZF+AC_{\omega}({}^{\omega} \omega)
\end{equation*}
существуют недетерминированные \ $\Pi_{1}^{1} $-игры. Этот
результат вряд ли можно улучшить, потому что все \
$\bigtriangleup_{1}^{1}$-игры детерминированы (Мартин~
\cite{Martin}).

\paragraph{IV. Континуальные множества}

\hfill {} \newline \hspace*{1em} Хорошо известен  ряд выдающихся
результатов, устанавливающих связи между регулярными свойствами
континуальных множеств и большими кардиналами (по относительной
непротиворечивости). Например, Шелах~\cite{Shelar} установил
необходимость недостижимых кардиналов для утверждения, что все
множества действительных чисел измеримы по Лебегу; измеримость \
$\Sigma _{3}^{1}$-множеств влечёт недостижимость \ $\omega _{1} $
\ в \ $L $ \ (см. также Raisonnier~\cite {Raisonnier}). Отсюда
следует существование неизмеримого \ $\Sigma_{3}^{1} $-множества
действительных чисел. Аналогично, непротиворечивость
\[
    ZF+DC+ \mbox{\it каждое несчётное множество } \qquad\qquad\qquad\qquad\qquad
\]
\vspace{-18pt}
\[
    \qquad \mbox{\it действительных чисел имеет совершенное ядро }
\]
не может быть доказана без Гипотезы измеримых кардиналов (см.
Мыциельский~ \cite{Mycielsky}). Кроме того, следующие гипотезы
равнонепротиворечивы над $ZF$:
\\
\quad \\
1) \hspace{\stretch{0.4}} $AC+\exists \mbox{\it \ недостижимый \
кардинал}$; \hspace{\stretch{0.6}}
\\
\quad \\
2) \hspace{\stretch{0.4}} $DC +$ ~ {\it каждое несчётное множество
действительных чисел имеет совершенное ядро};
\hspace{\stretch{0.6}}
\\
\quad \\
3) \hspace{\stretch{0.4}} $\omega_{1} \mbox{\it \ \ регулярен \ }
+ \ \forall a \in {}^{\omega} \omega \left ( \omega_{1}^{L[a]}<
\omega_{1} \right )$ \hspace{\stretch{0.6}}
\\
\quad \\
(см. Соловай~\cite{Solovay}, Шпекер~\cite{Specer},
Леви~\cite{Levy}). Поэтому $DC$ влечёт существование множества
действительных чисел без совершенного ядра и регулярность \
$\omega_{1} $ \ влечёт

\begin{equation*}
  \omega_{1}^{L[a]} \nless  \omega_{1}
\end{equation*}
\vspace{0pt}

\noindent для некоторого \ $a \in {}^{\omega} \omega$. \ Также
хорошо известно, что Гипотеза слабо компактного кардинала
равнонепротиворечива (над \ $ZFC~$) с утверждением о регулярных
свойствах всех континуальных несчётных множеств в \ $ZF+MA$ \
(Харингтон, Шелах~\cite{Harrington}). Следовательно, основная
теорема влечёт существование несчётных множеств \ $\subseteq
\omega_{\omega}$ \ без регулярных свойств.

\noindent Эти результаты можно уточнить; например, Соловай~
\cite{Solovay1} установил, что для всякого \ $a \in
{}^{\omega}\omega$ \quad $ \omega_{1}^{L[a]} < \omega_{1}$ \
эквивалентно свойству совершенного ядра всякого \ $\Pi_{1}^{1}(a)
$-множества действительных чисел. Отсюда следует, что регулярность
\ $\omega_{1} $ \ влечёт существование \ $\Pi_{1}^{1}(a)
$-множества действительных чисел без совершенного ядра для
некоторого \ $a \in {}^{\omega}\omega $.

\paragraph{V. Аксиома Мартина}

\hfill {} \newline \hspace*{1em} Следствия из результатов
Харингтона, Шелаха~\cite {Harrington} и основной теоремы,
упомянутые выше, можно сформулировать более точно, например,
следующим образом:
\\
\noindent $MA$ \ влечёт существование \ $\Delta
_{3}^{1}$-множеств, лишённых свойства Бэра, и неизмеримых \
$\Sigma _{3}^{1}$-множеств \ $\subseteq {}^{\omega}\omega $.

\noindent Кановей~\cite{Kanovei} установил в
\begin{equation*}
ZFC+MA+\left| R\right| >\omega _{1}+\forall x\subseteq  \omega _{1}~~\omega
_{1}^{L\left[ x\right] }=\omega _{0}
\end{equation*}
совместность
\begin{equation*}
ZFC+\exists \mbox{\it\;недостижимый кардинал}.
\end{equation*}
Следовательно, \ $MA$ \ несовместна с
\begin{equation*}
ZFC+\left| R\right| >\omega _{1}+\forall x\subseteq  \omega _{1}~~\omega
_{1}^{L\left[ x\right] }=\omega _{0}.
\end{equation*}

Ряд других следствий основной теоремы  слишком велик, чтобы здесь
их перечислить, поэтому автор намеревается представить более
детальный анализ таких следствий в дальнейших публикациях.

\theendnotes

\newpage

\thispagestyle{empty}

\begin{center}
\quad \\
\quad \\
\quad \\
\quad \\
\quad \\
\quad \\
{Научное издание} \\
\quad \\
{\bf Киселев} Александр Анатольевич\\
\quad \\
{\bf НЕДОСТИЖИМОСТЬ \\
И \\
СУБНЕДОСТИЖИМОСТЬ}\\
\quad \\
{В двух частях}\\
{Часть II}\\
\quad \\
{\small Ответственный за выпуск {\em Т. Е. Янчук}} \\
\quad \\
{\footnotesize Подписано в печать 10.07.2011. Формат
60$\times$84$^{1}/_{16}$. Бумага офсетная. Ризография. Усл. печ.
л. 9,07. Уч.-изд. л. 8,2.
\\
Тираж 100 экз. Заказ 473.} \\
\quad \\
{\footnotesize Республиканское унитарное предприятие} \\
{\footnotesize ``Издательский центр Беларуского Государственного Университета''} \\
{\footnotesize \selectlanguage{russian} ЛИ \No
\selectlanguage{russian}
02330/0494361 от 16.03.2009.} \\
{\footnotesize Ул. Красноармейская, 6, 220030, Минск.}\\
\quad \\
{\footnotesize Отпечатано с оригинал-макета заказчика
\\
в
Республиканском унитарном предприятии} \\
{\footnotesize ``Издательский центр Беларуского Государственного Университета''} \\
{\footnotesize \selectlanguage{russian} ЛП \No
\selectlanguage{russian} 02330/0494178 от 03.04.2009.} \\
{\footnotesize Ул. Красноармейская, 6, 220030, Минск.}\\
\end{center}
\label{end}

\end{document}